\documentclass[a4paper,10pt,leqno]{article} \usepackage[english]{babel}
\usepackage{german, umlaute, amsmath, amssymb , latexsym, theorem, epic,
epsfig, rotating}

\newcommand{\alphlist}{\begin{list}{(\alph{enumi})}{\usecounter{enumi}}}
\newcommand{\romanlist}{\begin{list}{(\roman{enumi})}{\usecounter{enumi}}}
\newcommand{\listend}{\end{list}}

\newcommand{\ld}{\ensuremath{,\ldots,}}
\newcommand{\ssq}{\ensuremath{\subseteq}}

\newcommand{\supnorm}[1]{\ensuremath{\parallel #1 \parallel_{\infty}}}

\newcommand{\ra}{\ensuremath{\rightarrow}}
\newcommand{\equi}{\ensuremath{\Leftrightarrow}}
\newcommand{\follows}{\ensuremath{\Rightarrow}}

\newcommand{\N}{\ensuremath{\mathbb{N}}} 
\newcommand{\R}{\ensuremath{\mathbb{R}}}
\newcommand{\Z}{\ensuremath{\mathbb{Z}}}
\newcommand{\Q}{\ensuremath{\mathbb{Q}}}

\newcommand{\kreis}{\ensuremath{\mathbb{T}^{1}}}

\newcommand{\nLim}{\ensuremath{\lim_{n\rightarrow\infty}}}

\newcommand{\nKonv}{\ensuremath{\stackrel{n\rightarrow
      \infty}{\longrightarrow}}}

\newcommand{\pKonv}{\ensuremath{\stackrel{p\rightarrow
      \infty}{\longrightarrow}}}

\newcommand{\solidqed}{{\raggedleft \large $\blacksquare$ \\}}

\newcommand{\proof}{\textit{Proof.}}
\newcommand{\qed}{{\raggedleft $\Box$ \\}}

\newcommand{\ind}{\ensuremath{\mathbf{1}}}

\newcommand{\twoscriptcup}[2]{ \ensuremath{ \bigcup_{ \begin{array}{c}
\scriptstyle #1 \\  \scriptstyle #2 \end{array} }}}

\newcommand{\insum}{\ensuremath{\sum_{i=1}^n}}

\newcommand{\inergsum}{\ensuremath{\sum_{i=0}^{n-1}}}

\newcommand{\dtel}{\ensuremath{\frac{1}{d}}}

\newcommand{\ntel}{\ensuremath{\frac{1}{n}}}

\newcommand{\halb}{\ensuremath{\frac{1}{2}}}
\newcommand{\drittel}{\ensuremath{\frac{1}{3}}}
\newcommand{\viertel}{\ensuremath{\frac{1}{4}}}

\newcommand{\thx}{\ensuremath{(\theta,x)}}
\newcommand{\thom}{\ensuremath{\theta + \omega}}

\newcommand{\phiplus}{\ensuremath{\varphi^+}}
\newcommand{\phimin}{\ensuremath{\varphi^-}}

\newcommand{\Tth}{\ensuremath{T_{\theta}}}

\newcommand{\Tthx}{\ensuremath{T_{\theta}(x)}}

\newcommand{\That}{\ensuremath{\widehat{T}}}

\title{\textsc{The creation of strange non-chaotic attractors in non-smooth
  saddle-node bifurcations}} \author{Tobias Henrik J\"ager \\
  Friedrich-Alexander-Universit\"at Erlangen-N\"urnberg}

\newtheorem{definition}{Definition}[section] 
\newtheorem{thm}[definition]{Theorem}
 
\newtheorem{lem}[definition]{Lemma} 
\newtheorem{cor}[definition]{Corollary}

\newtheorem{indscheme}[definition]{Induction scheme}
\newtheorem{adde}[definition]{Addendum}

\theorembodyfont{\rmfamily}
{\theoremstyle{plain} \newtheorem{claim}[definition]{Claim} }
\newtheorem{bem}[definition]{Remark}

\numberwithin{equation}{section}
\numberwithin{figure}{section}

\newcommand{\walpha}{\ensuremath{\sqrt{\alpha}}}
\newcommand{\alphtel}{\ensuremath{\frac{1}{\alpha}}}
\newcommand{\walphtel}{\ensuremath{\frac{1}{\sqrt{\alpha}}}}

\newcommand{\Balphcl}{\ensuremath{\overline{B_\frac{1}{\alpha}(0)}}}
\newcommand{\Bdalph}{\ensuremath{B_\frac{2}{\alpha}(0)}}
\newcommand{\omn}{\ensuremath{\omega_n}}
\newcommand{\omm}{\ensuremath{\omega_m}}
\newcommand{\omj}{\ensuremath{\omega_j}}

\newcommand{\omk}{\ensuremath{\omega_k}}

\newcommand{\JN}{\ensuremath{{\cal J}_N}} \newcommand{\J}{\ensuremath{{\cal
J}}} \newcommand{\Omtil}{\ensuremath{\tilde{\Omega}}}
\newcommand{\nutil}{\ensuremath{\tilde{\nu}}}

\newcommand{\Ttil}{\ensuremath{\tilde{T}}}
\newcommand{\Ftil}{\ensuremath{\tilde{F}}}

\newcommand{\betil}{\ensuremath{\tilde{\beta}}}
\newcommand{\util}{\ensuremath{\tilde{u}}}
\newcommand{\vtil}{\ensuremath{\tilde{v}}}
\newcommand{\thxo}{\ensuremath{(\theta_0,x_0)}}
\newcommand{\esscl}[1]{\ensuremath{\overline{#1}^\mathit{ess}}}
\newcommand{\err}{\ensuremath{\textrm{\upshape err}}}
\newcommand{\smin}{\ensuremath{\setminus}}

\begin{document}  

\setlength{\oddsidemargin}{0.08\textwidth}
\setlength{\evensidemargin}{0.12\textwidth}

\maketitle 

\begin{abstract}   
  We propose a general mechanism by which strange non-chaotic attractors (SNA)
  are created during the collision of invariant curves in quasiperiodically
  forced systems. This mechanism, and its implementation in
  different models, is first discussed on an heuristic level and by means of
  simulations. In the considered examples, a stable and an unstable
  invariant circle undergo a saddle-node bifurcation, but instead of a neutral
  invariant curve there exists a strange non-chaotic attractor-repeller pair
  at the bifurcation point. This process is accompanied by a very
  characteristic behaviour of the invariant curves prior to their collision,
  which we call \emph{`exponential evolution of peaks'}.

  This observation is then used to give a rigorous description of non-smooth
  saddle-node bifurcations and to prove the existence of SNA in certain
  parameter families of quasiperiodically forced interval maps. The
  non-smoothness of the bifurcations and the occurrence of SNA is established
  via the existence of \textit{`sink-source-orbits'}, meaning orbits with
  positive Lyapunov exponent both forwards and backwards in time. 
 
  The important fact is that the presented approach allows for a certain
  amount of flexibility, which makes it possible to treat different models at
  the same time - even if the results presented here are still subject to a
  number of technical constraints. This is unlike previous proofs for the
  existence of SNA, which are all restricted to very specific classes and
  depend on very particular properties of the considered systems. In order to
  demonstrate this flexibility, we also discuss the application of the results
  to the Harper map, an example which is well-known from the study of discrete
  Schr\"odinger operators with quasiperiodic potentials. Further, we prove the
  existence of strange non-chaotic attractors with a certain inherent
  symmetry, as they occur in non-smooth pitchfork bifurcations.
\end{abstract} 
 
\pagebreak 

\tableofcontents 
 
\pagebreak  


\section{Introduction} 

In the early 1980's, Herman \cite{herman:1983} and Grebogi et al.\
\cite{grebogi/ott/pelikan/yorke:1984} independently discovered the existence
of strange non-chaotic attractors (SNA's) in quasiperiodically forced (qpf)
systems. These objects combine a complicated geometry%
\footnote{This means in particular that they are not a piecewise
  differentiable (or even continuous) sub-manifold of the phase space.}
 with non-chaotic dynamics, a combination which is rather
  unusual and has only been observed in a few very particular cases before
  (the most prominent example is the Feigenbaum map, see \cite{milnor:1985}
  for a discussion and further references). In quasiperiodically forced
  systems, however, they seem to occur quite frequently and even over whole
  intervals in parameter space
  \cite{grebogi/ott/pelikan/yorke:1984,keller:1996,%
  osinga/wiersig/glendinning/feudel:2001}.  As a novel phenomenon this evoked
  considerable interest in theoretical physics, and in the sequel a large
  number of numerical studies explored the surprisingly rich dynamics of these
  relatively simple maps. In particular, the widespread existence of SNA's was
  confirmed both numerically (see
  \cite{bondeson/ott/antonsen:1985}--\nocite{romeiras/ott:1987}%
  \nocite{romeiras/etal:1987}%
  \nocite{ding/grebogi/ott:1989}\nocite{heagy/hammel:1994}%
  \nocite{pikovski/feudel:1995}%
  \nocite{feudel/kurths/pikovsky:1995}\nocite{witt/feudel/pikovsky:1997}%
  \nocite{chastell/glendinning/stark:1995}%
  \nocite{glendinning:1998}\nocite{sturman:1999}%
  \nocite{glendinning/feudel/pikovsky/stark:2000}%
  \nocite{negi/prasad/ramaswamy:2000}\cite{prasad/negi/ramaswamy:2001}, just
  to give a selection) and even experimentally
  \cite{Ditto/etal:1990,heagy/ditto:1991,zhou/moss/bulsara:1992}. Further, it
  turned out that SNA play an important role in the bifurcations of invariant
  circles \cite{osinga/wiersig/glendinning/feudel:2001,%
  chastell/glendinning/stark:1995,negi/prasad/ramaswamy:2000,haro/delallave:2006}.

 The studied systems were either discrete time maps, such as the qpf logistic
map \cite{heagy/hammel:1994,witt/feudel/pikovsky:1997,%
negi/prasad/ramaswamy:2000} and the qpf Arnold circle map
\cite{osinga/wiersig/glendinning/feudel:2001,ding/grebogi/ott:1989%
,feudel/kurths/pikovsky:1995,chastell/glendinning/stark:1995}, or skew product
flows which are forced at two or more incommensurate frequencies.  Especially
the latter underline the significance of qpf systems for understanding
real-world phenomena, as most of them were derived from models for different
physical systems (e.g.\ quasiperiodically driven damped pendula and Josephson
junctions \cite{bondeson/ott/antonsen:1985,romeiras/ott:1987,%
romeiras/etal:1987} or Duffing oscillators \cite{heagy/ditto:1991}. Their
Poincar\'e maps again give rise to discrete-time qpf systems, on which the
present article will focus.
\medskip

However, despite all efforts there are still only very few mathematically
rigorous results about the subject, with the only exception of qpf
Schr\"odinger cocycles (see below). There are results concerning the
regularity of invariant curves (\cite{stark:1999}, see also
\cite{sturman/stark:2000}), and there has been some progress in carrying over
basic results from one-dimensional dynamics
\cite{jaeger/keller:2005,fabbri/jaeger/johnson/keller:2004,%
jaeger/stark:2005}. But so far, the two original examples in
\cite{herman:1983} and \cite{grebogi/ott/pelikan/yorke:1984} remain the only
ones for which the existence of SNA's has been proved rigorously. In both
cases, the arguments used were highly specific for the respective class of
maps and did not allow for much further generalisation, nor did they give very
much insight into the geometrical and structural properties of the attractors.

The systems Herman studied in \cite{herman:1983} were matrix cocycles, with
quasiperiodic Schr\"odin\-ger cocycles as a special case. The linear structure
of these systems and their intimate relation to Schr\"odinger operators with
quasiperiodic potential made it possible to use a fruitful blend of techniques
from operator theory, dynamical systems and complex analysis, such that by now
the mathematical theory is well-developed and deep results have been obtained
(see \cite{avila/krikorian:2004} and \cite{avila/jitomirskaya:2005} for recent
advances and further reference).  However, as soon as the particular class of
matrix cocycles is left, it seems hard to recover most of these arguments.  One
of the rare exceptions is the work of Bjerkl\"ov in \cite{bjerkloev:2005a}
(taken from \cite{bjerkloev:2003}) and \cite{bjerkloev:2005}, which is based
on a purely dynamical approach and should generalise to other types of
systems, such as the ones considered here. (In fact, although implemented in a
different way the underlying idea in \cite{bjerkloev:2005} is very similar to
the one presented here, such that despite their independence the two articles
are closely related.)

On the other hand, for the so-called `\textit{pinched skew products}'
introduced in \cite{grebogi/ott/pelikan/yorke:1984}, establishing the
existence of SNA is surprisingly simple and straightforward (see
\cite{keller:1996} for a rigorous treatment and also \cite{glendinning:2002}
and \cite{jaeger:2004a}).  But one has to say that these maps were introduced
especially for this purpose and are rather artificial in some aspects. For
example, it is crucial for the argument that there exists at least one fibre
which is mapped to a single point. But this means that the maps are not
invertible and can therefore not be the Poincar\'e maps of any flow.
\medskip

\begin{figure}[h!]  
\noindent 
\begin{minipage}[t]{\linewidth}  
  \epsfig{file=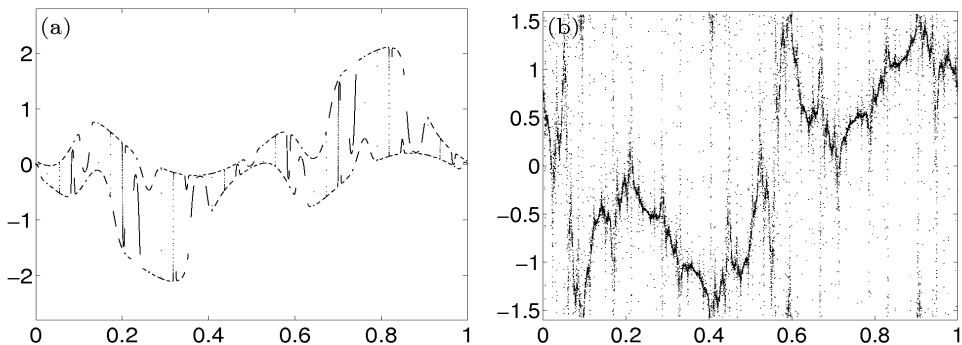, clip=, width=\linewidth,
    height=0.23\textheight}
         \caption{\small Two different types of strange non-chaotic
           attractors: The left picture shows a `strip-like' SNA in the system
           $\thx \mapsto (\thom,\tanh(5x) + 1.2015\cdot
           \sin(2\pi\theta))$. The topological closure of this object is
           bounded above and below by semi-continuous invariant graphs
           (compare (\ref{eq:boundinggraphs})). This is the type of SNA's that
           will be studied in the present work.  The right picture shows a
           different type that occurs for example in the critical Harper map
           (Equation~(\ref{eq:projharper}) with $\lambda=2$ and $E=0$; more
           details can be found in~\cite{datta/jaeger/keller/ramaswamy:2004}),
           where no such boundaries exist. In both cases $\omega$ is the
           golden mean.}
         \label{fig:twotypes}
\end{minipage} 
\end{figure}

The main goal of this article is to prove the existence of SNA in certain
parameter families of qpf systems where this has not been possible
previously. Thereby, we will concentrate on a particular type of SNA, namely
`strip-like' ones, which occur in saddle-node and pitchfork bifurcations of
invariant circles (see Figure~\ref{fig:twotypes}, for a more precise
formulation consider the definition of invariant strips in
\cite{jaeger/keller:2005} and \cite{jaeger/stark:2005}).
In such a saddle-node bifurcation, a stable and an unstable invariant circle
approach each other, until they finally collide and then vanish.  However,
there are two different possibilities. In the first case, which is similar to
the one-dimensional one, the two circles merge together uniformly to form one
single and neutral invariant circle at the bifurcation point. But it may also
happen that the two circles approach each other only on a dense, but
(Lebesgue) measure zero set of points. In this case, instead of a single
invariant circle, a strange non-chaotic attractor-repeller-pair is created at
the bifurcation point. Attractor and repeller are interwoven in such a way,
that they have the same topological closure. This particular route for the
creation of SNA's has been observed quite frequently
(\cite{feudel/kurths/pikovsky:1995,chastell/glendinning/stark:1995,%
glendinning:1998,prasad/negi/ramaswamy:2001}, see also
\cite{heagy/hammel:1994}) and was named \textit{`non-smooth saddle-node
bifurcation'} or \textit{`creation of SNA via torus collision'}. The only
rigorous description of this process so far was given by Herman in
\cite{herman:1983}. In a similar way, the simultaneous collision of two stable
and one unstable invariant circle may lead to the creation of two SNA's
embracing one strange non-chaotic repeller
\cite{osinga/wiersig/glendinning/feudel:2001,sturman:1999}.
\medskip

\textbf{Acknowledgments.} The results presented here were part of my thesis,
and I would like to thank Gerhard Keller for his invaluable advice and support
during all the years of my PhD-studies. I am also greatful to an anonymous
referee, whose thoughtful remarks greatly improved the manuscript. This work
was supported by the German Research Foundation (DFG), grant Ke 514/6-1.


\subsection{Overview} \label{Overview}

As mentioned above, the main objective of this article is to provide new
examples of SNA, by describing a general mechanism which is responsible for
the creation of SNA in non-smooth saddle-node bifurcations. While this
mechanism might not be the only one which exists, it seems to be common in a
variety of different models, including well-known examples like the Harper map
or the qpf Arnold circle map. The evidence we present will be two-fold: In the
remainder of this introduction we will explain the basic idea, and discuss on
an heuristic level and by means of numerical simulations how it is implemented
in the two examples just mentioned and a third parameter family, which we call
$\arctan$-family. An analogous phenomenom is also observed in so-called
Pinched skew products, first introduced in
\cite{grebogi/ott/pelikan/yorke:1984}, even if no bifurcation takes place in
these systems.

The heuristic arguments given in the introduction are then backed up by
Theorem~\ref{thm:snaexistence}, which provides a rigorous criterium for the
non-smoothness of saddle-node bifurcations in qpf interval maps. This leads to
new examples of strange non-chaotic attractors, and the result is flexible
enough to apply to different parameter families at the same time, provided
they have similar qualitative features and share a certain scaling
behaviour. Nevertheless, it must be said that there is still an apparent gap
between what can be expected from the numerical observations and what can be
derived from Theorem~\ref{thm:snaexistence}~. For instance, the latter does
not apply to the forced version of the Arnold circle map, and for the
application to the $\arctan$-family and the Harper map we have to make some
quite specific assumptions on the forcing function and the potential,
respectively. (Namely that these have a unique maximum and decay linearly in a
neighbourhood of it). However, our main concern here is just to show that the
general approach we present does lead to rigorous results at all, even if
these are still far from being optimal. The present work should therefore be
seen rather as a first step in this direction, which will hopefully inspire
further research, and not as an ultimate solution.
\medskip

The article is organised as follows: After we have given some basic
definitions, we will introduce our main examples in
Section~\ref{Examples}~. As mentioned, these are the arctan-family with
additive forcing, the Harper map, the qpf Arnold circle map and Pinched skew
products. The simulations we present mostly show the evolution of stable
invariant curves as the system parameters are varied. The crucial observation
is the fact that the behaviour of these curves prior to the bifurcation
follows a very characteristic pattern, which we call \emph{`exponential
evolution of peaks'}. In all the first three examples the qualitative features
of this process are similar, and even in Pinched skew products, where no
saddle-node bifurcation occurs, an analogue behaviour can be
observed. Finally, a slight modification of the arctan-family is used to
illustrate that the phenomenom is also present in non-smooth pitchfork
bifurcations.

On an heuristic level it is not difficult to give an explanation for this
behaviour, and this will be done in Section~\ref{Mechanism}~. The simple
geometric intuition obtained there will then determine the strategy for the
rigorous proof of the non-smoothness of the bifurcations in the later
sections.  More precisely, the heuristics indicate why the existence of SNA
should be linked to the appearance of sink-source-orbits in these situations,
and this will be one of the main ingredients of the proof.

Section~\ref{Main results} then contains the statement of our main results and
discusses their application to the examples from the introduction (or why such
an application is not possible, in the case of the qpf Arnold circle
map). Before we can turn to the existence of SNA and the
non-smoothness of bifurcations, we need to state two preliminary results. The
first, Theorem~\ref{thm:saddlenode}, provides a general framework in which
saddle-node bifurcations in qpf interval maps take place (smooth or
non-smooth). The second, Theorem~\ref{thm:sinksourcesna}, states that the
existence of sink-source-orbits%
\footnote{Orbits with positive Lyapunov exponent both forwards and backwards
  in time, see Definition~\ref{def:sinksource}~.}
 implies the existence of SNA's (although the converse is not true).  After
these statements and some related concepts have been introduced in Sections
\ref{GeneralSettingI} and \ref{SinkSourceSNA}, we can turn to the main result,
namely Theorem~\ref{thm:snaexistence}, which provides a criterium for the
existence SNA's created in non-smooth saddle-node bifurcations. The
counterpart for non-smooth pitchfork bifurcations is
Theorem~\ref{thm:symmetricsna}, which gives a criterium for the existence of
symmetric SNA's. More precisely, under the assertions of this theorem there
exists a triple consisting of two SNA, symmetric to each other, which embrace
a self-symmetric strange non-chaotic repeller. These objects are presumably
created by the simultaneous collision of two stable and one unstable invariant
curve. However, as the considered parameter families lack a certain
monotonicity property which is present in the situation of
Theorem~\ref{thm:snaexistence}, we cannot describe the bifurcation pattern in
a rigorous way as for the saddle-node bifurcations, such that the existence of
SNA is the only conclusion we draw in the symmetric setting.  The application
of these results to the arctan-family and the Harper map is then discussed in
detail in Section~\ref{Applications}, which resumes the structure of
Section~\ref{Examples} where these examples are introduced. As we have
mentioned before, the statement of Theorem~\ref{thm:snaexistence} is too
restricted to apply to the qpf Arnold circle map. However, in
Section~\ref{App:Remarks} we discuss some possible modifications, which might
allow to treat this example in a similar way, at least for particular forcing
functions.

Section~\ref{Generalsetting} provides the proofs for the more elementary
results (namely Theorems \ref{thm:saddlenode} and
\ref{thm:sinksourcesna}). All the remaining sections are then dedicated to the
proof of Theorems~\ref{thm:snaexistence} and \ref{thm:symmetricsna}, starting
with an outline of the construction in Section~\ref{Strategy}.


\subsection{Basic definitions and notations}  \label{BasicDefinitions}
A \textit{quasiperiodically forced (qpf) system} is a continuous map of the
form
\begin{equation} \label{eq:qpfs} 
    T : \kreis \times X \ra \kreis \times X \ \ \ , \ \ \ \thx \mapsto
    (\thom,\Tthx)
\end{equation}  
with irrational driving frequency $\omega$. At most times, we will restrict to
the case where the driving space $X=[a,b]$ is a compact interval and the
\textit{fibre maps} \Tth\ are all monotonically increasing on $X$. In this
case we say $T$ is a {\em qpf monotone interval map}. Some of the introductory
examples will also be qpf circle homeomorphisms, but there the situation can
often be reduced to the case of interval maps as well, for example when there
exists a closed annulus which is mapped into itself.

Two notations which will be used frequently are the following: Given any set
    $A \ssq \kreis \times X$ and $\theta \in \kreis$, we let $A_\theta :=
    \{x\in X\mid \thx \in A\}$. If $X=\R$ and $\varphi,\psi:\kreis \ra \R$ are
    two measurable functions, then we use the notation
\begin{equation}
  \label{eq:graphintervals}
  [\psi,\varphi] \ := \ \{ \thx \mid \psi(\theta) \leq x \leq
  \varphi(\theta) \} 
\end{equation}
similarly for $(\psi,\varphi)$, $(\psi,\varphi]$,
$[\psi,\varphi)$. 

Due to the minimality of the irrational rotation on the base there are
no fixed or periodic points for $T$, and one finds that the simplest
invariant objects are invariant curves over the driving space (also
invariant circles or invariant tori). More generally, a
($T$-)\textit{invariant graph} is a measurable function $\varphi :
\kreis \ra X$ which satisfies
\begin{equation} \label{eq:invgraph}
    \Tth(\varphi(\theta)) \ = \ \varphi(\thom) \quad \quad \forall \theta \in
    \kreis \  .
\end{equation}
This equation implies that the point set $\Phi:= \{(\theta,\varphi(\theta))
\mid \theta \in \kreis\}$ is forward invariant under $T$. As long as no
ambiguities can arise, we will refer to $\Phi$ as an invariant graph as well.

There is a simple way of obtaining invariant graphs from compact
invariant sets: Suppose $A \ssq \kreis \times X$ is $T$-invariant.
Then
\begin{equation} \label{eq:boundinggraphs}
    \phiplus_A(\theta) \ := \ \sup\{x\in X \mid \thx \in A \} 
\end{equation} 
defines an invariant graph (invariance following from the monotonicity
of the fibre maps). Furthermore, the compactness of $A$ implies that
$\phiplus_A$ is upper semi-continuous (see \cite{stark:2003}). In a
similar way we can define a lower semi-continuous graph $\phimin_A$ by
taking the infimum in (\ref{eq:boundinggraphs}). Particularly
interesting is the case where $A = \cap_{n\in\N}T^n(\kreis \times X)$
(the so-called global attractor, see \cite{glendinning:2002}). Then we
call $\phiplus_A$ ($\phimin_A)$ the \textit{upper (lower) bounding
  graph of the system}.

There is also an intimate relation between invariant graphs and
ergodic measures. On the one hand, to each invariant graph $\varphi$
we can associate an invariant ergodic measure by
\begin{equation} \label{eq:associatedmeasure}
  \mu_\varphi(A) \ := \ m(\pi_1(A \cap \Phi)) \ , 
\end{equation} 
where $m$ denotes the Lebesgue measure on \kreis\ and $\pi_1$ is the
projection to the first coordinate. On the other hand, if $f$ is a qpf
monotone interval maps then the converse is true as well: In this case, for
each invariant ergodic measure $\mu$ there exists an invariant graph
$\varphi$, such that $\mu=\mu_\varphi$ in the sense of
(\ref{eq:associatedmeasure}).  (This can be found in \cite{arnold:1998},
Theorem 1.8.4~. Although the statement is formulated for continuous-time
dynamical systems there, the proof literally stays the same.)

If all fibre maps are differentiable and we denote their derivatives
by $D\Tth$, then the stability of an invariant graph $\varphi$ is
measured by its \textit{Lyapunov exponent}%
\begin{equation} \label{eq:lyapexponent} 
  \lambda(\varphi) \ := \ \int_{\kreis} \log D\Tth(\varphi(\theta)) \
  d\theta \ . 
\end{equation} 
An invariant graph is called \textit{stable} when its Lyapunov
exponent is negative, \textit{unstable} when it is positive and
\textit{neutral} when it is zero.

Obviously, even if its Lyapunov exponent is negative an invariant graph does
not necessarily have to be continuous. This is exactly the case that has been
the subject of so much interest:
\begin{definition}[Strange non-chaotic attractors and repellers] 
  \label{def:sna} 
  A \textbf{strange non-chaotic attractor (SNA)} in a quasiperiodically forced
  system $T$ is a $T$-invariant graph which has negative Lyapunov exponent and
  is not continuous. Similarly, a \textbf{strange non-chaotic repeller (SNR)}
  is a non-continuous $T$-invariant graph with positive Lyapunov exponent.
\end{definition} 
This terminology, which was coined in theoretical physics, may need a little
bit of explanation. For example, the point set corresponding to a
non-continuous invariant graph is not a compact invariant set, which is
usually required in the definition of `attractor'. However, a SNA attracts
and determines the behaviour of a set of initial conditions of positive
Lebesgue measure (e.g.\ \cite{jaeger:2003}, Proposition~3.3), i.e.\ it carries
a `physical measure'.  Moreover, it is easy to see that the essential closure%
\footnote{The support of the measure $\mu_\varphi$ given by
  (\ref{eq:associatedmeasure}), where $\varphi$ denotes the SNA. See also
  Section~\ref{Essentialclosure}.} 
of a SNA is an attractor in the sense of Milnor \cite{milnor:1985}. `Strange'
just refers to the non-continuity and the resulting complicated structure of the
graph. The term `non-chaotic' is often motivated by the negative Lyapunov
exponent in the above definition \cite{grebogi/ott/pelikan/yorke:1984},
but actually we prefer a slightly different point of view: At least in the
case where the fibre maps are monotone interval maps or circle homeomorphisms,
the topological entropy of a quasiperiodically forced system is always zero,%
\footnote{For monotone interval maps this follows simply from the fact that
  every invariant ergodic measure is the projection of the Lebesgue measure on
  \kreis\ onto an invariant graph, such that the dynamics are isomorphic in
  the measure-theoretic sense to the irrational rotation on the
  base. Therefore all measure-theoretic entropies are zero, and so is the
  topological entropy as their supremum. In the case of circle homeomorphisms,
  the same result can be derived from a statement by Bowen (\cite{bowen:1971},
  Theorem~17).}
such that the system and its invariant objects should not be
considered as `chaotic'. This explains why we also speak of
\textit{strange non-chaotic repellers}. In fact, in invertible systems
an attracting invariant graph becomes a repelling invariant graph for
the inverse and vice versa, while the dynamics on them hardly changes.
Thus, it seems reasonable to say that `non-chaotic' should either
apply to both or to none of these objects.


\subsection{Examples of non-smooth saddle-node bifurcations} \label{Examples}

As mentioned, the crucial observation which starts our investigation here is
the fact that the invariant circles in a non-smooth bifurcation do not
approach each other arbitrarily. Instead, their behaviour follows a very
distinctive pattern, which we call \emph{exponential evolution of peaks}. In
this section we present some simulations which demonstrate this phenomenom in
the different parameter families mentioned in Section~\ref{Overview}~.
Although it seems difficult to give a precise mathematical definition of this
process, and we refrain from doing so here, this observation provides the
necessary intuition and determines the strategy of the proofs for the rigorous
results in the later chapters. (The same underlying idea can be found in
\cite{bjerkloev:2005} and \cite{jaeger:2004a}.)


\subsubsection{The arctan-family with additive forcing} \label{AtanFamily}

Typical representatives of the class of systems we will study in the later
sections are given by the family
\begin{equation} 
  \label{eq:atanfamily}  
  \thx \ \mapsto \ \left(\thom,\frac{\arctan(\alpha 
      x)}{\arctan(\alpha)}  - \beta\cdot (1-\sin(\pi\theta))\right) \ 
  . 
\end{equation}   
\begin{figure}[h!]  
\noindent  
\begin{minipage}[t]{\linewidth}  
  \epsfig{file=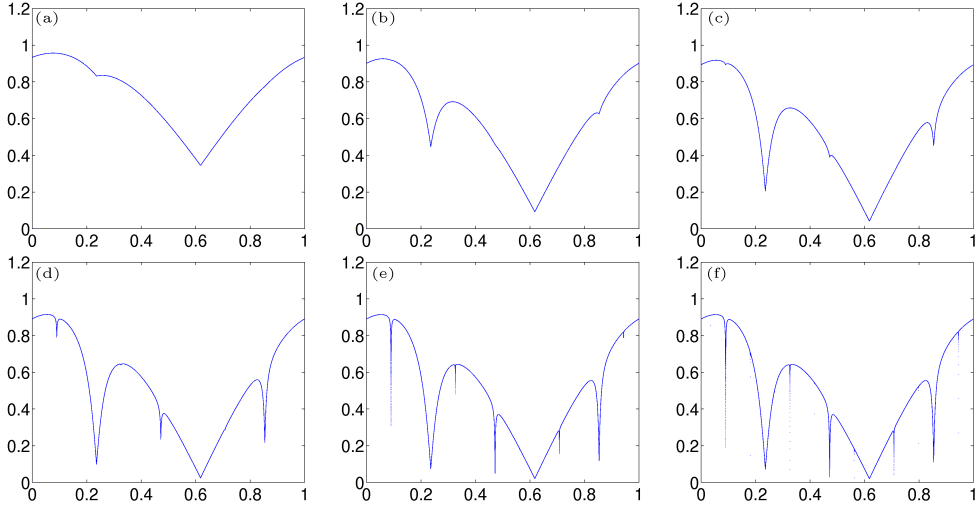, clip=, width=\linewidth, }
         \caption{\small Upper bounding graphs in the parameter family given
           by (\ref{eq:atanfamily}) with $\omega$ the golden mean and
           $\alpha=10$. The values for $\beta$ are: (a) $\beta=0.65$,
           (b) $\beta=0.9$, (c) $\beta=0.95$, (d) $\beta=0.967$, (e)
           $\beta=0.9708$, (f) $\beta=0.9710325$.}
         \label{fig:atanfamily}
\end{minipage}  
\end{figure}  
As we will see later on, these maps provide a perfect model for the
mechanism which is responsible for the exponential evolution of peaks
and the creation of SNA's in saddle-node bifurcations. The map $x
\mapsto \frac{\arctan(\alpha x)}{\arctan(\alpha)}$ has three fixed
points at $0$ and $\pm 1$, and for $\beta=0$ these correspond to three
(constant) invariant curves for (\ref{eq:atanfamily}). As the
parameter $\beta$ is increased, a saddle-node bifurcation between the
two upper invariant curves takes place: Only the lower of the three
curves persists, while the other two collide and cancel each other
out. In fact, it will not be very hard to describe this bifurcation
pattern in general (see Theorem~\ref{thm:saddlenode}), whereas proving
that this bifurcation is indeed \textit{`non-smooth'} will require a
substantial amount of work.

Figure~\ref{fig:atanfamily} shows the behaviour of the upper bounding
graph as the parameter $\beta$ is increased and reveals a very
characteristic pattern. The overall shape of the curves hardly
changes, apart from the fact that when the bifurcation is approached
they have more and more \textit{`peaks'} (as we will see there are
infinitely many in the end, but most of them are too small to be
seen). The point is that these peaks do not appear arbitrarily, but
one after each other in a very ordered way: In (a), only the first
peak is fully developed while the second just starts to appear. In (b)
the second peak has grown out and a third one is just visible, in (c)
and (d) the third one grows out and a fourth and fifth start to appear \ldots\ .
Further, each peak is exactly the image of the preceding one, and the
peaks become steeper and thinner at an exponential rate (which
explains the term \textit{`exponential evolution'} and the fact that
the peaks soon become too thin to be detected numerically).

As far as simulations are concerned, the pictures obtained with smooth
forcing functions in (\ref{eq:atanfamily}) instead of
$(1-\sin(\pi\theta))$, which is only Lipschitz-continuous and decays
linearly off its maximum at $\theta=0$, show exactly the same
behaviour. However, the rigorous results from the later sections only
apply to this later type of forcing. In Section~\ref{Mechanism} we
will discuss why this simplifies the proof of the non-smoothness of
the bifurcation to some extent.

Finally, it should mentioned that the phenomenom we just described does not at
all depend on any particular properties of the arcus tangent. Any strictly
monotone and bounded map of the real line with the same qualitative features,
which can vaguely described as being ``s-shaped'', can be used to replace the
arcus tangent in the above definitions without changing the observed behaviour
(e.g.\ $x \mapsto \tanh(x)$). If in addition this map has similar scaling
properties as the arcus tangent, as for example $x \mapsto \frac{x}{1+|x|}$,
then even the rigorous results we present in the later sections apply. We will
not prove this in detail, but it will be evident that the arguments which
we use in Section~\ref{App:atanfamily} to treat (\ref{eq:atanfamily}) can
be easily adjusted to this end. 

 
\subsubsection{The Harper map} \label{Harper}

The Harper map with continuous \emph{potential} $V:\kreis \ra \R$,
\emph{energy} $E$ and \emph{coupling constant} $\lambda$ is given by
\begin{equation}
  \label{eq:harper}
\thx \mapsto \left(\thom,\arctan\left(\frac{-1}{\tan(x)-E+\lambda V(\theta)}\right)\right) \ .
\end{equation}
It is probably the most studied example, and the reason for this is the fact that
its dynamics are intimately related to the spectral properties of discrete
Schr\"odinger operators with quasiperiodic potential (the so-called
almost-Mathieu operator in the case $V(\theta) = \cos(2\pi\theta)$). Before we
turn to the simulations, we briefly want to discuss this relation and the
arguments by which the existence of SNA in the Harper map is established in
\cite{herman:1983}. A more detailed discussion can be found in
\cite{haro/puig:2006}.

The map (\ref{eq:harper}) describes the projective action of the
$\textrm{SL}(2,\R)$-cocycle (or Schr\"odinger cocycle)
\begin{equation} \label{eq:cocycle} 
  (\theta,v) \mapsto (\thom,A_{\lambda,E}(\theta) \cdot v) \ , 
\end{equation}
where $$A_{\lambda,E}(\theta) \ = \ \left(\begin{array}{cc} E - \lambda
  V(\theta) & -1 \\ 1 & 0 \end{array}\right)$$ and $v=(v_1,v_2) \in
  \R^2$. This means that (\ref{eq:harper}) can be derived from
  (\ref{eq:cocycle}) by letting $x := \arctan(v_2/v_1)$. The Schr\"odinger
  cocycle in (\ref{eq:cocycle}) is in turn associated to the almost-Mathieu
  operator
 \begin{equation}
   \label{eq:mathieu}
  H_{\lambda,\theta} : \ell^2 \ra \ell^2 \ , \ (H_{\lambda,\theta}u)_n = u_{n+1} +
  u_{n-1} + \lambda V(\theta+n\omega)u_n \ ,
 \end{equation}
as each formal solution of the eigenvalue equation $H_{\lambda,\theta}u=Eu$
 satisfies $$\left(\begin{array}{c}u_{n+1} \\ u_n \end{array}\right) =
 A_{\lambda,E}(\theta+n\omega) \cdot \left(\begin{array}{c}u_n \\ u_{n-1}
 \end{array}\right)\ .$$ 
The existence of SNA for the Harper map is equivalent to \emph{non-uniform
 hyperbolicity} of the cocycle (\ref{eq:cocycle})
 \cite{herman:1983,haro/puig:2006}, a concept which is of fundamental
 importance in this context.

In order to explain it, recall that a $\textrm{SL}(2,\R)$-cocycle over an
irrational rotation is a mapping $\kreis \times \R^2 \ra \kreis
\times \R^2$ of the form $(\theta,v) \mapsto (\thom,A(\theta) \cdot v)$, where
$A : \kreis \ra \textrm{SL}(2,\R)$ is a continuous function. The Lyapunov
exponent of such a cocycle is given by
\begin{equation}
  \label{eq:cocycle-lyap}
\lambda(\omega,A) \ = \ \nLim \ntel \int_{\kreis} \log \| A_n(\theta)\| \ d\theta \ , 
\end{equation}
where $A_n(\theta) = A(\theta+(n-1)\omega) \circ \ldots \circ A(\theta)$. If
$\lambda(\omega,A)>0$, then Oseledets Multiplicative Ergodic Theorem implies
the existence of an invariant splitting $\R^2 = \mathbb{W}^s_\theta \oplus
\mathbb{W}^u_\theta$ (invariance meaning $A(\theta)(\mathbb{W}^i_\theta) =
\mathbb{W}^i_{\thom} \ (i=s,u)$), such that vectors in $\mathbb{W}^u_\theta$
are exponentially expanded and vectors in $\mathbb{W}^s_\theta$ are
exponentially contracted with rate $\lambda(\omega,A)$ by the action of
$A_n(\theta)$.  The cocycle $(\omega,A)$ is called {\em uniformly hyperbolic}
if the subspaces $\mathbb{W}^i_\theta$ depend continuously on $\theta$. If
they depend only measurably on $\theta$, but not continuously, then the
cocycle is called {\em non-uniformly hyperbolic}.

In order to see why the latter notion is equivalent to the existence of SNA,
note that the invariant subspaces can be written as
$$\mathbb{W}^i_\theta \ = \ \R \cdot\left(\begin{array}{c} 1\\ \tilde{\varphi}^i(\theta)
\end{array} \right) $$
with measurable functions $\tilde{\varphi}^i : \kreis \ra \R \cup \{\infty\}$, and
it follows immediately that by $\varphi^i:= \arctan(\tilde{\varphi^i})$ we can
define invariant graphs for the projective action of the cocycle (obtained by
letting $x=\arctan(v_2/v_1)$ as above). Moreover, it is not difficult to show
that $\lambda(\varphi^s) = 2\lambda(\omega,A)$ and $\lambda(\varphi^u) =
-2\lambda(\omega,A)$ in this case,%
\footnote{\label{foot:harper1} In the case of the Harper map, the crucial
  computation is the following: Fix $\theta \in \kreis$ and $v \in \R^2
  \setminus \{0\}$ and define vectors $v^n$ by $v^0 := v$ and
  $v^{n+1}:=A(\theta+n\omega)\cdot v^n$. Further, let $\theta_n :=
  \theta+n\omega \bmod 1$ and $x_n := \arctan(v^n_2/v^n_1)$, and denote the
  Harper map (\ref{eq:harper}) by $T$. Then,
  \begin{eqnarray*}
    DT_{\theta_k}(x_k) & = & \frac{1}{1+(\tan(x_k)-E+\lambda V(\theta_k))^{-2}}
    \cdot \frac{1+\tan(x_k)^2}{(\tan(x_k)-E+\lambda V(\theta_k))^2}  \\
    & = & \frac{1+\tan(x_k)^2}{1+(\tan(x_k)-E+\lambda V(\theta_k))^2} \ = \
    \frac{1+\tan(x_k)^2}{1+\tan(x_{k+1})^{-2}}  \ = \ \frac{\|
    v^k\|^2}{\|v^{k+1}\|^2} \ ,
  \end{eqnarray*}
  where we used $v^{k+1}_2 = v^k_1$ in the last step. Consequently, we obtain 
  \begin{eqnarray*}
    DT^n_\theta(x_0) \ = \ \prod_{k=0}^{n-1} DT_{\theta_k}(x_k) \ = \
    \frac{\|v^0\|^2}{\|v^n\|^2} \ ,
  \end{eqnarray*}
  and this establishes the asserted relation between the different Lyapunov
 exponents. The case of a general $\textrm{SL}(2,\R)$-cocycle can be treated
 in more or less the same way. }
 and conversely the existence of invariant
 graphs with non-zero Lyapunov exponent implies the existence of an invariant
 splitting with the mentioned properties. As the graphs $\varphi^i$ depend
 continuously on $\theta$ if and only if this is true for the subspaces
 $\mathbb{W}^i_\theta$, we obtain the claimed equivalence. \medskip

The crucial observation which was made by Herman is the fact that, using a
result from sub-harmonic analysis, lower bounds on the Lyapunov exponent can
be obtained for suitable choices of the $\textrm{SL}(2,\R)$-valued function
$A$. In the case of (\ref{eq:cocycle}) with potential $V(\theta) =
\cos(2\pi\theta)$, this bound is $\lambda(\omega,A_{\lambda,E}) \geq
\max\{0,\log(|\lambda|/2)\}$ \cite[Section 4.7]{herman:1983}. Consequently, if
$|\lambda| >2$ then the Lyapunov exponent of $(\omega,A_{\lambda,E})$ will be
strictly positive for all values of $E$.  On the other hand, it is well-known
that there cannot be a continuous splitting for all $E \in \R$, and
consequently for some $E$ the respective cocycle has to be non-uniformly
hyperbolic.

The simplest way to see this is probably to consider the rotation
number. Suppose $\omega \in \kreis \setminus \Q$ and $\lambda>2$ are
fixed. Then (\ref{eq:harper}) defines a skew-product map $T_E$ on the
two-torus, and for such maps a fibred rotation number $\rho(T_E)$ can be
defined, much in the way this is done for homeomorphisms of the circle. The
dependence of $\rho(T_E)$ on $E$ is continuous \cite[Section 5]{herman:1983},
and further it is easy to see that the existence of continuous invariant
graphs forces the fibred rotation number to be rationally related to $\omega$,
more precisely to take values in the module ${\cal M}_\omega := \{\frac{k}{q} \omega
\bmod 1 \mid k\in\Z, q\in \N\}$ (compare \cite[Section
5.17]{herman:1983}). However, if $E$ runs through the real line from $-\infty$
to $\infty$, then the rotation number $\rho(T)$ runs exactly once around the
circle \cite[Section 4.17(b)]{herman:1983}.  For all $E \in \R$ with
$\rho(T_E) \notin {\cal M}_E$, the existence of a SNA in (\ref{eq:harper})
follows. Refined results can be obtained by using the fact that the invariant
splitting cannot be continuous whenever $E$ belongs to the spectrum of the
almost-Mathieu operator. This is discussed in detail in
\cite{haro/puig:2006}. In particular, it allows to use lower bounds on the
measure of the spectrum to establish the existence of SNA for a set of
positive measure in parameter space. \medskip

For the simulations presented here we use a reflection
of (\ref{eq:harper}) w.r.t.\ the $\theta$-axis,
\begin{equation} \label{eq:projharper} 
  \thx \ \mapsto \ 
  \left(\thom,\arctan\left(\frac{1}{\tan(-x)- 
        E+\lambda V(\theta)}\right) \right)    
  \ , 
\end{equation} 
as this makes it easier to compare the pictures with the other examples. The
potential function which is used is $V(\theta) = \cos(2\pi\theta)$. Later, in
Section~\ref{App:Harper}, we have to make a different choice in order to
obtain rigorous results with the methods presented here.
\begin{figure}[h!] 
\noindent 
\begin{minipage}[t]{\linewidth} 
  \epsfig{file=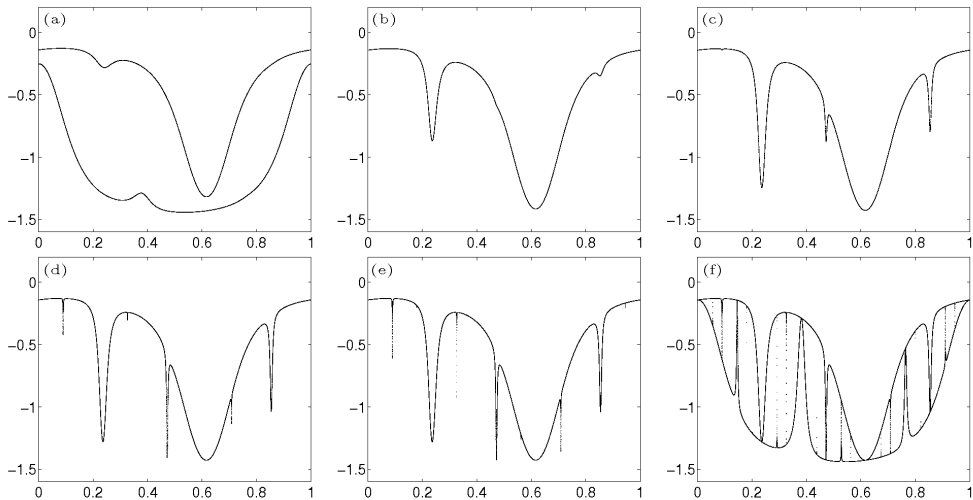, clip=, width=\linewidth, 
    }
         \caption{\small The stable invariant curves for the projected Harper
           map given by (\ref{eq:projharper}) with $\omega$ the golden
           mean, $\lambda=4$ and different values for $E$. (a)
           At $E=4.4$ the first peak is clearly visible, while
           the second just starts to appear. The repeller is close,
           but still a certain distance away. (b) At $E=4.3$
           the second peak has grown and the third starts to appear.
           This pattern continues, and more and more peaks can be seen
           in pictures (c) $E=4.289$, (d) $E=4.28822$
           and (e) $E=4.288208$. (f) finally shows attractor
           and repeller for $E=4.288207478$ just prior to
           collision.  }
         \label{fig:harpergraphs}
\end{minipage} 
\end{figure} 

As described by Herman in \cite[Section 4.14]{herman:1983}, when the parameter
$E$ approaches the spectrum of the almost-Mathieu operator from above, a
stable and an unstable invariant circle collide in a saddle-node
bifurcation. Even if the rigorous arguments used by Herman \cite{herman:1983}
are very specific for cocycles (as described above), the process seems to be
the same as in the arctan-family before: Figure~\ref{fig:harpergraphs} shows
the behaviour of the attractor before it collides with the repeller (the latter
is only depicted in Fig.~\ref{fig:harpergraphs}(a) and (f)). The pattern is
already familiar, the exponential evolution of peaks can be seen quite clearly
again.

Based on this observation, Bjerkl\"ov recently addressed a problem raised by
Herman \cite[Section 4.14]{herman:1983} about the structure of the minimal set
which is created in this bifurcation. Upon their collision, the stable and
unstable invariant circles are replaced by an upper, respectively lower
semi-continuous invariant graph. The region between the two graphs is a
compact and invariant set, but it is not at all obvious whether this set is
also minimal and coincides with the topological closures of the two graphs. In
\cite{bjerkloev:2005} Bjerkl\"ov gives a positive answer to this question,
provided the rotation number $\omega$ on the base is Diophantine and the
parameter $\lambda$ is sufficiently large. As his approach is purely dynamical
and does not depend on any particular properties of cocycles, it should be
possible to apply it to more general systems. This might allow to prove the
existence of SNA and to describe their structure, in the above sense, at the
same time.


\subsubsection{The quasiperiodically forced Arnold circle map} \label{ArnoldMap}

The most obvious physical motivation for studying qpf systems are probably
oscillators which are forced at two or more incommensurate frequencies.
\begin{figure}[h!]
\noindent
\begin{minipage}[t]{\linewidth}
  \epsfig{file=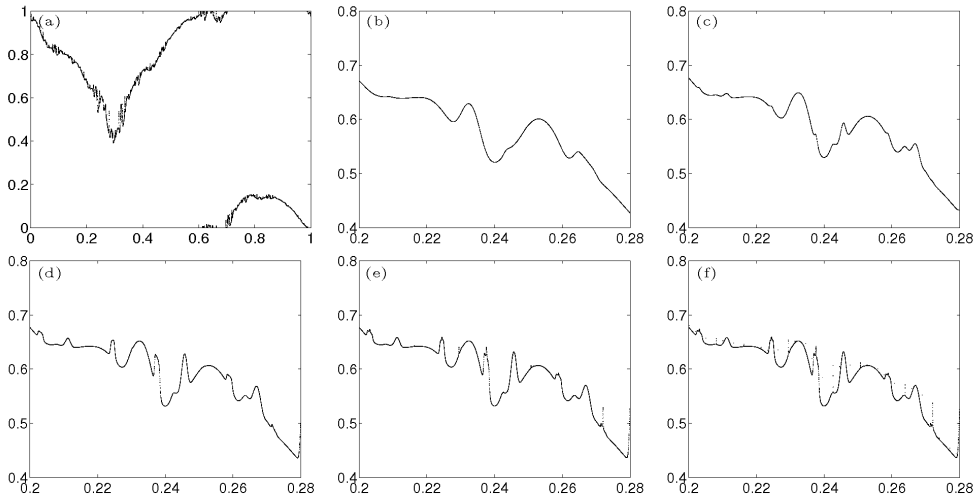, clip=, width=\linewidth,
    }
         \caption{\small Pictures obtained from the qpf Arnold circle map
           (\ref{eq:qpfarnold}) with $\alpha=0.99$ and $\beta=0.6$, 
           $\omega$ is the golden mean.
           (a) shows the attracting invariant curve for
           $\tau=0.3373547962$. As the exponential evolution of peaks
           takes place on a rather microscopic level, it is difficult
           to recognise any details. Therefore, the other pictures
           show the attractors only over the interval $[0.2,0.28]$.
           The $\tau$-values are (b) $\tau=0.337$, (c)
           $\tau=0.3373$, (d) $\tau=0.3373547$, (e)
           $\tau=0.33735479$, (f) $\tau=0.337357962$.}
         \label{fig:arnold}
\end{minipage}
\end{figure}
If these are modelled by differential equations, the Poincar\'e maps will be
of the form (\ref{eq:qpfs}). The qpf Arnold circle map, given by
\begin{equation}
  \label{eq:qpfarnold} \thx \mapsto
  \left(\thom,x+\tau+\frac{\alpha}{2\pi}\sin(2\pi x) 
  + \beta \sin(2\pi\theta) \right) 
\end{equation}
with real parameters $\alpha,\tau$ and $\beta$, is often studied
as a basic example (see \cite{ding/grebogi/ott:1989}). There are
several interesting phenomena which can be found in this family, such
as different bifurcation patterns, mode-locking or the transition to
chaos as the map becomes non-invertible \cite{ding/grebogi/ott:1989,feudel/kurths/pikovsky:1995,%
osinga/wiersig/glendinning/feudel:2001}.  Similar to the
unforced Arnold circle map
\cite{katok/hasselblatt:1995,demelo/vanstrien:1993}, there exist
so-called Arnold tongues -- regions in the parameter space on which
the rotation number stays constant. The reason for this is usually the
existence of (at least) one stable invariant circle inside of the
tongue. On the boundaries of the tongue this attractor collides with
an unstable invariant circle in a saddle-node bifurcation (see
\cite{osinga/wiersig/glendinning/feudel:2001,chastell/glendinning/stark:1995}
or \cite{glendinning/feudel/pikovsky/stark:2000} for a more detailed
discussion and numerical results).

\begin{figure}[h!] 
\noindent 
\begin{minipage}[t]{\linewidth} 
  \epsfig{file=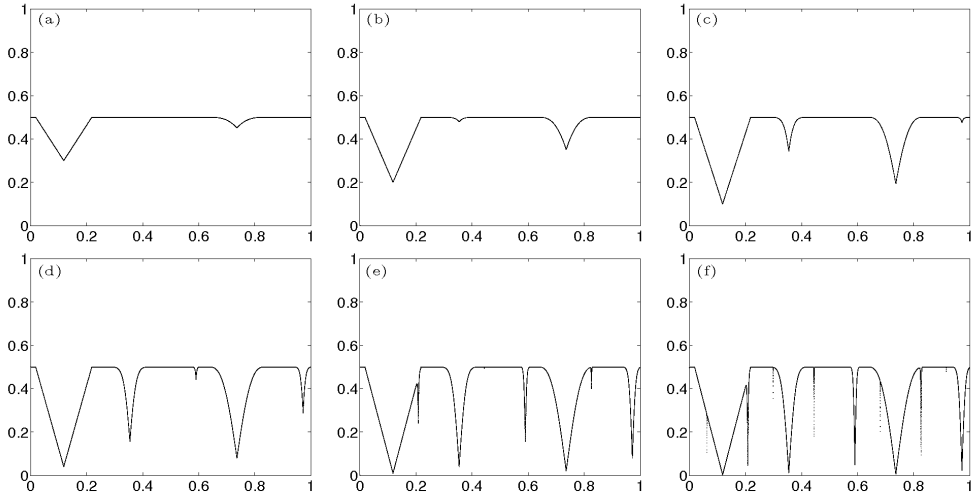, clip=, width=\linewidth, }
         \caption{\small The stable invariant curves in the system $\thx
           \mapsto
           \left(\thom,x+\tau+\frac{\alpha}{2\pi}\sin(2\pi\theta)\right.$
           $\left.-\beta \cdot \max\{0,1-10\cdot
           d(\theta,\halb)\}\right)$. This time the parameters $\alpha=0.99$
           and $\tau=0$ are fixed, while $\beta$ varies: (a) $\beta=-0.2$, (b)
           $\beta=-0.3$, (c) $\beta=-0.4$, (d) $\beta=-0.45$, (e)
           $\beta=-0.49$, (f) $\beta=-0.497$. Again, $\omega$ is the golden mean. The
           exponential evolution of peaks is clearly visible. }
         \label{fig:arnoldpeak}
\end{minipage} 
\end{figure} 

For our purpose it is convenient to study only those bifurcations which take
place on the boundary of the Arnold tongue with rotation number zero. In order
to do so, we fix the parameters $\alpha \in [0,1]$ and $\beta>0$, thus
obtaining a one-parameter family depending on $\tau$. As long as $\beta$ is
not too large, there exist a stable and an unstable invariant curve at
$\tau=0$.  Increasing or decreasing $\tau$ leads to the disappearance of the
two curves after their collision in a saddle-node bifurcation.  When $\alpha$
is close enough to 1 (where the map becomes non-invertible) this bifurcation
seems to be non-smooth \cite{osinga/wiersig/glendinning/feudel:2001}. The
problem here is the fact that the curves are already extremely `wrinkled'
before the exponential evolution of peaks really starts. Therefore, it is hard
to recognise any details in the global picture (see
Figure~\ref{fig:arnold}(a)).  This becomes different if we `zoom in' and only
look at the curves over a small interval. On this microscopic level, we
discover the more or less the same behaviour as before
(Figure~\ref{fig:arnold}(b)--(f)). Of course, this time we can not really
determine the order in which the peaks are generated, as we only see those
peaks which lie in our small interval. But we clearly see that more and more
peaks appear, and those appearing at a later time are smaller and steeper than
those before.

On the other hand, we can also use a more `\textit{peak-shaped}' forcing
function instead of the sine. In this case, the pictures we obtain look
exactly the same as the ones from the arctan-family above (see
Figure~\ref{fig:arnoldpeak}(a)-(f)). This effect will be discussed in more
detail in Section~\ref{Mechanism}~. Nevertheless, we should mention that, in
contrast to the two preceding examples, we do not provide any rigorous results
on the qpf Arnold circle map (see also Section~\ref{App:Remarks} for a
discussion).


\subsubsection{Pinched skew products} \label{PinchedSystems}

As for the Harper map, we refer to the original literature
\cite{grebogi/ott/pelikan/yorke:1984,keller:1996} for a more detailed
discussion of these systems.  Here, we will just have a look at the map
\begin{equation}
  \label{eq:pinchedskews} 
 \thx \mapsto (\thom,\tanh(\alpha x)\cdot \sin(\pi\theta)) \ , 
\end{equation} 
with real positive parameter $\alpha$, which is a typical representative of
this class of systems. Note that due to the multiplicative nature of the
forcing, the 0-line is \textit{a priori} invariant, and due to the zero of the
sine function there is one fibre which is mapped to a single point (hence
`pinched'). These are the essential features that are needed to prove the
existence of SNA in pinched skew products (see
\cite{keller:1996,glendinning:2002}).

\begin{figure}[h]  
\noindent 
\begin{minipage}[b]{\linewidth}  
  \epsfig{file=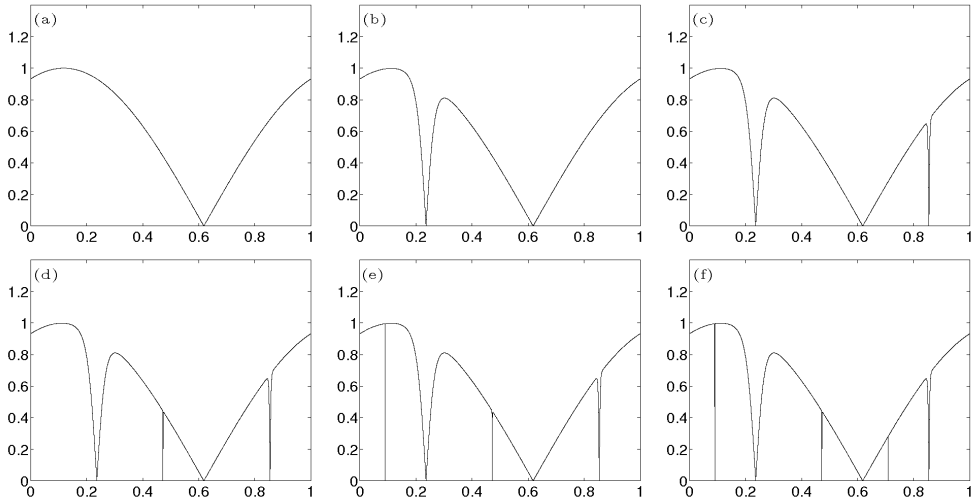, clip=, width=\linewidth, }
         \caption{\small The first six iterates of the upper boundary
           line for the pinched skew product given by
           (\ref{eq:pinchedskews}) with $\omega$ the golden mean and
           $\alpha=10$. In each step of the iteration one more peak
           appears, while apart from that the curves seem to stay the
           same. Further, the peaks become steeper and thinner at an
           exponential rate.}
         \label{fig:pinchedlines} 
\end{minipage}  
\end{figure} 

Figure~\ref{fig:pinchedlines} differs from the preceding ones insofar as it
does not show a sequence of invariant graphs as the systems parameters are
varied, but the first images of a constant line that is iterated with a fixed
map. Nevertheless, the behaviour is very much the same as before. The
exponential evolution of peaks can followed even easier here, as this time
each iterate produces exactly one further peak.

For Pinched skew products this process was quantified \cite{jaeger:2004a} in
order to describe the structure of the SNA's in more detail. The question
addressed there is basically the same as the one studied by Bjerkl\"ov in
\cite{bjerkloev:2005}, and the result is similar: The SNA, which is an upper
semi-continuous invariant graph above the 0-line in this situation, lies dense
in the region below itself and above the 0-line, provided the rotation number
$\omega$ on the base is Diophantine and the parameter $\alpha$ is large
enough.


\subsubsection{Non-smooth pitchfork bifurcations} \label{Pitchfork}

\begin{figure}[h!]
\noindent
\begin{minipage}[t]{\linewidth}
  \epsfig{file=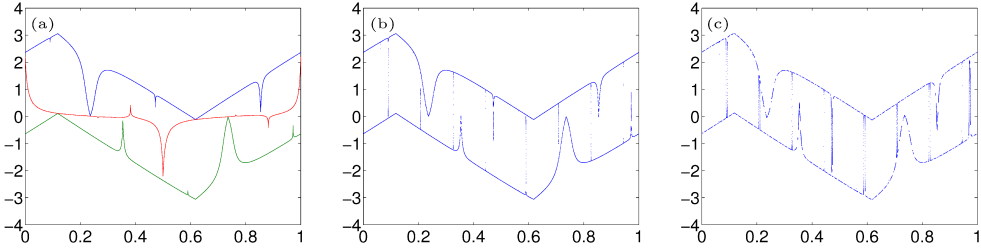, clip=, width=\linewidth,
    }
         \caption{\small A pitchfork bifurcation in the parameter
           family (\ref{eq:symmetricfamily}). (a) shows the upper and
           lower bounding graphs just prior to the collision. Note
           that here the two objects are still distinct, and three
           different trajectories (a backwards trajectory for the 
           repeller) are plotted to produce this picture.
           In contrast to this, (b) and (c) only show one single
           trajectory. There still exist two distinct SNA's, but these
           are interwoven in such a way that they cannot be
           distinguished anymore. Each of them seems to trace out a
           picture of both attractors before collision. The parameter
           values are $\alpha=10$ and (a) $\tau=1.64$, (b)
           $\tau=1.645$ and (c) $\tau=1.66$. $\omega$ is the golden mean.}
         \label{fig:pitchfork}
\end{minipage}
\end{figure}

Compared to saddle-nodes, pitchfork bifurcations are degenerate.
Usually they only occur if the system has some inherent symmetry that
forces three invariant circles to collide exactly at the same time.
Nevertheless, they have been described in the literature about SNA's
quite often (e.g.\ 
\cite{osinga/wiersig/glendinning/feudel:2001},\cite{sturman:1999}).
The reason for this is the fact that unlike in saddle-node
bifurcations, where the SNA's only occur at one single parameter, SNA's
which are created in pitchfork bifurcations seem to persist over a
small parameter interval. In addition, the
transition from continuous to non-continuous invariant graphs at the
collision point is much more distinct, as the SNA which is created
seems to trace out a picture of both stable invariant curves just
prior to the collision (see Figure~\ref{fig:pitchfork}).

We were not able to give a rigorous proof for this stabilising effect, or any
other details of a non-smooth pitchfork bifurcation.  However, by a slight
modification of the methods used for the non-smooth saddle-node bifurcation,
we can at least prove the existence of SNA's in systems with the mentioned
inherent symmetry (see Theorem~\ref{thm:symmetricsna} and
Section~\ref{SNAwithsymmetry}). For suitable parameters these systems have two
SNA's which are symmetric to each other and enclose a self-symmetric SNR, and
the three objects are interwoven in such a way that they all have the same
(essential) topological closure. As an example, we consider the parameter
family
\begin{equation} \label{eq:symmetricfamily}
\thx \mapsto \arctan(\alpha x) - \beta \cdot (1-4d(\theta,0)) \ .
\end{equation}
For Diophantine $\omega$ and sufficiently large $\alpha$ we will
obtain the existence of a SNA-SNR triple as described above for at least one
suitable parameter $\beta(\alpha)$.


\subsection{The mechanism: Exponential evolution of peaks} \label{Mechanism}
In the following, we will try to give a simple heuristic explanation for the
mechanism which is responsible for the exponential evolution of peaks.
Generally, one could say that it consists of a subtle interplay of an
\textit{`expanding region'} $\cal E$ and a \textit{`contracting region'} $\cal C$, which
communicate with each other only via a small \textit{`critical region'}
$\cal S$. In order to give meaning to this, we concentrate first on the
$\arctan$-family given by (\ref{eq:atanfamily}). \medskip

\begin{figure}[h!]
\noindent
\begin{minipage}[t]{\linewidth}
  \epsfig{file=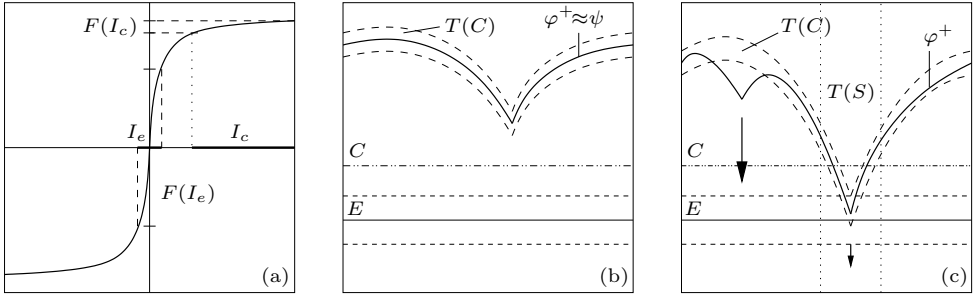, clip=, width=\linewidth,
    }
         \caption{\small As the fibre maps are expanding and contracting on $I_e$ and $I_c$,
           respectively, $T$ will be expanding in the vertical direction on
           ${\cal E} = \kreis \times I_e$ and contracting on ${\cal C} =
           \kreis \times I_c$. (b) As long as $\beta$ is not too large, $\cal
           C$ is mapped into itself. Thus, there exists a stable invariant
           circle inside of $T({\cal C})$ (in fact, as a point set this circle
           coincides with $\bigcap_{n\in\N} T^n({\cal C})$), which has approximately
           the shape of the forcing function. (c) When the first peak enters
           the expanding region it induces a second peak, which moves faster
           than the first one by the expansion factor in $\cal E$. The first
           peak is generated in the critical region $\cal S$, where the
           forcing achieves its maximum. Therefore, it is located in $T({\cal S})$. }
         \label{fig:mechanism}
\end{minipage}
\end{figure}

If we restrict to $\alpha \geq \tan(1)$ and $\beta \leq \pi$ in
(\ref{eq:atanfamily}), then we can choose $X=\left[-\frac{3}{2}\pi,\frac{3}{2}\pi\right]$
as the driven space, because in this case $\kreis \times
\left[-\frac{3}{2}\pi,\frac{3}{2}\pi\right]$ is always mapped into itself. Further, we
fix $\alpha$ sufficiently large, such that the map $F : x \mapsto
\arctan(\alpha x)$ has three fixed points $x^- < 0 < x^+$. As $0$ will be
repelling and $x^+$ attracting, we can choose a small interval $I_e$ around 0
which is expanded and an interval $I_c$ around $x^+$ which is contracted, and
define the expanding and contraction regions as ${\cal E}:= \kreis \times I_e$
and ${\cal C}:= \kreis \times I_c$ (see Figure~\ref{fig:mechanism}).  Of
course, there exists a second contracting region ${\cal C}^-$, corresponding
to $x^-$, but this does not take part in the bifurcation: Due to the one-sided
nature of the forcing, ${\cal C}^-$ is always a trapping region, independent
of the parameter $\beta$. Thus there always exists a stable invariant circle
inside of ${\cal C}^-$, and the saddle-node bifurcation only takes place
between the two invariant circles above.

By the choice of the intervals, the fibre maps $\Tth$ are contracting
on $I_c$ and expanding on $I_e$. Further, as long as $\beta$ is small
there holds
\begin{equation} \label{eq:expandinginterval}
    \Tth(I_c) \ssq I_c \ \ \textrm{ and } \ \ I_e \subseteq \Tth(I_e)
    \ 
\end{equation}  
for all $\theta \in \kreis$. Consequently, 
\begin{equation} \label{eq:expandingregion} 
    T({\cal C}) \ssq {\cal C} \ \ \textrm{ and } \ \ {\cal E} \subseteq
    T({\cal E}) \ . 
\end{equation}   
This means that $\cal C$ and $\cal E$ cannot interact, and there will be
exactly one invariant circle (stable and unstable, respectively) in each of
the two regions.  However, when $\beta$ is increased and approaches the
bifurcation point, (\ref{eq:expandingregion}) does not hold anymore.
Nevertheless, the relation (\ref{eq:expandinginterval}) will still be true for
`most' $\theta$, namely whenever the forcing function $(1-\sin(\pi\theta))$ in
(\ref{eq:atanfamily}) is not close to its maximum (see
Figure~\ref{fig:mechanism}(c)). Thus, even when $\cal E$ and $\cal C$ start to
interact, they will only do so in a vertical strip ${\cal S}:=W \times X$,
where $W \ssq \kreis$ is a small interval around $0$.

 This strip $\cal S$ is
the `\textit{critical region}' we referred to above and in which the first
peak is generated: As long as $T(\cal C) \ssq C$, the upper bounding graph
will be contained in $T(\cal C)$.  But this set is just a very small strip
around the first iterate of the line $\kreis \times \{x^+\}$, which is a curve
$\psi$ given by
\begin{equation} \label{eq:firstpeak}
    \psi(\theta) := x^+ - \beta\cdot (1-\sin(\pi(\theta-\omega))) \
\end{equation}
(see Figure~\ref{fig:mechanism}(b)). Consequently, the upper bounding
graph $\varphi^+$ will have approximately the same shape as $\psi$,
which means that it has a first peak centred around $\omega$, i.e.\ 
in $T(\cal S)$.  From that point on, the further behaviour is explained
quite easily. As soon as the first peak enters the expanding region,
its movement will be amplified due to the strong expansion in $\cal E$.
Thus a second peak will be generated at $2\omega \bmod 1$. It will be
steeper than the first one, and when $\beta$ is increased it also
grows faster by a factor which is more or less the expansion factor
inside $\cal E$. As soon as the second peak is large enough to enter the
expanding region, it generates a third one, which in turn induces a
fourth and so on $\ldots$ . \medskip

The picture we have drawn so far already gives a first idea about what
happens, although converting it into a rigorous proof for the existence of SNA
will still require a substantial amount of work. As we will see, it is not too
hard to give a good quantitative description of the behaviour of the peaks up
to a certain point, namely as long as the peaks do enter the critical region
(corresponding to the returns of the underlying rotation to the interval
$W$). But as soon as this happens, things will start to become
difficult. However, by assuming that the rotation number $\omega$ satisfies a
Diophantine condition we can ensure that such returns are not too frequent,
and that very close returns do not happen too soon. This will be sufficient to
ensure that the exponential evolution of peaks also carries on afterwards.
\medskip

\begin{figure}[h!]
\noindent
\begin{minipage}[t]{\linewidth}
  \epsfig{file=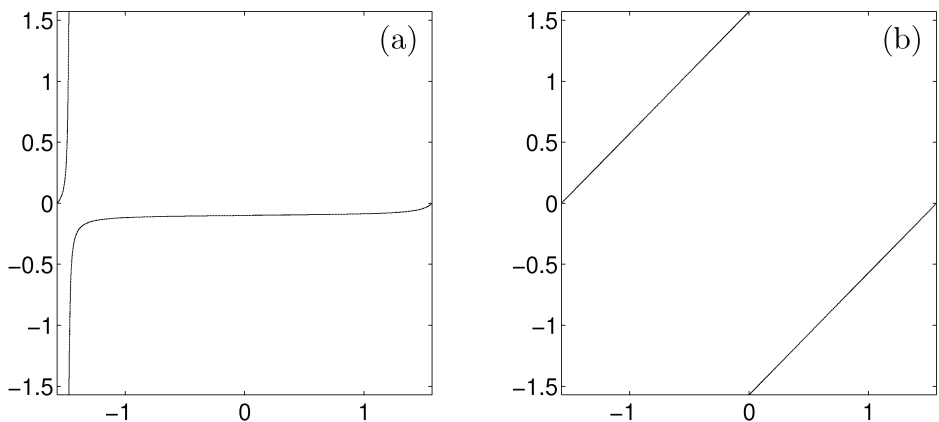, clip=, width=\linewidth,
    }
         \caption{\small Graphs of the projected
           M\"obius-transformations $x \mapsto
           \arctan\left(\frac{1}{\tan(-x)-10}\right)$ in (a) and $x
           \mapsto \arctan\left(\frac{1}{\tan(-x)}\right)$. }
         \label{fig:moebius}
\end{minipage}
\end{figure}

In principle, the mechanism is not different in the other parameter families
discussed in the last section. For the Harper map given by
(\ref{eq:projharper}), Figure~\ref{fig:moebius}(a) shows the graph of a
projected M\"obius-transformation $x \mapsto \arctan(\frac{1}{\tan(-x) - c})$
for large $c$. As long as $E \gg \lambda$, the fibre maps will all have
approximately this shape.
As we can see, there will be a repelling fixed point slightly above
$-\frac{\pi}{2}$ and an attracting one slightly below 0. This means
that if we choose $I_e$ and $I_c$ to be sufficiently small intervals around
these fixed points, then we have uniform expansion on $\cal E$, uniform
contraction on $\cal C$ and (\ref{eq:expandinginterval}) will be satisfied.
When $E \approx \lambda$, this will still be true on most
fibres. Only where the potential $\cos(2\pi\theta)$ is close to its
maximum at $\theta=0$, the picture changes
(Figure~\ref{fig:moebius}(b)). Here $-\frac{\pi}{2} \in I_c$ is mapped
close to 0, which means again that the expanding and contracting
region start to interact and a first peak is produced. (Thus, the
critical region $\cal S$ is again a vertical strip around 0.) As before,
this peak is amplified as soon as it enters the expanding region $\cal E$
and thus induces all others.

In some sense, the situation for the qpf Arnold circle map is even more
similar to the case of the arctan-family, as the forcing is additive again and
the fibre maps are clearly $s$-shaped as before. However, the difference is
the fact that while the derivative at the stable fixed point indeed vanishes,
such that the contraction becomes arbitrarily strong, the maximal expansion
factor is at most 2 (at least in the realm of invertibility $\alpha \leq 1$).
This explains why the resulting pictures in Figure~\ref{fig:arnold} are much
less clear. Roughly speaking, in combination with the limited expansion the
peak of the forcing function $\theta \mapsto \sin(\pi\theta)$ is just `too
blunt' to trigger the exponential evolution of peaks as easily as before. When
it finally does take place - as the simulations in Figure~\ref{fig:arnold}
suggest - the graphs are already too `wrinkled' to give a good picture. But of
course, if the shape of the forcing function is a second factor that decides
whether the exponential evolution of peaks takes place, then we can also
trigger this pattern by choosing one with a very sharp peak. This is exactly
what happened in Figure~\ref{fig:arnoldpeak}.

Finally, for Pinched skew products we refer to \cite{jaeger:2004a} for a more
detailed discussion.

\begin{figure}[h!] 
\noindent 
\begin{minipage}[t]{\linewidth} 
  \epsfig{file=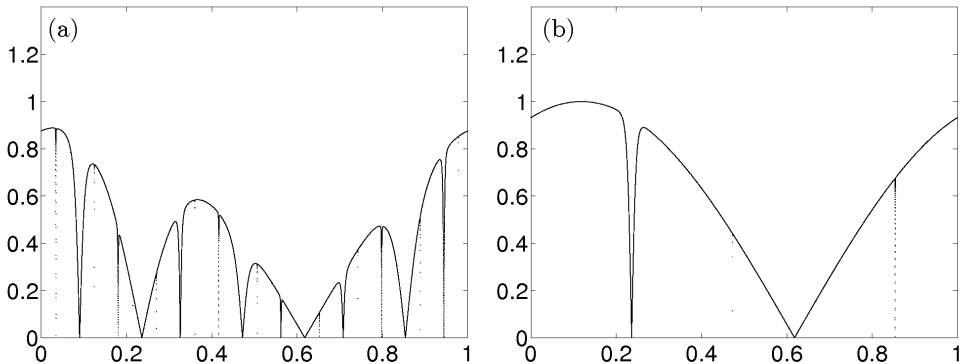, clip=, width=\linewidth, 
    }
         \caption{\small Upper bounding graphs in the pinched systems given by
           (\ref{eq:pinchedskews}). $\omega$ is the golden mean, the parameter
           values are (a) $\alpha=3$ and (b) $\alpha=32$. In (b), where the
           expansion is stronger, there seems to be less structure in
           comparison to (a). However, this is not a qualitative difference,
           but can be easily explained by the exponential evolution of peaks.
           If the expansion is stronger, the peaks of higher order are just
           not visible anymore, such that only the first few peaks can be
           seen.}
         \label{fig:structure}
\end{minipage} 
\end{figure} 

\begin{bem}
The preceding discussion gives a basic understanding of how SNA's are created
in the above examples. Although it might be very rudimentary, it can already
be used to anticipate a number of observations. Without trying to make things
very precise, we want to mention a few:
\alphlist
\item First, it is not hard
to guess in which parameter range the expanding and contracting regions start
to interact and the torus collision takes place in the above families, e.g.\
$E \approx \beta$ for the Harper map or $\beta \approx \frac{\pi}{2}$ for the
arctan-family.
\item  Another phenomenom which can be explained is the following:
The stronger the expansion and contraction are, i.e.\ the larger the
respective parameter is chosen, the less `structure' can be seen in the
pictures (see Figure~\ref{fig:structure}). However, obviously this `structure'
corresponds exactly to the peaks which are generated.  These can only be
detected numerically as long as they do not become too small, but of course
this happens faster if the expansion and contraction are
stronger. Figure~\ref{fig:structure} shows this effect for pinched
systems, but it can be observed similarly in all the examples we
treated. In particular, it is also present in the qpf Arnold circle map
(\ref{eq:qpfarnold}), which indicates again that the mechanism there is not
different from the other examples. 
\item As already mentioned, the exponential evolution of peaks is easier to
trigger if the forcing function has a very distinctive and sharp peak. 
Figures~\ref{fig:arnold} and \ref{fig:arnoldpeak} illustrate this in the
context of the qpf Arnold circle map.
\item In \cite{haro/delallave:2006}, the authors study (amongst other things)
  the parameter-dependence of the minimal distance $\Delta_\beta$ between the
  stable and instable invariant curve in a non-smooth bifurcation. Their
  situation is slightly different to the one considered here, since the
  dynamics take place on a torus and the attractor touches the repeller from
  above and below at the same time. Nevertheless, the pictures indicate that a
  process similar to the one described above takes place. The observation
  which was made by the authors is that the asymptotic dependence of
  $\Delta_\beta$ on $|\beta-\beta_c|$ seems to be a power law with exponent 1
  as $\beta \ra \beta_c$, i.e.\ $\Delta_\beta \sim |\beta-\beta_c|$ (where
  $\beta_c$ is the bifurcation parameter). Furthermore, this exponent seems to
  be universal for a certain class of models.

 At least in the situations we discussed, e.g.\ for (\ref{eq:atanfamily}), the
  exponential evolution of peaks offers a reasonable explanation for such a
  scaling behaviour: Since all peaks of the attractor have to touch the
  repeller at the same time and, according to our heuristics, all
  further peaks move much faster than the first one, it is the latter which
  determines the minimal distance of the two curves. However, as this first
  peak has approximately the shape of the forcing function (see
  (\ref{eq:firstpeak})), the position of its tip depends linearly on $\beta$.
  \listend
\end{bem}
Admittedly, some of the above remarks remain rather speculative unless they
are confirmed either by careful numerical studies or rigorous
proofs. Nevertheless, what we want to point out is that the mechanism we
described offers at least an heuristic explanation for a number of
observations which have sometimes been found to be puzzling or ever
confusing. Further, an intuitive understanding of the process should make it
easier to come up with reasonable conjectures, which can then (in the better
case) either be proved or at least be confirmed numerically. As already
mentioned, the issue we want to concentrate on in this article is a rigorous
proof for the existence of SNA. \medskip


Concerning the latter, the main problem we will
encounter is that we do not a priori know where the tips of the peaks are
located. If there is any chance of rigorously describing the exponential
evolution of peaks in a quantitative way, they must be located in the
expanding region at least most of the times. Otherwise, there would be no
plausible mechanism which forces the peaks to become steeper and steeper. But
the horizontal position is not the only problem. When we use a forcing
function with a quadratic maximum, then we do not even know the exact vertical
position: If the tip of one peak is on the fibre $\theta$, then the tip of the
next will be close to $\theta+\omega$, but it may be slightly shifted due to
the influence of the forcing. In order to explain this, suppose that a upper
bounding graph $\varphi^+$ of $T$ is differentiable and has a local minimum at
$\theta_0$. The derivative of $\varphi$ at $\theta_0+\omega$ is then given by
\begin{equation} \label{eq:peakimage}
\varphi'(\theta_0+\omega) \ = \ \frac{\partial}{\partial_\theta} (\pi_2
\circ T)(\theta_0,\varphi(\theta_0)) + DF_\theta(\varphi(\theta)) \cdot
\varphi'(\theta_0) \ = \ -\beta g'(\theta_0) \ .
\end{equation}
(Here we suppose that $T$ has fibre maps of the form $T_\theta = F(x) - \beta
g(\theta)$ as in (\ref{eq:atanfamily}).) Consequently, if $g'(\theta_0) \neq
0$, then $\theta_0+\omega$ is not a local minimum. 

This becomes different if the local minima, which we call `peaks', are
sufficiently `sharp and steep'. By this, we mean that that both $\lim_{\theta \nearrow
  \theta_0} -\varphi'(\theta)$ and $\lim_{\theta \searrow \theta_0}
\varphi'(\theta)$ are greater than a sufficiently large constant $M$
(depending on the ${\cal C}^1$-norms of $F$ and $g$). Then it can easily be
seen from (\ref{eq:peakimage}) that $\theta_0+\omega$ will be a local minimum
as well. If in addition $(\theta_0,\varphi(\theta_0))$ is located in the
expanding region and the expansion constant is sufficiently large, then the
peak at $\theta_0+\omega$ will again be sufficiently sharp and steep (in the
above sense). 

Our claim is now that we can produce such sharp peaks by choosing a forcing
function that, like $1-\sin(\pi\theta)$, is only Lipschitz-continuous at its
maximum and decays linearly in a neighbourhood. At least for the first peak
this is plausible, since we have argued above that at the onset of the
exponential evolution of peaks the invariant graph has approximately the shape
of the forcing function (see the discussion around (\ref{eq:firstpeak})). For
all further peaks we can expect the same, provided that the exponential
evolution of peaks is really caused by the mechanism described above, because
then the tips of the peaks are located in the expanding region (at least most
of the time).

However, we will not give a rigorous proof for this claim, since this would
require to describe the global structure of the invariant graphs.  In fact, we
argue that it is exactly this `localisation' of the tips of the peaks which
helps to overcome the need for such a global description (which would probably
be much more complicated on a technical level). In order to understand this,
note that (in case our claim holds), the tips of the peaks just correspond to
a single orbit, since then one minimum is mapped to another. Further, as
mentioned, we expect that this orbit spends most of the time in the expanding
region, and in fact this will already turn out to be sufficient to prove the
existence of a SNA: In this case there exists an orbit on the upper bounding
graph which has a positive vertical Lyapunov exponent, and this is not
compatible with the continuity of the upper bounding graph (the Lyapunov
exponent of the upper bounding graph is always non-positive, e.g.\ Lemma~3.5
in \cite{jaeger:2003}, and due to uniform convergence of the ergodic limits
this is true for any of its points).

 However, during the proof we will obtain even more information about this
particular orbit: It does not only have a positive Lyapunov exponent forwards,
but also backwards in time. Thus, concerning it Lyapunov exponents the orbit
behaves as if it was moving from a sink to a source (and referring to this we
will call it a \textit{`sink-source-orbit'}). As it will turn out, it is
contained in the intersection of the SNA and the SNR. The existence of such
atypical orbits is also well-known for the Harper map, where it is equivalent
to the existence of exponentially decaying eigenfunctions for the associated
Schr\"odinger operators
\footnote{\label{foot:localisation}Suppose $u \in \ell^2$ is a non-zero
  solution of the eigenvalue equation (\ref{eq:mathieu}) and let $x_n =
  \arctan(u_{n-1}/u_n)$ (see Section~\ref{Harper}).  Further, denote the
  Harper map given by (\ref{eq:harper}) by $T$. Then, using the formula
  derived in Footnote~\ref{foot:harper1}, we obtain that
  \begin{eqnarray*}
    DT^n_\theta(x_0) & = & \frac{u_0^2+u_{-1}^2}{u_n^2+u_{n-1}^2} \ .
  \end{eqnarray*}
  (Note that $u_{n-1}=u_n=0$ is not possible, as otherwise $u=0$.)
  Consequently, sink-source-orbits correspond to exponentially decaying
  eigenfunctions. The existence of such `localised' eigenfunctions for the
  almost-Mathieu operator was shown by Jitomirskaya in
  \cite{jitomirskaya:1999} (so-called {\em Anderson localisation}).}
 and indicates an intersection of the stable and unstable subspaces of the
matrix cocycle (see \cite{haro/puig:2006} for a more detailed discussion).

Summarising we can say that the `sharp' peak makes it possible to
concentrate on a single orbit instead of a whole sequence of graphs,
and the information about this orbit will already be sufficient to
establish the existence of a SNA. The fact that the construction in the proof
of our main results (Theorems~\ref{thm:snaexistence} and
\ref{thm:symmetricsna}), which is based of this idea, works fine can be seen
as an indirect `proof' of the claim we made in the above discussion.



\section{Statement of the main results and applications}
\label{Main results}

In this section we state and discuss the main results of this article and
their application to the examples from the introduction. The proofs are
postponed until the later sections, unless they can be given in a few
lines. In particular, this concerns the construction of sink-source-orbits,
which is carried out in Sections \ref{Strategy} to \ref{Symmetricsetting}.

Before we turn to results on the non-smoothness of bifurcations in
Section~\ref{NonSmoothSaddleNodes}, we provide a general setting in which
saddle-node bifurcations in qpf interval maps take place
(Section~\ref{GeneralSettingI}), and introduce sink-source-orbits as a
criterium for the existence of SNA (Section~\ref{SinkSourceSNA}).

\subsection{A general setting for saddle-node bifurcations in qpf interval
  maps} \label{GeneralSettingI}

Obviously, before we can study the non-smoothness of saddle-node bifurcations,
we have to provide a setting in which such bifurcations occur (smooth or
non-smooth). In order to do so, we will consider parameter families of maps
$T=T_\beta$ which are given by
\begin{equation}
  \label{eq:generalsystem}
  T_\beta\thx \ := \ (\thom,F(x) - \beta \cdot g(\theta)) \ ,
\end{equation}
where we suppose that, given a constant $C>0$, the functions $F$ and $g$
satisfy the following assumptions:
\begin{eqnarray}
  \label{cond:g1}
& & g:\kreis \ra [0,1] \textrm{ is continuous and takes the value 1 at least
once;} \\
\label{cond:F1} \textstyle
& &F : [-2C,2C] \ra [-C,C] \textrm{ is continuously differentiable with } F' > 0; \\
& &F \textrm{ has exactly three fixed points } x_-< 0, 0 \ \textrm{and } x_+>0
. \label{cond:F2}
\end{eqnarray}
Note that if we restrict to parameters $\beta \in [0,C]$, then we can choose
$X = [-2C,2C]$ as the driven space, because then $\kreis \times X$ is always
mapped into itself. Of course, this choice is somewhat arbitrary, the only
thing which is important is to fix some driven space $X$ independent of the
parameter $\beta$.  We also remark that the in the situations we will consider
later, $F$ is usually a bounded function which is defined on the whole real
line. In this case, we will only consider its restriction $F_{|[-2C,2C]}$,
where $C$ is any constant larger that $\sup_{x\in\R} |F(x)|$. This has the
advantage that we obtain a compact phase space in this way. In particular, it
allows to define the global attractor and the bounding graphs as it was done
in Section~\ref{BasicDefinitions}~.

As we chose the function $g$ to be
non-negative, the forcing only `acts downwards'. We will refer to this case as
\textit{`one-sided forcing'}.

The first problem we will encounter is to restrict the number of invariant
graphs which can occur. If there are too many, it will be hard to describe a
saddle-node bifurcation in detail. Fortunately, there exist general results
which allow this, without placing to restrictive conditions on the system. We
will discuss these in Section~\ref{SaddleNode} (see
Theorems~\ref{thm:schwarzian} and \ref{thm:concave}, taken from
\cite{jaeger:2003} and \cite{keller:1996}), before giving the proof of
Theorem~\ref{thm:saddlenode}~. The most convenient of these criteria is
to require $F$ to have negative Schwarzian derivative, which ensures that
there can be at most three different invariant graphs
(Theorem~\ref{thm:schwarzian}).%
\footnote{\label{foot:schwarzian}The Schwarzian derivative of a ${\cal C}^3$
    interval map $F$ is defined as
\[
     SF \ := \ \frac{F'''}{F'} -
    \frac{3}{2}\left(\frac{F''}{F'}\right)^2 \ .
\]
    It is intimately related to the cross ratio distortion of the map
    (see \cite{demelo/vanstrien:1993}), and this relation is exploited in
    \cite{jaeger:2003} to derive the mentioned statement. This is very similar
    to the proof of Theorem~\ref{thm:concave} given in
    Section~\ref{SaddleNode} (see Remark~\ref{bem:schwarzianproof}).}
However, in the particular situation we consider it will also be sufficient if
$F$ is convex on one suitable interval and uniformly contracting on
another. More precisely, we will use the following assumption:
\begin{equation} \label{cond:F3}
  \begin{array}{l} \textrm{Suppose that either } F \textrm{ is } {\cal C}^3
    \textrm{ and has negative Schwarzian derivative,} \\ \textrm{or there
      exists } c \in (x_-,0] \textrm{, such that } F_{|[-2C,x_-)} \textrm{ is
      uniformly contracting } \\\textrm{and } F_{|[c,2C)}\textrm{ is strictly concave.}
  \end{array}
\end{equation}
Now we can state the following result on the existence of saddle-node
bifurcations. 
\begin{thm}[Saddle-node bifurcation]
  \label{thm:saddlenode}
  Suppose $F$ and $g$ satisfy (\ref{cond:g1})--(\ref{cond:F3}) and let
  $X=[-2C,2C]$ and $\beta \in [0,C]$ as above. 
  Then the lower bounding graph of the system (\ref{eq:generalsystem}), which
  we denote by  $\varphi^-$,%
\footnote{We keep the dependence of $\varphi^-$ on $\beta$ implicit, same for
  $\psi$ and $\varphi^+$ below.}
is continuous and has negative Lyapunov exponent. Its dependence on $\beta$ is
  continuous (in ${\cal C}^0$-norm) and monotone: If $\beta$ is increased then
  $\varphi^-$ moves downwards, uniformly on all fibres. 

Further, there exists
  a critical parameter $\beta_c \in (0,C)$, such that the following holds:
\romanlist
\item If $\beta < \beta_c$, then there exist exactly two more invariant graphs
  above $\varphi^-$, both of which are continuous. We denote the upper one by
  $\varphi^+$ and the middle one by $\psi$, such that $\varphi^- < \psi <
  \varphi^+$. There holds $\lambda(\psi) > 0$ and $\lambda(\varphi^+) < 0$,
  and the dependence of the graphs on $\beta$ is continuous and monotone: If
  $\beta$ is increased then $\varphi^+$ moves downwards, whereas $\psi$ moves
  upwards, uniformly on all fibres.
\item If $\beta=\beta_c$, there exist either one or two more invariant graphs
  above $\varphi^-$. We denote them by $\psi$ and $\varphi^+$ again (allowing
  $\psi=\varphi^+$), where $\psi \leq \varphi^+$. Further, one of the two
  following holds:
  \begin{itemize}
    \item $\psi$ equals $\varphi^+$ $m$-a.s.\ and
    $\lambda(\psi)=\lambda(\varphi^+)=0$ {\em (Smooth Bifurcation)}.%
  \footnote{Of course, the natural possibility here is that $\psi$ and
    $\varphi^+$ are continuous and coincide everywhere. However, there is also
    a second, rather pathological alternative, which cannot be excluded: It
    might happen that there exists no continuous invariant graph apart from
    $\varphi^-$, but two semi-continuous invariant graphs $\psi$ and
    $\varphi^+$ which are $m$-a.s.\ equal. This is discussed in more detail in
    Section~\ref{Essentialclosure}~. Whether the bifurcation should really be
    called smooth in this case is certainly debatable. However, as the
    non-smooth bifurcations we prove to exist later on all involve non-zero
    Lyapunov exponents, we prefer this as a working definition in the context
    of this paper.}
    \item $\psi\neq\varphi^+$ $m$-a.s., $\lambda(\psi)>0$,
  $\lambda(\varphi^+)<0$ and both invariant graphs are non-continuous
  \emph{(Non-smooth Bifurcation)}.
  \end{itemize}
  In any case, the set $B := [\psi,\varphi^+]$ is compact and the set %
  $\{\theta \in \kreis \mid \psi(\theta) = \varphi^+(\theta) \}$ is dense in
  $\kreis$.%
  \footnote{A compact set $B \ssq \kreis \times X$ is called
    \textit{pinched}, if for a dense set of $\theta$ the set $B_\theta
    := \{ x \in X \mid \thx \in B\}$ consists of a single point. Thus,
    the last property could also be stated as `the set $B$ is
    pinched'.  }
\item If $\beta > \beta_c$, then $\varphi^-$ is the only invariant
  graph. 
\listend 
\end{thm}
The proof of Theorem~\ref{thm:saddlenode}, together with some preliminary
results which are needed, is given in Section~\ref{SaddleNode}~.
\medskip

When $F$ depends on an additional parameter, it is also natural to study
the dependence of the critical parameter $\beta_0$ on this parameter. We
refrain from producing a general statement and just concentrate on the
$\arctan$-family (\ref{eq:atanfamily}) given in the introduction. Let
\[
    F_\alpha(x) \ := \ \frac{\arctan(\alpha x)}{\arctan(\alpha)} \ .
\]
\begin{lem}
  \label{lem:parameterdependence}
  Let $\beta_0(\alpha)$ denote the critical parameter of the
  saddle-node bifurcation in Theorem~\ref{thm:saddlenode} with $F =
  F_\alpha$ in (\ref{eq:generalsystem}). Then $\alpha \mapsto
  \beta_0(\alpha)$ is continuous and strictly monotonically increasing
  in $\alpha$.
\end{lem}
Again, the proof is given in Section~\ref{SaddleNode}~. We note that while
continuity follows under much more general assumptions, the monotonicity
depends on the right scaling of the parameter family, namely on the fact that
the fixed points of $F_\alpha$ do not depend on $\alpha$.


\subsection{Sink-source-orbits and the existence of SNA} \label{SinkSourceSNA}

In this subsection we consider a slightly more general situation than in the last,
and suppose that
\begin{eqnarray}\label{eq:T1}
  && T \textrm{ is a qpf monotone interval map};\\\label{eq:T2}
  && \textrm{All fibre maps } T_\theta \textrm{ are differentiable with
  derivative } DT_\theta;\hspace{3eM}\\\label{eq:T3}
  && \thx \mapsto DT_\theta(x) \textrm{ is continuous and strictly positive.}
\end{eqnarray}
 In particular, this applies to parameter families which satisfy
 (\ref{cond:g1})--(\ref{cond:F2}). \medskip

In order to formulate the statements of this section, we have to introduce
different Lyapunov exponents. Let $\thx \in \kreis \times X$. Then the {\em
(vertical) finite-time forward} and {\em backward Lyapunov exponents} are
defined as
\begin{equation}
  \label{eq:finitelyap}
  \lambda^+(\theta,x,n)\ := \ \ntel\inergsum
  \log(DT_{\theta+i\omega}(\Tth^i(x)))
\end{equation}
and 
\begin{equation}
  \label{eq:finitebacklyap}
  \lambda^-(\theta,x,n)\ := \  -\ntel\insum
  \log(DT_{\theta-i\omega}(\Tth^{-i}(x))) \ .
\end{equation}
 When dealing with parameter families as in (\ref{eq:generalsystem}), we will
write $\lambda^\pm(\beta,\theta,x,n)$ for the pointwise finite-time Lyapunov exponents
with respect to the map $T_\beta$ if we want to keep the dependence on the
parameter $\beta$ explicit. 

As it is not always possible to ensure that the finite-time exponents
converge as $n\to\infty$, we distinguish between upper and lower Lyapunov
exponents: The \textit{(vertical) upper forward Lyapunov exponent} of a point
$\thx \in \kreis \times X$ is defined as
\begin{equation}
  \label{eq:pointwiselyap}
  \lambda^+(\theta,x)\ := \  \limsup_{n\to\infty} \lambda^+(\theta,x,n) \ .
\end{equation}
Similarly, the \textit{upper backward Lyapunov exponent} is defined as
\begin{equation}
  \label{eq:backpointwiselyap}
  \lambda^-(\theta,x)\ := \ \limsup_{n\to\infty} \lambda^-(\theta,x,n) \ .
\end{equation}
In the same way, we define the {\em lower forward} and {\em backward Lyapunov
 exponents}, replacing $\limsup$ by $\liminf$:
\begin{eqnarray}
  \label{eq:lowerpointwiselyap}
  \lambda_{\textrm{low}}^+(\theta,x)\ &= &  \liminf_{n\to\infty}
  \lambda^+(\theta,x,n) \  ; \\
  \label{eq:lowerbackpointwiselyap}
  \lambda_{\textrm{low}}^-(\theta,x)& := &\liminf_{n\to\infty}
  \lambda^-(\theta,x,n) \ .
\end{eqnarray}
Again, we write
$\lambda^\pm(\beta,\theta,x),\lambda^\pm_{\textrm{low}}(\beta,\theta,x)$ if we
want to keep the dependence on a parameter $\beta$ explicit. 

For any invariant graph $\varphi$, the Birkhoff ergodic theorem implies that
for $m$-a.e.\ $\theta \in \kreis$ the $\limsup$ and the $\liminf$ coincide
(i.e.\ the respective limits exists and we do not have to distinguish between
$\lambda^\pm$ and $\lambda^\pm_{\textrm{low}}$) and there holds
$\lambda^+(\theta,\varphi(\theta)) = -\lambda^-(\theta,\varphi(\theta)) =
\lambda(\varphi)$.  Further, when $\varphi$ is continuous the Uniform Ergodic
Theorem (e.g.\ \cite{katok/hasselblatt:1995}) implies that this holds for all
$\theta\in\kreis$ and the convergence is uniform on $\kreis$. Now, consider
the situation where $\psi$ is an unstable and $\varphi$ is a stable continuous
invariant graph, and there is no other invariant graph in between. Then points
on the repeller (or \textit{source}) $\psi$ will have a positive forward and a
negative backward Lyapunov exponent, and for points on the attractor (or
\textit{sink}) $\varphi$ it is just the other way around. Further, all points
between $\psi$ and $\varphi$ will converge to $\varphi$ forwards and to $\psi$
backwards in time, thus moving from source to sink, and consequently both
their exponents will be negative. These three cases should be considered as
more or less typical. In contrast to this, the remaining possibility of both
Lyapunov exponents being positive is rather strange, as it would suggest that
the orbit somehow moves from a sink to a source. This motivates the following
definition:
\begin{definition}[Sink-source-orbits]\label{def:sinksource}
  Suppose $T$ satisfies the assumptions (\ref{eq:T1})--(\ref{eq:T3}). Then we
  call an orbit of $T$ which has both positive forward and backward lower
  Lyapunov exponent a $\textbf{sink-source-orbit}$. If an orbit has both
  positive forward and backward upper Lyapunov exponent then we call it a
  $\textbf{weak sink-source-orbit}$.
\end{definition}
Obviously, every sink-source-orbit is also a weak sink-source-orbit. 

As mentioned in the introduction, the existence of sink-source-orbits is
already known for the Harper map (see Footnote~\ref{foot:localisation}), where
they only occur together with SNA (i.e.\ in the non-uniformly hyperbolic case,
as discussed in Section~\ref{Harper}).  This is not a mere coincidence:
\begin{thm}
  \label{thm:sinksourcesna}
  Suppose $T$ satisfies the assumptions (\ref{eq:T1})--(\ref{eq:T3}). Then the
  existence of a weak sink-source-orbit implies the existence of a SNA (and
  similarly of a SNR).
\end{thm}
The proof is given in Section~\ref{Sinksourceorbits}~. 
\begin{bem}\label{bem:sinksource}
\alphlist
\item In the proofs of Theorems~\ref{thm:snaexistence} and
  \ref{thm:symmetricsna} below, we actually construct
  sink-source-orbits. Thus, for the main purpose of this paper it would not
  have been necessary to introduce weak sink-source-orbits. However, since the
  existence of the latter is a much weaker assumption than the existence of a
  sink-source-orbit (see also (b) and (c) below), it seemed appropriate to state
  Theorem~\ref{thm:sinksourcesna} in this way. 
\item In some situations, it is also possible to obtain results in the
  opposite direction. For example, if $M$ is a minimal set which contains both
  a SNA and a SNR, then weak sink-source-orbits are dense (even residual) in
  $M$. In order to see this, note that, in
  the above situation, for some constant $c>0$ the set $M$ contains a point
  $(\theta_1,x_1)$ with $\lambda^+(\theta_1,x_1) > c$ and a point
  $(\theta_2,x_2)$ with $\lambda^-(\theta_2,x_2) > c$. Due to minimality, it
  follows that the open sets
  $$A_n := \{ \thx \in M \mid \exists m\geq n : \lambda^+(\theta,x,m) > c\}$$
  and $$B_n := \{ \thx \in M \mid \exists m \geq n: \lambda^-(\theta,x,m) >
  c\}$$ are both dense in $M$. By Baire's Theorem, their intersection
  $S:=\bigcap_{n\in\N} A_n\cap B_n$ is residual, and obviously every point in
  $S$ belongs to a weak sink-source-orbit.
\item The preceding remark becomes false if `weak sink-source-orbit' is
  replaced by sink-source-orbit. In fact, it is well-known that SNA may exist
  in the absence of sink-source-orbits, even if there is a minimal set which
  contains both an SNA and a SNR. Examples are provided by the Harper map: As
  we have discussed in Section~\ref{Harper}, the existence of a
  sink-source-orbit is equivalent to the existence of an exponentially
  decaying eigenfunction for the corresponding Schr\"odinger
  operator. However, there are situations in which, for certain energies in
  the spectrum of $H_{\lambda,\theta}$ (which does not depend on $\theta$),
  there exist no such `localised' eigenfunctions, independent of
  $\theta$. This follows for example from Theorem 5 in \cite{puig:2006},
  together with the concept of Aubry-duality, which is explained in Section 2
  of the same paper (the original source is \cite{aubry/andre:1980}).  The
  fact that there is a (unique) minimal set which contains both an SNA and an
  SNR in these examples is shown in \cite[Section 4.17]{herman:1983}.

 For a more detailed discussion of the implications of spectral-theoretic
  results for the Harper map, we also refer to \cite{haro/puig:2006} (this
  particular issue is addressed in Section V(C)).  \listend
\end{bem}

An observation which was made frequently in numerical studies of SNA is a very
unusual distribution of the finite-time Lyapunov exponents. The interesting
fact is that although in the limit all observed Lyapunov exponents were
negative, the distribution of the finite-time Lyapunov exponents still showed
a rather large proportion of positive values, even at very large times (see
\cite{pikovski/feudel:1995},\cite{prasad/negi/ramaswamy:2001}).  Of course,
the existence of a sink-source orbit could be a possible explanation for such
a behaviour. On the other hand, we can also use information about the
finite-time Lyapunov exponents to establish the existence of a
sink-source-orbit, and this will play a key role in the proof of our main
results:
\begin{lem}
  \label{lem:sinksourceshadowing}
  Let $I$ be a compact metric space $\R$ and $(T_\beta)_{\beta\in I}$ be a
  parameter family of qpf monotone interval maps which all satisfy the
  assumptions (\ref{eq:T1})--(\ref{eq:T3}) above. Further, assume that the
  dependence of the maps $T_\beta$ and $\thx \mapsto DT_{\beta,\theta}(x)$ on
  $\beta$ is continuous (w.r.t.\ the topology of uniform convergence).

Suppose there exist sequences of integers $l^-_1,l^-_2,\ldots \nearrow
\infty$ and $l^+_1,l^+_2,\ldots \nearrow \infty$, a sequence
$(\theta_p,x_p)_{p\geq 1}$ of points in $\kreis \times X$ and a
sequence of parameters $(\beta_p)_{p\geq 1}$, such that for all
$p\in\N$ there holds
\[
    \lambda^+(\beta_p,\theta_p,x_p,j) \ > \ c \ \ \ \forall j=1\ld l^+_p
\]
and
\[
    \lambda^-(\beta_p,\theta_p,x_p,j) \ > \ c \ \ \ \forall j=1\ld
    l^-_p \ 
\]
for some constant $c>0$.  Then there is at least one $\beta_0 \in I$,
such that there exists a sink-source-orbit (and thus a SNA-SNR-pair) for the
map $T_{\beta_0}$.
\end{lem}
\proof\ In fact, the statement is a simple consequence of compactness and
continuity: By going over to suitable subsequences if necessary, we can assume
that the sequences $(\theta_p)_{p\geq 1},(x_p)_{p\geq 1}$ and
$(\beta_p)_{p\geq 1}$ converge. Denote the limits by $\theta_0,\ x_0$ and
$\beta_0$, respectively.

Now, due to the assumptions on $T_\beta$ and $DT_{\beta,\theta}(x)$
the functions $(\beta,\theta,x) \mapsto \lambda^\pm(\beta,\theta,x,j)$
are continuous for each fixed $j\in\N$. Thus, we obtain
\[
    \lambda^\pm(\beta_0,\theta_0,x_0,j) \ = \ \lim_{p\ra\infty}
    \lambda^\pm(\beta_p,\theta_p,x_p,j) \ \geq \ c \ \quad \forall j \in \N\ ,
\]
such that 
\[
     \lambda_{\textrm{low}}^\pm(\beta_0,\theta_0,x_0) \ = \ \liminf_{j \ra \infty}
     \lambda^\pm(\beta_0,\theta_0,x_0,j) \ \geq \ c  \ > \ 0 \ .
\]
Hence, the orbit of $(\theta_0,x_0)$ is a
sink-source-orbit for the map $T_{\beta_0}$.

\qed

%


\subsection{Non-smooth bifurcations} \label{NonSmoothSaddleNodes}

In order to formulate the results concerning the non-smoothness of
bifurcations and the existence of SNA, we first have to quantify the
qualitative features of the functions $F$ and $g$ which were used in the
discussion in Section~\ref{Mechanism}~. Some of the assumptions we will make
below are quite specific and could in principle be formulated in a more
general way. However, as the proofs of Theorems~\ref{thm:snaexistence} and
\ref{thm:symmetricsna} are quite involved anyway, we refrain from introducing
any more additional parameters, even if this could lead to slightly more
flexible results. As we have mentioned before, our main goal here is just to
show that the presented approach does lead to rigorous results at all, we do
not aim for the greatest possible generality. Hence, we content ourself here
to provide a statement which applies, after suitable rescaling and
reparametrisation, to at least two of the main examples from the introduction
(see Sections \ref{App:atanfamily} and \ref{App:Harper}). \medskip

First of all, we will suppose that $\gamma$ and $\alpha$ are positive constants which
satisfy
\begin{eqnarray}
  \gamma & \leq & 1/16 \ ; \label{cond:gamma0} \\  \walpha & > & 4/\gamma \
  \  \geq \ 64 \
      \label{cond:alphagamma0}  \ .
\end{eqnarray}
Further, we will assume (in addition to (\ref{cond:g1})--(\ref{cond:F3})), that
\begin{eqnarray}
      && F([-3,3]) \ssq [-3/2,3/2] \ \textrm{ (in other words\ } C=3/2 \textrm{ in
      (\ref{cond:F1}))}; \\ && F(0) = 0\ \textrm{ and }\ F(\pm x_\alpha) = \pm
      x_\alpha \ \textrm{ where } \ {\textstyle x_\alpha \ := \ 1 +
      \frac{2}{\sqrt{\alpha}} } \ ;
      \label{cond:Ffixedpoints} \\
      && 2\alpha^{-2} \ \leq \ F'(x) \ \leq \ \alpha^2 \hspace{2eM} \forall x
       \in [-3,3] \ ;
       \label{cond:Funiformbounds} \\
       && F'(x)  \ \geq \ 2\alpha^{\halb} \hspace{5.5eM} \forall x \in
       \overline{B_\frac{2}{\alpha}(0)} \ ;
       \label{cond:Fexpansion} \\
       && F'(x) \ \leq \  {\textstyle \halb\alpha^{-\halb} } \hspace{4.7eM}
        \forall x : |x|  \geq  \gamma \ ;
       \label{cond:Fcontraction}\\
       &&F({\textstyle\alphtel}) \ \geq \ 1-\gamma \ \textrm{ and }
       \ F(-{\textstyle \alphtel}) \ \leq \ -(1-\gamma) \ .
       \label{cond:Fmapsover}
\end{eqnarray}
Finally, we will require that
\begin{equation}
  g : \kreis \to [0,1] \ \textrm{ has the unique maximum } g(0) = 1
\end{equation}
and for some constants $L_1,L_2>0$ there holds
\begin{eqnarray}
 && g \ \textrm{ is Lipschitz-continuous with Lipschitz constant } L_1\ \label{eq:g-lipschitz} ;\\ &&
 g(\theta) \ \leq \ \max\{ 1-3\gamma , 1-L_2\cdot d(\theta,0) \} \ ,
      \label{cond:sharppeak}
\end{eqnarray}
where $d$ denotes the usual Euclidean distance on the circle.
Essentially, this quantifies the properties which we have already
mentioned in Section~\ref{Mechanism}: $F$ has three fixed points
(\ref{cond:Ffixedpoints}), acts highly expanding close to 0
(\ref{cond:Fexpansion}) and highly contracting further away
(\ref{cond:Fcontraction}). Thus, the expanding region $\cal E$ from
Section~\ref{Mechanism} corresponds to $\kreis \times \overline{B_\frac{2}{\alpha}(0)}$,
whereas the contracting region $\cal C$ corresponds to $\kreis \times
[\gamma,3]$. Further, (\ref{cond:Fmapsover}) ensures that $\kreis \times \Balphcl$ is
mapped over itself in a very strong sense, and finally condition
(\ref{cond:sharppeak}) makes precise what we meant when speaking of a
`sharp peak' before.

The last assumption we need is a Diophantine condition on the rotation
number $\omega$. We use the notation
\begin{equation}
          \omn := n\omega \bmod 1 \ 
\end{equation}
and suppose that there exist constants $c,d > 0$, such that
\begin{equation}
   \label{cond:diophantine} 
    d(\omn,0) \ \  \geq \ \ c \cdot n^{-d} \ \ \ \forall n \in \N \ . 
\end{equation}
(Here $d(\theta,\theta')$ denotes the usual Euclidean distance of two points
$\theta,\theta'\in \kreis$.)
\begin{thm}
  \label{thm:snaexistence}
  Suppose $\alpha,\gamma,F$ and $g$ are chosen such that
  (\ref{cond:g1})--(\ref{cond:F3}) and
  (\ref{cond:gamma0})--(\ref{cond:sharppeak}) hold. Further, assume that $\omega$
  satisfies the Diophantine condition (\ref{cond:diophantine}) and let
  $$T_\beta\thx \ = \ (\thom,F(x)-\beta g(\theta))$$ as in
  (\ref{eq:generalsystem}). Let $\beta_c\in(0,3/2)$ be the critical
  parameter of the saddle-node bifurcation described in
  Theorem~\ref{thm:saddlenode}~. Then there exist constants $\gamma_0 =
  \gamma_0(L_1,L_2,c,d)>0$ and $\alpha_0 = \alpha_0(L_1,L_2,c,d)>0$ with the
  following property:

  If $\gamma < \gamma_0$ and $\alpha > \alpha_0$, then there exists a
  sink-source-orbit for the system $T_{\beta_c}$. Consequently, there exists a
  SNA (the invariant graph $\varphi^+$ in Theorem~\ref{thm:saddlenode}(ii))
  and a SNR ($\psi$ in Theorem~\ref{thm:saddlenode}(ii)), and both objects
  have the same essential closure.%
\footnote{See Section~\ref{Essentialclosure} for the definition of the
  essential closure.}
\end{thm}
The proof of this theorem is given in the Sections
\ref{Strategy}--\ref{Construction}, an outline of the strategy is given at the
beginning of Section~\ref{Strategy}. 

\begin{bem}\label{bem:snaexistence}
  \alphlist
\item We remark that the existence of a sink-source-orbit in the parameter
  family $T_\beta$ in the above theorem does not depend on the statement of
  Theorem~\ref{thm:saddlenode}. Even if the assumptions
  (\ref{cond:g1})--(\ref{cond:F3}) are dropped and
  Theorem~\ref{thm:saddlenode} no longer applies, we still obtain the
  existence of a parameter $\beta_c$ for which $T_{\beta_c}$ has a
  sink-source-orbit and a SNA-SNR-pair, provided $\gamma$ is sufficiently
  small and $\alpha$ sufficiently large. However, in this case $\beta_c$ is
  not necessarily unique anymore. Further, it is not possible to say whether
  it is a bifurcation parameter, nor to control the number of invariant graphs
  which might occur.
\item The dependence of $\gamma_0$ and $\alpha_0$ on $L_1,L_2,c$ and $d$ can
  be made explicit. More precisely, the conditions which have to be satisfied
  are (\ref{cond:alpha1}), (\ref{cond:hfunctions1}), (\ref{cond:hfunctions2})
  and (\ref{cond:u})--(\ref{cond:alpha3}). Conditions (\ref{cond:hfunctions1})
  and (\ref{cond:hfunctions2}) are somewhat implicit, but once the parameters
  $u$ and $v$ are fixed according to (\ref{cond:u})--(\ref{cond:sigma}),
  explicit formulas can be derived from the proof of
  Lemma~\ref{lem:omegaestimates}~.
\item Numerical observations (as well as the statement of the above theorem)
  suggest that there might be a critical parameter $\alpha^*$, such that the
  saddle-node bifurcation in the family $\Ttil_{\alpha,\beta}$ with fixed
  $\alpha$ is smooth whenever $\alpha<\alpha^*$ and non-smooth whenever
  $\alpha>\alpha^*$.  However, whether this is really the case is completely
  open.
\listend
\end{bem}

As we have mentioned in Section~\ref{Mechanism}, the sharp peak of the forcing
function leads to a localisation of the sink-source-orbit. In fact, its
construction in the later sections yields enough information to determine it
precisely:
\begin{adde}\label{adde}
  In the situation of Theorem~\ref{thm:snaexistence} denote the SNA by
  $\varphi^+$ and the SNR by $\psi$. Then the point $(\omega,\varphi^+(\omega))$ belongs
  to a sink-source-orbit.%
\footnote{There might be more than one sink-source-orbit, but this is the
  particular one which we will construct in the later sections.}
Further, this sink-source-orbit is contained in the intersection
  $\Phi^+\cap\Psi$ (which means $\varphi^+(\omega)=\psi(\omega)$).
\end{adde}
The proof is given in Section~\ref{Construction}.\medskip

Next, we turn to the existence of SNA in symmetric systems. In order to do so,
we have to modify the assumptions on $F$ and $g$. First of all, instead of
(\ref{cond:g1}) and (\ref{cond:F3}) we will assume that
\begin{eqnarray} && \label{cond:g-pitchfork}
  g:\kreis \to [-1,1] \ \textrm{ is continuous and has the unique maximum }
  g(0)=1 ; \\
  &&\label{eq:F-has-schwarzian} F \textrm{ is } {\cal C}^3 \textrm{ and has negative Schwarzian
    derivative.}
\end{eqnarray}
(The alternative in (\ref{cond:F3}) only works for one-sided
forcing). Further, we will require the following symmetry conditions
\begin{eqnarray}
  && \label{cond:Fsymmetry} F(-x) \ = \ -F(x) \ ; \\&& \label{cond:gsymmetry}
  \textstyle g(\theta+\halb) \ = \ -g(\theta) \ .
  \end{eqnarray}
Finally, (\ref{cond:sharppeak}) will be replaced by
\begin{equation} \textstyle
  \label{cond:symmetricpeak}
  |g(\theta)| \ \leq \ \max\left\{1-3\gamma,1-L_2\cdot
    d(\theta,\{0,\halb\})\right\}   \ . 
\end{equation}

Note that (\ref{cond:Fsymmetry}) and  (\ref{cond:gsymmetry})  together imply that
the map $T=T_\beta$ given by (\ref{eq:generalsystem}) has the following
symmetry property:
\begin{equation}
  \label{eq:systemsymmetry}
  -\Tth(x) \ = \ T_{\theta+\halb}(-x) 
\end{equation}
Now suppose that $\varphi$ is a $T$-invariant graph. Then due to
(\ref{eq:systemsymmetry}) the graph given by
\begin{equation}
  \label{eq:symmetricgraph} \textstyle
  \overline{\varphi}(\theta) \ := \ -\varphi(\theta+\halb)
\end{equation}
is invariant as well. In particular, this implies that the upper and lower
bounding graphs satisfy $\varphi^+(\theta) = -\varphi^-(\theta+\halb)$, and if
one of these graphs undergoes a bifurcation, then the same must be true for
the second one as well. As the negative Schwarzian derivative of $F$ will
allow us to conclude that there is only one other invariant graph $\psi$ apart
from the bounding graphs $\varphi^\pm$, this implies that any possible
collision between invariant graphs has to involve all three invariant graphs at
the same time and must therefore be a pitchfork bifurcation. However, as we
have mentioned before, due to the lack of monotonicity in the symmetric
setting we cannot ensure that there is a unique bifurcation
point. Nevertheless, we obtain the following result concerning the existence
of SNA with symmetry:
\begin{thm}
  \label{thm:symmetricsna}
  Suppose $\gamma,\alpha,F$ and $g$ are chosen, such that
  (\ref{cond:gamma0})--(\ref{cond:Fmapsover}), (\ref{eq:g-lipschitz}) and
  (\ref{cond:g-pitchfork})--(\ref{cond:symmetricpeak}) hold. Further, assume
  $\omega$ satisfies the Diophantine condition (\ref{cond:diophantine}) and
  let
  $$T_\beta\thx \ = \ (\thom,F(x)-\beta g(\theta))$$ as in
  (\ref{eq:generalsystem}). Then
  there exist constants $\gamma_0 = \gamma_0(L_1,L_2,c,d)>0$ and $\alpha_0 =
  \alpha_0(L_1,L_2,c,d)>0$ with the following property:

  If $\gamma < \gamma_0$ and $\alpha > \alpha_0$, then there is a parameter
  $\beta_c$ such that there exist two SNA $\varphi^-$ and $\varphi^+$ and a
  SNR $\psi$ with $\varphi^- \leq \psi \leq \varphi^+$ for $T_{\beta_c}$.
  Further, there holds $\esscl{\Phi^-} = \esscl{\Psi} = \esscl{\Phi^+}$, and
  the invariant graphs satisfy the symmetry equations \[ \textstyle
  \varphi^-(\theta) = -\varphi^+(\theta + \halb) \quad \textrm{and} \quad
  \psi(\theta) = -\psi(\theta+\halb) \ .
\]
\end{thm}

\begin{bem}
  As in Theorem~\ref{thm:snaexistence}, the dependence of $\gamma_0$ and
  $\alpha_0$ on $L_1,L_2,c$ and $d$ can be made explicit (compare
  Remark~\ref{bem:snaexistence}(b)). The
  conditions which have to be satisfied are (\ref{cond:alpha1}),
  (\ref{cond:hfunctions1}), (\ref{cond:hfunctions2}),
  (\ref{cond:u})--(\ref{cond:alpha3}), (\ref{cond:alphagammasymmetric}) and
  (\ref{cond:gvalues}).
\end{bem}

\subsection{Application to the parameter families} \label{ex:atanexample} \label{Applications}

The assumptions on $F$ and $g$ used in Theorems~\ref{thm:snaexistence} and
\ref{thm:symmetricsna} are somewhat technical and might seem very
restrictive. However, in this subsection we will see that they are more
flexible that they might look like on first sight (although there are surely
some constraints). In particular, after performing some surgery we can apply
them at least to two of the parameter families from the introduction, namely
the $\arctan$-family with additive forcing and the Harper map. In both cases,
the respective parameters have to chosen sufficiently large, but of course
this goes perfectly well with the statement of
Theorem~\ref{thm:snaexistence}~. As a consequence, the respective corollaries
become a lot easier to formulate.

The qpf Arnold map then demonstrates the limits of
Theorem~\ref{thm:snaexistence}, since it is not possible to apply the result
in this case. This is briefly discussed in Subsection~\ref{App:Remarks}~.

\subsubsection{Application to the arctan-family.} \label{App:atanfamily}

Applied to the $\arctan$-family with additive forcing,
Theorem~\ref{thm:snaexistence} yields the following:

\begin{cor} \label{cor:atanfamily}
  Suppose $\omega$ satisfies the Diophantine condition
  (\ref{cond:diophantine}). Then there exists $\alpha_0 = \alpha_0(c,d)$ such
  that for all $\alpha>\alpha_0$ the system $T_{\alpha,\beta}$ given by
  (\ref{eq:atanfamily}) undergoes a non-smooth saddle-node bifurcation as the
  parameter $\beta$ is increased from 0 to 1.
\end{cor}
\begin{bem}
  As already mentioned in Section~\ref{AtanFamily}, the above statement
  remains true if the $\arctan$ in (\ref{eq:atanfamily}) is replaced by the map $x
  \mapsto\frac{x}{1+|x|}$, or any other function which has similar scaling
  properties. This will become obvious in the following proof, but we refrain
  from producing a more general statement here. 
\end{bem}

\textit{Proof of Corollary~\ref{cor:atanfamily}~.} Since the system (\ref{eq:atanfamily}) does not satisfy
(\ref{cond:Ffixedpoints}), we cannot apply Theorem~\ref{thm:snaexistence}
directly. Therefore, we start by considering a slightly rescaled version of
(\ref{eq:atanfamily}). Let
\[
     \Ftil_\alpha(x) \ := \ C(\alpha)\cdot \arctan(\alpha^{\frac{4}{3}}x) \ \
     \textrm{ where } \ \ C(\alpha) \ := \
     \frac{ 1+\frac{2}{ \sqrt{\alpha} } }{\arctan(\alpha^\frac{4}{3} +
     2\alpha^\frac{5}{6})} \ 
\]
and
\[
     g(\theta) \ := \ 1 - \sin(\pi\theta) \ .
\]
Note that $\Ftil_\alpha$ always satisfies (\ref{cond:Ffixedpoints}). The
important thing we have to ensure is that whenever we fix a suitably small
$\gamma$, such that (\ref{cond:gamma0}), (\ref{cond:sharppeak}) and any
additional smallness conditions on $\gamma$ which appear later on are satisfied, then
(\ref{cond:Fmapsover}) holds for all sufficiently large values of
$\alpha$. This means that we can first fix $\gamma$, and then ensure that all
inequalities involving $\alpha$ alone or both $\alpha$ and $\gamma$, such as
(\ref{cond:alphagamma0}), hold by choosing $\alpha$ sufficiently large,
without worrying about (\ref{cond:Fmapsover}).  However, in this particular
case it is easy to see that $\Ftil_\alpha(\alphtel) = C(\alpha) \cdot
\arctan(\alpha^{\drittel}) \ra 1$ as $\alpha \ra \infty$ (note that
$\lim_{\alpha\ra\infty} C(\alpha) = \frac{2}{\pi}$), which is exactly what we
need.

Now, if $\gamma$ is chosen small enough $g(\theta) = |1-\sin(\pi\theta)|$
clearly satisfies $(\ref{cond:sharppeak})$, for example with $L_2 := 2$. The
Lipschitz-constant $L_1$ is $\pi$. Thus, it remains to check
the assumptions on the derivative of $\Ftil_\alpha$. To that end, note that
\[
    \Ftil_\alpha'(x) \ = \ C(\alpha) \cdot
    \frac{\alpha^\frac{4}{3}}{1+\alpha^\frac{8}{3}x^2} 
\]
We have $\Ftil_\alpha'(0) \sim \alpha^\frac{4}{3},\
\Ftil_\alpha'(\frac{2}{\alpha}) \sim \alpha^\frac{2}{3}$ and
$\Ftil_\alpha'(\gamma) \sim \alpha^{-\frac{4}{3}}$ for each fixed $\gamma > 0$
as $\alpha \ra \infty$. Therefore, the conditions
(\ref{cond:Funiformbounds}),(\ref{cond:Fexpansion}) and
(\ref{cond:Fcontraction}) will always be satisfied when $\alpha$ is large
enough. Consequently we can apply Theorem~\ref{thm:snaexistence} and obtain
that there exists some $\tilde{\alpha}_0$ such that for all $\alpha\geq
\tilde{\alpha}_0$ the parameter family $$\Ttil_{\alpha,\beta} \ : \ \thx \mapsto
(\thom,\Ftil_\alpha(x) - \beta g(\theta))$$ undergoes a non-smooth pitchfork
bifurcation (in the sense of Theorem~\ref{thm:saddlenode}) as $\beta$ is
increased from $0$ to $3/2$. \medskip

Now denote the map given by (\ref{eq:atanfamily}) by $T_{\alpha,\beta}$. We
claim that there exists a monotonically increasing function $\sigma : \R^+ \ra
\R^+$ and a function $\tau : \R^+ \ra \R^+$ such that $T_{\alpha,\beta}$ is
smoothly conjugate to
$\Ttil_{\sigma(\alpha),\tau(\alpha)\beta}$. Consequently, the parameter family
$T_{\alpha,\beta}$ equally exhibits non-smooth saddle-node bifurcations if
$\alpha$ is chosen sufficiently large (larger than
$\sigma^{-1}(\tilde{\alpha}_0)$).

In order to define $\sigma$, it is convenient to introduce an
intermediate parameter family $\That_{\alpha,\beta}$ with fibre maps
\[
  \That_{\alpha,\beta,\theta}(x) \ = \ \arctan(\alpha x) - \beta
  g(\theta) \ .
\]
We let $h_1(\theta,x) = (\theta,\arctan(\alpha)x)$,
$\sigma_1(\alpha) = \arctan(\alpha)^{-1} \alpha$ and $\tau_1(\alpha) = 
\arctan(\alpha)$. Then 
\[
   T_{\alpha,\beta} \ = \ h_1^{-1} \circ
   \That_{\sigma_1(\alpha),\tau_1(\alpha)\beta} \circ h_1 \ ,
\]
such that $T_{\alpha,\beta} \sim \That
  _{\sigma_1(\alpha),\tau_1(\alpha)\beta}$, where $\sim$ denotes
  the existence of a smooth conjugacy. 

On the other hand, let $h_2(\theta,x) = (\theta, C(\alpha)^{-1}x)$,
$\sigma_2(\alpha) = C(\alpha)\alpha^{\frac{4}{3}}$ and $\tau_2(\alpha) 
= C(\alpha)^{-1}$. Again, a simple computation yields
\[
  \Ttil_{\alpha,\beta} \ = \ h^{-1}_2 \circ
  \That_{\sigma_2(\alpha),\tau_2(\alpha)\beta} \circ h_2 \ .
\]
As $\sigma_1$ and $\sigma_2$ are both strictly monotonically
increasing and therefore invertible, this implies
$\That_{\alpha,\beta} \sim
\Ttil_{\sigma_2^{-1}(\alpha),\tau_2(\sigma_2^{-1}(\alpha))^{-1}\beta}$
and consequently
\[
    T_{\alpha,\beta} \ \sim \
    \That_{\sigma_1(\alpha),\tau_1(\alpha)\beta} \ \sim \
    \Ttil_{\sigma_2^{-1} \circ \sigma_1(\alpha),\tau_2(\sigma_2^{-1}
      \circ \sigma_1(\alpha))^{-1} \tau_1(\alpha)\beta} \ .
\]
Hence, we can define $\sigma = \sigma_2^{-1} \circ \sigma^1$ and $\tau 
= \frac{\tau_1}{\tau_2\circ \sigma_2^{-1} \circ \sigma_1}$ as
claimed. 

Finally, since $F_\alpha$ has the fixed point $x_+=1$, $g(0)=1$ and we are in
the case of one-sided forcing, it can easily be seen that the bifurcation must take place
before $\beta=1$ (meaning that the critical parameter $\beta_c$ given
Theorem~\ref{thm:saddlenode} is strictly smaller than one). For larger
$\beta$-values all orbits eventually end up below the $0$-line and
consequently converge to the lower bounding graph $\varphi^-$, such that this
is the only invariant graph. (Compare with the proof of
Theorem~\ref{thm:saddlenode} in Section~\ref{SaddleNode}.) This completes the
proof.

\qed

%

\subsubsection{Application to the Harper map} \label{App:Harper}

We want to emphasise that we do not claim any originality for the presented
results on the Harper map. Our aim is merely to demonstrate the flexibility of
our general statements by applying them to this well-known family. For the
particular case of the Harper map there surely exist more direct and elegant
ways to produce such results, starting with \cite{herman:1983}. Usually such
results require more regularity than we use here (the potentials we can treat
are only Lipschitz-continuous), but from the physics point of view this is
surely the more interesting case anyway. Further, although potentials which
are only Lipschitz-continuous are not explicitly treated in
\cite{bjerkloev:2005}, the methods developed there surely allow to do this as
well. Thus, the real achievement here is rather to show that the underlying
mechanism for non-smooth bifurcations is in principle the same in the Harper
map as in other parameter families, like the $\arctan$-family with additive
forcing, despite the very particular structures which distinguish
Schr\"odinger cocycles from other models. \medskip

As Theorem~\ref{thm:snaexistence} is tailored-made for qpf interval maps, its
application to the Harper map is somewhat indirect. This means that we have to
perform a number of modifications before the system in (\ref{eq:projharper})
is in a form which meets the assumptions of the theorem. First of all, we
remark that the dynamics of (\ref{eq:projharper}) are equivalent to those of
the map
\begin{equation}\label{eq:riccati} \thx \ 
    \mapsto \ \left(\thom,\frac{-1}{x} + E - \lambda V(\thom)\right) \
  \end{equation}
defined on $\kreis \times \overline{\R}$.  In order to see this, note that, by
taking the inverse and replacing $\omega$ by $-\omega$ in (\ref{eq:riccati}),
we obtain the system
\begin{equation}
\thx \mapsto \left(\thom,\frac{-1}{x-E+\lambda V(\theta)} \right)
\end{equation}

Using the change of variables $x \mapsto \tan(-x)$, this yields
(\ref{eq:projharper}). The proof of Corollary~\ref{cor:harperI} below will
mainly consist in showing that there exists a parameter family of qpf interval
maps which satisfies the assumptions of Theorem \ref{thm:snaexistence}, such
that it exhibits a non-smooth saddle-node bifurcation, and which is at the
same time conjugated to (\ref{eq:riccati}), provided that both systems are
restricted to the relevant part of the phase space in which the bifurcation
takes place.

In order to proceed, we will now first take the rather atypical viewpoint of
fixing $E$ and considering $\lambda$ as the bifurcation parameter (whereas
usually in the study of Schr\"odinger cocycles the coupling constant $\lambda$
is fixed and the spectral parameter $E$ is varied). However, in this particular
situation the two viewpoints are actually equivalent and the analogous result
from the standard viewpoint can be recovered afterwards. More precisely, we
first show that Theorem~\ref{thm:snaexistence} implies the following:
\begin{cor} \label{cor:harperI}
  Suppose $\omega$ satisfies (\ref{cond:diophantine}) and the potential $V$ is
  non-negative, Lipschitz-continuous and decays linearly in a neighbourhood of
  its unique maximum. Then there exists a constant $E_0 = E_0(V,c,d)$ with the
  following property:

  If $E \geq E_0$, then there exists a unique parameter
  $\lambda_c=\lambda_c(E)$, such that for all $\lambda \in [0,\lambda_c]$
  there exist exactly two invariant graphs for the system
  (\ref{eq:projharper}) (and likewise for (\ref{eq:riccati})), one with positive and one with
  negative Lyapunov exponent. If $\lambda< \lambda_c$ then both these graphs
  are continuous, if $\lambda = \lambda_c$ they are non-continuous (i.e.\ a
  SNA and a SNR) and have the same topological closure.  Furthermore, the
  mapping $E \mapsto \lambda_c(E)$ is strictly monotonically increasing.
\end{cor}
Due to the monotone dependence of $\lambda_c(E)$ on $E$, Corollary
\ref{cor:harperI} immediately
implies
\begin{cor} \label{cor:harperII} 
  Suppose $\omega$, $V$ and $E_0$ are chosen as in Corollary~\ref{cor:harperI}
  and let $\lambda_0 := \lambda_c(E_0)$. Then the following holds:

  If $\lambda \geq \lambda_0$, then there exists a unique parameter
  $E_c=E_c(\lambda)\geq E_0$, such that for all $E\geq E_c$ there exist
  exactly two invariant graphs for the system (\ref{eq:projharper}) (and
  likewise for (\ref{eq:riccati})), one with positive and one with negative
  Lyapunov exponent. If $E>E_c$ then both these graphs are continuous, if $E =
  E_c$ they are non-continuous (i.e.\ a SNA and a SNR) and have the same
  topological closure.  The mapping $\lambda \mapsto E_c(\lambda)$ is the
  inverse of the mapping $E \mapsto \lambda_c(E)$.

\end{cor}
We remark that the Harper map can be viewed as a qpf circle homeomorphism (by
identifying $\R \cup \{\infty\}$ with $\kreis$). Since we do not want to
introduce rotation numbers for such systems here, we do not speak more
precisely about what happens if $E$ is decreased beyond $E_c$ (or $\lambda$ is
increased beyond $\lambda_c$) and just mention that in this case the rotation
number starts to increase and becomes non-zero for $E < E_c$. Invariant graphs
and even continuous invariant curves may exist in this situation, but they
will have a different homotopy type (i.e.\ they `wind around the torus' in the
vertical direction).

\ \\ \textit{Proof of Corollary~\ref{cor:harperI}.}  In the following we
always assume that the parameter $E$ is chosen sufficiently large, without
further mentioning. (In particular, most of the statements below are only true
for large $E$.)

As the two systems (\ref{eq:projharper}) and (\ref{eq:riccati}) are equivalent
(as mentioned above), it suffices to show that the statement is true for
(\ref{eq:riccati}).  Further, for the sake of simplicity we assume that $V$ is
normalised, i.e.\ $\sup_{\theta\in\kreis} V(\theta) = 1$. Let $\alpha :=
E^{3/2}$ and
\[
  F_1(x) \ := \ -1/x + E \ .
\]
Then $$F_1([1/E,2/E]) \ \supseteq \ [1/E,2/E] \ =: \ I_1$$ and $$F_1([3E/4,E]) \
\ssq \
[3E/4,E] \ =: \ I_2 \ .$$ Further $F_1$ is uniformly expanding on $I_1$ and
uniformly contracting on $I_2$. As $F_1$ is strictly concave on $(0,\infty)$,
it follows that $F_{1 \mid (0,\infty)}$ has exactly two fixed points $x_1 \in
I_1$ and $x_2 \in I_2$. \medskip

Let $s:= \frac{1+2/\sqrt{\alpha}}{x_2-x_1}$ and $h(x) := (x-x_1) \cdot
s$. Note that we have
\[
     s \ \in \ [1/E,2/E] \ .
\]
As
$h$ sends $x_1$ to 0 and $x_2$ to $1+2/\sqrt{\alpha}$, the map
\begin{equation} \label{eq:F2-def}
F_2(x) \ := \ h \circ F_1 \circ h^{-1}(x) \ = \ \frac{-s}{x/s+x_1} + s\cdot (E-x_1) 
\end{equation}
has fixed points $0$ and $1+2/\sqrt{\alpha}$. In addition, if $\gamma \in
(0,1)$ is fixed, then it is easy to check  that on the one hand 
$$F_2'(x) \ \in \ [1/4E^2,4/\gamma^2E^2] \ \ssq \ \left[\alpha^{-2},\halb\alpha^{-\halb}\right]
\quad \quad \forall x \in [\gamma,1+2/\sqrt{\alpha}] $$ and on the other hand
$$ F_2'(x) \ \in \ [E/16,E^2] \ \ssq \ [2\alpha^{\halb},\alpha^2] \quad
\quad \forall x \in [0,2/\alpha] \ .
$$ (Always assuming that $E$ is sufficiently large.)  Further, there holds
\[
F_2(1/\alpha) \ = \ \frac{-s}{1/\alpha s+x_1} + s\cdot (E-x_1) \ \geq \ -4/\sqrt{E} +
1-4/E^2 \ \stackrel{E \ra \infty}{\longrightarrow} \ 1 \ ,
\]
such that we can assume $F_2(1/\alpha) \geq 1-\gamma$. \medskip
 
 Due to the definition of $F_2$ in (\ref{eq:F2-def}) and as $h$ is affine with
slope $s$, the map $H:\thx \mapsto (\theta,h(x))$ smoothly conjugates
(\ref{eq:riccati}) with
\begin{equation}
  \thx \ \mapsto \ \left(\thom,F_2(x) - s\lambda V(\thom)\right) \ .
\end{equation}
Now we choose a ${\cal C}^1$-map $F : [-3,3] \ra [-\frac{3}{2},\frac{3}{2}]$,
such that $F_{|[0,1+2/\sqrt{\alpha}]} = F_{2\mid[0,1+2/\sqrt{\alpha}]}$ and
which satisfies the requirements
(\ref{cond:Ffixedpoints})--(\ref{cond:Fmapsover}). This is possible, since we
have shown above that $F_{2\mid[0,1+2/\sqrt{\alpha}]}$ has all the required
properties. In addition, $F$ can be chosen such that it is strictly concave on
$(0,3]$, has a unique fixed point $x_-$ in $[-3,0)$ and is uniformly
contracting on $[-3,x_-]$. Consequently, it satisfies the second alternative
of (\ref{cond:F3}). Further, if we let $g(\theta) := V(\thom)$ and $\beta :=
s\lambda$ and define
\[
T_\beta\thx \ = \ (\thom,F(x) - \beta g(\theta)
\]
 as in Theorem~\ref{thm:snaexistence}, then $H$ conjugates $T_\beta$
restricted to $ \kreis \times [0,1+2/\sqrt{\alpha}]$ and $S_{E,\lambda}$
restricted to $\kreis \times [x_1,x_2]$, where $S_{E,\lambda}$ denotes the map
given by (\ref{eq:riccati}).

The parameter family $T_\beta$ satisfies all the requirements of
Theorem~\ref{thm:snaexistence}. Hence, we obtain the existence of a critical
parameter $\beta_c$, which is the bifurcation parameter in a non-smooth
saddle-node bifurcation. Further, due to the monotonicity described in
Theorem~\ref{thm:saddlenode}(i), the $T_\beta$-invariant graphs $\psi$ and
$\varphi^+$ are always contained in \mbox{$\kreis \times
[0,1+2/\sqrt{\alpha}]$}. Consequently their preimages under $H$, which we
denote by $\widehat{\psi}$ and $\widehat{\varphi}^+$, are
$S_{E,\lambda}$-invariant and contained in $\kreis \times
(0,+\infty)$. Therefore, the parameter family $S_{E,\lambda}$ equally
undergoes a non-smooth saddle-node bifurcation with critical parameter
$\lambda_c = \beta_c/s$.\medskip

In order to complete the proof only two things remain to be shown: The
monotonicity of $E \mapsto \lambda_c(E)$ and the fact that for all $\lambda
\leq \lambda_c$ the two graphs $\widehat{\psi}$ and $\widehat{\varphi}^+$ are
indeed the only ones for the system $S_{E,\lambda}$. In order to see the
latter, we note that restricted to $[0,+\infty)$ all the fibre maps of
$S_{E,\lambda}$ are strictly concave, such that there can be only two
invariant graphs in $\kreis \times [0,+\infty)$. However, as $S_{E,\lambda}$
maps $\kreis \times [-\infty,0)$ into the forward invariant set
$[\varphi^+,+\infty)$, there cannot be any other invariant graphs in $\kreis
\times [-\infty,0)$ either. (Of course, the same conclusion also follows by
considering the associated $\textrm{SL}(2,\R)$-cocycle: Due to the non-zero
Lyapunov exponent there exists an invariant splitting into stable and unstable
subspaces. These correspond exactly to the two invariant graphs above and
there will be no others (compare Section~\ref{Harper}).)

In order to see the strict monotonicity of $E \mapsto \lambda_c(E)$, fix
$\epsilon > 0$ and suppose that $E_2 = E_1 + 3\epsilon$ and
$\lambda'<\lambda_c(E_1)+\epsilon$. Then $S_{E_1,\lambda_c(E_1)-\epsilon}$ has
a continuous invariant graph $\varphi^+_1$ contained in $\kreis \times
(0,+\infty)$, and some iterate of $S_{E_1,\lambda_c(E_1)-\epsilon}$ acts
uniformly contracting in the vertical direction on $[\varphi^+_1,\infty)$.
(This follows from the Uniform Ergodic Theorem in combination with the fact
that $DS_{E,\lambda,\theta}(x) = 1/x^2$ is decreasing in $x$.) However, since
\[
S_{E_2,\lambda',\theta}(x) \ \geq \
S_{E_1,\lambda_c(E_1)-\epsilon,\theta}(x) \quad \forall x \in
[\varphi^+_1,+\infty) \ ,
\] 
(recall that we assumed $V$ to be normalised) this implies that
$S_{E_2,\lambda'}$ maps $[\varphi^+_1,\infty)$ into itself and the respective
iterate of $S_{E_2,\lambda'}$ also acts uniformly contracting in the vertical
direction on this set. (Note that $DS_{E,\lambda,\theta}$ does not depend on
the parameters $E$ and $\lambda$.) Consequently $S_{E_2,\lambda'}$ has an
attracting and continuous invariant graph contained in
$[\varphi^+_1,+\infty)$. As this is true for all $\lambda' \leq
\lambda_c(E_1)+\epsilon$, this implies $\lambda_c(E_2) >
\lambda_c(E_1)+\epsilon$.

\qed

\subsubsection{Remarks on the qpf Arnold circle map and Pinched systems} \label{App:Remarks}

We neither apply the results from Section~\ref{NonSmoothSaddleNodes} to
Pinched skew products nor the the qpf Arnold map, but for very different
reasons. In the case of Pinched skew products, this would seem like using a
sledgehammer to crack a nut. In these systems the existence of SNA can be
established by a few short and elegant arguments, making use of their
particular structure (see \cite{grebogi/ott/pelikan/yorke:1984} and
\cite{keller:1996}). Even the exponential evolution of peaks can be described
in a much more simple way in this setting, a fact which was used in
\cite{jaeger:2004a} to study the topological structure of SNA in Pinched skew
products. (In fact, this preceding result and the striking similarities
between the pictures in Figures \ref{fig:atanfamily} and
\ref{fig:pinchedlines}, which strongly suggested some common underlying
pattern, were the starting point for the work presented here).  In principle
it is possible to view the SNA in these systems as
being created in non-smooth bifurcations, as this is done in
\cite{glendinning:2004}. However, as treating them with the methods presented
here would even need some additional modifications, we refrain from doing so.

For the case of the qpf Arnold circle map, the situation is completely
different. Here it is just not possible to apply our results. The reason for
this is the fact that no matter how the parameter $\alpha \in [0,1]$ in
(\ref{eq:qpfarnold}) is chosen, the maximal expansion rate is always at most
two. Further, for any interval of fixed length the uniform contraction rate
also remains bounded. Although the derivative goes to zero at $\theta=\halb$
if $\alpha$ is close to $1$, a strong contraction only takes place
locally. This means that the expansion and contraction rates one can work with
will always be moderate and cannot be chosen arbitrarily large by adjusting
the parameters. However, this is exactly what would be necessary for the
application of Theorem~\ref{thm:snaexistence}~.  In the case of the forcing
function $\theta \mapsto \sin(2\pi\theta)$ used in (\ref{eq:qpfarnold}), there
is also not much hope that a refinement of our methods would yield results. As
the simulations in Figure~\ref{fig:arnold} indicate, the exponential evolution
of peaks is only present in a very weak form in this case. Therefore, it
should be doubted that this process can be described in a rigorous way with
approximative methods as the ones we use in the proof of
Theorem~~\ref{thm:snaexistence}, which necessarily involve a lot of rough
estimates.

However, as already indicated in Section~\ref{ArnoldMap}, this might become
different if one chooses a more suitable forcing function, and considers for
example the parameter family
\begin{equation}
  \label{eq:arnoldmodified}
\thx \ \mapsto \ \left(\thom,x+\tau+\frac{\alpha}{2\pi} \sin(2\pi x) - \beta
\cdot \max\left\{0,1-\sigma \cdot d(\theta,0)\right\}\right) 
\end{equation}
with sufficiently large parameter $\sigma$. In this case the exponential
evolution of peaks is very distinct again, as one can see in
Figure~\ref{fig:arnoldpeak}. Consequently, it should also be possible to treat
this situation rigorously. Nevertheless, Theorem~\ref{thm:snaexistence} is not
sufficient for this purpose. Changing the forcing function does not have any
influence on the expansion and contraction rates, such that these will still
be too weak to meet our assumptions. Yet, there is an additional fact which we
do not make use of in the proof of Theorem~\ref{thm:snaexistence}: In the
situation of (\ref{eq:arnoldmodified}) with large $\sigma$, the forcing
function vanishes almost everywhere, apart from a small neighbourhood of
0. This means that after every visit in this neighbourhood, the expansion,
respectively contraction, has a long time to work, without any quasiperiodic
influence, before the next return. It seems reasonable to expect that this
could be used to make up for the weak expansion and contraction rates, for
example by regarding a renormalisation of the original system after a
sufficiently large finite time. However, the implementation of this idea is
left for the future \ldots\ .
 

\subsubsection{SNA's with symmetry} \label{SNAwithsymmetry}

Similar to the proof of Corollary~\ref{cor:atanfamily}, it is possible to show
that for sufficiently large parameters $\alpha$ the parameter family
(\ref{eq:symmetricfamily}) satisfies the assumptions of
Theorem~\ref{thm:symmetricsna}~. This leads to the following
\begin{cor}
  Suppose $\omega$ satisfies the Diophantine condition
  (\ref{cond:diophantine}). Then there exists $\alpha_0 = \alpha_0(c,d)$ such
  that for all $\alpha>\alpha_0$ there is a parameter $\beta_c =
  \beta_c(\alpha)$ such that the system (\ref{eq:symmetricfamily}) with
  parameters $\alpha$ and $\beta_c$ has two SNA and one SNR, with the
  properties described in Theorem~\ref{thm:symmetricsna}, and no other
  invariant graphs.
\end{cor}
As the details are more or less the same as in Section~\ref{App:atanfamily},
we omit the proof. 

To the knowledge of the author, this is the first situation where existence of
such a triple of intermingled invariant graphs can be described
rigorously. Similarly, it is the first example of a qpf monotone interval map without
continuous invariant graphs.


\section{Saddle-node bifurcations and sink-source-orbits} \label{Generalsetting}

The aim of this section is threefold: First, it is to introduce a
general setting where a (not necessarily non-smooth) saddle-node
bifurcation occurs and can be described rigorously. Secondly, we
will show that the presence of a \textit{`sink-source-orbit'} implies
the non-smoothness of the bifurcation, and how the existence of such
an orbit can be established by approximation with finite trajectories.
The construction of such trajectories with the required properties
will then be carried out in the succeeding Sections~\ref{Strategy} to
\ref{Construction}. Finally, before we can start we have to address a
subtle issue concerning the definition of invariant graphs:


\subsection{Equivalence classes of invariant graphs and the essential 
  closure} \label{Essentialclosure}

The problem we want to discuss is the following: Any invariant graph
$\varphi$ can be modified on a set of measure zero to yield another
invariant graph $\tilde{\varphi}$, equal to $\varphi$ $m$-a.s.\ (where
$m$ denotes the Lebesgue measure on $\kreis$). We usually do not want
to distinguish between such graphs.  On the other hand, especially
when topology is concerned we sometimes need objects which are
well-defined everywhere.  So far, this has not been a problem. The
bounding graphs of invariant sets defined by (\ref{eq:boundinggraphs})
are well-defined everywhere, and for the definition of the associated
measure (\ref{eq:associatedmeasure}) it does not matter. But in
general, some care has to be taken. We will therefore use the
following convention:

We will consider two invariant graphs as equivalent if they are $m$-a.s.\
equal and implicitly speak about equivalence classes of invariant graphs (just
as functions in ${\cal L}^{1}_{\textrm{Leb}}(\R)$ are identified if they are
Lebesgue-a.s.\ equal). Whenever any further assumptions about invariant graphs
such as continuity, semi-continuity or inequalities between invariant graphs
are made, we will understand it in the way that there is at least one
representative in each of the respective equivalence classes such that the
assumptions are met. All conclusions which are then drawn from the assumed
properties will be true for all such representatives.

There is one case where this terminology might cause confusion: It is
possible that an equivalence class contains both an upper and a lower
semi-continuous graph, but no continuous graph.%
\footnote{To get an idea of what could happen, consider the function
  $f : x \mapsto \sin \frac{1}{x} \ \forall x \neq 0$. By choosing
  $f(0) = 1$ we can extend it to an upper semi-continuous function, by
  choosing $f(0) = -1$ to a lower semi-continuous function, but there
  is no continuous function in the equivalence class.}
This rather degenerate case cannot occur when the Lyapunov exponent of
the invariant graph is negative (see \cite{stark:2003},
Proposition~4.1), but when the exponent is zero it must be taken into
account. To avoid ambiguities, we will explicitly mention this case
whenever it can occur.

In order to assign a well defined point set to an equivalence class of
invariant graphs, we introduce the \textit{essential closure}:

\begin{definition}                    \label{def:essclosure}
  Let $T$ be a qpf monotone interval map. If $\varphi$ is an invariant
  graph, we define its essential closure as
\begin{equation}                           \label{essclosure}
     \esscl{\Phi} := \{ (\theta,x) :
     \mu_{\varphi}(U) > 0 \ \forall \textit{open
       neighbourhoods U of } (\theta,x) \} \ ,
\end{equation}
where the associated measure $\mu_\varphi$ is given by
(\ref{eq:associatedmeasure}).

\end{definition}
Several facts follow immediately from this
definition:
\begin{itemize}
\item 
       \esscl{\Phi}\ is a compact set.
     \item $\esscl{\Phi}$ is equal to the topological support
       $\textrm{supp}(\mu_{\varphi})$ of the measure $\mu_\varphi$,
       which in turn implies \ $\mu_{\varphi}(\esscl{\Phi}) = 1$ (see
       e.g.  \cite{katok/hasselblatt:1995}).
\item
        Invariant graphs from the same equivalence class have the same
        essential closure (as they have the same associated measure).
      \item $\esscl{\Phi}$ is contained in every other compact set
        which contains $\mu_{\varphi}$-a.e.\ point of $\Phi$, in
        particular in $\overline{\Phi}$.
\item
        \esscl{\Phi}\ is forward invariant under $T$.%
        \footnote{This can be seen as follows: Suppose $x \in
          \esscl{\Phi}$ and $U$ is an open neighbourhood of $T(x)$.
          Then $T^{-1}(U)$ is an open neighbourhood of $x$, and
          therefore $\mu_{\varphi}(U) = \mu_{\varphi} \circ T^{-1}(U)
          > 0$. This means $T(x) \in \esscl{\Phi}$, and as $x \in
          \esscl{\Phi}$ was arbitrary we can conclude that
          $T(\esscl{\Phi}) \subseteq \esscl{\Phi}$. On the other hand
          $T(\esscl{\Phi})$ is a compact set which contains
          $\mu_{\varphi}$-a.e.\ point in $\Phi$, therefore
          $\esscl{\Phi} \subseteq T(\esscl{\Phi})$.  }
\end{itemize}


\subsection{Saddle-node bifurcations: Proof of Theorem~\ref{thm:saddlenode}}
\label{SaddleNode}

As mentioned, the first problem we have to deal with is to restrict the number
of invariant graphs which can occur. If there are too many, it will be hard to
describe a saddle-node bifurcation in detail. However, there is a result which
is very convenient in this situation:
\begin{thm}[Theorem 4.2 in \cite{jaeger:2003}]    \label{thm:schwarzian}
  Suppose $T$ is a qpf monotone interval map and all fibre maps $\Tth$ are
  ${\cal C}^3$. Further assume $\thx \mapsto D\Tth(x)$ is continuous and all
  fibre maps have strictly positive derivative and strictly negative
  Schwarzian derivative (see Footnote~\ref{foot:schwarzian}).  Then there are
  three possible cases:
\romanlist
\item There exists one invariant graph $\varphi$ with
  $\lambda(\varphi) \leq 0$.
\item There exist two invariant graphs $\varphi$ and $\psi$ with
  $\lambda(\varphi) < 0$ and $\lambda(\psi) = 0$.
\item There exist three invariant graphs $\varphi^- \leq \psi \leq
    \varphi^+$ with $\lambda(\varphi^-) < 0$, $\lambda(\psi)
    > 0$ and $\lambda(\varphi^+) < 0$. 
\listend
Regarding the topology of the invariant graphs, there are the
following possibilities: 
\begin{list}{(\roman{enumi})'}{\usecounter{enumi}}
\item   If the single invariant graph has negative
        Lyapunov exponent, it is always continuous. Otherwise the
        equivalence class contains at least an upper and a lower
        semi-continuous representative.
\item The upper invariant graph is upper semi-continuous, the
        lower invariant graph lower semi-continuous. If $\varphi$ is
        not continuous and $\psi$ (as an equivalence class) is only
        semi-continuous in one direction, then $\esscl{\Phi} =
        \esscl{\Psi}$.
\item $\psi$ is continuous if and only if $\varphi^+$ and
        $\varphi^-$ are continuous.  Otherwise $\varphi^-$ is at least
        lower semi-continuous and $\varphi^+$ is at least upper
        semi-continuous. If $\psi$ not lower semi-continuous then
        $\esscl{\Phi^-} = \esscl{\Psi}$, if $\psi$ is not upper
        semi-continuous then $\esscl{\Psi} = \esscl{\Phi^+}$.
\listend
Finally, as long as $\lambda(\varphi^-) < 0$ the graph $\psi$ can be defined by
\begin{equation}
  \label{eq:psidef}
  \psi(\theta) \ := \ \sup\{ x \in X \mid \lim_{n\ra\infty} |\Tth^n(x)
  - \varphi^-(\theta+n\omega)| = 0 \} \ .
\end{equation}
\end{thm}

In order to use the alternative assumption in (\ref{cond:F3}), we need a
similar result for concave fibre maps, which is due to Keller. The
main idea of the argument is contained in \cite{keller:1996}. However, as the statement
was never published in this form,  we include a proof.
\begin{thm}[G.\ Keller]
  \label{thm:concave} Suppose $T$ is a qpf monotone interval map, all fibre
  maps $T_\theta$ are differentiable and $\thx \mapsto DT_\theta(x)$ is
  continuous. Further, assume that there exist measurable functions
  $\gamma^\pm : \kreis \to X$, such that for all $\theta \in \kreis$ the fibre
  maps $T_\theta$ are strictly concave on $I(\theta)=
  [\gamma^-(\theta),\gamma^+(\theta)] \ssq X$. Then there exist at most two
  invariant graphs taking their values in $I(\theta)$, i.e.\ satisfying
 \begin{equation} \label{eq:concavity-interval} \varphi(\theta) \ \in \ I(\theta)
  \quad \forall \theta \in \kreis \ .
 \end{equation}
 If there exist two invariant graphs $\varphi_1\leq \varphi_2$ which both
 satisfy (\ref{eq:concavity-interval}), then $\lambda(\varphi_1) > 0$ and
 $\lambda(\varphi_2) < 0$. 

Further, if the graphs $\gamma^\pm$ are continuous and are mapped below
themselves, meaning that there holds
\begin{equation} 
T_\theta(\gamma^\pm(\theta)) \ \leq \  \gamma^\pm(\theta) \quad \forall \theta\in\kreis
 \ ,
\end{equation}
then either $\varphi_1,\varphi_2$ are both continuous, or $\varphi_1$ is lower
semi-continuous, $\varphi_2$ is upper semi-continuous and
$\esscl{\Phi_1}=\esscl{\Phi_2}$. (If there is only one invariant graph
which satisfies (\ref{eq:concavity-interval}), then it always
contains an upper and a lower semi-continuous representative in its
equivalence class.)
\end{thm}
\proof \ Suppose for a contradiction that there exist three different
invariant graphs $\varphi_1 \leq \varphi_2 \leq \varphi_3$ which all satisfy
(\ref{eq:concavity-interval}). As we identify invariant graphs which belong to
the same equivalence class, we have $\varphi_1(\theta) < \varphi_2(\theta) <
\varphi_3(\theta)$ $m$-almost surely. Due to the strict concavity of the fibre
maps and the invariance of the three graphs we obtain
\begin{equation}
  \label{eq:difference-quotients}
\log\left(\frac{\varphi_2(\thom)-\varphi_1(\thom)}{\varphi_2(\theta)-\varphi_1(\theta)}\right)
\ > \ \log\left(
\frac{\varphi_3(\thom)-\varphi_2(\thom)}{\varphi_3(\theta)-\varphi_2(\theta)} \right)
\quad m\textrm{-a.s.} \ .
\end{equation}
However, the following Lemma~\ref{lem:keller} applied to $Y=\kreis$,
$S(\theta)=\thom$, $\nu=m$ and $f = \log(\varphi_{i+1}-\varphi_i)$ ($i=1,2$)
yields that the integral with respect to $m$ on both sides equals zero, thus leading
to a contradiction.  Note that $f\circ S- f$ has the constant majorant
$\log(\max_{\thx \in \kreis \times X} DT_\theta(x))$.

\begin{lem}[Lemma 2 in \cite{keller:1996}] \label{lem:keller}
  Suppose $(Y,S,\nu)$ is a measure-preserving dynamical system, $f:Y \ra \R$
  is measurable and $f\circ S - f$ has an integrable majorant or
  minorant. Then $\int_Y f\circ S - f \ d\nu = 0$.
\end{lem}

For the estimates on the Lyapunov exponents, note that due to the strict
concavity there holds
\begin{eqnarray*}
\lambda(\varphi_1) & = & \int_{\kreis} \log\left( \lim_{t \to 0} 
\frac{T_\theta(\varphi_1(\theta)+t) - \varphi_1(\thom)}{t}\right) \ d\theta \\ & > &
\int_{\kreis}\log\left(
\frac{\varphi_2(\thom)-\varphi_1(\thom)}{\varphi_2(\theta)-\varphi_1(\theta)}\right)
\ d\theta \ = \ 0 \ .
\end{eqnarray*}
(The last equality follows again from Lemma~\ref{lem:keller}~.) Similarly, we
obtain $\lambda(\varphi_2) < 0$. \medskip

Now suppose $\gamma^+$ is continuous and $T_\theta(\gamma^+(\theta)) \leq
\gamma^+(\thom) \ \forall \theta \in \kreis$. Then we can define a sequence of
monotonically decreasing continuous curves by
\[
\gamma^+_n(\theta) \ := \ T^n_{\theta-n\omega}(\gamma^+(\theta-n\omega)) \ .
\]
As this sequence is bounded below by the invariant graph $\varphi_2$ it
converges pointwise, and the limit has to be an invariant graph. Since there
are no other invariant graphs between $\varphi_2$ and $\gamma^+$, we must have
$\varphi_2 = \lim_{n\to\infty} \gamma^+_n$. Consequently $\varphi_2$ is upper
semi-continuous as the monotone limit of a sequence of continuous curves. In
the same way one can see that $\varphi_1$ must be lower semi-continuous.

If $\varphi_1$ is not continuous, then the upper bounding graph of the compact
invariant set $\esscl{\Phi_1}$ must be an upper semi-continuous invariant
graph which lies between $\gamma^-$ and $\gamma^+$. The only candidate for
this is $\varphi_2$, such that $\Phi_2 \ssq \esscl{\Phi_1}$. However, this is
only possible if $\varphi_2$ is not continuous. Otherwise, as
$\lambda(\varphi_2) < 0$ and due to the Uniform Ergodic Theorem, some iterate
of $T$ would act uniformly contracting in the fibres on some neighbourhood $U$
of $\varphi_2$. In this case no other invariant graph could intersect $U$ on a
set of positive measure, contradicting $\Phi_2 \ssq \esscl{\Phi_1}$. Replacing
$T$ by $T^{-1}$ we can repeat the same argument for $\varphi_2$, such that
either both graphs are continuous or both are only semi-continuous and have the
same essential closure. This completes the proof.

\qed
\medskip

\begin{bem} \label{bem:schwarzianproof}
  \alphlist
\item The proof of Theorem~\ref{thm:schwarzian} in \cite{jaeger:2003}
  basically relies on the same idea as the above proof of
  Theorem~\ref{thm:concave}~. It depends on the fact that negative Schwarzian
  derivative of a ${\cal C}^3$-map $F: X \ra X$ is equivalent to strictly
  negative cross ratio distortion. The latter is defined as
  \[
    {\cal D}_F(w,x,y,z) \ = \ \frac{\frac{F(y)-F(y)}{y-x} \cdot
    \frac{F(z)-F(w)}{z-w}}{\frac{F(x)-F(w)}{x-w} \cdot \frac{F(z)-F(y)}{z-y}}
    \ ,
  \]
  where $w<x<y<z \in X$.  Applying the resulting inequality to four invariant
  graphs and integrating over the circle leads to a contradiction, similar to
  the argument after (\ref{eq:difference-quotients}). This excludes the
  existence of more than three invariant graphs in the situation of
  Theorem~\ref{thm:schwarzian}, and in more or less the same way one obtains
  the inequalities for the Lyapunov exponents.
\item We remark that the first part of Theorem~\ref{thm:schwarzian} (meaning
  statements (i)--(iii)) still holds if the dependence of $T_\theta$ on $\theta$ is only
  measurable, provided all other assumptions of the theorem are met and
  $\theta \mapsto \log(\max_{x\in X} DT_\theta(x))$ has an integrable majorant
  or minorant. Similarly, in Theorem~\ref{thm:concave} the statement about the
  number of the invariant graphs and the Lyapunov exponents remain true in the
  analogous case. 
\listend
\end{bem}

The preceding statements now allow to prove Theorem~\ref{thm:saddlenode}:
\bigskip

\textit{\bf Proof of Theorem~\ref{thm:saddlenode}~.} We start with the case
where all fibre maps have negative Schwarzian derivative
(see~(\ref{cond:F3})). Then due to Theorem~\ref{thm:schwarzian}, the number of
graphs which can exist is at most three . In order to show that the lower
bounding graph $\varphi^-$ is always continuous, let us first collect some
facts about the map $F$: As $F$ has three fixed points, there must exist some
$c\in[-2C,2C]$ with $F''(c) = 0$. However, the negative Schwarzian derivative
implies that $F'''(x) < 0$ whenever $F''(x) = 0$ for some $x \in [-2C,2C]$. Thus
there can be only one $c$ with $F''(c)=0$, and in addition $F''(x)$ will be
strictly positive for $x<c$ and strictly negative for $x>c$. Therefore
$F_{|[-2C,c)}$ will be strictly convex and $F_{|(c,2C]}$ strictly concave, and
this in turn implies that 0 is an unstable fixed point whereas $x^-$ and $x^+$
are stable. Further $F-\textrm{Id}$ is strictly positive on $(0,x^+)$ and
strictly negative on $(x^-,0)$, and finally $F$ is a uniform contraction on
$[-2C,x^-]$.%
\footnote{Note that we do not know whether $c \in (x^-,0]$, such that this
does not imply the second alternative in (\ref{cond:F3}).}

As we are in the case of one-sided forcing, for any %
$\epsilon$ with $-\epsilon \in (x^-,0)$ the set $\kreis \times
[-2C,-\epsilon]$ is mapped into itself, independent of $\beta$. Further,
as $g$ does not
vanish identically, there exist $\epsilon > 0$ and $n \in \N$ such that %
$T^n(M) \ssq \kreis \times [-2C,-\epsilon]$, where %
$M:=\kreis \times [-2C,0]$. Consequently  
\begin{eqnarray*}
\bigcap_{n\in\N} T^n(M) & \ssq 
& \bigcap_{n\in\N} T^n(\kreis \times
[-2C,-\epsilon]) \\
& \ssq & \bigcap_{n\in\N} \kreis \times
[-2C,F^n(-\epsilon)] \ = \ \kreis \times [-2C,x^-] \ =: \ N \ .
\end{eqnarray*}
Now $T$ acts uniformly contracting on $N$ in the vertical direction.
This means that there will be exactly one invariant graph contained in
$N \ssq M$, which is stable and continuous, and this is of course the
lower bounding graph $\varphi^-$. In particular $\varphi^- < 0$
independent of $\beta$. Furthermore, no other invariant graph can intersect
$N$. 
\romanlist
\item On the one hand, there obviously exist three invariant graphs at
  $\beta=0$, namely the constant lines corresponding to the three fixed
  points.  As these are not neutral, they will also persist for small values
  of $\beta$. On the other hand consider $\beta = C$. As we assumed
  that $g$ takes the maximum value of 1 at least for one $\theta_0 \in
  \kreis$, the point $(\theta_0,C)$ is mapped into $M$. (Recall that
  $F:[-2C,2C] \ra [-C,C]$.) But as we have seen, any point
  in $M$ is attracted to $\varphi^-$ independent of $\beta$. Thus there exists
  an orbit which starts above the upper bounding graph and ends up converging
  to $\varphi^-$. This means that there can be no other invariant graph apart
  from $\varphi^-$, and as this situation is stable the same will also be true for
  all $\beta$ sufficiently close to $C$. 
  
  Consequently, if we define $\beta_0$ as the infimum of all $\beta
  \in (0,C)$ for which there do not exist three continuous
  invariant graphs, then $\beta_0 \in (0,C)$ and statement
  (i) holds by definition.
  
  It remains to show that the graphs $\varphi^\pm$ and $\psi$ depend
  continuously and monotonically on $\beta$. Continuity simply follows from
  the fact that invariant curves with non-zero Lyapunov exponents depend
  continuously on ${\cal C}^1$-distortions of the system.  For the
  monotonicity of $\varphi^+$, note that since there is no other invariant
  graph above, $\varphi^+$ is the limit of the iterated upper boundary lines
  $\varphi_n$, which are defined by $\varphi_n(\theta) :=
  T_{\theta-n\omega}^n(2C)$. Due to the one-sided forcing, each of these
  curves will decrease monotonically as $\beta$ is increased, and this carries
  over to $\varphi^+$ in the limit. The same argument applies to $\varphi^-$,
  as this is the pointwise limit of the iterated \textit{lower} boundary
  lines.  Finally, note that $\psi$ can be defined as the upper boundary of
  the set
  \begin{eqnarray*}
     \lefteqn{\{ (\theta,x) \mid \lim_{n\ra\infty} |\Tth^n(x) -
    \varphi^-(\theta+n\omega)| = 0 \} } \\ & = & \{ \thx \mid \exists n\in\N :
    T^n\thx \in M \} \ .
  \end{eqnarray*}
  This set increases with $\beta$, and thus the graph $\psi$ will move
  upwards.
  
\item As all points in $M$ are attracted to $\varphi^-$, the two upper
  invariant graphs for $\beta<\beta_0$ must be contained in $M^c$.  Simply due
  to continuity, for $\beta \ra \beta_0$ the pointwise limits of these curves
  will be invariant graphs for $T_{\beta_0}$, although not necessarily
  continuous. By compactness, they will be contained in $\overline{M^c}$ and
  can therefore not coincide with $\varphi^-$.  Further, they cannot be both
  distinct and continuous: Due to the non-zero Lyapunov exponents given by
  Theorem~\ref{thm:schwarzian}(iii), this is a stable situation, contradicting
  the definition of $\beta_0$. Thus there only remain the two stated
  possibilities: Either the two graphs are distinct and not continuous, or
  they coincide $m$-a.s.\ and are neutral (see
  Theorem~\ref{thm:schwarzian}). The compactness of $B$ simply follows from
  the semi-continuity of the graphs $\psi$ and $\varphi^+$.
  
  In the case where $\psi$ equals $\varphi^+$ $m$-a.s., the fact that
  $B$ is pinched is obvious. Otherwise, it follows from
  Theorem~\ref{thm:schwarzian} that the two graphs have the same
  essential closure, which we denote by $A$. Now all invariant ergodic
  measures supported on $B$ (namely $\mu_\psi$ and $\mu_{\varphi^+}$)
  have the same topological closure $A$, which means that $A$ is
  minimal and there is no other minimal subset of $B$. Therefore
  Theorem~4.6 in \cite{stark:2003} implies that $B$ is pinched.
\item Suppose $\betil = \beta_0 + 2\epsilon$ for any $\epsilon > 0$.
  We have to show that there is no other invariant graph apart from
  the lower bounding graph $\varphi^-$. For this, it suffices to find
  an orbit which starts on the upper boundary line and ends up in $M$:
  This means that it finally converges to $\varphi^-$, which is
  impossible if there exists another invariant graph above.
  
  First, consider $\beta=\beta_0$ and let $\theta_1$ be chosen such
  that $\psi(\theta_1) = \varphi^+(\theta_1)$. As the pinched fibres
  are dense in $\kreis$ and $g(\theta_0)=1$, we can assume w.l.o.g.\ 
  that $g(\theta_1-\omega) \geq \halb$. Further, as the upper boundary lines
  converge pointwise to $\varphi^+$, there exists some $n\in\N$ such
  that
  \[
      \varphi_{n}(\theta_1) \  = \ T_{\beta_0,\theta_1-n\omega}^{n}(2C) 
     \ \leq \ \varphi^+(\theta_1) +
      \frac{\epsilon}{2} \ .
  \]
  Now, as the forcing is one-sided (i.e.\ $g\geq0$) we have
  $T^{n-1}_{\betil,\theta_1-n\omega}(2C) \leq
  T^{n-1}_{\beta_0,\theta_1-n\omega}(2C)$ and consequently
  \begin{eqnarray*}
      T^n_{\betil,\theta_1-n\omega}(2C) & = &
      T_{\betil,\theta_1-\omega}(T^{n-1}_{\betil,\theta_1-n\omega}(2C)) 
      \\ & \leq &
      T_{\betil,\theta_1-\omega}(T^{n-1}_{\beta_0,\theta_1-n\omega}(2C)) \\  
      & = & F(T^{n-1}_{\beta_0,\theta_1-n\omega}(2C)) - \betil \cdot
      g(\theta_1-\omega) \\
      & = & T^n_{\beta_0,\theta_1-n\omega}(2C) - (\betil - \beta_0)
      \cdot g(\theta_1-\omega) \\ 
      & \leq & \varphi^+(\theta_1) + \frac{\epsilon}{2} - \epsilon 
      \ < \ \varphi^+(\theta_1) \ = \
      \psi(\theta_1) \ .
  \end{eqnarray*}   
  However, already for $T_{\beta_0}$ the orbits of all points below
  $\psi$ eventually enter $M$, and again due to the one-sided nature
  of the forcing this will surely stay true for the respective orbits
  generated with $T_{\betil}$. Thus, for $\beta=\betil$ the orbit
  starting at $(\theta_1,2C)$ ends up in $M$ and therefore converges to
  the lower bounding graph. As $\epsilon > 0$ was arbitrary, this
  proves statement (iii).
\listend
Now assume the second alternative in (\ref{cond:F3}) holds, i.e.\ for some $c
\in (x^-,0]$ the map $F_{|[c,2C]}$ is strictly concave and $F_{[-2C,x_-]}$ is
uniformly contracting. Then the above proof basically remains the same, the
only difficulty is to see that for any $\beta \in [0,C]$ there cannot be more
than three invariant graphs. However, on the one hand it can be seen as above
that the lower bounding graph $\varphi^-$ is the only invariant graph in $M$
and no other invariant graph intersects $M$, since all orbits in this set
converge to $\varphi^-$. On the other hand we can apply
Theorem~\ref{thm:concave} with $I(\theta) = [0,2C]$ to see that there can be
at most two invariant graphs in $M^c$. 

Apart from this, the above arguments work in exactly the same way, replacing
Theorem~\ref{thm:schwarzian} by Theorem~\ref{thm:concave} where necessary. 

\qed
\medskip

\textit{\bf Proof of Lemma \ref{lem:parameterdependence}~.}  The continuity
  simply follows from the fact that both the situations above and below the
  bifurcation are stable, due to the non-zero Lyapunov
  exponents. Consequently, the sets $\{ (\alpha,\beta) \mid \beta <
  \beta_0(\alpha) \}$ and $\{(\alpha,\beta) \mid \beta> \beta_0(\alpha)\}$ are
  open, which means that $\alpha \mapsto \beta_0(\alpha)$ must be continuous.

In order to see the monotonicity, let $T_{\alpha,\beta}$ be the system
given by (\ref{eq:generalsystem}) with $F = F_\alpha$. Suppose that
$\tilde{\alpha} > \alpha$. Denote the upper bounding graph of the
system $T_{\alpha,\beta_0(\alpha)}$ by $\varphi^+$, the invariant
graph in the middle by $\psi$. As all points on or below the 0-line eventually
converge to the lower bounding graph (see the proof of
Theorem~\ref{thm:saddlenode}), the invariant graphs $\psi$ and
$\varphi^+$ must be strictly positive. As $\psi$ is lower
semi-continuous and $\varphi^+ \geq \psi$, both graphs are uniformly
bounded away from 0. Thus, there exists some $\delta > 0$ such that
$\delta \leq \varphi^+ \leq 1-\delta$.

For any $x \in [\delta,1-\delta]$ the map $F_\alpha(x)$ is
strictly increasing in $\alpha$.%
\footnote{We have 
\[
\frac{\partial}{\partial \alpha} F_\alpha(x) \ = \ 
\frac{\partial}{\partial \alpha} \left( \frac{\arctan(\alpha
    x)}{\arctan(\alpha)} \right) \ = \ \left( \frac{x\cdot
    \arctan(\alpha)}{1+\alpha^2x^2} - \frac{\arctan(\alpha
    x)}{1+\alpha^2} \right) \cdot \arctan(\alpha)^{-2} \ .
\]
This is positive if and only if
\[
   G_\alpha(x) \ := \ x \cdot \arctan(\alpha) \cdot
(1+\alpha^2) - \arctan(\alpha x) \cdot (1 + \alpha^2x^2)
\]
is positive. However, it is easy to verify that $G_\alpha(0) =
G_\alpha(1) = 0$ and $G_\alpha$ is strictly concave on $[0,1]$, i.e.\ 
$\frac{\partial^2}{\partial^2 x} G_\alpha(x) < 0 \ \forall x \in
[0,1]$, such that $G_\alpha(x) > 0 \ \forall x \in (0,1)$.  }
 Due to compactness this means that
there exists $\epsilon > 0$, such that $F_{\tilde{\alpha}} >
F_{\alpha}+\epsilon$ on $[\delta,1-\delta]$. Let $\betil :=
\beta_0(\alpha) + \epsilon$. Then
\[
    T_{\tilde{\alpha},\betil,\theta}(x) \ > \
    T_{\alpha,\beta_0(\alpha),\theta}(x) \ \ \ \ \ \forall \thx \in 
    \kreis \times [\delta,1-\delta] \ .
\]
Consequently $T_{\tilde{\alpha},\betil}$ maps the graph
$\varphi^+$ strictly above itself, which means that the upper bounding 
graph $\tilde{\varphi}^+$ of this system must be above
$\varphi^+$. It can therefore not coincide with the lower bounding
graph, which lies below the 0-line. Hence $\beta_0(\tilde{\alpha})
\geq \betil > \beta_0(\alpha)$.

\qed


\subsection{Sink-source-orbits and SNA: Proof of
  Theorem~\ref{thm:sinksourcesna}} \label{Sinksourceorbits} 

Suppose that $T$ satisfies the assumptions of Theorem~\ref{thm:sinksourcesna}
and denote the upper and lower bounding graph by $\varphi^+$ and $\varphi^-$,
respectively. Suppose there exists no non-continuous invariant graph with
negative Lyapunov exponent, but a point $(\theta_0,x_0) \in \kreis \times X$
with $\lambda^+\thxo > 0$ and $\lambda^-\thxo >0$ (i.e.\ a
sink-source-orbit). Let
\[
    \psi^+(\theta) \ := \ \inf \{ \varphi(\theta) \mid \varphi
    \textrm{ is a continuous } T \textit{-invariant graph with }
    \varphi(\theta_0) \geq x_0 \} \ ,
\]
with $\psi^+ :\equiv \varphi^+$ if no such graph $\varphi$
exists. Similarly, define
\[
    \psi^-(\theta) \ := \ \sup \{ \varphi(\theta) \mid \varphi
    \textrm{ is a continuous } T \textit{-invariant graph with }
    \varphi(\theta_0) \leq x_0 \} \ ,
\]
with $\psi^- :\equiv \varphi^-$ if there is no such graph $\varphi$.
By the continuity and monotonicity of the fibre maps, $\psi^+$ and
$\psi^-$ will be invariant graphs again. In addition, $\psi^+$ will be
upper and $\psi^-$ lower semi-continuous and $\psi^- \leq \psi^+$.
Thus, the set $A := [\psi^-,\psi^+]$ is compact. By a semi-uniform
ergodic theorem contained in \cite{sturman/stark:2000} (Theorem~1.9),
both $\lambda^+\thxo$ and $-\lambda^-\thxo$ must be contained in the
convex hull of the set
\[
    \left\{ \int_A \log D\Tthx \ d\mu\thx \mid \mu \textit{ is a }
    T_{|A}\textit{-invariant and ergodic probability measure} \right\} \ .
\]
As all ergodic measures are associated to invariant graphs (see
(\ref{eq:associatedmeasure})), this means that there must exist
invariant graphs with positive and negative Lyapunov exponents in $A$.
However, as we assumed that all stable invariant graphs are continuous
and there are no continuous invariant graphs contained in the interior
of $A$ by the definition of $\psi^\pm$, the only possible candidates
for a negative Lyapunov exponent are $\psi^+$ and $\psi^-$. We
consider the case where only $\lambda(\psi^-) < 0$, if $\psi^+$ or
both invariant graphs are stable this can be dealt with similarly.
Note that by the assumption we made at the beginning, the negative
Lyapunov exponent ensures that $\psi^-$ must be continuous.

Consequently, the convergence of the Lyapunov exponents is uniform on
$\psi^-$, such that there there is and open neighbourhood of this
curve which is uniformly contracted in the vertical direction by some
iterate of $T$. Therefore, if we define 
\[
\tilde{\psi}^-(\theta) \ := \ \inf\{x \geq \psi^-(\theta) \mid
\limsup_{n\ra \infty}
|\Tth^n(x) - \psi^-(\theta+n\omega)| > 0 \} \ .
\]
then $\tilde{\psi}^- > \psi^-$, and in addition $\tilde{\psi}^-$ is
lower semi-continuous. Note that 
\[
  \nLim |\Tth^n(x) - \psi^-(\theta+n\omega)| \ = \  0 \ \ \ \ \ 
  \forall \thx \in [\psi^-,\tilde{\psi}^- )
\]
by definition. The forward orbit of \thxo\ cannot converge to
$\psi^-$ as this contradicts $\lambda^+\thxo > 0$. Therefore $x_0 \geq
\tilde{\psi}^-(\theta_0)$. Further, there holds $\tilde{\psi}^- \leq
\psi^+$. This means that \thxo\ is contained in the compact set
$\tilde{A} := [\tilde{\psi^-},\psi^+]$. But as $\tilde{A}$ does not
contain an invariant graph with negative Lyapunov exponent anymore,
this contradicts $\lambda^-\thxo > 0$, again by Theorem~1.9 in
\cite{sturman/stark:2000}.

The existence of a strange non-chaotic repeller follows in the same way
by regarding the inverse of $T$ restricted to the global attractor.

\qed

%

\section{The strategy for the construction of the sink-source-orbits} \label{Strategy}

The inductive construction of longer and longer trajectories which are
expanding in the forwards and contracting in the backwards direction (compare
Lemma~\ref{lem:sinksourceshadowing}) will be a rather complicated inductive
procedure. On the one hand, a substantial amount of effort will have to be put
into introducing the right objects and providing a number of preliminary
estimates and technical statements in Section~\ref{Tools}.  On the other hand,
it will sometimes be quite hard to see the motivation for all this until the
actual construction is carried out in Section~\ref{Construction}. In order to
give some guidance to the reader in the meanwhile, we will try to sketch a
rough outline of the overall strategy in this section, and discuss at least
some of the main problems we will encounter. In particular, we will try to
indicate how a recursive structure appears in the construction, induced by the
recurrence behaviour of the underlying irrational rotation.

To this end, we will start by deriving some first (easy) estimates, which will
make it much easier to talk about what happens further. This will show that up
to a certain point the construction is absolutely straightforward. The further
strategy will then only be outlined, as the tools developed in Section
\ref{Tools} are needed before it can finally be converted into a rigorous
proof in Section~\ref{Construction}.


\subsection{The first stage of the construction}

As mentioned in Section~\ref{Mechanism}, for a suitable choice of the
functions $F$ and $g$ in (\ref{eq:generalsystem}) we can expect that
the tips of the peaks correspond to a sink-source-orbit.  However, as
we do not know the bifurcation parameter exactly, we can only
approximate it and show that in each step of the approximation there
is a longer finite trajectory with the required behaviour. The
existence of the sink-source-orbit at the bifurcation point will then
follow from Lemma~\ref{lem:sinksourceshadowing}~.

As we will concentrate only on trajectories in the orbit of the
0-fibre, the following notation will be very convenient:
\begin{definition} \label{def:parameterfamily}
For the map $T_\beta$ defined in Theorem~\ref{thm:snaexistence} with fibre
maps $T_{\beta,\theta}$, let
\[
T_{\beta,\theta,n} := T_{\beta,\theta+\omega_{n-1}} \circ \ldots \circ
T_{\beta,\theta} 
\]
if $n > 0$ and $T_{\beta,\theta,0} := \textrm{Id}$. Further, for any pair
$l \leq n$ of integers let
\[
   \xi_n(\beta,l) :=  T_{\beta,\omega_{-l},n+l}(3) \ .
\]
In other words, $\xi_n(\beta,l)$ is the $x$-value of that point from
the $T_\beta$-forward orbit of $(\omega_{-l},3)$, which lies on the
$\omn$-fibre. Thus, the lower index always indicates the fibre on
which the respective point is located. 
\end{definition}
Slightly abusing language, we
will refer to $(\xi_j(\beta,l))_{n\geq-l}$ as the forward orbit of the 
point $(\omega_{-l},3)$, suppressing the $\theta$-coordinates. 
\medskip

Note that under the assumptions of Theorem~\ref{thm:snaexistence} (which imply
in particular that we are in the case of one-sided forcing, i.e.\ $g\geq 0$)
the mapping $\beta \mapsto \xi_n(\beta,l)$ is monotonically decreasing for any
fixed numbers $l$ and $n$, with strict monotonicity if $l\geq 0$ and $n\geq 1$
since $g(0) =1$. In addition, we claim that when $n\geq 1$ and $l \geq 0$, the
interval $\Balphcl$ is covered as $\beta$ increases from $0$ to $\frac{3}{2}$,
i.e.\
\begin{equation} \label{eq:betacovers} \textstyle
    \xi_n(\frac{3}{2},l) \ < \ -\alphtel   \ .
\end{equation}
In order to see this, note that $\xi_0(\beta,l)$ is always smaller
than 3, such that $\xi_0(\beta,l)-x_\alpha \leq 2 -
\frac{2}{\walpha}$. Therefore, using $F(x_\alpha)=x_\alpha$,
(\ref{cond:Fcontraction}) and $g(0)=1$ we obtain
\[  \textstyle
\xi_1(\frac{3}{2},l) \ = \ F(\xi_0(\beta,l))-\frac{3}{2}\cdot g(0) \ 
\leq \ x_\alpha + \frac{2-\frac{2}{\walpha}}{2\walpha} - \frac{3}{2} \ 
= \ \frac{3}{\walpha} - \alphtel - \frac{1}{2} \ .
\]
By (\ref{cond:alphagamma0}) the
right side is smaller than $-\alphtel$, and %
as $\kreis \times [-3,-\alphtel)$ is always mapped into itself
this proves our claim.

\ \\
From now on, we use the following notation: For any pair $k,n$ of
integers with $k \leq n$ let
\begin{equation} \label{eq:kninterval}
    [k,n] \ := \ \{ k \ld n \} \ .
\end{equation}
What we want to derive is a statement of the following kind
\begin{quotation} \emph{
  If $\xi_N(\beta,l) \in \Balphcl$ for `suitable' integers $l \leq 0 $
  and $N \geq 1$, then $\xi_j(\beta,l) \in \Balphcl$ for `most' $j \in
  [1,N]$ and $\xi_j(\beta,l) \geq \gamma$ for `most' $j \in [-l,0]$. }
\end{quotation}
Of course, we have to specify what `suitable' and `most' mean, but as
this will be rather complicated we postpone it for a while. As
(\ref{eq:betacovers}) implies that there always exist values of $\beta
\in [0,\frac{3}{2}]$ with $\xi_n(\beta,l) \in \Balphcl$, such a
statement would ensure the existence of trajectories which spend most
of the backward time in the contracting region and most of the forward
time in the expanding region. This is exactly what is needed for the
application of Lemma~\ref{lem:sinksourceshadowing}~. As mentioned, up
to a certain point things are quite straightforward:

\begin{lem}
  \label{lem:inductionstart}
Suppose that the assumptions of Theorem~\ref{thm:snaexistence} hold. Further, let
$n\geq 1$, $l \geq 0$ and assume that
\[
d(\omj,0) \ \geq \ \frac{3\gamma}{L_2} \ \ \ \forall j \in [-l,-1]
\cup [1,n-1] \ 
.
\]
Then
$\xi_n(\beta,l) \in \Balphcl$ implies %
$\beta \in [1+\walphtel,1+\frac{3}{\walpha}]$,
\begin{equation}
    \xi_j(\beta,l) \ \in \ \Balphcl \hspace{3.8eM} \forall j \in [1,n] \ 
\end{equation}
and 
\begin{equation}
    \xi_j(\beta,l) \ \geq \ \gamma \hspace{6eM} \forall j \in [-l,0] \ 
    .
\end{equation}
\end{lem}
The proof relies on the following basic estimate:
\begin{lem}
  \label{lem:basicestimate}
  Suppose that the assumptions of Theorem~\ref{thm:snaexistence}
  hold. Further, assume that $\beta \leq 1 + \frac{4}{\walpha}$, $j \geq -l$
  and $d(\omj,0) \geq \frac{3\gamma}{L_2}$. Then $\xi_j(\beta,l) \geq
  \alphtel$ implies $\xi_{j+1}(\beta,l) \geq \gamma$ and $\xi_j(\beta,l) \leq
  - \alphtel$ implies $\xi_{j+1}(\beta,l) \leq -\gamma$. Consequently,
  $\xi_{j+1}(\beta,l) \in \Balphcl$ implies $\xi_j(\beta,l) \in \Balphcl$.
\end{lem}
\proof \
Suppose that $\xi_j(\beta,l) \geq \alphtel$. Using $d(\omj,0) \geq
\frac{3\gamma}{L_2}$ and (\ref{cond:sharppeak}) we obtain that
$g(\omj) \leq 1-3\gamma$. Therefore
\begin{eqnarray*}
  \xi_{j+1}(\beta,l) & = & F(\xi_j(\beta,l)) - \beta\cdot g(\omj) \\
   & \stackrel{(\ref{cond:Fmapsover})}{\geq} &  \textstyle 1 - \gamma -
   (1+\frac{4}{\walpha})(1-3\gamma) \ \geq \ 2\gamma -  
  \frac{4}{\walpha} \ \stackrel{(\ref{cond:alphagamma0})}{\geq} \
  \gamma \ .
  \end{eqnarray*}
  As $g \geq 0$, we also see that $\xi_j(\beta,l) \leq -\alphtel$
  implies 
  \[
  \xi_{j+1}(\beta,l) \ \leq \ F(\xi_j(\beta,l)) \
  \stackrel{(\ref{cond:Fmapsover})}{\leq} \
  -(1-\gamma) \ \stackrel{(\ref{cond:gamma0})}{\leq} \  -\gamma \ .
  \]

\qed

\ \\
\textit{Proof of Lemma \ref{lem:inductionstart}:} \\
Suppose that $\xi_n(\beta,l) \in \Balphcl$. We first show that $\beta
\leq 1+\frac{3}{\walpha}$: As $\xi_0(\beta,l) \leq 3$ we can use
$F(x_\alpha) = x_\alpha$ and (\ref{cond:Fcontraction}) to see that
$F(\xi_0(\beta,l)) \leq 1+\frac{3}{\walpha}-\alphtel$. As $g(0)=1$
this gives
\[
    \xi_1(\beta,l) \ = \ F(\xi_0(\beta,l)) - \beta \ \leq \
    \left(1+\frac{3}{\walpha} - \beta\right) - \alphtel \ .
\]
Thus, for $\beta > 1 + \frac{3}{\walpha}$ we have $\xi_1(\beta,l) <
-\alphtel$, and as $\kreis \times [-3,-\alphtel)$ is mapped into
itself this would yield $\xi_n(\beta,l) < - \alphtel$, contradicting
our assumption. Therefore $\xi_n(\beta,l) \in \Balphcl$ implies $\beta
\leq 1 + \frac{3}{\walpha}$.

Now we can apply Lemma~\ref{lem:basicestimate} to all $j \in [1,n-1]$
and obtain $\xi_j(\beta,l) \in \Balphcl \ \forall j \in [1,n]$ by
backwards induction on $j$, starting at $j=n$. Similarly,
$\xi_j(\beta,l) \geq \gamma \ \forall j = -l \ld 0$ follows from
$\xi_{-l}(\beta,l) = 3 \geq \gamma$ by forwards induction, as we can
again apply Lemma~\ref{lem:basicestimate} to all $j \in [-l,-1]$.

It remains to prove that $\beta \geq 1 + \walphtel$. We already showed
that $\xi_0(\beta,l) \geq \gamma \geq x_\alpha-1$, such that we can
use $F(x_\alpha) = x_\alpha$ and (\ref{cond:Fcontraction}) again to
see that
\[
\xi_1(\beta,l) \ \geq \ x_\alpha - \frac{1}{2\walpha} - \beta \ = \ 1
+ \frac{3}{2\walpha} - \beta \ 
\stackrel{(\ref{cond:alphagamma0})}{\geq} \ \left(1 + \walphtel -
  \beta \right) + \alphtel \ .
\]
As we also showed above that $\xi_1(\beta,l) \leq \alphtel$, the
required estimate follows.

\qed


\subsection{Dealing with the first close return} \label{Firstreturn}

As we have seen above, everything works fine as long as the $\omj$ do not
enter the interval $B_{\frac{3\gamma}{L_2}}(0)$ again. Thus, in the context of
Section~\ref{Mechanism} the critical region ${\cal C}$ corresponds to the
vertical strip $B_{\frac{3\gamma}{L_2}}(0) \times [-3,3]$. We will now sketch
the argument by which the construction can be continued even beyond the first
return to this critical region:

Suppose $m\in\N$ is the first time such that $d(\omm,0) <
\frac{3\gamma}{L_2}$ and fix some $l \leq m-1$.  Then
Lemma~\ref{lem:inductionstart} yields information up to time $m$,
meaning that we can apply it whenever $n \leq m$. But
we cannot ensure that $\xi_{m+1}(\beta,l) \in \Balphcl$ implies
$\xi_m(\beta,l)\in \Balphcl$ as before. In fact, this will surely be
wrong when $\omm$ is too close to 0, such that $g(\omm) \approx 1$.
In order to deal with this, we will define a certain `exceptional'
interval $J(m) = [ m-l^- \ld m + l^+ ]$. The integers $l^-$ and
$l^+$ will have to be chosen very carefully later on, but for now the
reader should just assume that they are quite small in comparison to
both $m$ and $l$.  Then, instead of showing that $\xi_{m+1}(\beta,l)
\in \Balphcl$ implies $\xi_m(\beta,l) \in \Balphcl$ as before, we will
prove that
\begin{equation}
  \label{eq:secondstage}
  \xi_{m+l^++1}(\beta,l) \in \Balphcl \ \ \textit{ implies } \ \
  \xi_{m-l^--1}(\beta,l) \in \Balphcl \ .
\end{equation}
Using Lemma~\ref{lem:inductionstart}, the latter then ensures that
$\xi_j(\beta,l) \in \Balphcl \ \forall j \in [1,m-l^--1]$.

Recall that as we are in the case of one-sided forcing, the dependence of
$\xi_n(\beta,l)$ on $\beta$ is strictly monotone. Thus, in order to prove
(\ref{eq:secondstage}), it will suffice to consider the two unique parameters
$\beta^+$ and $\beta^-$ which satisfy
\begin{equation}
  \label{eq:betaplusdef}
  \xi_{m-l^--1}(\beta^+,l) \ = \ \alphtel
\end{equation}
and 
\begin{equation}
  \label{eq:betaminusdef}
  \xi_{m-l^--1}(\beta^-,l) \ = \ -\alphtel \ .
\end{equation}
If we can then show the two inequalities
\begin{equation} \label{eq:betaplus}
  \xi_{m+l^++1}(\beta^+,l) \ > \ \alphtel
\end{equation}
and
\begin{equation} \label{eq:betaminus} 
  \xi_{m+l^++1}(\beta^-,l) \ < \ -\alphtel \ ,
\end{equation}
this immediately implies (\ref{eq:secondstage}). 

Now, first of all the fact that (\ref{eq:betaminus}) follows from
(\ref{eq:betaminusdef}) is obvious, as $\kreis \times [-3,-\alphtel]$
is mapped into $\kreis \times [-3,-(1-\gamma)]$ by
(\ref{cond:Fmapsover}), independent of the parameter $\beta$. Thus, it
remains to show (\ref{eq:betaplus}). This will be done by comparing the orbit%
\footnote{Recall that we suppress the $\theta$-coordinate $\omj$ of
  points $(\omj,\xi_j(\beta,l))$ from the forward orbit of
  $(\omega_{-l},3)$.}
\begin{equation}
\xi_{m-l^--1}(\beta^+,l) \ld \xi_{m+l^++1}(\beta^+,l)
\end{equation}
with suitable \textit{`reference orbits'}, on which information is
already available by Lemma~\ref{lem:inductionstart}~. In order to
make such comparison arguments precise (as sketched in
Figure~\ref{fig:firstreturn} below), we will need the following
concept:
\begin{definition}
   \label{def:errorterm}
   For any $\beta_1,\beta_2 \in [0,\frac{3}{2}]$ and
   $\theta_1,\theta_2 \in \kreis$, the \textbf{error term} is defined
   as
\begin{eqnarray*}
    \err(\beta_1,\beta_2,\theta_1,\theta_2) & := & \ \sup_{n \in \Z}
    |\beta_1\cdot g(\theta_1+\omn) -  \beta_2\cdot g(\theta_2+\omn)| \ 
    .
\end{eqnarray*}
Note that $\err(\beta_1,\beta_2,\theta_1,\theta_2) = 
    \sup_{n\in\Z} \supnorm{T_{\beta_1,\theta_1+\omn} -
      T_{\beta_2,\theta_2+\omn}}$. 
\end{definition}

The next remark gives a basic estimate:
\begin{bem}
  \label{bem:errorterm}
Suppose that $g$ has Lipschitz-constant $L_1$ (as in
(\ref{eq:g-lipschitz})). Further, assume that $\theta_1 = \omk,\ \theta_2 =
\omega_{k+m}$ for some $k,m\in \Z$, $d(\omm,0) \leq \frac{2\epsilon}{L_2}$,
and $\beta_1,\beta_2 \in [1,\frac{3}{2}]$ satisfy $|\beta_1-\beta_2| <
2\epsilon$. Then
\[ 
    \err(\beta_1,\beta_2,\theta_1,\theta_2) \leq K \cdot \epsilon
\] 
where $K := 3\cdot \frac{L_1}{L_2} + 2$. 
\end{bem}  
\proof \
For any $n\in \N$, let $j:= k+n$. Then $\omk+\omn = \omj$ and
$\omega_{k+m} + \omn = \omega_{j+m}$. Thus, the above estimate follows from
\begin{eqnarray*}   
    \lefteqn{|\beta_1\cdot g(\omj) - \beta_2 \cdot g(\omega_{j+m})| \leq} \\ & &  
    \beta_1\cdot |g(\omj)-g(\omega_{j+m})| +  
     g(\omega_{j+m})\cdot |\beta_1-\beta_2| \ \leq \ \beta_1  \cdot
     \frac{2\epsilon}{L_2} \cdot 
           L_1 + 2\epsilon \ \leq K \cdot \epsilon
\end{eqnarray*}    

\qed

\ \\ Thus, even if two finite trajectories are generated with slightly
different parameters and are not located on the same but only on nearby
fibres, the fibre maps which produce them will still be almost the same. This
makes it possible to compare two such orbits, at least up to a certain
extent. For the remainder of this section, the reader should just assume that
the remaining differences between the fibre maps can always be neglected. Of
course, when the construction is made rigorous later on it will be a main
issue to show that this is indeed the case.

Let us now turn to Figure~\ref{fig:firstreturn}, which illustrates the
argument used to derive (\ref{eq:betaplus}). The first
reference orbit, shown as crosses, is generated with the unique
parameter $\beta^*$ that satisfies $\xi_m(\beta^*,l) = 0$. Due to
Lemma~\ref{lem:inductionstart} (with $n=m$), we know that this orbit
always stays in the expanding region before, i.e.\ 
\begin{equation}
  \label{eq:reforbitone}
  \xi_j(\beta^*,l) \  \in \ \Balphcl \ \ \ \forall j = 1 \ld m-1 \ .
\end{equation}
Recall that $\beta^+$ was defined by $\xi_{m-l^--1}(\beta^+,l) =
\alphtel$. This implies $\xi_{m-l^-}(\beta^+,l) \geq \gamma$ by
Lemma~\ref{lem:basicestimate}~.  Thus, the `new' orbit
$\xi_{m-l^--1}(\beta^+,l),$ $\ldots,$ $\xi_{m+l^++1}(\beta^+,l)$
(corresponding to the black squares in Figure~\ref{fig:firstreturn})
leaves the expanding region and enters the contracting region (A),
whereas the reference orbit (crosses) stays in the expanding region at the same
time, i.e.\ $\xi_{m-l^-}(\beta^*,l) \in \Balphcl$, by
(\ref{eq:reforbitone}). 
\begin{figure}[t]
\noindent
\begin{minipage}[t]{\linewidth}
  \epsfig{file=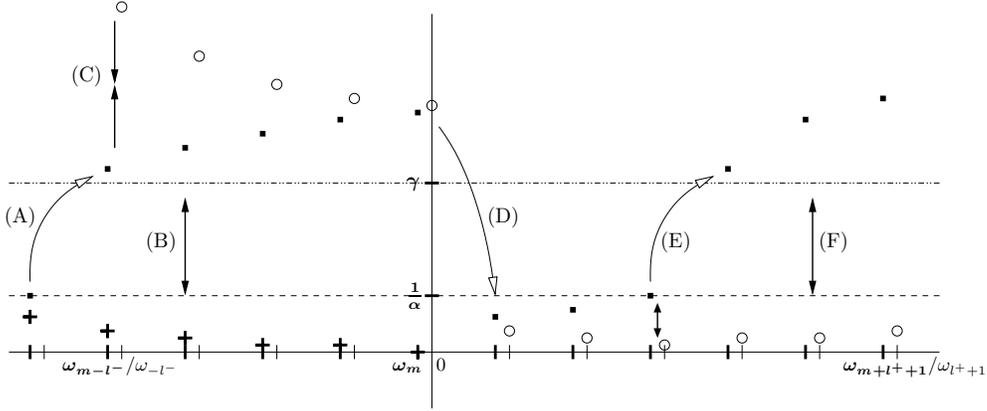, clip=, width=\linewidth,
    }
         \caption{\small The above diagram shows three finite
           trajectories: The `new' orbit $\xi_{m-l^--1}(\beta^+,l),$
           $\ldots,$ $\xi_{m+l^++1}(\beta^+,l)$ (black squares), the
           first reference orbit $\xi_{m-l^--1}(\beta^*,l),$ $\ldots,$
           $\xi_m(\beta^*,l)$ (crosses) and the second reference orbit
           $\xi_{-l^-}(\beta^+,l),$ $\ldots,$ $\xi_{l^++1}(\beta^+,l)$
           (circles). For convenience, successive iterates on the
           circle are drawn in straight order. (This corresponds to the
           situation where either the rotation number $\omega$ is very
           small, or where we consider a $q$-fold cover of the circle
           $\kreis$.)  After the first iterate, the new orbit leaves
           the expanding and enters the contracting region (A).
           Afterwards, the first reference orbit together with the
           strong expansion on $\kreis \times \overline{B_{\frac{2}{\alpha}}(0)}$
           ensure that the new orbit stays in the contracting region
           as the $\omm$-fibre is approached (B). Consequently, it
           gets attracted to the second reference orbit, which also
           lies in the contracting region (C). When the 0-fibre is
           passed, the forcing acts stronger on the second reference
           orbit (which passes exactly through the 0-fibre) than on
           the new orbit (which only passes through the $\omm$-fibre).
           Therefore, the new orbit will be slightly above the second
           reference orbit afterwards (D).  From now on, the expansion
           on $\kreis \times \overline{B_{\frac{2}{\alpha}}(0)}$ ensures that the
             new orbit eventually gets pushed out of the expanding
             region (E), and stays in the contracting region
             afterwards (F). }
         \label{fig:firstreturn}
\end{minipage}
\end{figure}
Afterwards, due to the strong expansion on
$\kreis \times \overline{B_\frac{2}{\alpha}(0)}$ it is not possible for the new
orbit to approach the reference orbit anymore, such that it will stay
`trapped' in the contracting region (B). In this way, we will
obtain%
\footnote{We should mention that in this particular situation
  (\ref{eq:trapped}) could still be derived directly from
  Lemma~\ref{lem:basicestimate} . However, the advantage of the
  described comparison argument is that it is more flexible and
  will also work for later stages of the construction.}
%
\begin{equation}
  \label{eq:trapped}
  \xi_{j}(\beta^+,l) \ \geq \ \gamma \ \ \ \forall j = m- l^- \ld m \ .
\end{equation}
Now we start to use a second reference orbit, namely
$\xi_{-l^-}(\beta^+,l),$ $\ldots,$ $\xi_{l^++1}(\beta^+,l)$, shown by
the circles in Figure~\ref{fig:mechanism}. Note that this time it will
be generated with exactly the same parameter $\beta^+$ as the new
orbit, but located on slightly different fibres.  By
Lemma~\ref{lem:inductionstart} (with $n=m-l^--1$, note that
$\xi_{m-l^--1}(\beta^+,l) \in \Balphcl$ by definition), we know that
\begin{equation}
  \label{eq:reforbittwoa}
  \xi_j(\beta^+,l) \ \geq \ \gamma \ \ \ \forall j = -l^- \ld 0 
\end{equation}
and
\begin{equation}
  \label{eq:reforbittwob}
  \xi_j(\beta^+,l) \ \in \ \Balphcl \ \ \ \forall j = 1 \ld l^++1 \ .
\end{equation}
Combining (\ref{eq:trapped}) and (\ref{eq:reforbittwoa}), we see that
the two orbits we want to compare both spend the first $l^-$ iterates
in the contracting region. Thus they are attracted to each other, and
consequently $|\xi_0(\beta^+,l) - \xi_m(\beta^+,l)|$ will be very small (C).
In fact, if $l^-$ has been chosen large enough, then this difference
will be of the same magnitude as $\epsilon := L_2 \cdot d(\omm,0)$,
i.e.\ 
\begin{equation} \label{eq:orbitcontraction}
  |\xi_0(\beta^+,l) - \xi_m(\beta^+,l)| \ \leq \ \kappa \cdot \epsilon
\end{equation}
for a suitable constant $\kappa > 0$. 

The next step is crucial: When going from $\xi_0(\beta^+,l)$ to
$\xi_1(\beta^+,l)$, the downward forcing takes its maximum (i.e.\ 
$g(0)=1$). In contrast to this, in the transition from
$\xi_m(\beta^+,l)$ to $\xi_{m+1}(\beta^+,l)$ the forcing function
$g(\omm)$ is only close to 1. More precisely, (\ref{cond:sharppeak})
yields $g(\omm) \leq 1-\epsilon$. Therefore
\begin{eqnarray*} \label{eq:crucialstep}
    \lefteqn{\xi_{m+1}(\beta^+,l) - \xi_1(\beta^+,l) \ \geq} \\
    & \geq & \beta^+ \cdot
    \epsilon  - |F(\xi_m(\beta^+,l)) - F(\xi_0(\beta^+,l))|  \
    \stackrel{(\ref{cond:Fcontraction})}{\geq}  \ \epsilon
    -\frac{\kappa\cdot \epsilon}{\walpha}  \ \geq \
    \frac{\epsilon}{2} \ ,    
\end{eqnarray*}
where we have assumed that $\walpha$ will be larger than $2\kappa$ and
$\beta^+ \geq 1$.  Thus, when the orbits pass the 0- and $\omm$-fibre,
respectively, a difference is created and the new orbit will be
slightly above the reference orbit afterwards (D). But from that point
on, the reference orbit stays in the expanding region by
(\ref{eq:reforbittwob}). Therefore, the small difference will be
expanded until finally the new orbit is `thrown out' upwards (E) and
gets trapped in the contracting region again (F). This will complete the
proof of (\ref{eq:betaplus}).

The crucial point now is the fact that the scheme in
Figure~\ref{fig:firstreturn} offers a lot of flexibility. We have
described the argument for the particular case of the first close
return, but in fact all close returns will be treated in a similar
way. The only difference will be the fact that the reference orbits we use
in the later stages of the construction may not stay in the expanding
(or respectively contracting) region all of the considered
times. However, this will still be true for most times, and that is
sufficient to ensure that on average the expansion (or contraction)
overweights and the new orbit shows the required behaviour.


\subsection{Admissible and regular times} \label{Admissiblepoints}

The picture we have drawn so far is already sufficient to motivate
some further terminology. As we have seen above, not all times
$N\in\N$ are suitable for the construction, in the sense of the
statement given below (\ref{eq:kninterval}). Thus, we will distinguish
between times which are \textit{`admissible'} and others which are
not. Only for admissible $N$ we will show that $\xi_N(\beta,l)\in
\Balphcl$ allows to draw conclusions about previous times $j < N$. To
be more precise, for any given admissible $N$ we will define a set %
$R_N \ssq [1,N]$ and show that $\xi_N(\beta,l) \in \Balphcl$ implies
$\xi_j(\beta,l) \in \Balphcl \ \forall j \in R_N$. The integers $j \in
R_N$ will then be called \textit{`regular with respect to N'}. The
precise definitions of admissible and regular times will be given in
Sections~\ref{Admissibletimes} and \ref{Regulartimes}~.

In order to give an example, consider the situation of the previous   
section: There, all points $N\leq m$ are admissible, and so is 
$m+l^++1$, but $m+1 \ld m+l^+$ are
not admissible. Further, for any $N \leq m$ we can choose %
$R_N = [1,N]$, and the set $R_{m+l^++1}$ contains at least all points from
$[1,m+l^++1] \setminus J(m)$. However, it will turn out that we have to
define even more times as regular w.r.t. $m+l^++1$, and thus derive
information about them, as this will be needed in the later stages of
the construction. Namely, the additional points we need to be regular
are $m+1 \ld m+l^+$. The reason why this is necessary is explained
in Section~\ref{Outline} and Figure~\ref{fig:rightside}. However, in this particular
situation it is not difficult to achieve this:

As $\omm$ is a close return, we can expect (and also ensure by using
the diophantine condition and suitable assumptions on $\gamma$) that
$\omega_{m+1} \ld \omega_{m+l^+}$ are rather far away from 0, in
particular not contained in $B_{\frac{3\gamma}{L_2}}(0)$. But this
means that we can apply Lemma~\ref{lem:basicestimate} to $m+1 \ld
m+l^+$ and obtain that $\xi_{m+l^++1}(\beta,l) \in \Balphcl$ implies
$\xi_j(\beta,l) \in \Balphcl \ \forall j = m+1 \ld m+l^+$ by backwards
induction on $j$.
Thus, if we divide the interval $J(m)$ into two parts %
$J^-(m) := [m-l^-,m]$ and $J^+(m) := [m+1,m+l^+]$, then we can also
define all points in the right part $J^+(m)$ as regular, such that %
$R_{m+l^++1} = [1,m+l^++1] \setminus J^-(m)$.

The reader should keep in mind that although most points will be both
regular and admissible, the difference between the two notions is
absolutely crucial. For example, for the argument in the previous
section it was vitally important that $m$ itself is admissible, as the 
first reference orbit ended exactly on the $\omm$-fibre. But on the
other hand, $m$ will not be regular w.r.t.\ any $N\geq m$, as it is a
close return itself and certainly contained in $J^-(m)$.


\subsection{Outline of the further strategy} \label{Outline}

For a certain while the arguments from Section~\ref{Firstreturn} will allow to
continue the construction as described. When there is another close return at
time $m'>m$ and $d(\omega_{m'},0)$ is approximately of the same size as
$d(\omm,0)$, then the diophantine condition will ensure that $m$ and $m'$ are
far apart. Thus, if we define an exceptional interval $J(m')$ again, this will
be far away from $J(m)$ and we can proceed more or less as before.  However,
we have also seen that the minimal lengths of $l^-$ and $l^+$ depend on how
close $\omm$ is to 0, as there must be enough time for the contraction to work
until (\ref{eq:orbitcontraction}) is ensured, and similarly for the expansion
until the new orbit is pushed out of the expanding region.  To be more
precise, let $p \in \N_0$ such that $\epsilon = L_2 \cdot d(\omm,0) \in
[\alpha^{-(p+1)},\alpha^{-p})$.  Then the minimal lengths of $l^-$ and $l^+$
will depend linearly on $p$, as the expansion and contraction rates are always
between $\alpha^{\pm\halb}$ and $\alpha^{\pm 2}$ by (\ref{cond:Fexpansion})
and (\ref{cond:Fcontraction}). Thus, at some stage we will encounter a close
return at time $\hat{m}$, for which the quantities $\hat{l}^-$ and $\hat{l}^+$
needed to define a suitable interval  $J(\hat{m}) =
[\hat{m}-\hat{l}^-,\hat{m}+\hat{l}^+]$ are larger than $l$ and $m$.
 
At first, assume that only $\hat{l}^+ > m$, whereas $\hat{l}^-$ is still
smaller that $l$. As mentioned, we will be able to show that 
\begin{equation}
  \xi_{\hat{m}+\hat{l}^++1}(\beta,l) \in \Balphcl \ \ \textit{ implies } \ \
  \xi_{\hat{m}-\hat{l}^--1}(\beta,l) \in \Balphcl
\end{equation}
by a slight modification of the argument sketched in
Figure~\ref{fig:firstreturn}~. In fact, for the left side there is no
difference: If $\beta^+$ and $\beta^*$ are again chosen such that
$\xi_{\hat{m}-\hat{l}^--1}(\beta^+,l) = \alphtel$ and
$\xi_{\hat{m}}(\beta^*,l) = 0$, then the first reference orbit
$\xi_{\hat{m}-\hat{l}^--1}(\beta^*,l),$ $\ldots,$
$\xi_{\hat{m}}(\beta^*,l)$ will again stay in the expanding region all
the time. Therefore we can use it to control the first part
$\xi_{\hat{m}-\hat{l}^--1}(\beta^+,l),$ $\ldots,$
$\xi_{\hat{m}}(\beta^+,l)$ of the new orbit as before, and conclude
that it always stays in the contraction region. As the same will be
true for the first part $\xi_{-\hat{l}^-}(\beta^+,l),$ $\ldots,$
$\xi_0(\beta^+,l)$ of the second reference orbit, the contraction
ensures again that $|\xi_{\hat{m}}(\beta^+,l) - \xi_0(\beta^+,l)|$ is
small enough (compare (\ref{eq:orbitcontraction})), and consequently
$\xi_{\hat{m}+1}(\beta^+,l)$ will be slightly above $\xi_1(\beta^+,l)$
after the 0-fibre is passed (compare (\ref{eq:crucialstep})).

\begin{figure}[t]
\noindent
\begin{minipage}[t]{\linewidth}
  \epsfig{file=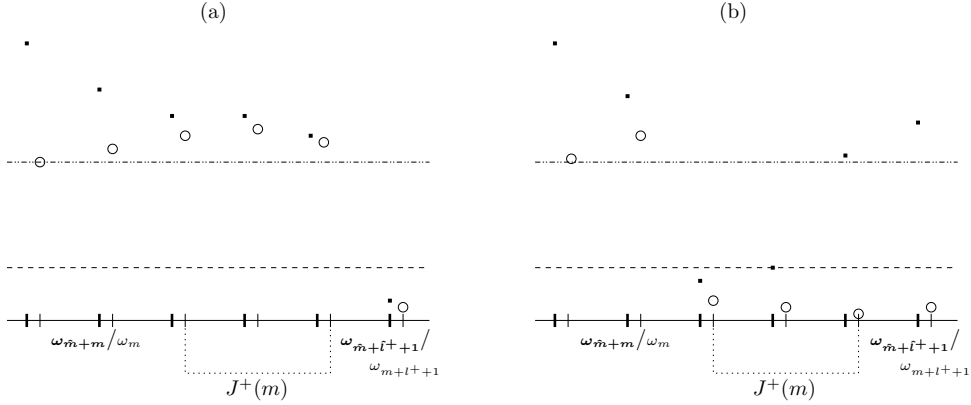, clip=, width=\linewidth,
    }
         \caption{\small In the above diagram, $J(m)$ is located at
           the end of $[1,\hat{l}^+]$, such that $m+l^+ = \hat{l}^+$.
           At first, the new orbit $\xi_{\hat{m}+1}(\beta^+,l),$
           $\ldots,$ $\xi_{\hat{m}+\hat{l}^++1}(\beta^+,l)$ will be
           pushed out of the expanding region (not shown). But at the
           end of the interval $[1,\hat{l}^+]$ the reference orbit
           $\xi_1(\beta^*,l),$ $\ldots,$
           $\xi_{\hat{l}^++1}(\beta^*,l)$ leaves the expanding region
           for a few iterates. Thus, the new orbit may approach the
           reference orbit during this time and enter the expanding
           region again afterwards. Now we consider two different
           situations: In (a) we assume that the reference orbit
           spends all times $j\in J(m)$ outside of the expanding
           region. This is what we have to take into account if we do
           not define the points in $J^+(m)$ as regular, and
           consequently do not derive any information about them. Then
           the new orbit may still be close to the reference orbit
           until the very last step, and thus lie in the expanding
           region at the end. (b) On the other hand, if we can obtain
           information about the $j\in J^+(m)$ and thus define them as
           regular, then we know that the reference orbit stays in the
           expanding region at these times. Therefore the new orbit
           may enter the expanding region after time $\hat{m}+m$, but
           it will be pushed out again before the end of the interval
           $J(\hat{m})$ is reached.}
         \label{fig:rightside}
\end{minipage}
\end{figure}

But afterwards, the second part $\xi_1(\beta^+,l) \ld
\xi_{\hat{l}^++1}(\beta^+,l)$ of the reference orbit will not stay in
the expanding region all the time, as the exceptional interval $J(m)$
is contained in $[1,\hat{l}^+]$ and the points in $J^-(m)$ will not be
regular w.r.t.\ $\hat{m}-\hat{l}^--1$. However, as all other points in
$[1,\hat{l}^+]$ are regular, it is still possible to show that the new
orbit is eventually pushed out of the expanding region again, but this
needs a little bit more care than before. Figure~\ref{fig:rightside}
shows one of the problems we will encounter, and thereby explains why
it is so vitally important that we have information about the points
in $J^+(m)$ as well, i.e.\ define them as regular before.

\begin{figure}[t]
\noindent
\begin{minipage}[h!]{\linewidth}
  \epsfig{file=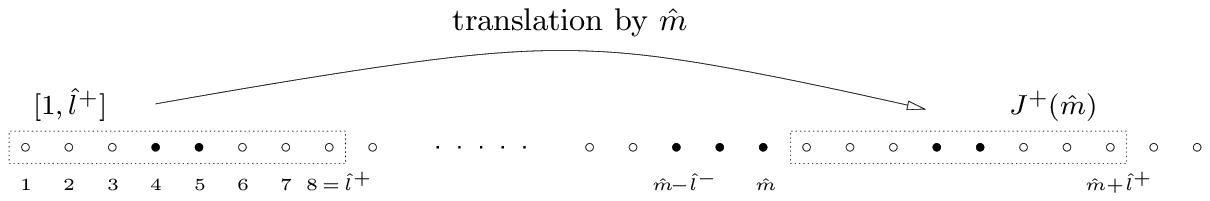, clip=, width=\linewidth,
    }
         \caption{\small Recursive structure of the sets
           $R_N$. Regular points are shown in white, exceptional ones
           in black. The set $R_N \cap J^+(\hat{m})$ is a
           translate of the set $R_N \cap [1,\hat{l}^+]$. }
         \label{fig:recursive}
\end{minipage}
\end{figure}

Now, we can begin to see how a recursive structure in the definition of
the sets $R_N$ appears: In order to have enough information for even
later stages in the construction, we will again have to define at
least most points in $J^+(\hat{m}) = [\hat{m}+1,\hat{m}+\hat{l}^+]$ as
regular. As it will turn out, we will be able to
show that $\xi_{\tilde{m}+\tilde{l}^++1} \in \Balphcl$ implies
$\xi_{\tilde{m}+j}(\beta,l) \in \Balphcl$ exactly whenever the
respective point $\xi_j(\beta^+,l)$ of the reference orbit lies in the
expanding region as well. In other words, a point $\hat{m}+j \in
J^+(\hat{m})$ will be regular if and only if %
$j \in [0,\hat{l}^+]$ was regular before. This leads to a kind of
self-similar structure in the sets of regular points, which will
express itself in relations of the following form:
\begin{equation} \label{eq:selfsimilarity}
     R_N \cap J^+(\hat{m}) \ = \left(R_N \cap [1,\hat{l}^+]\right) +
       \hat{m} \ = \ R_{\hat{l}^+} + \hat{m} 
\end{equation}
In other words, the structure of the sets $R_N$ after a close return,
i.e.\ in the right part $J^+$ of an exceptional interval, is the same
as their structure at the origin (see Figure~\ref{fig:recursive}).

\ \\
What remains is to extend the construction not only forwards, but also
backwards in time. As we have mentioned above, for some close return
$\tilde{m}$ we will eventually have to choose $\tilde{l}^-$ larger
than $l$. In this case, it is not sufficient anymore to have reference
orbits starting on the $\omega_{-l}$-fibre. However, we can still
carry out the construction exactly up to $\tilde{m}$. Thus, if
$\beta^*$ is chosen such that $\xi_{\tilde{m}}(\beta^*,l) = 0$, then
we will know that $\xi_{\tilde{m}-\tilde{l}^-}(\beta^*,l),$ $\ldots,$
$\xi_{\tilde{m}}(\beta^*,l)$ spends `most' of the time in the
expanding region. Therefore, we can use it as a reference orbit in
order to show that $\xi_{-\tilde{l}^-}(\beta,\tilde{l}^-),$ $\ldots,$
$\xi_0(\beta,\tilde{l}^-)$ stays in the contracting region `most' of
the time, at least for parameters $\beta$ which are close enough to
$\beta^*$. (Recall that this orbit starts on the upper boundary line,
i.e.\ $\xi_{-\tilde{l}^-}(\beta,\tilde{l}^-) = 3$ by definition.) It
will then turn out that it suffices to consider such parameter values.

In this way, the construction will be extended backwards and we can
then start to look at the forward part of the trajectories starting on
the $\omega_{-\tilde{l}}$-fibre.  Consequently, when we reach
$\tilde{m}$ again the backwards part of the trajectories is long
enough to carry on beyond this point, again using the same comparison
arguments as above. The only difference to
Figure~\ref{fig:firstreturn} will be that now the reference orbits
only stay most and not all of the time in the expanding or contracting
region, respectively. Nevertheless, this will still be sufficient to proceed
more or less in the same way. Hence we can continue the
construction, until we reach some even closer return. Then the
trajectories have to be extended further in the backwards direction
again and so on \ldots \ .


\section{Tools for the construction} \label{Tools}

In this section, we will provide the the necessary tools for the
construction of the sink-source-orbits in Sections~\ref{Construction}
and \ref{Symmetricsetting}.  As we have seen, there are mainly two
things which have to be done: First, we need some statements about the
comparison of orbits, namely one about expansion and one about
contraction. These will be derived in Section~\ref{Comparingorbits}.
Secondly, we have to define the sets of admissible and regular times,
which will be done in Sections~\ref{Admissibletimes} and
\ref{Regulartimes}~. However, before this we will have to introduce
yet another collection of sets $\Omega_p$ ($p\in\N_0$) in
Section~\ref{Approximatingsets}. These sets $\Omega_p$ will be used as
an approximation for the sets of \textit{non-regular} times and will
make it possible to control the frequency with which these can occur.


\subsection{Comparing orbits} \label{Comparingorbits}

The two statements we aim at proving here are Lemma
\ref{lem:orbitscontraction} and Lemma \ref{lem:orbitsthrowout}. They will
allow to compare two different orbit-segments which (i) start on nearby fibres
and (ii) result from systems $T_{\beta_1}, T_{\beta_2}$ with parameters
$\beta_1,\beta_2$ close together (compare Definition~\ref{def:errorterm} and
Remark~\ref{bem:errorterm}). The reader should note that throughout this
wubsection we only use assumptions (\ref{cond:gamma0}),(\ref{cond:alphagamma0}),%
(\ref{cond:Funiformbounds})--(\ref{cond:Fcontraction}) and the
Lipschitz-continuity of $g$. In particular, we neither use the fact that $g$
is non-negative, nor (\ref{cond:sharppeak}). Therefore, we will also be able
to use the results for the case of symmetric forcing in
Section~\ref{Symmetricsetting}. The diophantine condition on $\omega$ as well
as (\ref{cond:Fmapsover}) and (\ref{cond:sharppeak}) will not be needed until
the next section. Before we start, we make one more assumption on the
parameter $\alpha$: We suppose that $K$ is chosen as in Remark
\ref{bem:errorterm} and assume 
\begin{equation} 
   \label{cond:alpha1}
    \walpha \ \geq \ 2K \ .
\end{equation}
The following notation is tailored to our purpose of
comparing two orbits:
\begin{definition} \label{def:orbits}
Suppose $T_\beta$ is defined as in Theorem~\ref{thm:snaexistence}~.  If
  $\theta_1,\theta_2 \in \kreis$, $x^1_1,x^2_1 \in [-3,3]$ and
  $\beta_1,\beta_2 \in [0,\frac{3}{2}]$ are given, let
  \begin{equation}
    \label{eq:orbits}
    x^1_n \ := \ T_{\beta_1,\theta_1,n-1}(x^1_1) \ \ \ , \ \ \ x^2_n
    \ :=  \ T_{\beta_2,\theta_2,n-1}(x^2_1)
  \end{equation}
and 
\begin{equation}
  \label{eq:taun}
  \tau(n) \ := \ \# \{ j \in [1,n] \mid x^1_j \notin \Balphcl \} \ .
\end{equation}

\end{definition}
We start with a lemma about orbit-contraction. Essentially, the
statement is that if two orbits spend most of the time in the
contracting region above the line $\kreis \times \{\gamma\}$, then
their distance in the vertical direction gets contracted up to the
magnitude of the error term:
\begin{lem}
   \label{lem:orbitscontraction}
   Suppose than conditions (\ref{cond:Funiformbounds}) and
   (\ref{cond:Fcontraction}) hold and
   $$\err(\beta_1,\beta_2,\theta_1,\theta_2) \leq K\cdot \epsilon$$ for some
   $\epsilon> 0$. Let further
\begin{equation}
\eta(k,n) := \# \{ 
   j \in [k,n] \mid x^1_j \textrm{ or } x^2_j < \gamma \} 
\end{equation}
and assume that $\eta(j,n) \leq
\frac{n+1-j}{10} \ \forall j = 1 \ld n$ and $\alpha^{-\frac{n}{4}} \leq
\epsilon$. Then
\begin{equation}
    |x^1_{n+1} -x^2_{n+1}| \ \leq \ \epsilon \cdot \left( 6 + K \cdot 
     \sum_{j=0}^\infty \alpha^{-\viertel j} \right) \ .
\end{equation}
A similar statement holds for $\tilde{\eta}(k,n) := \# \{ j \in [k,n]
\mid x^1_j \textrm{ or } x^2_j > -\gamma \}$.
\end{lem} 
\proof\  
We prove the following statement by backwards induction on $k$:
For all $k = 1 \ld n+1$ there holds
\begin{eqnarray} \nonumber
    | x^1_{n+1}-x^2_{n+1}| & \leq & \label{eq:contractionformula} |x^1_k - x^2_k| \cdot
      \alpha^{-\halb(n+1-k-5\eta(k,n))} \\ &&+ \ K\cdot
      \epsilon \cdot  \sum_{j=k+1}^{n+1} 
      \alpha^{-\halb(n+1-j-5\eta(j,n))}  \quad
\end{eqnarray}
The case $k=n+1$ is obvious. For the induction step, first suppose
$x^1_k$ or $x^2_k < \gamma$, such that $\eta(k,n) = \eta(k+1,n) + 1$.
Then, by (\ref{cond:Funiformbounds}) we have
\[
|x^1_{k+1} - x^2_{k+1}| \leq |x^1_k-x^2_k| \cdot \alpha^2 + K\cdot
\epsilon \ ,
\]
and by applying the statement for $k+1$ we get
\begin{eqnarray*}
    \lefteqn{|x^1_{n+1}-x^2_{n+1}| \ \leq } \\
     &\leq & (|x^1_k - x^2_k| \cdot \alpha^2 + K \cdot \epsilon)
     \cdot \alpha^{-\halb(n-k-5\eta(k+1,n))} + K \cdot \epsilon  \sum_{j=k+2}^{n+1}
     \alpha^{-\halb(n+1-j-5\eta(j,n))} \\ 
     & = & |x^1_k - x^2_k| \cdot \alpha^{-\halb(n+1-k-5\eta(k,n))} +
     K \cdot \epsilon \cdot \sum_{j=k+1}^{n+1} \alpha^{-\halb(n+1-j-5\eta(j,n))} \ .
\end{eqnarray*}
On the other hand, suppose $x^1_k,x^2_k \geq \gamma$, such that
$\eta(k,n) = \eta(k+1,n)$. In this case we can use
(\ref{cond:Fcontraction}) to obtain
\[
     |x^1_{k+1} - x^2_{k+1}| \leq
     |x^1_k-x^2_k| \cdot  \alpha^{-\halb}+ K \cdot \epsilon
\]
and thus 
\begin{eqnarray*}
    |x^1_{n+1}-x^2_{n+1}| & \leq & (|x^1_k - x^2_k| \cdot \alpha^{-\halb} +
     K\cdot \epsilon) \cdot \alpha^{-\halb(n-k-5\eta(k+1,n))} \\ &&  + \ K\cdot
     \epsilon \cdot \sum_{j=k+2}^{n+1} \alpha^{-\halb(n+1-j-5\eta(j,n))} \\ &
     = & |x^1_k - x^2_k| \cdot \alpha^{\halb(n+1-k-5\eta(k,n))} \\ && + \ K \cdot
     \epsilon \cdot \sum_{j=k+1}^{n+1} \alpha^{-\halb(n+1-j-5\eta(j,n))} \ .
\end{eqnarray*}
The statement of the lemma is now just an application of
(\ref{eq:contractionformula}).  Note that $|x^1_1-x^2_1|$ is always
bounded by 6.

\qed

\ \\
The result about orbit-expansion we will need is a little bit more
intricate. The problem is the following: We have one reference orbit,
which spends most of the time well inside of the expanding region
$\kreis \times \overline{\Bdalph}$. A second orbit starts a certain distance
above, and we want to conclude that at some point it has to leave the
expanding region while the first orbit remains inside at the same
time. The following case is still quite simple:
\begin{lem}
  \label{lem:orbitsbasic}
  Suppose that conditions (\ref{cond:alphagamma0}), (\ref{cond:Fexpansion})
  and (\ref{cond:alpha1}) hold and further
  $$\err(\beta_1,\beta_2,\theta_1,\theta_2) \leq K \cdot \alpha^{-1}$$ and
  $x^2_1 \geq x^1_1 + \alphtel$.  Then as long as $\tau(n) = 0$ there holds $
  x^2_{n+1} \geq x^1_{n+1} + \frac{3}{\alpha}$. Thus $x^2_{n+1} \geq
  \frac{2}{\alpha}$ if $x^1_{n+1} \in \Balphcl$. A similar statement holds if
  $x^2_1 \leq x^1_1 - \alphtel$.
\end{lem}
\proof\ 
This follows from
      \[ \textstyle
           x^2_{n+1} \stackrel{(\ref{cond:Fexpansion})}{\geq} x^1_{n+1} +
           2\walpha\cdot \alphtel - K \cdot \alpha^{-1} \geq x^1_{n+1} +
           \alphtel(2\walpha - K) \stackrel{(\ref{cond:alpha1})}{\geq}
           x^1_{n+1} + \walphtel
      \]
      as long as $x^2_n - x^1_n \geq \alphtel$ and $x^1_n \in \Balphcl$.
      Note that $\walphtel \geq \frac{3}{\alpha}$ by
      (\ref{cond:alphagamma0}).

\qed

\ \\
However, it is not always that easy, because we also need to address
the case where the first orbit does not stay in the contracting region
all but only `most' of the times. This needs a little bit more care,
and there are some natural limits: For example, $x^1_j$ must not spend
to many iterates in the contracting region, even if these only make up
a very small proportion of the length of the whole orbit segment.
Otherwise the vertical distance between the two orbits may be
contracted until it is of the same magnitude of the error term, and
then the order of the orbits might get reversed. Another requirement
is that $x^1_j$ does not leave the expanding region too often towards
the end of the considered time interval. The reason for this was
already demonstrated in Figure~\ref{fig:rightside}~.

In the end we aim at proving Lemma \ref{lem:orbitsthrowout}, which is
the statement that will be used later on. However, in order to do so we
need two intermediate lemmas first.

\begin{lem}
  \label{lem:orbitsexpansion}
  Suppose that conditions (\ref{cond:Funiformbounds}), (\ref{cond:Fexpansion})
  and (\ref{cond:alpha1}) hold and further
  $$\err(\beta_1,\beta_2,\theta_1,\theta_2) \leq K\cdot\epsilon$$ with $\epsilon
  \leq \alpha^{-q}$ for some $q \geq 1$ and
\begin{equation} 
    x^2_1 \ \geq \ x^1_1 + \frac{\epsilon}{2} \cdot \alpha^r
\end{equation}
with $0 \leq r < q$. Suppose further that for all $j = 1 \ld n$ there
holds 
\begin{equation}
    x^1_j \in \Balphcl \ \follows \ x^2_j \in \overline{\Bdalph}
\end{equation}
and 
\begin{equation} \textstyle
  \label{eq:taucond1}
               r + \halb(j - 5\tau(j) ) \ \geq \ \halb \ .
\end{equation}
Then
\begin{equation} \label{eq:orbitexpansion} 
    x^2_{n+1} \ \geq \ x^1_{n+1} + \frac{\epsilon}{2} \cdot
    \alpha^{r+\halb(n-5\tau(n))} \ .
\end{equation}
A similar statement holds if $x^2_1 \leq x^1_q -
\frac{\epsilon}{2}\cdot \alpha^r$.  
\end{lem}
Note that (\ref{eq:taucond1}) is always guaranteed if either $\tau(n) \leq
\max\{0,\frac{2r-1}{4}\}$ (as $5\tau(j) -j \leq 4\tau(j) \leq 4\tau(n)$), or
if $\tau(j) \leq \frac{j}{8} \ \forall j = 1 \ld n$.

\ \\
\proof\ We prove
(\ref{eq:orbitexpansion}) by induction on $n$. The case $n=0$ is obvious. For
the induction step, we have to distinguish two cases:

\ \\
\underline{Case 1:} \hspace{3eM} $x^1_n \in \Balphcl$ , \ \ i.e.\ $\tau(n) =
      \tau(n-1)$ \\
      \begin{eqnarray*}
           x^2_{n+1} & \stackrel{(\ref{cond:Fexpansion})}{\geq} &
           x^1_{n+1} + 2\walpha \cdot \frac{\epsilon}{2} \cdot 
           \alpha^{r+\halb(n-1-5\tau(n-1))} - K\cdot \epsilon \\
           & = & x^1_{n+1} + \epsilon \cdot \left(
             \alpha^{r+\halb(n-5\tau(n))} - K \right) \
           \geq  \
           x^1_{n+1} + \frac{\epsilon}{2} \cdot
           \alpha^{r+\halb(n-5\tau(n))} 
      \end{eqnarray*}
      where we used $\alpha^{r+\halb(n-5\tau(n))} \geq \sqrt{\alpha}
      \geq 2K$ by (\ref{eq:taucond1}) and (\ref{cond:alpha1}) in the
      last step.
    
     \ \\
     \underline{Case 2:} \hspace{3eM} $x^1_n \notin \Balphcl$ , \ \
     i.e.\ $\tau(n) = \tau(n-1) + 1$  
     \begin{eqnarray*}
          x^2_{n+1} & \stackrel{(\ref{cond:Funiformbounds})}{\geq}  &
          x^1_{n+1} + 2 \alpha^{-2} \cdot
          \frac{\epsilon}{2} \cdot \alpha^{r+\halb(n-1-5\tau(n-1))} -
          K \cdot \epsilon \\
          & = & x^1_{n+1} + \epsilon \cdot \left(
             \alpha^{r+\halb(n-5\tau(n))} - K \right) \
           \geq  \
           x^1_{n+1} + \frac{\epsilon}{2} \cdot
           \alpha^{r+\halb(n-5\tau(n))} 
     \end{eqnarray*}
     where we used $\alpha^{r+\halb(n-5\tau(n))} \geq 2K$ again in the
     step to the last line.

\qed
     
\begin{lem}
   \label{lem:orbitstimebelow}
Suppose that conditions (\ref{cond:alphagamma0}), (\ref{cond:Fexpansion}) and
  (\ref{cond:alpha1}) hold and $$\err(\beta_1,\beta_2,\theta_1,\theta_2) \leq K
  \cdot \alpha^{-q}$$ for some $q \geq 1$. Further, assume that $x^1_1,
  x^1_{n+1} \in \Balphcl, \ x^2_1 \geq \frac{2}{\alpha}$ and $\tau(n) \leq
  \max\{0,\frac{2q-3}{4}\}$. Then $x^2_j \geq x^1_j \ \forall j = 1 \ld n$ and
  there holds
\[
\# \underbrace{\{ j \in [2,n+1] \mid x^1_j \notin \Balphcl \textrm{ or
  } x^2_j \in \overline{\Bdalph} \}}_{{\displaystyle =: \Upsilon}} \ 
\leq 5\tau(n) \ .
\]
A similar statement holds if $x^2_1 \leq -\frac{2}{\alpha}$. 
\end{lem}
\proof\   
It suffices to obtain a suitable upper bound on $\# \tilde{\Upsilon}$ where
\[
    \tilde{\Upsilon} := \{ j \in [2,n+1] \mid x^1_j \in \Balphcl
    \textrm{ and } x^2_j \in \overline{\Bdalph} \} \ ,
\]
as obviously $\# \Upsilon = \# \tilde{\Upsilon} + \tau(n+1) = \#
\tilde{\Upsilon} + \tau(n)$. (Note that $\tau(n+1) = \tau(n)$ as
$x^1_{n+1} \in \Balphcl$ by assumption.) To that end, we can separately
consider blocks $[k+1,l]$ where $k,l$ are chosen such that
\romanlist
  \item $1 \leq k < l \leq n+1$
  \item $x^1_j \in \Balphcl \follows x^2_j \in \overline{\Bdalph} \ \forall j \in [k+1,l]$
  \item $x^1_k \in \Balphcl$ and  $x^2_k \notin \overline{\Bdalph}$
  \item $x^1_l \in \Balphcl$ and $x^2_l \in \overline{\Bdalph}$
  \item $l$ is the maximal integer in $[k+1,n+1]$ with the
        above properties (ii) and (iv).  
\listend
Note that $\tilde{\Upsilon}$ is contained in the disjoint union of all such
blocks $[k+1,l]$. 

\ \\
We now want to apply Lemma~\ref{lem:orbitsexpansion}, but starting
with $x_k^i$ instead of $x^i_1\ (i=1,2)$. Therefore, let
$\tilde{\theta_i}=\theta_i+\omega_{k-1}$, \ %
$\tilde{x}^i_1 = x^i_k $ and $\tilde{n}=l-k$ in Definition
\ref{def:orbits}. Note that %
$\tilde{\tau}(\tilde{n}) = \tau(l-1) - \tau(k-1)$, but as we assumed
that $x^1_k,x^1_l \in \Balphcl$ in (iii) and (iv) we equally have %
$\tilde{\tau}(\tilde{n}) = \tau(l)-\tau(k)$. As %
$x^2_k \geq x^1_k + \alphtel$ by (iii), we can apply
Lemma~\ref{lem:orbitsexpansion} with %
$\epsilon = \alpha^{-q}$ and $r = q-1$ to obtain
\[   %
    x^2_{l} \ = \ \tilde{x}^2_{\tilde{n}+1} \ \geq \ \tilde{x}^1_{\tilde{n}+1} +
    \frac{\alpha^{-1}}{2} \cdot \alpha^{\halb(\tilde{n} -
      \tilde{\tau}(\tilde{n}))} \ = \ 
    \ x^1_l + \frac{\alpha^{-1}}{2} \cdot
    \alpha^{^\halb(l-k-5(\tau(l)-\tau(k)))} \ .  
\]
As $|x^2_l - x^1_l| \leq \frac{3}{\alpha} < \frac{1}{2\walpha}$ by (iv) and
(\ref{cond:alphagamma0}), we
must therefore have $l-k-5(\tau(l)-\tau(k)) \leq 0$ or equivalently
$l-k \ \leq \ 5(\tau(l)-\tau(k))$. Thus
\[
\# \left(\tilde{\Upsilon} \cap [k+1,l]\right) \ = \ l - k - (\tau(l) -
\tau(k)) \ \leq \ 4(\tau(l)-\tau(k)) \ .
\]
Summing over all such blocks $[k+1,l]$ we obtain $\# \tilde{\Upsilon}
\leq 4\tau(n)$, and this completes the proof.

\qed

\begin{lem}
  \label{lem:orbitsthrowout}
  Suppose that conditions (\ref{cond:alphagamma0}), (\ref{cond:Funiformbounds}),
  (\ref{cond:Fexpansion}) and (\ref{cond:alpha1}) hold and
  $$\err(\beta_1,\beta_2,\theta_1,\theta_2) \leq K \cdot \epsilon$$ for some
  $\epsilon \in [\alpha^{-(q+1)},\alpha^{-q})$, $ q \geq 1$. Further, assume
  that for some $n \in \N$ with $x^1_{n+1} \in \Balphcl$ there holds $\tau(n)
  \leq \max\{0,\frac{2q-3}{4}\}$ and
  \begin{equation} \label{eq:orbitsthrowout}
  \tau(n)-\tau(j) \leq \frac{n-j}{6} \ \forall j \in [1,n] \ .
  \end{equation}
  Then
  \alphlist
  \item $x^1_1 \in \Balphcl$ but $x^2_1 \geq \frac{2}{\alpha}$ implies
      $x^2_{n+1} \geq \frac{2}{\alpha}$.
  \item If $n \geq 5q$ and $\tau(j) \leq \frac{j}{8} \ \forall j=1 \ld n$,
      then $x^2_1 \geq x^1_1 + \frac{\epsilon}{2}$ implies $x^2_{n+1}
      \geq \frac{2}{\alpha}$.
  \listend
  Again, similar statements hold for the reverse inequalities.
\end{lem}
\proof\ 
\alphlist
\item Note that $\tau(1) = 0$ as $x^1_1 \in \Balphcl$ by assumption. By
  Lemma \ref{lem:orbitstimebelow} we have
  \[ 
     \# \Upsilon \ \leq \  5\tau(n) \
     \stackrel{(\ref{eq:orbitsthrowout})}{\leq} \ \frac{5(n-1)}{6} \ \leq \ n-1 \ ,
  \]
  Thus there exists $j_0 \in [2,n+1]$ such that $x^1_{j_0} \in \Balphcl$
  but $x^2_{j_0} \geq \frac{2}{\alpha}$.

  If we shift the starting points in Definition \ref{def:orbits} to
  $\tilde{\theta}_i := \theta_i +\omega_{j_0-1}$ and $\tilde{x}_1^i =
  x^i_{j_0}  \ (i=1,2)$ and denote the resulting sequences by
  $\tilde{x}^1_j, \tilde{x}^2_j$, then $\tilde{n} := n - j_0 +1$
  satisfies the same assumptions as before. As $\tilde{n} < n$ we can
  repeat this procedure until $\tilde{n} < 6$. But then
  $\tilde{\tau}(\tilde{n}) = 0$, such that $\tilde{x}^1_1 \in \Balphcl$
  and $\tilde{x}^2_1 \geq \frac{2}{\alpha}$ implies $\tilde{x}^2_{\tilde{n}+1}
  = x^2_{n+1} \geq \frac{2}{\alpha}$ by Lemma \ref{lem:orbitsbasic},
  proving statement (a).
\item
   We claim that there exists $j_1 \in [1,n+1]$ such that $x^1_{j_1} \in
  \Balphcl$ but $x^2_{j_1} \notin \overline{\Bdalph}$.
  
  Suppose there exists no such $j_1$ and let $k$ be the largest
  integer in $[1,n]$ such that $x^1_{k+1} \in \Balphcl$. As $x^1_j \in
  \Balphcl \follows x^2_j \in \overline{\Bdalph}$ holds for all $j=1 \ld k$, we can
  apply Lemma \ref{lem:orbitsexpansion} with $r=0$ to obtain
  \begin{equation} \label{eq:throwout}
    x^2_{k+1} \ \geq  \ x^1_{k+1} + \frac{\epsilon}{2}
        \cdot \alpha^{\halb(k-5\tau(k))} 
        \ \geq \ x^1_{k+1} + \halb \alpha^{\halb(k-5\tau(k))-q-1} \ .
  \end{equation}
  Now $\tau(n) = \tau(k+1) + n-k-1$ by definition of $k$. Further
  $\tau(k) = \tau(k+1)$, as $x^1_{k+1} \in \Balphcl$ by the choice of
  $k$. Therefore
  \begin{eqnarray*} \textstyle
     \lefteqn{ \halb(k-5\tau(k)) \ =}   \\
   & = & \ \halb(k+1-5\tau(k+1)) - \halb  
      \ \geq \ \textstyle
      \halb(n-5\tau(n)) - \halb \\
      & \geq &  \halb(5q-5\cdot \frac{2q-2}{4}) - \halb \ = \ \frac{5}{2}q -
    \frac{5}{4}q + \frac{12}{4} - \halb \  > \ q \ .
  \end{eqnarray*}
  Plugged into (\ref{eq:throwout}) this yields $x^2_{k+1} \geq
  x^1_{k+1} + \halb$, contradicting $x^1_{k+1} \in \Balphcl$ and
  $x^2_{k+1} \in \overline{\Bdalph}$.
  
  Thus we can choose $j_1$ with $x^1_{j_1} \in \Balphcl$ and $x^2_{j_1}
  \geq \frac{2}{\alpha}$ as claimed. Shifting the starting points as
  before we can now apply (a) to complete the proof.
\listend

\qed


\subsection{Approximating sets} \label{Approximatingsets}

As mentioned in Section~\ref{Strategy}, for each close return $m \in
\N$ with $d(\omm,0) \leq \frac{3\gamma}{L_2}$ we will introduce an
exceptional interval $J(m)$. However, before we can do so we first
have to define some intermediate intervals $\Omega_p(m)$. These will
contain the intervals $J(m)$, such that they can be used to obtain
estimates on the `density' of the union of exceptional intervals. As
we need a certain amount of flexibility, we have to introduce a whole
sequence of such approximating sets $(\Omega_p(m))_{p \in \N_0}$,
which will be increasing in $p$.

The statements of this as well as the two following subsections do not involve
the dynamics of $T$, they are only related to the underlying irrational
rotation by $\omega$. Therefore, the only assumptions which are used are the
Diophantine condition~(\ref{cond:diophantine}) as well as (\ref{cond:gamma0})
and (\ref{cond:alphagamma0}).
\begin{definition}
   \label{def:omegasets}
\alphlist
\item
Let 
\[
    S_n(\alpha) \ := \ \left\{ \begin{array}{cll} \sum_{i=0}^{n-1}
      \alpha^{-i} & \textrm{if} & n \in \N \cup \{ \infty \} \\ \\
        1 &  \textrm{if} & n \leq 0 \end{array} \right.  \ .
\]
\item
For $p \in \N_0 \cup \{\infty\}$ let $Q_p : \Z \ra \N_0$ be defined by
\[
    Q_p(j)  :=  \left\{ \begin{array}{cll} q  & \mathrm{if} &
    d(\omj,0) \in \left[ S_{p-q+1}(\alpha) \cdot
    \frac{\alpha^{-q}}{L_2},S_{p-q+2}(\alpha)\cdot \frac{\alpha^{-(q-1)}}{L_2}
    \right) \ \textrm{for  } q \geq 2
    \\ \\
     1 & \mathrm{if} & d(\omj,0) \in \left[ S_p(\alpha) \cdot 
    \frac{\alpha^{-1}}{L_2},  \frac{4\gamma}{L_2} +
    S_{p}(\alpha)\cdot \frac{\alpha^{-1}}{L_2} \cdot (1-\ind_{\{0\}}(p))  \right)
     
     \\ \\
     0 & \mathrm{if} & d(\omj,0) \geq \frac{4\gamma}{L_2} +
    S_{p}(\alpha)\cdot \frac{\alpha^{-1}}{L_2} \cdot (1-\ind_{\{0\}}(p)) 
     \end{array} \right. 
\]
if $j \in \Z \setminus \{ 0 \}$ and $Q_p(0) := 0$. Further let
\[
        p(j) \ := \ Q_0(j) \ .
\]
\item
For fixed $u,v \in \N$ let $\util:=u+2$ and $\vtil:=v+2$. Then, for any $j \in \Z$ define
\[
     \Omega^-_p(j) := [j-\util\cdot Q_p(j),j] \ , \ \ \Omega^+_p(j)
     :=  [j+1,j+\vtil\cdot Q_p(j)] 
\]
and 
\[ 
        \Omega_p(j)  :=  \Omega^-_p(j) \cup \Omega^+_p(j) 
\] 
if $Q_p(j) > 0$, with all sets being defined as empty if $Q_p(j) =
0$. Further let 
\[
    \Omega^{(\pm)}_p := \bigcup_{j \in \Z} \Omega^{(\pm)}_p(j) \ \ \ \
    \textrm{ and } \ \ \ \ 
    \Omtil^{(\pm)}_p := \! \! \! \!\!\twoscriptcup{j \in \Z}{Q_p(j) \leq p}
    \Omega_p^{(\pm)}(j) \ .
\]
\item
Finally,  let 

\[
    \begin{array}{llll} \textstyle
    \nu(q) & := & \min\left\{ j \in \N \mid p(j) \geq q \right\}
    & \forall q \in \N \\ \\
    \nutil(q) & := & 
    \min\left\{ j \in \N \mid d(\omj,0) < 
    3S_\infty(\alpha)\cdot \frac{\alpha^{-(q-1)}}{L_2} \right\} & \textrm{if } q \geq 2
    \   \textrm{ and } \\ \\
    \nutil(1) & := & \min\left\{ j
    \in \N \mid d(\omj,0) < 3\left(\frac{4\gamma}{L_2} +
    S_\infty(\alpha)\cdot\frac{\alpha^{-1}}{L_2} \right) \right\}  .
\end{array}
\]
\listend

\end{definition}

\begin{bem}
   \label{bem:alphagamma2}
Suppose that (\ref{cond:gamma0}) and (\ref{cond:alphagamma0}) hold, such that
we have $\walpha \geq \frac{4}{\gamma} \geq 64$. As $S_\infty(\alpha) =
\frac{\alpha}{\alpha - 1}$, the following estimates can be deduced easily from
this:
\begin{eqnarray}
    \label{cond:alpha2}
     \alpha & \geq & S_\infty(\alpha) + 1 \\
     \gamma & \geq & \frac{S_\infty(\alpha)+1}{\alpha}
     \label{cond:alphagamma1} \ .
\end{eqnarray}
\end{bem}

\begin{bem}
As in the prededing remark, we suppose that (\ref{cond:gamma0}) and
(\ref{cond:alphagamma0}) hold.
  \label{bem:omegasets}
\alphlist
\item By definition, we have $Q_{p'}(j) \leq Q_p(j) \ \forall j \in
  \Z$ whenever $p' \leq p $.  Further, there holds $ Q_\infty(j) \leq
  p(j) + 1 \ \forall j \in \N \ .  $ For $p(j) \geq 1$ this follows
  from (\ref{cond:alpha2}), which implies
  $\frac{S_\infty(\alpha)}{\alpha} \leq 1$. In the case $p(j) = 0$
  this is true by (\ref{cond:alphagamma1}). Altogether, this yields
  \begin{equation} \label{eq:pjestimates}
     p(j) \ \leq \ Q_p(j) \ \leq \ Q_\infty(j) \ \leq \ p(j) + 1 \ \ \ \forall j
     \in \Z,\ p \in \N
  \end{equation}
\item As a direct consequence of (a) we have $\Omega^{(\pm)}_{p'}(j)
  \subseteq \Omega^{(\pm)}_p(j) \ \forall j \in \N$ whenever $p' \leq
  p$. The same holds for the sets $\Omega^{(\pm)}_p$ and
  $\Omtil^{(\pm)}_p$.  
\listend
\end{bem}
The following two lemmas collect a few basic properties of the sets
$\Omega^{(\pm)}_p$ and $\Omtil^{(\pm)}_{p}$. The first one is a
certain `almost invariance' property under translations with $m$ if
$\omm$ is close to $0$. This is closely related to the recursive
structure of the sets $R_N$ of regular points mentioned in the last
section (see (\ref{eq:selfsimilarity})), and explains why we had to
introduce a whole family $(\Omega_p)_{p\in\N_0}$ of approximating
sets. 

Lemma \ref{lem:omegaestimates} then contains the estimates which can
be obtained from the diophantine condition. These will allow us to
control the ``density'' the sets of $\Omega^{(\pm)}_\infty$ (and thus
of the sets $R_N$ defined later on) by making suitable assumptions on
the parameters.
\begin{lem}
  \label{lem:omegatransition} Suppose that conditions (\ref{cond:gamma0}) and
(\ref{cond:alphagamma0}) hold.  Let $p \geq 2$ and suppose $p(m) \geq p$ and
  $Q_{p-2}(k) \leq p-2$. Then
\alphlist
\item \ \  $Q_{p-2}(k) \ \leq \ Q_{p-1}(k \pm m) \ \leq \ Q_{p-2}(k) + 1$
\item \ \ $\Omtil_{p-2}^{(\pm)} \pm m \ \subseteq \
  \Omtil_{p-1}^{(\pm)}$. Using $\Omtil_{-1} := \emptyset = \Omtil_0$,
  this also holds if $p=1$. 
\listend
\end{lem}
\proof\ 
\alphlist
\item 
     Let $q := Q_{p-2}(k)$, so that $p-q \geq 2$ by assumption. We first
     show that 
     \begin{equation} \label{eq:omegatransitionA}
         Q_{p-1}(k+m) \ \geq \ q \ .
     \end{equation}
     To that end, first suppose $q \geq 2$, such that
     $d(\omk,0) < S_{p-q}(\alpha) \cdot \frac{\alpha^{-(q-1)}}{L_2}$. Then
     \begin{eqnarray*}
           \lefteqn{d(\omega_{k \pm m},0) \ \leq \ d(\omk,0) + d(\omm,0) 
           \ < \ S_{p-q}(\alpha)\cdot\frac{\alpha^{-(q-1)}}{L_2} +
             \frac{\alpha^{-(p-1)}}{L_2} }
           \\
           & = & \left(S_{p-q}(\alpha) +
             \alpha^{-(p-q)}\right)\cdot \frac{\alpha^{-(q-1)}}{L_2} \ = \
           S_{p-q+1}(\alpha)\cdot\frac{\alpha^{-(q-1)}}{L_2} \ . 
     \end{eqnarray*}
     This proves (\ref{eq:omegatransitionA}) in case $q \geq 1$. The case $q=1$ is treated
     similarly, if $q=0$ there is nothing to show.

     It remains to prove that 
     \begin{equation} \label{eq:omegatransitionB}
             Q_{p-1}(k+m) \ \leq \ q + 1 \ .
     \end{equation}
     This time, first assume $q \geq 1$, such that $d(\omk,0) \geq
     S_{p-q-1}(\alpha)\cdot \frac{\alpha^{-q}}{L_2}$. Then 
     \begin{eqnarray*}
           \lefteqn{d(\omega_{k \pm m},0) \ \geq \ d(\omk,0) - d(\omm,0) 
           \ \geq \ S_{p-q-1}(\alpha)\cdot \frac{\alpha^{-q}}{L_2} -
             \frac{\alpha^{-(p-1)}}{L_2} }
           \\
           & = & (\underbrace{\alpha \cdot S_{p-q-1}(\alpha) -
            \alpha^{-(p-q-2)}}_{ \geq \ \alpha - 1 \ \geq \ S_\infty(\alpha)
           \textrm{ by (\ref{cond:alpha2}) } }) \cdot \frac{\alpha^{-(q+1)}}{L_2} \ \geq \
           S_{p-q-1}(\alpha)\cdot \frac{\alpha^{-(q+1)}}{L_2} \ . 
     \end{eqnarray*}
     This implies (\ref{eq:omegatransitionB}). Again, the case $q=0$
     is treated similarly, using (\ref{cond:alphagamma1}) instead of
     (\ref{cond:alpha2}).
     
   \item Now suppose $j \in \Omtil_{p-2}^{(\pm)}$. Then $\exists k \in \Z$
     such that $Q_{p-2}(k) \leq p -2$ and $j \in \Omega^{(\pm)}_{p-2}(k)$. As
     $Q_{p-1}(k \pm m) \geq Q_{p-2}(k)$ by (a), this implies $j \pm m \in
     \Omega^{(\pm)}_{p-1}(k \pm m)$, and as $Q_{p-1}(k+m) \leq Q_{p-2}(k)
     +1 \leq p-1$ this set is contained in $\Omtil_{p-1}^{(\pm)}$. 
\listend

\qed

\begin{lem} \label{lem:omegaestimates}
  Let $u,v \in \N$ be fixed and suppose $\omega$ satisfies the
  diophantine condition (\ref{cond:diophantine}). Then there exist
  functions $h,H : \R_+^2 \ra \R_+$ with $h(\gamma,\alpha) \nearrow
  \infty$ and $H(\gamma,\alpha) \searrow 0$ as $(\gamma + \alpha^{-1})
  \searrow 0$, such that
\alphlist
\item
 \[ 
     \nu(q) \ \geq \ \nutil(q) \  \geq  \ h(\gamma,\alpha) \cdot (q+2) 
     \cdot w \ \ \  \forall q \in \N
 \]
where $w:= \tilde{u}+\tilde{v}+1 = u+v+5$.
\item
  \[ 
      \# ([1,N] \cap \Omega_\infty ) \ \leq \
      H(\gamma,\alpha) \cdot N \ \ \textrm{and} \ \ \ 
      \# ([-N,-1] \cap \Omega_\infty ) \ \leq \  H(\gamma,\alpha)
      \cdot N   \ \ \ \forall N \in \N .
  \]
\listend
\end{lem}
\proof\   
\alphlist
\item
  The diophantine condition implies that $c\cdot \nutil(q)^{-d} \leq
  2S_\infty(\alpha)\cdot \frac{\alpha^{-(q-1)}}{L_2}$ (if $q \geq 2$). Thus, a simple
  calculation yields
  \[
       \frac{\nutil(q)}{w\cdot(q+2)} \ \geq \
       \left( \frac{c\cdot L_2}{2S_\infty(\alpha)} \right)^\dtel \cdot
       \frac{\alpha^{\frac{q-1}{d}}}{w\cdot(q+2)} 
  \]
  and the right hand side converges to $\infty$ uniformly in $q$ as $\alpha
  \ra \infty$. Similarly, 
  \[
       \frac{\nutil(1)}{3w} \ \geq \ \frac{1}{3w} \cdot \left(
         \frac{(c\cdot L_2)}{(8\gamma +
       2S_\infty(\alpha)\cdot\alpha^{-1})}\right)^\dtel 
  \]
  and again the right hand side converges to $\infty$ as $\gamma +
  \alpha^{-1} \searrow 0$. Thus we can define the minimum of both
  estimates as $h(\gamma,\alpha)$. 
\item 
  As we have seen in (a) we have $\nutil(q) \geq
  \left( \frac{c\cdot L_2}{2S_\infty(\alpha)} \right)^\dtel \cdot \alpha^{\frac{q-1}{d}}$
  if $q \geq 2$. Now $[1,N] \cap \Omega_\infty(j) =
  \emptyset$ if $j > N + w\cdot Q_\infty(j)$. Therefore
  \begin{eqnarray*} 
     \lefteqn{\textstyle \frac{1}{N} \cdot \# ( [1,N] \cap \Omega_\infty )} \\
     & \leq & \frac{1}{N} \sum_{q=1}^{\infty} q \cdot w \cdot \# \{ 1 \leq j
     \leq N+q\cdot w \mid Q_\infty(j) = q \}  
     \\  
     & \leq &  \frac{1}{N} \left( \frac{N+w}{\nutil(1)}\cdot w + \sum_{q=2}^{\infty}
    q\cdot w \cdot \frac{N+q\cdot w}{\nutil(q)} \right) \\
     & \leq & 
    \frac{w+\frac{w^2}{N}}{\nutil(1)} + \sum_{q=2}^{\infty}
     \frac{q\cdot w+q^2 \cdot
    \frac{w^2}{N}}{\left( \frac{c\cdot L_2}{2S_\infty(\alpha)} \right)^\dtel \cdot
     \alpha^{\frac{q-1}{d}} } \ .
  \end{eqnarray*}
  The right hand side converges to 0 uniformly in $N$ as $\gamma + \alpha^{-1} \ra
  0$ and we can use it as the definition of $H(\gamma,\alpha)$.
\listend

\qed


\subsection{Exceptional intervals and admissible times} \label{Admissibletimes}

In order to decide whether a time $N \in \N$ is admissible, in the sense of
Section~\ref{Admissiblepoints}, we will first have to introduce exceptional
intervals $J(m)$ corresponding to close returns $m\in \N$ with $d(\omm,0) \leq
\frac{3\gamma}{L_2}$.  For the sets $\Omega_p$ defined above, two different
intervals $\Omega_p(m)$ and $\Omega_p(n) \ (m \neq n)$ can overlap, without
one of them being contained in the other. This is something we want to exclude
for the exceptional intervals, and we can do so by carefully choosing their
lengths. To this end, we have to introduce two more assumptions on the
parameters:

We let $h$ and $H$ be as in Lemma~\ref{lem:omegaestimates} and suppose that
$\gamma$ and $\alpha$ are chosen such that $h(\gamma,\alpha) \geq 1$ and
$H(\gamma,\alpha) \leq \frac{1}{12w}$. In other words, we will assume
that for all $q,N \in \N$ there holds
\begin{eqnarray}
    \nutil(q) & \geq & (q+2)\cdot w \ , \label{cond:hfunctions1} \\
    \#([-N,-1]\cap \Omega_\infty) & \leq & \frac{N}{12w} \ \
    \textit{ and } \ \    \#([1,N]\cap \Omega_\infty) \ \  \leq \ \
    \frac{N}{12w}  \ .
    \label{cond:hfunctions2}  
\end{eqnarray}

\begin{bem}
  \label{bem:omegaestimates}
Suppose that (\ref{cond:gamma0}), (\ref{cond:alphagamma0}) and
  (\ref{cond:hfunctions1}) hold. Assumption~(\ref{cond:hfunctions1}) ensures
  that on the one hand the sets $\Omega_\infty(j)$ never contain the origin
  (and are, in fact, a certain distance away from it), and on the other hand
  two such sets of approximately equal size do not interfere with each other.
  This will be very convenient later on. To be more precise:
\alphlist
\item There holds
  \begin{equation}\label{eq:origin}
        -2,-1,0,1,2 \ \notin \ \Omega_\infty  \ .
  \end{equation} 
\item
  If $Q_\infty(j) \geq q$ for some $j \in \Z$ then
  \begin{equation} \label{eq:ominftybounds}
      [-\util \cdot (q+2),\vtil\cdot (q+2)] \cap \Omega_\infty(j) =  \emptyset \ . 
  \end{equation}
\item Let $m,n\in\Z$, $m\neq n$. Then $\Omega_\infty(m) \cap
  \Omega_\infty(n) = \emptyset$ whenever $|Q_\infty(m) - Q_\infty(n)|
  \leq 2$ or $|Q_p(m) - Q_p(n)| \leq 1$ for some $p \in \N_0$.  \listend
\end{bem}
\proof\ 
(a) and (b) follow immediately from (\ref{cond:hfunctions1}) and the
definition of the sets $\Omega_\infty(j)$. In order to prove (c), let
$q := \min\{ Q_\infty(m),Q_\infty(n) \}$. Then necessarily 
$d(\omega_{m-n},0) = d(\omm,\omn) < 2S_\infty(\alpha) \cdot 
\frac{\alpha^{-(q-1)}}{L_2}$ and thus %
$|m-n| \geq \nutil(q) \geq (q+2)\cdot w$ by (\ref{cond:hfunctions1}).
On the other hand both $Q_\infty(m)$ and $Q_\infty(n)$ are at most 
$q+2$, and thus the definition of the $\Omega_\infty(j)$ implies the 
disjointness of the two sets. Finally, note that $|Q_p(m) - Q_p(n)| 
\leq 1$ implies $|Q_\infty(m) - Q_\infty(n)| \leq 2$ by
(\ref{eq:pjestimates}).
  
\qed

\begin{bem}
 \label{bem:existenceofls}
Suppose that (\ref{cond:gamma0}), (\ref{cond:alphagamma0}),
  (\ref{cond:hfunctions1}) and (\ref{cond:hfunctions2}) hold.
  (\ref{cond:hfunctions2}) ensures that the ``density'' of the set
  $\Omega_\infty$ is small enough, and this will be very important for the
  construction later on. On the other hand, it also enables us now to choose
  suitable lengths for the exceptional intervals $J(m)$:
 
 We have $\# ([-\util\cdot q,-1] \cap \Omega_\infty) \leq
 \frac{q}{12}$. This implies that we can find at least two
 consecutive integers outside of $\Omega_\infty$ in the interval
 $[-\util\cdot q,-u\cdot q]$. In other words, for all $q \in \N$ there 
  exists $l^-_q \in \N$ such that
  \begin{eqnarray} \label{eq:lmin}
     u\cdot q \ \leq \ l^-_q \ <  \ \util \cdot q & \ \textrm{and} \ & -l^-_q, \
     -l^-_q-1  \notin \Omega_\infty \ .
  \end{eqnarray}
  Similarly, there exists $l^+_q\in\N$, such that
  \begin{eqnarray} \label{eq:lplus}
    v \cdot q \ \leq \ l^+_q \ < \ \vtil \cdot q & \ \textrm{and} \ &
    l^+_q, \ l^+_q + 1 \notin   \Omega_\infty \ .
  \end{eqnarray}
  In addition, we can assume that $l^\pm_p \geq l^\pm_q$ whenever $p
  \geq q$. (If $l^+_q,l^+_q+1$ are both contained in $[v\cdot
  (q+1),\vtil \cdot (q+1)]$, then we can just take $l^+_{q+1} =
  l^+_q$. Otherwise, we find a suitable $l^+_{q+1} > l^+_q$ in this
  interval.) Note also that (\ref{eq:lmin}),(\ref{eq:lplus}) and 
  (\ref{cond:hfunctions1}) together imply that
\begin{equation} \label{eq:lplusestimates}
     \min\{u,v\} \cdot q \ \leq \  l^\pm_q \ < \
     \nutil(\max\{1,q-2\}) \ \leq \ \nu(\max\{1,q-2\}) \ . 
\end{equation}
\end{bem}
Now we are able to define the exceptional intervals:
\begin{definition}[Exceptional intervals]
  \label{def:ansets}
Suppose that (\ref{cond:gamma0}), (\ref{cond:alphagamma0}),
  (\ref{cond:hfunctions1}) and (\ref{cond:hfunctions2}) hold.  Then for any $q
  \in \N$, choose $l^{\pm}_q$ as in Remark~\ref{bem:existenceofls} and define,
  for any $m \in\N$ with $p(m) \geq 0$,
\begin{eqnarray*}
    \lambda^-(m) & := & m - l^-_{p(m)} \ \ \ , \ \ \  \lambda^+(m) \ \ 
    := \ \ m+    l^+_{p(m)}  \\ \nonumber \\
    J^-(m) & := & [\lambda^-(m),m] \ \ \ , \ \ \ J^+(m) \ \ := \ \ 
    [m+1,\lambda^+(m)] 
\end{eqnarray*} 
and 
\[
    J(m) \ \ := \ \ J^-(m) \cup J^+(m) \ .
\]
If $p(m) = -1$, then $J^{(\pm)}:= \emptyset$. Further, let
\[
    A_N \ \ := \ \ [1,N] \setminus \bigcup_{1 \leq m < N} J(m)
\ \ \textrm{ and } \ \ \Lambda_N \ \ := \ \ [1,N] \setminus A_N
\]

\end{definition}

From now on, we will use conditions (\ref{cond:gamma0}),
  (\ref{cond:alphagamma0}), (\ref{cond:hfunctions1}) and
  (\ref{cond:hfunctions2}) as standing assumptions in the remainder of this
  subsection, as well as in Subsection~\ref{Regulartimes} (since all the
  statements concern the preceding definition, directly or indirectly).

\begin{bem}
  \label{bem:ansets}
\alphlist
\item As we have mentioned before, the exceptional intervals are
  contained in the approximating sets. To be more precise, for each $m 
  \in \N$ with $p(m) \geq 0$ there holds
\begin{eqnarray} 
    J(m) & \subset & [\lambda^-(m)-1,\lambda^+(m)+1] \nonumber \\
    & \ssq & \Omega_0(m) 
    \ \ssq \ \Omega_p(m) \ \ssq \ \Omega_\infty(m) \ , \label{eq:jminominfty}
    \end{eqnarray}
where $p\in\N$ is arbitrary. This follows from the choice of the
$l^\pm_q$ in Remark~\ref{bem:existenceofls} together with the
definition of the intervals $\Omega_p(m)$. As a consequence, we have that
\begin{equation} \label{eq:lambdaninominfty} 
\Lambda_N \ \ssq \ \Omega_0 \ \ssq \ \Omega_p \ \ssq \ \Omega_\infty\ \ 
\ \forall N,p \in \N \  .
\end{equation}
\item Further, suppose that $m\neq n$ and $|Q_\infty(m) - Q_\infty(n)|
  \leq 2 $ or $|Q_p(m) - Q_p(n)| \leq 1$ for some $p \in \N_0$. Then
  (a) together with Remark~\ref{bem:omegaestimates}(c) implies that
\begin{eqnarray} \nonumber
\lefteqn{J(m) \cap J(n) \ = \ \emptyset \ = } \\
& = & [\lambda^-(m)-1,\lambda^+(m)+1] \cap
[\lambda^-(n)-1,\lambda^+(n)+1]  \ .\label{eq:jmndisjointness}
\end{eqnarray}
In particular this is true if $|p(m) - p(n)| \leq 1$ (recall that
$p(j) = Q_0(j)$). 
\item The sets $A_N$ were defined as subsets of $[1,N]$, and it will
  turn out that they contain a very large proportion of the points
  from that interval. This could lead to this impression they form an
  increasing sequence of sets, but this is not true. For example,
  suppose that $N$ itself is a close return, such that $p(N) \geq 1$.
  In this case $N$ may still be contained in $A_N$, as the exceptional
  interval $J(N)$ is not taken into account in the definition of this
  set, but surely $N \notin A_{N+1}$. Thus, whenever we reach a close
  return, there may be a sudden decrease in the sets $A_N$ in the next
  step.  In general, we only have the two relations
  \begin{equation}
       A_{N_2} \setminus A_{N_1} \ \ssq \ [N_1+1,N_2] \ \ \textrm{ and } \ \ 
  \end{equation}
  and 
  \begin{equation} \label{eq:aninequality}
       A_{N_2} \cap [1,N_1] \ \ssq \ A_{N_1} \ .
  \end{equation}
  where $N_1 \leq N_2$. However, the fluctuations and sudden decreases
  will only take place at the end of the interval $[1,N]$, and the
  starting sequence of the sets $A_N$ will stabilize at some point:
  Suppose $N_0 \leq N_1 \leq N_2$ and $N_0 \in A_{N_2}$. Then
  \begin{equation} \label{eq:stabilizingeffect}
    A_{N_1} \cap [1,N_0] \ = \ A_{N_2} \cap [1,N_0] \ = \ A_{N_0} \ .
  \end{equation}
  This simply follows from the fact that when $N_0$ is contained in
  $A_{N_2}$ no exceptional interval $J(m)$ with $m \in [N_0+1,N_2-1]$ can
  reach into $[1,N_0]$, as it would then have to contain $N_0$. Thus %
  $A_{N_1} \cap [1,N_0]  =  A_{N_2} \cap [1,N_0]  = [1,N_0] \setminus
  \bigcup_{1 \leq m < N_0} J(m)$. 
  
  Note that (\ref{eq:stabilizingeffect}) is always true whenever $N_0$ is
  not contained in any exceptional interval, i.e.\ $N_0 \notin
  \bigcup_{m\in\N} J(m) \ssq \Omega_0 \ssq \Omega_\infty$. In this
  case we have
  \begin{equation}
    \label{eq:finallystabilized}
    A_N \cap [1,N_0] \ = \ A_{N_0} \ \ \ \forall N \geq N_0 \ .
  \end{equation}
  In particular, as $l^+_q \notin \Omega_\infty$, this implies
  $[1,l^+_q] \cap A_N = A_{l^+_q} \ \forall N \geq l^+_q$. 
\item Note also that it is not always true that %
  $\Lambda_N = \bigcup_{1\leq m < N} J(m)$, as one of the exceptional
  intervals might extend beyond $N$, whereas $\Lambda_N$ was defined
  as a subset of $[1,N]$. However, as we will see this
  relation holds as soon as we restrict to \textit{`admissible'} times
  (see below).
\listend
\end{bem}
The sets $A_N$ will serve three different aims: First of all, they
will play an important role in the construction of the
sink-source-orbits themselves. Secondly, they will also be intermediates for
the definition of the sets $R_N$ of regular points. And finally, we will
now use them to define admissible times:
\begin{definition}[Admissible times]
  \label{def:admissibletimes}
A time $N \in \N$ is called \textbf{admissible} if $N \in A_N$ (which
is equivalent to $N \notin \Lambda_N$). The set $\{N\in\N \mid N
\textrm{ is admissible } \}$ will be denoted by $A$.

\end{definition}
\begin{bem}
  \label{bem:admissibletimes}
\alphlist
\item Any $N \in \N \setminus \Omega_0$ is admissible (see
  Remark~\ref{bem:ansets}(a)). In particular, $l^+_q$ and $l^+_q+1$
  are admissible for any $q \geq 1$. 
\item As we mentioned above, for any admissible time $N$ there holds
     \begin{equation} \label{eq:lambdanequality}
            \Lambda_N \ = \ \bigcup_{1\leq m<N} J(m) \ ,
     \end{equation}
     as $N \in A_N$ ensures that none of the exceptional intervals
     $J(m)$ with $m < N$ extends further than $N-1$. 
   \item For any $N_1 \in \N$, all times $N_0 \in A_{N_1}$ are
     admissible. This is a direct consequence of the fact that
     $A_{N_1} \cap [1,N_0] \ssq A_{N_0}$ (see
     (\ref{eq:aninequality})). However, as already mentioned there
     might also be further admissible times contained in
     $[1,N_1]\setminus A_{N_1} = \Lambda_{N_1}$.
   \item Note that $A=\bigcup_{N\in\N} A_N$. The inclusion $\ssq$
     follows directly from the definition, whereas $\supseteq$ is a
     consequence of (c). 
\listend
\end{bem}
Now we can also verify the property of the exceptional intervals which
was mentioned at the beginning of this section: Whenever two such
intervals $J(m)$ and $J(n)$ intersect, one of them is contained in the
other. We do not prove this statement in full, but rather concentrate
on `maximal' intervals, as this will be sufficient for our purposes.
\begin{lem}
  \label{lem:centralpoints} Let $N\in\N$ be admissible and suppose $J$ is
  a non-empty maximal interval in $\Lambda_N = [1,N] \setminus A_N$. Then
  there exists a unique $m \in J$ with $p(m) = \max_{j\in J} p(j)$, and there
  holds $J = J(m)$. Furthermore, $p(j) < p-1 \ \forall j \in J \setminus
  \{m\}$.
\end{lem}
\proof \
Let $p:= \max_{j \in J} p(j)$ and $m \in J$ with $p(m) = p$. Obviously
there holds $J(m) \ssq J$. By definition, there cannot be any $j \in J
\supseteq J(m)$ with $p(j) > p$. Therefore, as
Remark~\ref{bem:ansets}(b) implies that $|p(j)-p| > 1 \ \forall j \in
J(m) \setminus \{m\}$, there holds $p(j) < p-1 \ \forall j \in J(m)
\setminus \{m\}$. Thus, it suffices to prove that $J=J(m)$. This will
in turn follow if we can show that $\lambda^-(m) - 1$ and
$\lambda^+(m)+1$ are not contained in $\Lambda_N$, because then $J(m)$
is a maximal interval in $\Lambda_N$ itself and must therefore be
equal to $J$. We will only treat the case of $\lambda^-(m)-1$, the other
one is similar. In order to show that $\lambda^-(m)-1$ is not
contained in $J(k)$ for any $k = 1 \ld N$ we distinguish three different
cases, according to the value of $Q_{p-2}(k)$:

First suppose $Q_{p-2}(k) > p+1$. Then $p(k) > p$ by
(\ref{eq:pjestimates}). If $\lambda^-(m)-1 \in J(k)$, then
 $J(k) \cup J(m)$ is an interval and therefore $k \in J(k) \ssq
J$. But this contradicts the definition of $p$. 

If $Q_{p-2}(k) \in \{ p-1,p,p+1\}$, then $|Q_\infty(k) -
Q_\infty(m)| \leq 2$ (again (\ref{eq:pjestimates})) and
therefore  $\lambda^-(m)-1 \notin \Lambda(k)$ by Remark~\ref{bem:ansets}(b).

This only leaves the possibility $Q_{p-2}(k) \leq p-2$. But in this
case $\lambda^-(m)-1 \in \Lambda(k)$ implies $\lambda^-(m)-1 \in
\Omega_{p-2}(k) \subseteq \Omtil_{p-2}$ (see
Remark~\ref{bem:ansets}(a)). As $p(m)=p$
we can apply Lemma \ref{lem:omegatransition}(b) to obtain that
$\lambda^-(m) - 1 - m = -l^-_p - 1 \in \Omtil_{p-1}$, contradicting
$-l^-_p - 1 \notin \Omega_\infty $ (by the choice of the $l^\pm_q$ in
Remark \ref{bem:existenceofls}).

As mentioned, the same arguments apply to $\lambda^+(m) +
1$, which completes the proof.

\qed

\ \\
This naturally leads to the following
\begin{definition} \label{def:centralpoints}
If $N$ is admissible and $A_N = \{ a_1 \ld a_n \}$ with $1=a_1  < \ldots < a_n=N$, let
\[
        \JN \ := \{[a_k+1,a_{k+1}-1] \mid k=1 \ld n-1 \} \setminus
        \{ \emptyset \}
\]
be the family of all maximal intervals in $\Lambda_N = [1,N] \setminus
A_N$ and $\J := \bigcup_{N \in \N} \JN$. For any $J \in \J$ let $p_J:=
\max_{j \in J} p(j)$ and define $m_J$ as the unique $m \in J$ with
$p(m) = p_J$. $m_J$ will be called the \textbf{central point} of the
interval $J$.

Further, let $J^- := J^-(m_J)$ and $J^+ := J^+(m_J)$ (note that
$J=J(m_J)$ by Lemma~\ref{lem:centralpoints}).
\end{definition}
Note that not for every $n \in \N$ with $p(n) > 0$ the interval $J(n)$ 
is contained in $\J$. In fact, this will be wrong whenever $J(n) \ssq
J^+(m)$ for some $m<n$. 

Among some other facts, the following lemma states that central points
are always admissible. In the light of the discussion in
Section~\ref{Admissiblepoints}, it is not surprising that this will
turn out to be crucial for the construction.
\begin{lem} \label{lem:endpoints}
\alphlist
\item Let $J \in \J$. Then $\lambda^-(m_J) -1 \in A_{m_J}$,
  $\lambda^-(m_J) \in A_{m_J}$ and $m_J \in A_{m_J}$. In particular,
  $\lambda^-(m_J)-1, \lambda^-(m_J)$ and $m_J$ are admissible.
  Further, there holds
  \begin{equation} \label{eq:pjinjmspecial}
  p(j) \leq Q_\infty(j) \leq \max\{0,p_J-2\} \ \ \ \forall j \in J
  \setminus \{m_J\} \ .
  \end{equation}
\item More generally, if $J \in \J$ and $q \leq p_J$, then
  $m_J-l^-_q-1,\ m_J-l^-_q,\ m_J\in A_{m_J}$. In particular,
  $m_J-l^-_q-1,\ m_J-l^-_q$ and $m_J$ are admissible.  Further there
  holds
  \begin{equation} \label{eq:pjinjm}
     p(j) \leq Q_\infty(j) \leq \max\{0,q-2\} \ \ \  \forall j \in [m_J - l^-_q,m_J+l^+_q]
     \setminus \{m_J\} \ .
  \end{equation}
\item If $J \in \J$, then $\lambda^+(m_J) + 1$ is admissible. 
\item For all $q \in \N$ there holds $\nu(q)-l^-_q -1,\
  \nu(q)-l^-_q,\ \nu(q) \in A_{\nu(q)}$ and 
  \[
     Q_\infty(j) \leq \max\{0,q-2\} \ \ \  \forall j \in [\nu(q) -
     l^-_q,\nu(q)+l^+_q] \setminus \{\nu(q)\} \ .
  \]
  In particular $\nu(q)-l^-_q-1,\nu(q)-l^-_q$ and $\nu(q)$ are admissible. 
 \listend
\end{lem}
\proof\ 
\alphlist
\item This is a special case of (b), which we prove below.
\item Let $m := m_J$ and $j \neq m$. 
  Suppose $Q_\infty(j) \geq q-2$. Then
  \[ 
  d(\omega_{m-j},0) \ = \ d(\omm,\omj) \ \leq \ 2 S_\infty(\alpha) \cdot
  \frac{\alpha^{-(q-3)}}{L_2} \ .
 \] 
 Therefore $|m-j| \geq \nutil(q-2) > l^\pm_q$ by
 (\ref{eq:lplusestimates}), which implies that $j \notin
 [m-l^-_q,m+l^+_q]$. This proves (\ref{eq:pjinjm}).
 
 As $J \in \J$, there exists some $N > m$ such that $J$ is a maximal
 interval in $\Lambda_N$ and consequently $\lambda^-(m) - 1$ is
 contained in $A_N$ (in particular $p(\lambda^-(m)-1)=0$). Hence, for
 any $n < \lambda^-(m)-1$ the interval $J(n)$ lies strictly to the
 left of $\lambda^-(m)-1$ and can therefore not intersect $J$.  Thus,
 in order to show that $m - l^-_q-1,l^-_q,m \in A_m$, it suffices to
 show that none of these points is contained in $U :=
 \bigcup_{n\in[\lambda^-(m),m-1]} J(n)$. However, by
 (\ref{eq:jminominfty}) and (\ref{eq:pjinjm}) there holds $U \ssq
 \Omtil_{q-2}$. As $p(m) \geq q$ by assumption,
 Lemma~\ref{lem:omegatransition}(b) implies $U-m \ssq \Omtil_{q-1}
 \ssq \Omega_\infty$ and the statement follows from $-l^-_q-1,l^-_q,0
 \notin \Omega_\infty$.
 
 Finally, note that $m-l^-_q \in A_m$ implies $m-l^-_q \in 
 A_{m-l^-_q}$ by (\ref{eq:aninequality}), similarly for $m-l^-_q-1$, 
 such that these points are both admissible. 
\item As $m$ is admissible, $\lambda^+(m)+1$ cannot be contained in 
  $J(n)$ for any $n < m$ (as all of these intervals must be contained 
  in $[1,m-1]$). Thus, it suffices to show that $\lambda^+(m)+1$ is
  not contained in $\tilde{U} := \bigcup_{n\in[m+1,\lambda^+(m)]}
  J(n)$. But this set is again contained in $\Omtil_{p_J - 2}$ by
  (\ref{eq:pjinjm}). Therefore $\tilde{U} -m \ssq \Omega_\infty$ by
  Lemma~\ref{lem:omegatransition}(b), and the statement follows
  from $l^+_{p_J}+1 \notin \Omega_\infty$.
\item We show that $\nu(q)$ is admissible.
  Lemma~\ref{lem:admissiblecentral} below then implies that $\nu(q)$ is a
  central point, and we can therefore apply (b) in order to prove (d).
  
  Suppose $n< \nu(q)$. We have to show that $\nu(q) \notin
  \Omega_\infty(n) \supseteq J(n)$. In order to see this, note that
  $p(j) < q$ by definition of $\nu(q)$. Thus $d(\omega_{\nu(q)-n},0) =
  d(\omega_{\nu(q)},\omn) \geq \nutil(q) \geq (q+2) \cdot w$ by
  (\ref{cond:hfunctions1}), and consequently $\nu(q) \notin
  \Omega_\infty(n)$. As $n<\nu(q)$ was arbitrary, this implies $\nu(q)
  \in A_{\nu(q)}$, such that $\nu(q)$ is admissible.

\listend

\qed
 
\ \\
For Part (a) of the preceding lemma, the inverse is true as well:
\begin{lem}
  \label{lem:admissiblecentral}
  Suppose $m\in \N$ is admissible and $p(m) > 0$. Then $J(m) \in
  \J_{\lambda^+(m)+1} \ssq \J$ and $\lambda^-(m)-1,\lambda^-(m)$ and
  $\lambda^+(m)+1$ are admissible.
\end{lem}
\proof\ We start by proving that $\lambda^+(m)+1$ is admissible, i.e.\
contained in $A_{\lambda^+(m)+1}$. First of all, the fact that $m$ is
admissible ensures that none of the intervals $J(n)$ with $n < m$ intersects
$[m+1,\lambda^+(m)+1]$. Therefore, none of these intervals can contain
$\lambda^+(m)+1$, and for $J(m)$ the same is true by definition. Now suppose
$n \in [m+1,\lambda^+(m)]$. Then, similar as in the proof of
Lemma~\ref{lem:endpoints}(b) we obtain $p(n) \leq p(m)-2$ and therefore $J(n)
\ssq \Omtil_{p(m)-2}$. Thus $J(n)-m$ is contained in $\Omtil_{p(m)-1} \ssq
\Omega_\infty$ by Lemma~\ref{lem:omegatransition}(b) and can therefore not
contain $l^+_{p(m)}+1 \notin \Omega_\infty$. Thus
$\lambda^+(m)+1=m+l^+_{p(m)}+1$ is admissible.
  
By Lemma~\ref{lem:centralpoints}, for any maximal interval $J = J(n)
\in \J_{\lambda^+(m)+1}$ that intersects $J(m)$ there holds either
$J(n) = J(m)$, such that $n=m$, or $J(m) \subseteq J(n)$.  However,
the second case cannot occur if $n < m$ (as $m$ is admissible), and
for $n\in [m+1,\lambda^+(m)]$ it is ruled out as we have just argued that $p(n) < p(m)$
for such $n$. This proves $J(m) \in \J_{\lambda^+(m)+1}$.
  
Finally, we can apply Lemma~\ref{lem:endpoints}(a) to $J=J(m)$, which yields
that $\lambda^-(m)-1$ and $\lambda^-(m)$ are admissible as well.

\qed


\subsection{Regular times} \label{Regulartimes}

Now we can turn to defining the sets of regular points $R_N
\ssq [1,N]$. The sets $A_N$ already contain all points outside of the
exceptional intervals $J(m)\ (m \in [1,N-1])$. As described in
Section~\ref{Strategy}, we have to add certain points from the right
parts $J^+(m)$ of these intervals. In order to do so, for each $J \in
\J_N$ we will define a set $R(J) \ssq J^+$ and then let $R_N = A_N
\cup \bigcup_{J\in\J_N} R(J)$. Both $R_N$ and $R(J)$ will be defined by
induction on $p$. To be more precise, in the $p$-th step of the
induction we first define $R(J)$ for all $J\in\J$ with $p_J \leq p-1$,
and then $R_N$ for all admissible times $N \leq \nu(p)$.
\medskip

As in the preceding one, conditions (\ref{cond:gamma0}),
  (\ref{cond:alphagamma0}), (\ref{cond:hfunctions1}) and
  (\ref{cond:hfunctions2}) will be used as standing assumptions in this
  subsection (since all of the statements in this subsection directly or
  indirectly depend on Definition~\ref{def:ansets}).
\begin{definition}[Regular times]
  \label{def:regularsets}
  As mentioned, we proceed by induction on $p$. Note that the
  inclusions $R_N \ssq [1,N]$ and $R(J) \ssq J^+$ are preserved in
  every step of the induction.

\ \\
\underline{$p=1:$} \ In order to start the induction
  let
\[
    R_N \ := \ [1,N] \ 
\]
for any $N \leq \nu(1)$.  Note that by definition there is no $J\in\J$
with $p_J = 0$.

\ \\
\underline{$p \ra p+1$:} \ Suppose $R(J)$ has been defined for all $J
\in \J$ with $p_J \leq p-1$ and $R_N$ has been defined for all
admissible times $N \leq \nu(p)$. In particular, this means that
$R_{l^+_p}$ has defined already.%
\footnote{As $l^+_p \leq \nu(p)$ by (\ref{eq:lplusestimates}) and
  $l^+_p$ is admissible by Remark~\ref{bem:admissibletimes}(a).}
Then, for all $J \in \J$ with $p_J = p$ let
\begin{equation}
  \label{eq:recursivedef1}
    R(J)  \ = \ R_{l^+_p} + m_J \ .
\end{equation}
Note that as $J^+ = [m_J+1,m_J+l^+_p]$, the inclusion $R(J) \ssq J^+$
follows from $R_{l^+_p} \ssq [1,l^+_p]$. 
Further, for any admissible $N \in [\nu(p)+1,\nu(p+1)]$ let
\begin{equation}
   \label{eq:recursivedef2}
R_N \ := \ A_N \cup \bigcup_{J\in\J_N} R(J) \ .
\end{equation}
Here the inclusion $R_N \ssq [1,N]$ follows from $R(J) \ssq J^+  \ssq
 J \ \forall J \in \J_N$, as $J \ssq [1,N] \ \forall J \in \J_N$ by
definition (see Definition~\ref{def:centralpoints}). 
\medskip

\noindent
Finally, we call $j\leq N$ \textbf{regular with respect to $N$} if $j$ is
contained in $R_N$.
\end{definition}
\begin{bem}
  \label{bem:regularsets} 
\alphlist
\item Obviously any $j\in A_N$ is regular with respect to $N$. As
  $[1,N] \setminus A_N = \Lambda_N \ssq \Omega_\infty$ (see
  (\ref{eq:lambdaninominfty})), this implies that any $j \in \N \setminus
  \Omega_\infty$ is regular with respect to any $N \geq j$. In
  particular (see Remarks~\ref{bem:omegaestimates} and \ref{bem:existenceofls})
  \begin{equation} \label{eq:lpminrn}
       1,2, l^+_q , l^+_q+1 \ \in \ A_N \ \ssq \ R_N \ \ \ \forall q \in
        \N, \ N \geq l^+_q+1 \ .
  \end{equation}
\item Similar to the sets $A_N$, the sequence $(R_N)_{N\in\N}$ is not
  increasing (compare Remark~\ref{bem:ansets}(c)). However, if $N_0
  \leq N_1 \leq N_2$ are all admissible and $N_0 \in A_{N_2}$, then
  \begin{equation} \label{eq:regularstabilized}
      R_{N_1} \cap [1,N_0] \ = \ R_{N_2} \cap [1,N_0] \ = \ R_{N_0} .
  \end{equation}
  This can be seen as follows: $N_0 \in A_{N_2}$ implies that no interval
  $J(m) \ (N_0 \leq m < N_2)$ can reach into $[1,N_0]$, and in addition $N_0$
  is admissible (see (\ref{eq:stabilizingeffect})).  Therefore, since $R(J)
  \ssq J$, all three sets in (\ref{eq:regularstabilized}) coincide with
  $A_{N_0} \cup \bigcup_{J \in \J_{N_0}} R(J)$.

  In particular, by (\ref{eq:lpminrn}) this implies that 
  \begin{equation}
   \label{eq:lplusstabilized}
    R_N \cap [1,l^+_q]  = R_{l^+_q} 
    \ \textrm{ and } \  R_N \cap [1,l^+_q+1] = R_{l^+_q+1} \ \ \  \forall
    N \geq l^+_q+1 \ .  
  \end{equation}
\item Let $J \in \J$. As $R(J) = R_{l^+_{p_J}} +m_J$, statement (a)
  implies
    \begin{equation}
      \label{eq:mponeinrj}
       m_J+1,m_J+2,m_J+l^+_q,m_J+l^+_q+1  \in  R(J) \ \ \ \forall q \leq p_J \ .
    \end{equation}
\listend
\end{bem}
It will also be useful to have a notation for the sets of non-regular
points:
\begin{definition}
  \label{def:gammasets}
For each admissible time $N \in \N$ let
\[
    \Gamma_N \ := \ [1,N] \setminus R_N
\]
and for each $J \in \J$ let
\[
    \Gamma^+(J) \ := \ J^+ \setminus R(J) \ \ \textrm{ and } \ \
    \Gamma(J) \ := \ J^- \cup \Gamma^+(J) \ .
\]

\end{definition}
\begin{bem}
  \label{bem:gammasets}
\alphlist
\item Note that 
      \begin{equation} \label{eq:gammanequality}
            \Gamma_N = \bigcup_{J \in \J_N} \Gamma(J) \ = \
            \bigcup_{J\in\J_N} J^- \cup \Gamma^+(J) \   . 
      \end{equation} 
\item Similar to (\ref{eq:recursivedef1}), the sets $\Gamma^+(J)$ satisfy
  the recursive equation
  \begin{equation}
    \label{eq:recursivegammas} \Gamma^+(J) \ = \ \Gamma_{l^+_{p_J}} +
    m_J \ .
  \end{equation}
\item As $A_N \ssq R_N$, there holds $\Gamma_N \ssq \Lambda_N$. Thus,
  Remark~\ref{bem:ansets}(a) implies
  \begin{equation} \label{eq:gammaninominfty}
        \Gamma_N \ \ssq \ \Lambda_N \ \ssq \  \Omega_0 \ \ssq \
        \Omega_p \ \ssq \Omega_\infty 
  \end{equation}
  for all admissible times $N \in \N$. $p \in \N$ is arbitrary. 
\item Suppose both $N$ and $N+1$ are admissible. Then $p(N)=0$, such
  that $J(N) = \emptyset$, otherwise $N+1$ would be contained in
  $J(N)$ could therefore not be admissible. Thus there holds $\Lambda_N =
  \Lambda_{N+1}$ (see (\ref{eq:lambdanequality})). But this means that
  $J_N = J_{N+1}$ and consequently $\Gamma_N=\Gamma_{N+1}$ (see
  (\ref{eq:gammanequality})). In particular, this is true whenever
  $N,N+1 \notin \Omega_\infty$, such that we obtain
  \begin{equation}
    \label{eq:gammalqplus}
    \Gamma_{l^+_q} \ = \ \Gamma_{l^+_q+1} \ \ \ \forall q \in \N \ .
  \end{equation}
\listend
\end{bem}
Now we must gather some information about the sets $R_N$ and
$\Gamma_N$. First of all, the following lemma gives some basic
control. In order to state it,  let
\begin{equation} \label{eq:omtilmin}
\Omtil^{(\pm)}_{-1} \ :=  \ \emptyset \ 
\end{equation}
and note that $\Omtil^{(\pm)}_{0} \ = \emptyset$ as well.
\begin{lem}
  \label{lem:regularsetsbasic}
  \alphlist
\item For any $J \in \J$ there holds $\Gamma(J) \ssq \Omtil_{p_J-2}$. Further,
  for any admissible $N \leq \nutil(q)$ there holds $\Gamma_N \ssq
  \Omtil_{q-1}$.
\item If $j \in R(J)$ for any $J \in \J$, then 
\begin{equation} \label{eq:rjdistance}
d(\omj,0) \ \geq \ \frac{4\gamma}{L_2} - S_{p_J-1}(\alpha) \cdot 
\frac{\alpha^{-1}}{L_2} \ \geq \ \frac{3\gamma}{L_2} \ .
\end{equation}
Further, for any admissible $N \leq \nu(q)$ there holds
\begin{equation} \label{eq:rndistance}
      d(\omj,0) \ \geq \ \frac{4\gamma}{L_2} - S_{q-1}(\alpha) \cdot 
     \frac{\alpha^{-1}}{L_2} \ \geq \ \frac{3\gamma}{L_2} \ \ \ \forall j
     \in R_N \setminus \{N\} \ .
\end{equation}
\listend
\end{lem}
\proof\  
\alphlist
\item   We proceed by induction on $q$. More precisely, we prove the
  following induction statement:
  \begin{eqnarray}
    \label{eq:rsbinduction1}
    \Gamma^+(J) & \ssq & \Omtil^-_{p_J-2} \ \ \  \forall J \in \J :
    p_J \leq q \  \\ \label{eq:rsbinduction2} 
    \Gamma_N & \ssq & \Omtil_{q-1}^- \ \ \ \forall N \leq \nutil(q) \ .
  \end{eqnarray}
  For $q=1$ note that $\Gamma_N$ is empty for all $N \leq \nutil(1)$. In
  particular $\Gamma_{l^+_1}$ is empty, as $l^+_1 \leq \nutil(1)$ by
  (\ref{eq:lplusestimates}). But this means in turn that for any
  $J\in\J$ with $p_J = 1$ the set $\Gamma^+(J) = \Gamma_{l^+_1} + m_J$
  is empty as well (see (\ref{eq:recursivegammas})).
  
  Let $p \geq 1$ and suppose the above statements hold for all $q \leq
  p$. Further, let $J\in \J$ with $p_J = p+1$. Then $\Gamma_{l^+_{p+1}} \ssq
  \Omtil^-_{p-2}$ as $l^+_{p+1} < \nutil(p-1)$ by
  (\ref{eq:lplusestimates}). Therefore
  \[
       \Gamma^+(J) \ = \ \Gamma_{l^+_{p+1}} + m_J \ \ssq \
       \Omtil^-_{p-2} + m_J \ \ssq \ \Omtil^-_{p-1} \ .
  \] 
  by Lemma~\ref{lem:omegatransition}(b). Thus (\ref{eq:rsbinduction1}) 
  holds for $q=p+1$. 
  
  Now suppose $N \leq \nutil(p+1)$ and note that this implies $Q_p(m) \leq p \
  \forall m < N$.  Further, we have $\Gamma_N = \bigcup_{J\in\J_N} J^- \cup
  \Gamma^+(J)$ by (\ref{eq:gammanequality}). As $J^- \ssq \Omega^-_p(m_J) \
  \forall J\in \J$ and $m_J < N \ \forall J \in \J_N$, there holds $J^- \ssq
  \Omtil^-_p$ for any $J \in \J_N$, and for $\Gamma^+(J)$ the same follows
  from (\ref{eq:rsbinduction1}). This proves (\ref{eq:rsbinduction2}) for
  $q=p+1$.
\item Suppose (\ref{eq:rjdistance}) holds whenever $p_J \leq p$. This implies
  (\ref{eq:rndistance}) for all $q \leq p$: We have $d(\omj,0) \geq
  \frac{4\gamma}{L_2}$ whenever $j\in A_N\setminus \{N\}$ for some $N\in\N$,
  and further $p_J < q \ \forall J \in \J_N$ whenever $N \leq \nu(q)$.

  It remains to prove (\ref{eq:rjdistance}) by induction on $p_J$. If $p_J
  \leq 2$ the statement is obvious, because then $p(j) = 0 \ \forall j \in J
  \setminus \{ m_J\}$ by Lemma~\ref{lem:endpoints}(a).

  Suppose now that (\ref{eq:rjdistance}) holds whenever $p_J \leq p$. As
  mentioned above, (\ref{eq:rndistance}) then holds for all $q \leq p$. Let
  $p_I = p+1$ for some $I \in\J$ and $j \in R(I)$.  Then $j-m_I \in
  R_{l^+_{p+1}}$ (see (\ref{eq:recursivedef1})), and as $l^+_{p+1} \leq
  \nu(p)$ we can apply (\ref{eq:rndistance}) with $q=p$ to obtain that
  \[ 
  d(\omega_{j-m_I},0) \ \geq \ \frac{4\gamma}{L_2} - S_{p-1}(\alpha)
  \cdot \frac{\alpha^{-1}}{L_2} \ .
  \]
  Consequently
  \begin{eqnarray*}
      d(\omj,0) & \geq & d(\omega_{j-m_I},0) - d(\omega_{m_I},0) 
      \ \geq \ \frac{4\gamma}{L_2} - S_{p-1}(\alpha) \cdot
      \frac{\alpha^{-1}}{L_2} - \frac{ \alpha^{-p}}{L_2} 
      \\
      & = & \frac{4\gamma}{L_2} - \left(S_{p-1}(\alpha) +
      \alpha^{-(p-1)} \right) \cdot \frac{\alpha^{-1}}{L_2} \ = \
      \frac{4\gamma}{L_2} - S_p(\alpha) \cdot \frac{\alpha^{-1}}{L_2} \ .
  \end{eqnarray*}
\listend

\qed

\ \\ As a consequence of Lemma \ref{lem:endpoints} and the preceding lemma, we
obtain the following statements and estimates. In order to motivate these, the
reader should compare the statements with the assumptions of
Lemma~\ref{lem:orbitsthrowout}~. 
\begin{lem}
   \label{lem:regsets}
\alphlist
\item For any admissible $N\in\N$ there holds
  \begin{equation}
  \#\left([1,j] \setminus R_N \right) \ \leq \ \frac{j}{12w} \ \ \ \ \ \forall
  j \in [1,N] \ .
  \end{equation}
  In particular 
  \begin{equation} \label{eq:regsets3}
    \# \left( [1,l^+_q] \setminus R_N \right) \ \leq \ \left[\frac{q}{12}\right] \ \leq \
    \max\left\{0,\frac{2q-5}{4}\right\} \ \ \ \ \ \forall q \in \N \ ,
  \end{equation}
 where $[x]$ denotes the integer part of $x\in\R^+$. 
\item Let $q \geq 1$ and $\sigma := \frac{u+3}{u+v}$. Then
\begin{equation}
   \# ([j+1,l^+_q] \setminus R_{l^+_q}) \ \leq \ \sigma\cdot (l^+_q - j) \ \ \
   \forall j \in [0,l^+_q-1] \ .
\end{equation}
\item Let $N\in\N$ be admissible, $J\in\JN$ and $\lambda^+:=\lambda^+(m_J)$. Then 
\begin{equation}
   \# ([j+1,\lambda^+] \cap \Gamma_N) \ \leq \ \sigma\cdot (\lambda^+-j) \
   \ \ \forall j \in [0,\lambda^+-1] \ .
\end{equation}
\item
Suppose $m\in \N$ is admissible and $p(m) \geq 1$, such that $J := J(m) \in
\J$ by Lemma~\ref{lem:admissiblecentral}~. Then for all $q \leq p_J$ there holds
       \begin{equation} \label{eq:approachfraction}
             \#([m-l^-_q,m] \setminus R_m) \ \leq \ \frac{q}{12} \ \leq \
             \max\left\{0,\frac{2q-5}{4}\right\} \ . 
       \end{equation}
\listend
\end{lem}
We recall that we use conditions
(\ref{cond:gamma0}), (\ref{cond:alphagamma0}), (\ref{cond:hfunctions1}) and
(\ref{cond:hfunctions2}) as standing assumptions in this subsection.

\ \\
\proof\ 
Recall that $([1,j] \setminus R_N ) = ([1,j] \cap \Gamma_N)$.
\alphlist
\item This is a direct consequence of (\ref{eq:gammaninominfty}) and
  (\ref{cond:hfunctions2}). For the second inequality in (\ref{eq:regsets3}),
  note that $\# ( [1,l^+_q] \setminus R_N ) = 0$ whenever $q < 12$.
\item 
  We prove the following statement by induction on $q$:
  \begin{equation}
          \label{eq:gammasets1}
             \forall j \in [0,l^+_q-1]  \ \exists n \in [j+1,l^+_q]: \ \#
              [j+1,n] \cap \Gamma_{l^+_q}\ \leq \ \sigma\cdot (n-j) \ .
  \end{equation}
  This obviously implies the statement, as it ensures the existence of
  a partition of $[j+1,l^+_q]$ into disjoint intervals $I_i =[j_i+1,j_{i+1}]$
  with $j=j_1 < j_2 < \ldots < j_k=l^+_q$ which all satisfy
  \[ 
          \# \left(I_i \cap \Gamma_{l^+_q}\right) \ \leq \ \sigma
          \cdot (j_{i+1} - j_i) \ .
  \]
  If $q=1$, then (\ref{eq:gammasets1}) is obvious as
  $\Gamma_{l^+_1} \ssq \Lambda_{l^+_1} = \emptyset$ (see
  (\ref{eq:gammaninominfty}) and note that $l^+_1 \leq \nu(1)$ by
  (\ref{eq:lplusestimates})). Now suppose (\ref{eq:gammasets1}) holds
  for all $q \leq p$. In order to show (\ref{eq:gammasets1}) for
  $p+1$, we have to distinguish three possible cases. Recall that by
  (\ref{eq:gammanequality}) and (\ref{eq:gammaninominfty})
  \[
  \Gamma_{l^+_{p+1}} \ = \ \bigcup_{J\in\J_{l^+_{p+1}}} J^- \cup
  \Gamma^+(J) \ \ssq \ \Lambda_{l^+_{p+1}} \ .
  \]
  If $j+1 \notin \Gamma_{l^+_{p+1}}$ we can choose $n =
  j+1$.
  
  If $j+1 \in \Gamma^+(J)$ for some $J \in {\cal J}_{l^+_{p+1}}$ then
  $p_J \leq p$ as $l^+_{p+1} < \nu(p)$ by
  (\ref{eq:lplusestimates}). By (\ref{eq:recursivegammas}) there holds
  $j - m_J \in \Gamma_{l^+_{p_J}} \ssq [0,l^+_{p_J}-1]$. Thus we can
  apply the induction statement with $q=p_J$ to $j-m_J$ and obtain some
  $\tilde{n} \in [j - m_J +1,l^+_{p_J}]$ with
  \[
  \#([j-m_J+1,\tilde{n}] \cap \Gamma_{l^+_{p_J}}) \ \leq \ 
  \sigma \cdot (\tilde{n}-j+m_J) \ .
  \]
  As $\Gamma^+(J) = \Gamma_{l^+_{p_J}} + m_J$ (again by
  (\ref{eq:recursivegammas})), $n := \tilde{n} + m_J$ has the required
  property.
  
  Finally, if $j+1 \in J^-$ for some $J \in {\cal
    J}_{l^+_{p+1}}$ then $[\lambda^-(m_J),j+1] \ssq J^- \ssq
  \Gamma_{l^+_{p+1}}$. Therefore
  \begin{eqnarray} \label{eq:regsets3b}
             \lefteqn{ \frac{ \# \left( [j+1,\lambda^+(m_J)] \cap
             \Gamma_{l^+_{p+1}} \right) }{\lambda^+(m_J) - j}  \ \leq \  
             \frac{ \# \left(J \cap \Gamma_{l_{p+1}^+}
             \right) }{\# J } } \\
             & \leq & \frac{\#(J^- \cup \Gamma^+(J))}{(u+v)\cdot p_J}  \
             \leq \ \frac{(u+2) \cdot p_J
             + \# \Gamma_{l^+_{p_J}}}{(u+v)\cdot p_J} \ \leq \frac{u+3}{u+v} 
             \ . \nonumber
  \end{eqnarray}
  where we used part (a) of this lemma with $j=N=l^+_{p_J}$ to
  conclude that $\# \Gamma_{l^+_{p_J}} \leq p_J$.
\item Similar to (a), we prove
that
  \[
  \forall j \in [0,\lambda^+-1] \ \exists n \in [j+1,\lambda^+] : \ \#([j+1,n] \cap
  \Gamma_N) \ \leq \ \sigma(n-j) \ .
  \]
  Again, we have to distinguish three cases:
 
  If $j+1 \notin \Gamma_N$ we can choose $n = j+1$.
      
  If $j+1 \in \Gamma^+(I)$ for some $I \in \JN$, then we can choose $n =
       \lambda^+(m_I) = m_I + l^+_{p_I}$. Using that $\Gamma^+(I) =
       \Gamma_{l^+_{p_I}} + m_I$ by (\ref{eq:recursivegammas}), part (b)
       implies that $n$ has the required property.
       
       If $j+1 \in I^-$ for some $I \in \JN$ we can choose $n =
       \lambda^+(m_I)$ and proceed exactly as in (\ref{eq:regsets3b}),
       with $J$ being replaced by $I$.
\item By Lemma~\ref{lem:endpoints}(b) there holds $m-l^-_q \in
  A_m$. Therefore (\ref{eq:gammanequality}), $\Gamma(J) \ssq J$
  and (\ref{eq:jminominfty}) imply that
  \[
         [m-l^-_q,m] \cap \Gamma_m  \ \ssq  \bigcup_{j \in
         [m-l^-_q+1,m-1]} \hspace{-2eM}
         \Omega_{q-2}(j) \  =: \ U \ .
  \]
  Further, Lemma~\ref{lem:endpoints}(b) yields that $U \ssq \Omtil_{q-2}$,
  such that $U-m \ssq \Omega_\infty$ by
  Lemma~\ref{lem:omegatransition}(b). Consequently
  \begin{eqnarray*}
    \lefteqn{\#\left([m-l^-_q,m] \cap \Gamma_m\right) \ \leq \ \#U \ = \
    \#(U-m)} \\ & \leq &
    \#\left([-l^-_q,-1] \cap \Omega_\infty\right) \
    \stackrel{(\ref{cond:hfunctions2})}{\leq} \ \frac{l^-_q}{12w} \ \leq \
    \frac{q}{12} \ \leq \ \max\left\{0,\frac{2q-5}{12}\right\} \ . 
  \end{eqnarray*}
\listend

\qed


\section{Construction of the sink-source orbits: One-sided forcing} \label{Construction}

We now turn to the construction of the sink-source-orbits in the case of
one-sided forcing. Before we start with the core part of the proof, we have to
add some more assumptions on the parameters. Further, we restate two estimates
from the preceding section, together with a few other facts that will be used
frequently in the construction.
\medskip 

First of all, we choose $u$ and $v$ such that
\begin{eqnarray}
  u & \geq & 8  \label{cond:u} \ , \\  v & \geq & 8 \label{cond:v} \
  , \\
  \label{cond:sigma} 
  \sigma & \leq & \textstyle \frac{1}{6} \ .              
\end{eqnarray}
In addition, we assume that 
\begin{eqnarray} \textstyle
    \halb \walpha & \geq & 6 + K\cdot S_\infty(\alpha^{\viertel}) \ \
    . \label{cond:alpha3} 
\end{eqnarray}
Further, we remark that (\ref{cond:alphagamma0}) implies
\begin{eqnarray}
   \alpha & \geq & 4S_\infty(\alpha) \ . \label{cond:alpha4}
\end{eqnarray}
    
Now suppose that (\ref{cond:gamma0}), (\ref{cond:alphagamma0}),
  (\ref{cond:hfunctions1}) and (\ref{cond:hfunctions2}) hold. Together with
  the above assumptions and the respective results from the last section (see
  (\ref{eq:lplusestimates}), Lemma~\ref{lem:regsets}(b) and
  (\ref{eq:lplusstabilized})), this yields that for any $q \geq 1$ the
  following estimates hold:
\begin{eqnarray}
  \label{eq:lpmestimates}
   4(q+1) \ \leq \ 8 q \ \leq \  l^\pm_q &  <  &
     \nutil(\max\{1,q-2\}) \  \leq  \ \nu(\max\{1,q-2\}) \ \\
   \#([j+1,l^+_q] \smin R_N) & \leq & \frac{l^+_q-j}{6} \ \ \ \ \
   \forall N \geq l^+_q, \ j\in [0,l^+_q-1] \label{eq:etajnestimate} \ .
\end{eqnarray}
Recall that $(\xi_n(\beta,l))_{n\geq -l}$ corresponds the forward
orbit of the point $(\omega_{-l},3)$ under the transformation
$T_\beta$, where we suppress the $\theta$-coordinate (see Definition
\ref{def:parameterfamily}). As we are in the case of one-sided
forcing, we can use the fact that for all $l,n \in \Z, n \geq -l$ the
mapping $\beta \mapsto \xi_n(\beta,l)$ is monotonically decreasing in
$\beta$. For $l \geq 0$ and $n\geq 1$ the monotonicity is even strict
(as $g(0) = 1 > 0$ and $F$ is strictly increasing by
(\ref{cond:Funiformbounds})). This has some very convenient implications. First
of all, we can uniquely define parameters $\beta^+_{q,n}$ and
$\beta^-_{q,n}$ ($q,n \in \N$) by the equations
\begin{equation}
  \label{eq:betapluspndef}
  \xi_n(\beta_{q,n}^+,l^-_q) \ = \ \alphtel 
\end{equation}
and
\begin{equation}
  \label{eq:betaminpndef}
  \xi_n(\beta_{q,n}^-,l^-_q) \ = \ -\alphtel \ .
\end{equation}
In addition, we let
\begin{equation}
  \label{eq:lpmzero}
  l^-_0 \ := \ 0 \ \ \ \ \textrm{ and } \ \ \ \ l^+_0 \ := \ 0 \ 
\end{equation}
(note that so far the $l^\pm_q$ had only been defined for $q \geq 1$)
and extend the definitions of $\beta^\pm_{q,n}$ to $q=0$. If we now
want to show that $\xi_n(\beta,l^-_q) \in \Balphcl$ implies
$\xi_j(\beta,l^-_q) \in \Balphcl$ for some $j<n$, we can do so by
proving that
\begin{equation}
  \label{eq:betapluspn}
  \xi_n(\beta_{q,j}^+,l^-_q) \ \geq \ \alphtel 
\end{equation}
and
\begin{equation}
  \label{eq:betaminpn}
  \xi_n(\beta_{q,j}^-,l^-_q) \ \leq \ -\alphtel 
\end{equation}
(compare (\ref{eq:betaplusdef})--(\ref{eq:betaminus})). Furthermore,
(\ref{eq:betaminpn}) is a trivial consequence of the fact that $\kreis
\times [-3,-\alphtel]$ is mapped into $\kreis \times [-3,-(1-\gamma)]
\ssq \kreis \times [-3,-\alphtel)$ (see (\ref{cond:Fmapsover})). Thus,
it always suffices to consider (\ref{eq:betapluspn}).
\medskip

Now we can formulate the induction statement we want to prove: 
\begin{indscheme} \label{thm:indscheme}
Suppose the assumptions of Theorem~\ref{thm:snaexistence} are satisfied and
(\ref{cond:alpha1}), (\ref{cond:hfunctions1}), (\ref{cond:hfunctions2}) and
(\ref{cond:u})--(\ref{cond:alpha3}) hold. Then for any $q \in \N_0$ there holds
\begin{list}{\textbf{\Roman{enumi}.} \ \ }{\usecounter{enumi}}
\item
      If $\xi_{l^+_q+1}(\beta,l^-_q) \in \Balphcl$ then 
      \begin{equation} \label{eq:indstatementI}
                 \xi_j(\beta,l^-_q) \ \geq  \ \gamma \ \ \ \forall j \in
                 [-l^-_q,0] \setminus \Omega_\infty 
      \end{equation}
       and $\beta \in \left[1+\walphtel,1+\frac{3}{\walpha}\right]$.
     \item Suppose $n \in [l^+_q+1,\nu(q+1)]$ is admissible. Then
       $\xi_n(\beta,l^-_q) \in \Balphcl$ implies that (\ref{eq:indstatementI}) holds,
       \begin{equation} \label{eq:indstatementII}
          \ \xi_j(\beta,l^-_q)
       \in \Balphcl \ \ \ \ \ \forall j \in R_n  
       \end{equation}
       and $\beta \in \left[1+\walphtel,1+\frac{3}{\walpha}\right]$.
    \item Let $1 \leq q_1 \leq q$ and suppose $n_1 \in 
      [l^+_{q_1}+1,\nu(q_1+1)]$ and $n_2 \in [l^+_{q}+1,\nu(q+1)]$ 
      are both admissible.
      \alphlist
         \item  If $q_1 = q$ and $n_1 \in R_{n_2}$, then  
         \begin{equation} \label{eq:indstatementIIIa}
              |\beta^+_{q_1,n_1} - \beta^+_{q,n_2}| \ \leq \
              2\alpha^{-\frac{n_1}{4}} \ .
         \end{equation} 
       \item If $q_1 < q$ there holds
         \begin{equation} \label{eq:indstatementIIIb}
              |\beta^+_{q_1,n_1} - \beta^+_{q,n_2}| \ \leq \
              3\cdot \hspace{-0.7eM}\sum_{i=q_1+1}^{q} \alpha^{-i}\ \leq \ \alpha^{-q_1} \ .
         \end{equation} 
      \listend
\end{list}
\end{indscheme}
The proof is given in the next subsection. The statement of
Theorem~\ref{thm:snaexistence} now follows easily, with the help of Lemma~\ref{lem:sinksourceshadowing}:
\medskip

\noindent
\textit{\bf Proof of Theorem \ref{thm:snaexistence}~.}
In order to apply Lemma~\ref{lem:sinksourceshadowing} we can use the
same sequences $l^\pm_p$ as in Induction Scheme \ref{thm:indscheme}.
Further, let $\beta_p := \beta^+_{p,l^+_p+1}$, $\theta_p:=\omega$ and
$x_p := \xi_1(\beta_p,l^-_p)$. From Part II of the induction statement
with $q=p$ and $n=l^+_p+1$ we obtain that
\[
    \xi_j(\beta_p,l^-_p) \ \in \ \Balphcl \ \ \ \ \
    \forall j \in R_{l^+_p+1}  \ ,
\]
and
Lemma~\ref{lem:regsets}(a) implies that 
\[ 
\# \left([1,j] \cap R_{l^+_p+1}\right) \ \geq \
\frac{11}{12} \cdot j \ \ \ \ \ \forall j \in [1,l^+_p] \ . 
\]
Therefore it follows from (\ref{cond:Funiformbounds}) and
(\ref{cond:Fexpansion}) that
\begin{eqnarray*}
    \lefteqn{\lambda^+(\beta_p,\theta_p,x_p,j) \ = } \\ & = &
    \frac{1}{j} \sum_{i=1}^{j} 
    \log F'(\xi_i(\beta_p,l^-_p)) 
    \ \geq \ \frac{11}{12} \cdot \frac{\log \alpha}{2} - \frac{2\log 
      \alpha}{12} \ = \ \frac{7}{24} \cdot \log\alpha \ \ \ \ \forall
    j \in [1,l^+_p] \ .
\end{eqnarray*}
Likewise, we can conclude from Part I of the induction statement with
$q=p$ in combination with (\ref{cond:hfunctions2}),
(\ref{cond:Funiformbounds}) and (\ref{cond:Fcontraction}) that
\[
\lambda^-(\beta_p,\theta_p,x_p,j) \ \geq \ \frac{7}{24} \cdot
\log\alpha \ \ \ \ \forall 
    j \in [1,l^-_p] \ .
\]
Consequently, the assertions of Lemma~\ref{lem:sinksourceshadowing}
are satisfied, such that there is at least one parameter value at
which a sink-source-orbit and consequently an SNA and an SNR occur
(see Theorem~\ref{thm:sinksourcesna}). Due to
Theorem~\ref{thm:saddlenode}, the only parameter where this is
possible is the critical parameter $\beta_0$. Finally the statement
about the essential closure again follows from
Theorem~\ref{thm:saddlenode}~.

\qed
\bigskip

\noindent
\textit{\bf Proof of Addendum~\ref{adde}.}  Define $\beta_p$, $\theta_p$ and
$x_p$ as above.  From Part III of the induction statement it follows that
$(\beta_p)_{p\in\N}$ is a Cauchy-sequence and must therefore converge to
$\beta_0$ (instead of only having a convergent subsequence, as in the proof of
Lemma~\ref{lem:sinksourceshadowing}). To be more precise, if $p < q$ we have
$|\beta_p-\beta_q| \leq \alpha^{-p}$, such that
\begin{equation}
  \label{eq:betapdistance}
  |\beta_p - \beta_0| \ \leq \ \alpha^{-p} \ \ \ \ \ \forall p \in \N
  \ .
\end{equation}
Further, let
\[
\theta_0 \ := \ \omega 
\] 
and
\begin{equation} \label{eq:xplimit}
x_0 \ := \ \lim_{p\ra\infty} x_p \ .
\end{equation}
If the limit in (\ref{eq:xplimit}) does not exist,%
\footnote{In fact it is possible to show that $(x_p)_{p\in\N}$ is a
  Cauchy-sequence as well, by using Lemma~\ref{lem:orbitscontraction}
  and Part I of the induction statement. However, we refrain from
  doing so as this is not relevant for the further argument.}
we just go over to a suitable subsequence. From Part II of the
induction statement with $q=p$ and $n=l^+_p+1$, it follows that
\[
    T_{\beta_p,\omega,j-1}(x_p) \ = \ \xi_j(\beta_p,l^-_p) \ \in
    \Balphcl \ \ \ \ \ \forall j \in R_{l^+_p+1} \ .
\]
Using that $R_{l^+_p+1} \ssq R_{l^+_q+1} \ \forall q \geq p$ by
(\ref{eq:lplusstabilized}) and the continuity of the map $(\beta,x)
\mapsto T_{\beta,\omega,j-1}(x)$, we see that
\begin{equation}
  \label{eq:xoinBalpha}
  T_{\beta_0,\omega,j-1}(x_0) \ \in \ \Balphcl \ \ \ \ \ \forall j \in
  R_{l^+_p+1}, p \in \N \ .
\end{equation}

Now $\varphi^+$ and $\psi$ can be defined pointwise as the upper bounding
graph (see (\ref{eq:boundinggraphs}) of the system $T_{\beta_0}$ and by
equation (\ref{eq:psidef}), respectively. Then the fact that
\begin{equation}
\label{eq:xogeqpsiomega} \psi(\omega) \ \leq \  x_0
\end{equation}
is obvious, otherwise the forward orbit of $(\theta_0,x_0)$ would
converge to the lower bounding graph $\varphi^-$ and its forward
Lyapunov exponent would therefore be negative. On the other hand
suppose
\[
   \psi(\omega) \geq x_0 - \alpha^{-p}
\]
for some $p \geq 2$. Then we can compare the orbits 
\begin{equation}
  x^1_1 \ld x^1_n \ := \ x_0 \ld T_{\beta_0,\omega_1,l^+_p-1}(x_0)
\end{equation}
and 
\begin{equation}
   x^2_1 \ld x^2_n \ := \ \psi(\omega_1) \ld \psi(\omega_{l^+_p})
\end{equation}
via Lemma~\ref{lem:orbitsthrowout}(b)%
\footnote{We can choose $\epsilon=\frac{\alpha^{-p}}{2}$, such that
  $q=p$. Note that the error term is zero, as we consider orbits which
  are located on the same fibres and generated with the same
  parameter. As $l^+_p+1 \in R_{l^+_p+1}$, \ $x^1_{n+1} \in \Balphcl$
  follows from (\ref{eq:xoinBalpha}). $\tau(n) \leq \frac{2p-3}{4}$
  and $\tau(j) \leq \frac{j}{8}$ follow from
  Lemma~\ref{lem:regsets}(a), whereas $\tau(n)-\tau(j) \leq
  \frac{n-j}{6}$ follows from (\ref{eq:etajnestimate}). Finally
  $n=l^+_p \geq 5p$ by (\ref{eq:lpmestimates}). }
and obtain that $\psi(\omega_{l^+_p+1}) \leq - \frac{2}{\alpha}$. But
as we have seen in the proof of Theorem~\ref{thm:saddlenode} that all
points below the 0-line eventually converge to the lower bounding
graph, this contradicts the definition of $\psi$. Consequently 
\[
   x_0 \ \leq \ \psi(\omega) + \alpha^{-p} \ \ \ \forall p \in \N \ .
\]
Together with (\ref{eq:xogeqpsiomega}) this implies that $x_0 =
\psi(\omega)$. 

\ \\
As $\psi \leq \varphi^+$, we immediately obtain $x_0 \leq
\varphi^+(\omega)$, such that it remains to show
\begin{equation}
  \label{eq:xoabove}
  x_0 \ \geq \ \varphi^+(\omega)  \ .
\end{equation}
To that end, we denote the upper boundary lines of the system
(\ref{eq:generalsystem}) by $\varphi_n$ if $\beta = \beta_0$ and by
$\varphi_{p,n}$ if $\beta=\beta_p$. Now either infinitely many $\beta_p$
are below $\beta_0$, or infinitely many $\beta_p$ are above $\beta_0$.
Therefore, by going over to a suitable subsequence if necessary, we
can assume w.l.o.g.\ that either $\beta_p \leq \beta_0 \ \forall p
\in\N$ or $\beta_p \geq \beta_0 \ \forall p \in \N$. The first case is
treated rather easily: If $\beta_p \leq \beta_0$, then
\[
    x_p \ = \ \xi_1(\beta_p,l^-_p) \ = \ \varphi_{p,l^-_p+1}(\omega) \ 
    \geq \ \varphi_{l^-_p+1}(\omega) \ \geq \ \varphi^+(\omega) \ .
\]
Passing to the limit $p \ra \infty$, this proves (\ref{eq:xoabove}). 
 
\ \\
On the other hand, suppose $\beta^p \ \geq \beta_0$. In this case, we
will show that
\begin{equation} \label{eq:xodistance}
 |x_p - \varphi_{l^-_p+1}(\omega)| \ = \ |\xi_1(\beta_p,l^-_p) -
 \xi_1(\beta_0,l^-_p)| \ \leq \ \alpha^{-p} \cdot \left(6+K\cdot
   S_\infty(\alpha^{\viertel}) \right)  \ .
\end{equation}
As $\varphi_n(\omega) \nKonv \varphi^+(\omega)$ and $x_p \pKonv x_0$,
this again proves (\ref{eq:xoabove}). Note that as $\beta_p \geq
\beta_0$ we have $\xi_j(\beta_0,l^-_p) \geq \xi_j(\beta_p,l^-_p) \ 
\forall j \geq -l^-_p$, such that $\xi_j(\beta_p,l^-_p) \geq \gamma$
implies $\xi_j(\beta_0,l^-_p) \geq \gamma$. This allows to compare the 
orbits 
\begin{equation}
  x^1_1 \ld x^1_n \ := \ \xi_{-l^-_p}(\beta_p,l^-_p) \ld \xi_0(\beta_p,l^-_p)
\end{equation}
and 
\begin{equation}
  x^1_1 \ld x^1_n \ := \ \xi_{-l^-_p}(\beta_0,l^-_p) \ld \xi_0(\beta_0,l^-_p)
\end{equation}
via Lemma~\ref{lem:orbitscontraction},%
\footnote{With $\epsilon = \alpha^{-p}$. We have $|\beta_p-\beta_0|
  \leq \alpha^{-p}$ by (\ref{eq:betapdistance}), such that
  $\err(\ldots) \leq \epsilon$. $\eta(j,n) \leq \frac{n+1-j}{10}$
  follows from Part I of the induction statement with $q=p$ together
  with (\ref{cond:hfunctions2}) and $0\notin \Omega_\infty$. Finally
  $n = l^-_p+1 \geq 4p$ by (\ref{eq:lpmestimates}), such that
  $\alpha^{-\frac{n}{4}} \leq \epsilon$. } which yields
(\ref{eq:xodistance}).

\qed

%

\subsection{Proof of the induction scheme} \label{Onesidedproof}

\textbf{Standing assumption:} In this whole subsection, we always assume that the
assumptions of the Induction scheme \ref{thm:indscheme} are
satisfied. \medskip

 Before we start the proof of the Induction statement, we provide the
following lemma, which will be used in order to obtain estimates on the
parameters $\beta^+_{q,n}$:
\begin{lem}
\label{lem:parameterestimates} Suppose   Let $n$ be admissible and
$\xi_n(\beta_1,l),\xi_n(\beta_2,l) \in \Balphcl$. Further, suppose that
$\xi_n(\beta,l) \in \Balphcl$ implies $\xi_j(\beta,l) \in \Balphcl \ \forall j
\in R_n$. Then
\[
      |\beta_1-\beta_2| \ \leq \ 2\alpha^{-\frac{n}{4}} \ .
\] 
\end{lem}
\proof\ 
Note that 
\begin{eqnarray}
\lefteqn{\frac{\partial}{\partial \beta} \xi_{j+1}(\beta,l) \ = \ \label{eq:xiderivatives}
\frac{\partial}{\partial \beta} \left( F(\xi_j(\beta,l)) - \beta
  \cdot g(\omega_j) \right)} \\ \nonumber 
  & = & F'(\xi_j(\beta,l)) \cdot
\frac{\partial}{\partial \beta} \xi_j(\beta,l) - g(\omega_j)
\ \stackrel{(g \geq 0)}{\leq} \
 F'(\xi_j(\beta,l)) \cdot
\frac{\partial}{\partial \beta} \xi_j(\beta,l)\ . 
\end{eqnarray}
W.l.o.g.\ we can assume $\beta_1 < \beta_2$. As we have
$\frac{\partial}{\partial \beta} \xi_0(\beta,k) \leq -1$, the
inductive application of (\ref{eq:xiderivatives}) together with
(\ref{cond:Funiformbounds}) and (\ref{cond:Fexpansion}) yields
\[
\frac{\partial}{\partial \beta} \xi_n(\beta,l) \ \leq \ 
-(\alpha^{\halb})^{\#([1,n-1]\cap R_n])} \cdot (\alpha^{-2})^{\#(
  [1,n-1]\setminus R_n)} \ = \ -\alpha^{\halb(n-1-5\cdot \#\Gamma_n)}
\]
as long as $\xi_n(\beta,k) \in \Balphcl$. (Recall that $[1,n] \setminus
R_n = \Gamma_n$ by definition and $n \in R_n$ by assumption.) In particular this is
true for all %
$\beta \in [\beta_1,\beta_2]$. Lemma~\ref{lem:regsets}(a) yields %
$\# \Gamma_n = \# ([1,n-1] \smin R_n) \leq \frac{n-1}{10}$, such
that we obtain
\[
\frac{\partial}{\partial \beta} \xi_n(\beta,l) \ \leq \ -\alpha^{\frac{n-1}{4}}
\]
The required estimate now follows from $|\xi_n(\beta_1,l) -
\xi_n(\beta_2,l)| \leq \frac{2}{\alpha}$.

\qed
\bigskip

We prove the Induction scheme \ref{thm:indscheme} by induction on $q$, 
proceeding in six steps. The first one starts the induction:

\ \\
\underline{\textbf{Step 1:}} \ \ \ \textit{Proof of the statement for
  $q=0$.} \\ \ \\
As $d(\omj,0) \geq \frac{4\gamma}{L_2} \ \forall j \in [1,\nu(1)-1]$,
Part I and II of the induction statement are already contained in
Lemma~\ref{lem:inductionstart}, and Part III is still void.

\solidqed


\ \\ Now let $p \geq 1$ and assume that the statement of Induction
scheme \ref{thm:indscheme} holds for all $q \leq p-1$. We have to show
that the statement then holds for $p$ as well. The next two steps
will prove Part I of the induction statement for $p$. Note that for
$p=1$ Part I of the induction statement is still contained in
Lemma~\ref{lem:inductionstart} as $l^\pm_1 < \nu(1)$ by
(\ref{eq:lpmestimates}). Therefore, we can assume
\begin{equation}
  \label{eq:pgeqone}
  p \ \geq 2
\end{equation}
during Step 2 and Step 3.

\ \\ \underline{\textbf{Step 2:}} \ \ \ 
\textit{If
  $|\beta-\beta^+_{p-1,\nu(p)}| \leq \alpha^{-p}$, then
  $\xi_j(\beta,l^-_p) \geq \gamma \ \forall j \in [-l^-_p,0] \setminus
  \Omega_\infty$.}

\ \\
This is a direct consequence of the following lemma with $q=p$, $l^* =
l^-_{p-1}$, $l = l^-_p$, $\beta^* = \beta^+_{p-1,\nu(p)}$,
$m=\nu(p)$ and $k=-\nu(p)$,  .%
\footnote{Note that $\nu(p)$ is admissible by
  Lemma~\ref{lem:endpoints}(d). Therefore $\beta^* = \beta^+_{p-1,\nu(p)} \in
  \left[1+\walphtel,1+\frac{3}{\walpha}\right]$ and $\xi_j(\beta^*,l^-_{p-1}) \in
  \Balphcl \ \forall j \in R_{\nu(p)}$ follow from Part II of the
  induction statement with $q=p-1$ and $n=\nu(p)$. \label{foot:step2a} }
Note that $\Omtil_{p-2} - \nu(p) \ssq \Omtil_{p-1} \ssq
\Omega_\infty$ by Lemma \ref{lem:omegatransition}(b). The statement of
the lemma is slightly more general because we also want to use it in
similar situations later. Recall that $\Omtil_{-1} = \Omtil_0 =
\emptyset$, see (\ref{eq:omtilmin}). 
\begin{lem}
  \label{lem:approaching}
  Let $q\geq1$, $l^*,l \geq 0$, $\beta^* \in
  \left[1+\walphtel,1+\frac{3}{\walpha}\right]$ and $|\beta - \beta^*| \leq 2\alpha^{-q}$.
  Suppose that $m$ is admissible, $p(m) \geq q$ and either $k=0$ or
  $p(k) \geq q$. Further, suppose
\[
      \xi_j(\beta^*,l^*) \ \in \ \Balphcl \ \ \ \forall j \in R_m \
\] 
and $\xi_{m+k-l^-_q}(\beta,l) \geq \gamma$. Then 
\[
        \{ j \in [m-l^-_q,m] \mid \xi_{j+k}(\beta,l) < \gamma \} \ \ssq \
        \Omtil_{q-2} \ .
\]
\end{lem}
\ \\
\proof\ 
We have that $J(m) \in \J$ by Lemma~\ref{lem:admissiblecentral}, such
that we can apply Lemma~\ref{lem:endpoints}(b) to $J:=J(m)$. Note that
$m=m_J$ and $p_J = p(m) \geq q$ in this case. Let $t := m -
l^-_q$. We will show that
\[
    \{ j \in [t,m] \mid \xi_{j+k}(\beta,l) < \gamma \} \ \ssq \ \bigcup_{t
    \leq j < m} [\lambda^-(j),\lambda^+(j)+1] \ .
\]
As $[\lambda^-(j),\lambda^+(j)+1] \ssq \Omega_0(j) \ssq
\Omega_{q-2}(j)$ and $Q_{q-2}(j) \leq Q_\infty(j) \leq \max\{0,q-2\} \ \forall j
\in [t,m-1]$ (see Lemma \ref{lem:endpoints}(b)), this proves the
statement.

\ \\
Let $J_1 \ld J_r$ be the ordered sequence of intervals $J \in \J_m$ with $J
\ssq [t,m]$, such that 
\[
    [t,m] \setminus A_m \ = \ [t,m] \cap \Lambda_m \ = \
    \bigcup_{i=1}^r J_i  
\]
(recall that $[1,m] \setminus A_m = \Lambda_m$ by definition).
Further, define
\[
    j^-_i \ := \ \lambda^-(m_{J_i}) \ \ \ \textrm{and } \ \ j^+_i \ := \
    \lambda^+(m_{J_i})  \ ,
\]
such that $J_i = [j_i^-,j^+_i]$. We have to show that
$\xi_{j+k}(\beta,l) \geq \gamma$ whenever $j$ is contained in
$[j_i+2,j^-_{i+1}-1]$ for some $i = 1 \ld r$ or in %
$[t,j^-_1-1] \cup [j^+_r+2,m]$, and we will do so by induction on $i$.
The case where $j^+_i+1 = j^-_{i+1}-1$, such that
$[j^+_i+2,j^-_{i+1}-1]$ is empty, is somewhat special and will be
addressed later, so for now we always assume $j^+_i+1 < j^-_{i+1}-1$.
  
\ \\
Let us first see that $\xi_{j^+_i+2+k}(\beta,l) \geq \gamma$ implies %
$\xi_{j+k}(\beta,l) \geq \gamma$ %
$ \forall j \in [j^+_i+2,j^-_{i+1}-1]$ : If %
$j \in [j^+_i+2,j^-_{i+1}-1]$, then $j \in A_m$. Hence $d(\omega_j,0)
\geq \frac{4\gamma}{L_2}$ and therefore 
\[
    d(\omega_{j+k},0) \ \geq \ \frac{4\gamma-\alpha^{-q}}{L_2} \geq
    \frac{3\gamma}{L_2}
\] 
by (\ref{cond:alphagamma0}) if $q \geq 2$. In case $q=1$ we obtain the
same result, as $\nutil(1) > l^-_1$ by (\ref{eq:lpmestimates}) then
implies
that $d(\omj,0) \geq \frac{8\gamma}{L_2} \ \forall j \in [t,m]$. Further %
$\beta \in [1,1+\frac{4}{\walpha}]$ as $|\beta-\beta^*| \leq
2\alpha^{-q} \leq \frac{2}{\alpha}
\stackrel{(\ref{cond:alphagamma0})}{\leq} \walphtel$ and %
$\beta^* \in [1+\frac{1}{\walpha},1+\frac{3}{\walpha}]$. Inductive
application of Lemma \ref{lem:basicestimate} therefore yields
\[
    \xi_{j+k}(\beta,l) \ \geq \gamma \ \ \ \forall j \in
    [j^+_i+2,j^-_{i+1}-1] \ . 
\]
The same argument also starts and ends the induction: As
$\xi_{t+k}(\beta,l) \geq \gamma$ by
assumption we get %
$\xi_{j+k}(\beta,l) \geq \gamma \ \forall j \in [t,j^-_1-1]$, and for
$j \in [j^+_r+2,m]$ this follows from %
$\xi_{j^+_r+2+k}(\beta,l) \geq \gamma$. 

\ \\
If $q=1$, then Lemma~\ref{lem:endpoints}(b) yields that %
$p(j) = 0 \ \forall j \in [t,m]$ and consequently %
$[t,m] \setminus A_m = \emptyset$, such that we are already finished 
in this case. Therefore, we can assume from now on that $q \geq 2$. It
remains to prove that
\begin{equation} \label{eq:approaching}
  \xi_{j^-_i-1+k}(\beta,l) \geq \gamma \ \ \textrm{ implies } \ \
  \xi_{j^+_i+2+k}(\beta,l) \geq \gamma \ .
\end{equation} In order to do this,
we have to apply Lemma \ref{lem:orbitsthrowout}(a): Let %
$\epsilon := \alpha^{-(q-1)}$ and choose
\begin{equation} \label{step2:reference}
   x^1_1 \ld x^1_n \ := \ \xi_{j^-_i-1}(\beta^*,l^*) \ld
   \xi_{j^+_i}(\beta^*,l^*) 
\end{equation}
and 
\begin{equation} \label{step2:new}
x^2_1 \ld x^2_n \ := \ \xi_{j^-_i-1+k}(\beta,l) \ld
\xi_{j^+_i+k}(\beta,l)\ .
\end{equation}
As $d(\omega_k,0) \leq \frac{\alpha^{-(q-1)}}{L_2}$ and %
$|\beta-\beta^*| \leq 2\alpha^{-q}$ we have
$\err(\beta_1,\beta_2,\theta_1,\theta_2) \leq K\cdot \epsilon$ by
Remark \ref{bem:errorterm}~. Further $x^1_1 =
\xi_{j_i^--1}(\beta^*,l^*) \in \Balphcl$ and $x^1_{n+1} =
\xi_{j_i^++1}(\beta^*,l^*) \in \Balphcl$ by assumption (as
$j^-_i-1,j^+_i+1 \in A_m \ssq R_m$), whereas %
$x^2_1 = \xi_{j^-_i-1+k}(\beta,l)\geq \gamma$ $\geq \frac{2}{\alpha}$.
Applying Lemma~\ref{lem:regsets}(d) we obtain that
\[
    \tau(n)  =  \#([j^-_i-1,j^+_i] \setminus R_m)  \leq \#([t,m]
    \setminus R_m) \stackrel{(\ref{eq:approachfraction})}{\leq} \min\left\{
     0,\frac{2q-5}{4}\right\} 
\]
Finally, we have 
\[
  |\tau(n) - \tau(j)| \leq \#([j_i^+-(n-j)+1,j^+_i]) \setminus R_m \
  \leq \ -\sigma\cdot(n-j) \leq \frac{n-j}{6} 
\] 
by Lemma \ref{lem:regsets}(c) (with $N=m$, $J=J_i$ and $\lambda^+ =
\lambda^+(m_{J_i}) = j^+_i$). Thus all the assumptions of
Lemma~\ref{lem:orbitsthrowout} are satisfied and we can conclude that
$x_{n+1}^2 = \xi_{j^+_i+1+k}(\beta,l) \geq \frac{2}{\alpha}$. As we have
$d(\omega_{j^+_i+1+k},0) \geq \frac{3\gamma}{L_2}$ again, Lemma
\ref{lem:basicestimate} now implies %
$\xi_{j^+_i+2+k}(\beta,l)$ $\geq \gamma$.

As mentioned, we still have to address the case where $j^+_i+1 =
j^-_{i+1}-1$, such that $[j^+_i+2,j^-_{i+1}-1]$ is empty. In this case
we still obtain that $\xi_{j^+_i+1+k}(\beta,l) =
\xi_{j^-_{i+1}-1+k}(\beta,l) \geq \frac{2}{\alpha}$. But this is
sufficient in order to apply Lemma~\ref{lem:orbitsthrowout}(a) once more,
in exactly the same way as above, to conclude that
$\xi_{j^+_{i+1}+1+k}(\beta,l) \geq \frac{2}{\alpha}$.  Thus in the
next step we obtain $\xi_{j^+_{i+1}+2+k}(\beta,l) \geq \gamma$ as
before, unless again $j^+_{i+1}+1 = j^-_{i+2}-1$. In any case, the
induction can be continued.

\qed

\solidqed


\ \\ \underline{\textbf{Step 3:}} \ \ \ \textit{$\xi_{l^+_p+1}(\beta,l^-_p)
\in \Balphcl$ \ implies \ $|\beta - \beta^+_{p-1,\nu(p)}|  \leq  \alpha^{-p}$.} 

\ \\ Recall that we can assume $p \geq 2$, see (\ref{eq:pgeqone}). Let
$\beta^* := \beta^+_{p-1,\nu(p)}$, $\beta^+:=\beta^*-\alpha^{-p}$ and
$\beta^- := \beta^*+ \alpha^{-p}$. We prove 
\begin{claim}
  \label{claim:step3}
  \hspace{5eM} $\xi_{l^+_p+1}(\beta^+,l^-_p) > \alphtel$ \ .
\end{claim}
As $\xi_{l^+_p+1}(\beta^-,l^-_p) < -\alphtel$ follows in
exactly the same way, this implies the statement.

\ \\ \textit{Proof of the claim:} \\
Using Step 2, we see that 
\begin{equation} \label{step3:a}
\xi_j(\beta^+,l^-_p) \ \geq \ \gamma \ \ \ \ \ \forall j \in [-l^-_p,0]
\setminus \Omega_\infty \ .
\end{equation}
On the other hand, from Part II of the the induction statement with
$q=p-1$ and $n=\nu(p)$ it follows that%
\footnote{Note that $\nu(p)$ is admissible by
  Lemma~\ref{lem:endpoints}(d) and $\xi_{\nu(p)}(\beta^*,l^-_{p-1}) \in
  \Balphcl$ by definition of $\beta^*= \beta^+_{p-1,\nu(p)}$.
\label{foot:step3a} }
\begin{equation} \label{step3:b}
\xi_j(\beta^*,l^-_{p-1}) \ \geq \ \gamma \ \ \ \ \ \forall j \in
[-l^-_{p-1},0] \setminus \Omega_\infty 
\end{equation}
and
\begin{equation}
  \label{step3:c}
  \xi_j(\beta^*,l^-_{p-1}) \ \in \ \Balphcl \ \ \ \ \ \forall j \in
  R_{\nu(p)} \ .
\end{equation}
Thus we can use Lemma \ref{lem:orbitscontraction} with $\epsilon =
\alpha^{-p}$ to compare the sequences
\begin{equation}
  \label{step3:referenceone}  
   x^1_1 \ld x^1_n \ := \
\xi_{-l^-_{p-1}}(\beta^*,l^-_{p-1}) \ld
\xi_{-1}(\beta^*,l^-_{p-1})
\end{equation} 
and 
\begin{equation}
  \label{step3:new}
  x^2_1 \ld x^2_n :=
\xi_{-l^-_{p-1}}(\beta^+,l^-_p) \ld \xi_{-1}(\beta^+,l^-_p)
\end{equation}
and obtain that%
\footnote{As the two orbits lie on the same fibres and
  $\beta^*-\beta^+ = \alpha^{-p}$, we have $\err(\ldots) \leq K
  \cdot \epsilon$, see Remark~\ref{bem:errorterm}~. Further, by
  (\ref{step3:a}) and (\ref{step3:b}) we have $\eta(j,n) \leq
  \#([-(n-j),-1] \cap \Omega_\infty) \leq \frac{n+1-j}{10}$ by
  (\ref{cond:hfunctions2}) and $n = l^-_{p-1} \geq 
  4p$ by (\ref{eq:lpmestimates}).\label{foot:step3b} }
\[
|\xi_0(\beta^+,l^-_p) - \xi_0(\beta^*,l^-_{p-1})| \ \leq \ 
\alpha^{-p} \cdot (6 + K\cdot S_\infty(\alpha^{-\viertel})) \ .
\]
Note that (\ref{step3:a}) and (\ref{step3:b}) in particular imply that 
both $\xi_0(\beta^*,l^-_{p-1}) \geq \gamma$ and $\xi_0(\beta^+,l^-_p)
\geq \gamma$. Therefore we can use (\ref{cond:Fcontraction}) to obtain
\begin{eqnarray*}
    \xi_1(\beta^+,l^-_p) & \geq & \xi_1(\beta^*,l^-_{p-1}) +
    (\beta^*-\beta^+) - \alpha^{-p} \cdot \frac{6+K\cdot
    S_\infty(\alpha^{\viertel})}{2\walpha} \\
    & \stackrel{(\ref{cond:alpha3})}{\geq} & \xi_1(\beta^*,l^-_{p-1}) +
    \frac{\alpha^{-p}}{2} \ .
\end{eqnarray*}
Now we compare 
\begin{equation}
  \label{step3:reference2}
  x^1_1 \ld x^1_n \ := \
\xi_{1}(\beta^*,l^-_{p-1}) \ld \xi_{l^+_p}(\beta^*,l^-_{p-1})
\end{equation}
 and 
 \begin{equation}
   \label{step3:new2}
   x^2_1 \ld x^2_n :=
\xi_1(\beta^+,l^-_p) \ld \xi_{l^+_p}(\beta^+,l^-_p)
 \end{equation}
 via Lemma
 \ref{lem:orbitsthrowout}(b)%
 \footnote{\label{footnote2}Again, the assumptions of the lemma with
   $\epsilon = \alpha^{-p}$ are all satisfied: We have
   $\err(\ldots) \leq K\cdot\epsilon$ as before. (\ref{step3:c})
   together with Lemma~\ref{lem:regsets}(a) implies
\[
 \tau(n) \ 
   \leq \ \#([1,l^+_p] \setminus R_{\nu(p)}) \ \leq
   \frac{2p-3}{4} \ .
\]
and similarly $\tau(j) \leq \frac{j}{8})$.  As $l^+_p+1 \in
R_{\nu(p)}$ by (\ref{eq:lpminrn}) we also have $x^1_{n+1} =
\xi_{l^+_p+1}(\beta^*,l^-_{p-1}) \in \Balphcl$.  Further
$\tau(n)-\tau(j) \leq \frac{n-j}{6}$ follows from
(\ref{eq:etajnestimate}), and $n=l^+_p \geq 5p$ by
(\ref{eq:lpmestimates}). 
\label{foot:step3c}}
and obtain that $ \xi_{l^+_p}(\beta^+,l^-_p) = x^2_{n+1} \geq
\frac{2}{\alpha}$.

\qed

\solidqed

\ \\ Step 2 and 3 together prove Part I of the induction statement for
$q=p$, apart from $\beta \in \left[1+\walphtel,1+\frac{3}{\walpha}\right]$ whenever
$\xi_{l^+_p+1}(\beta,l^-_p) \in \Balphcl$. This will be postponed until after
the next step. However, Step~3 implies the
slightly weaker estimate 
\[
    \beta \in
\left[1+\walphtel-\alpha^{-p},1+\frac{3}{\walpha} + \alpha^{-p}\right] \ .
\]
(Note that as $\nu(p)$ is admissible the induction statement can be applied
to $q=p-1$ and $n=\nu(p)$, such that $\beta^+_{p-1,\nu(p)} \in
\left[1+\walphtel,1+\frac{3}{\walpha}\right]$.) This will be sufficient in the meanwhile.

  
The next three steps will prove Part II and III of the induction
statement for $q=p$. In order to do so we will proceed by induction on
$n \in [l^+_p+1,\nu(p+1)]$. The next step starts the induction, by
showing Part II for $n=l^+_p+1$.

\ \\ \underline{\textbf{Step 4:}} \ \ \ 
\textit{$\xi_{l^+_p+1}(\beta,l^-_p) \in \Balphcl$ \ implies \ 
  $\xi_j(\beta,l^-_p) \in \Balphcl \ \forall j \in R_{l^+_p+1} $}

\ \\ Assume that $\xi_{l^+_p+1}(\beta,l^-_p) \in \Balphcl$. As we are
in the case of one-sided forcing, $\xi_j(\beta,l^-_p) \leq -\alphtel$
for any $j \in [1,l^+_p]$ implies $\xi_{l^+_p+1}(\beta,l^-_p) \leq
-\alphtel$ (compare the discussion below (\ref{eq:betaminpn})).
Therefore, it suffices to show that for any $j \in R_{l^+_p+1}
\setminus \{l^+_p+1\}$
\begin{equation} \label{step4:a}
   \xi_j(\beta,l^-_p) \ \geq \ \alphtel \ \ \ \ \textit{ implies } \ \ \ \
   \xi_{l^+_p+1}(\beta,l^-_p) \ \geq \ \alphtel \ .
\end{equation}
Using the two claims below, this can be done as follows: Suppose $j
\in R_{l^+_p+1}$ and $\xi_j(\beta,l^-_p) \geq \alphtel$. Then
$d(\omj,0) \geq \frac{3\gamma}{L_2}$ by
Lemma~\ref{lem:regularsetsbasic}(b), such that
Lemma~\ref{lem:basicestimate} implies $\xi_{j+1}(\beta,l^-_p) \geq
\gamma \geq \frac{2}{\alpha}$. Therefore (\ref{step4:a}) follows
directly from Claim~\ref{claim:step4a} with $k=j+1$, provided $j+1 \in
R_{l^+_p+1}$. On the other hand, if $j+1 \in \Gamma_{l^+_p+1}$ then
Claim~\ref{claim:step4b} (with $k=j$) yields the existence of a
suitable $\tilde{k}$, such that (\ref{step4:a}) follows again from
Claim~\ref{claim:step4a}. As %
$R_{l^+_p+1} \cup \Gamma_{l^+_p+1} = [1,l^+_p+1]$, this covers all
possible cases.

\begin{claim}
   \label{claim:step4a} 
   Suppose $\xi_k(\beta,l^-_p) \geq \frac{2}{\alpha}$ for some $k \in
   R_{l^+_p+1}$. Then $\xi_{l^+_p+1}(\beta,l^-_p) > \alphtel$.
 \end{claim}
 \proof\  Let $\beta^* := \beta^+_{p-1,\nu(p)}$ as in Step~3. Note
 that $\xi_{l^+_p+1}(\beta,l^-_p) \in \Balphcl$ implies $|\beta -
 \beta^*| \leq \alpha^{-p}$ by Step 3. Further, we can again apply
 Part II of the induction statement to $q=p-1$ and $n=\nu(p)$. As
 $R_{l^+_p+1} \ssq R_{\nu(p)}$ (see (\ref{eq:lplusstabilized})) we
 obtain
 \begin{equation}
   \label{step4:b}
   \xi_j(\beta^*,l^-_{p-1}) \ \in \ \Balphcl \ \ \ \ \ \forall j \in
   R_{l^+_p+1} \ .
 \end{equation}
The claim now follows from Lemma~\ref{lem:orbitsthrowout}(a), which we apply
to compare the orbits%
\footnote{We choose $\epsilon = \alpha^{-p}$. $\err(\ldots) \leq K \cdot
  \epsilon$ follows from $|\beta-\beta^*| \leq \alpha^{-p}$. As $k\in
  R_{l^+_p+1}$ by assumption and $l^+_p+1 \in R_{l^+_p+1}$ by
  (\ref{eq:lpminrn}), we have $x^1_1,x^1_{n+1} \in \Balphcl$ by
  (\ref{step4:b}). Finally $\tau(n) \leq \min\{0,\frac{2p-3}{4}\}$ and
  $\tau(n) - \tau(j) \leq \frac{n-j}{6}$ follow from
  Lemma~\ref{lem:regsets}(a) and
  (\ref{eq:etajnestimate}).\label{foot:step4a}}
\begin{equation}
  \label{step4:referenceone}
x^1_1 \ld x^1_n \ := \ \xi_k(\beta^*,l^-_{p-1}) \ld \xi_{l^+_p}(\beta^*,l^-_{p-1})
\end{equation}
and  
\begin{equation}
  \label{step4:newone}
  x^2_1 \ld x^2_n \ := \ \xi_k(\beta,l^-_p) \ld \xi_{l^+_p}(\beta,l^-_p) \ .
\end{equation}
Thus we obtain $\xi_{l^+_p+1}(\beta,l^-_p) = x^2_{n+1} \geq
\frac{2}{\alpha}$.

\qed

\begin{claim}
  \label{claim:step4b} Suppose $k\in R_{l^+_p+1},\ k+1\in \Gamma_{l^+_p+1}$
  and $\xi_k(\beta,l^-_p) \geq \alphtel$. Then there exists some $\tilde{k} \in
  R_{l^+_p+1}$ with $\xi_{\tilde{k}}(\beta,l^-_p) \geq \frac{2}{\alpha}$. 
\end{claim}
\proof\  First of all, note that $\Gamma_{l^+_p+1} = \Gamma_{l^+_p}$
by (\ref{eq:gammalqplus}) and $\Gamma_{l^+_p} =
\bigcup_{J\in\J_{l^+_p}} \Gamma(J)$ by (\ref{eq:gammanequality}).
Therefore, there must be some $J_1 \in \J_{l^+_p}$ such that $k+1 \in
\Gamma(J_1)$. Let $m_1 := m_{J_1}$ and $p_1 := p(m_1)$. As
$\Gamma(J_1)=J^-_1 \cup \Gamma^+(J_1)$, we have two possibilities:

Either $k+1
\in J^-_1$, which means that $j+1=\lambda^-(m_1)$ (as %
$J^-_1 = [\lambda^-(m_1),m_1]$ is an interval and we assumed %
$k \in R_{l^+_p+1}$). In this case define $m=m_1$ and $t=0$.

The other alternative is that $k+1 \in \Gamma^+(J_1)$, and in this
case we have to `\textit{go backwards through the recursive structure
  of the set $R_{l^+_p+1}$'}, until we arrive at the first
alternative:

As $\Gamma^+(J_1) = \Gamma_{l^+_{p_1}} + m_1$ by
(\ref{eq:recursivegammas}), $k+1 \in \Gamma^+(J_1)$ means that
$k-m_1+1 \in \Gamma_{l^+_{p_1}}$. Hence, similar to before there
exists some $J_2 \in \J_{l^+_{p_1}}$ such that either $k-m_1+1 =
\lambda^-(m_{J_2})$ or $k-m_1+1 \in \Gamma^+(J_2)$. Let $m_2 :=
m_{J_2}$ and $p_2 := p(m_2)$.  If we are in the second case where
$j-m_1+1 \in \Gamma^+(J_2)$ we continue like this, but after finitely
many steps the procedure will stop and we arrive at the first
alternative. This is true because in each step the $p_i$ become
smaller, more precisely $p_{i+1} \leq
p_i-3$,%
\footnote{Note that there is no $J \in \J_{l^+_{p_i}}$ with $p_J \geq p_i-2$
by (\ref{eq:lpmestimates}). \label{foot:step4b} }
and finally $\Gamma_{l^+_1}$ is empty. Thus we obtain two sequences
$p_1 > \ldots > p_r \geq 0$ and $m_1 > \ldots > m_r$ with $p_i =
p(m_i)$, such that $k-\sum_{i=1}^{r-1} m_i + 1 = \lambda^-(m_r)$ for
some $r \in \N$. Let $m := m_r$ and $t := \sum_{i=1}^{r-1} m_i$, such
that $p_r=p(m)$ and $k+1-t = \lambda^-(m)$. Note that for $r=1$ this
coincides with the above definitions of $m$ and $t$ in the first case.
We have
\begin{eqnarray} 
    d(\omega_t,0) & \leq & \sum_{i=1}^{r-1} d(\omega_{m_i},0) \nonumber \ \leq \
    \sum_{i=1}^{r-1} \frac{\alpha^{-(p_i-1)}}{L_2} 
     \\ & \leq &
    \frac{\alpha^{-(p(m)+2)} \cdot S_\infty(\alpha)}{L_2} \
    \stackrel{(\ref{cond:alpha4})}{\leq} \ \viertel \cdot
    \frac{\alpha^{-(p(m)+1)}}{L_2} \ . \label{eq:domegat}
\end{eqnarray}
Now choose some $q' \geq p(m) \geq 1$ such that %
$l^+_{q'}+1 \leq m \leq \nu(q'+1)$. This is possible as %
$m \geq \nu(p(m)) \geq l^+_{p(m)}+1$, and because the intervals
$[l^+_q+1,\nu(q)]$ overlap by (\ref{eq:lpmestimates}). In addition, we
can choose $q' < p-1$ as $m \leq l^+_p+1 < \nu(p-1)$.

We now want to apply Lemma~\ref{lem:approaching} with $\beta^* :=
\beta^+_{q',m}$ $q=p(m)$, $l^* = l^-_{q'}$ $l=l^-_p$ and $k=t$. Note
that we can apply Part II of the induction statement with $q=q'$ and
$n=m$ to obtain that %
$\beta^* \in \left[1+\walphtel,1+\frac{3}{\walpha}\right]$,
\begin{equation}
  \label{step4:c} \xi_j(\beta^*,l^-_{q'}) \ \geq \ \gamma \ \ \ \ \ \forall j
  \in [l^-_{q'},0] \setminus \Omega_\infty 
\end{equation}
and
\begin{equation}
  \label{step4:d} \xi_j(\beta^*,l^-_{q'}) \ \in \Balphcl \ \ \ \ \ \forall j
  \in R_m \ .
\end{equation}
Further, Part III of the induction statement%
\footnote{With $q_1 = q'$, $q=p-1$, $n_1=m$ and
  $n_2=\nu(p-1)$. \label{foot:step4d} }
together with Step 3 imply that
\begin{equation} \label{step4:e}
|\beta-\beta^*| \ \leq \ |\beta - \beta_{p-1,\nu(p)}| +
  |\beta_{p-1,\nu(p)} - \beta^*| \ \leq \ \alpha^{-p} +
\alpha^{-q'} \ \leq \ 2\alpha^{-q'} \  
\end{equation}
and finally $\xi_{k+1}(\beta,l^-_p) \geq \gamma$ if
$\xi_k(\beta,l^-_p) \geq \alphtel$ by Lemma~\ref{lem:basicestimate}~.
Thus Lemma \ref{lem:approaching} yields
\begin{equation} \label{step4:f}
    \{ j \in [\lambda^-(m),m] \mid \xi_{j+t}(\beta,l^-_p) < \gamma \} \ \ssq \
    \Omtil_{p(m)-2} \ .
\end{equation}
Consequently (Lemma~\ref{lem:omegatransition}(b))
\begin{equation} \label{step4:g}
    \{ j \in [-l^-_{p(m)},0] \mid \xi_{j+m+t}(\beta,l^-_p) < \gamma \} \ \ssq \
    \Omega_\infty \ .
\end{equation}
This means that we can compare the two sequences 
\begin{equation}
  \label{step4:referencetwo}
x^1_1 \ld x^1_n \ := \ \xi_{-l^-_{p(m)}}(\beta^*,l^-_{q'}) \ld \xi_{-1}(\beta^*,l^-_{q'})
\end{equation}
and 
\begin{equation}
  \label{step4:newtwo}
x^2_1 \ld x^2_n \ := \ \xi_{m+t-l^-_{p(m)}}(\beta,l^-_p) \ld
\xi_{m+t-1}(\beta,l^-_p)
\end{equation}
 via Lemma~\ref{lem:orbitscontraction} with %
 $\epsilon \ := \ L_2 \cdot d(\omm,0) \ \in $ 
$(\alpha^{-p(m)},\alpha^{-(p(m)-1)}]$ to obtain that%
\footnote{We have $q=p(m)-1$. Note that $d(\omega_{m+t},0) \leq
  \frac{2\epsilon}{L_2}$ (see (\ref{eq:domegat})) and $|\beta-\beta^*|
  \leq 2\alpha^{-q'} \leq 2\epsilon$ by (\ref{step4:e}), such that
  $\err(\ldots) \leq K \cdot \epsilon$ by Remark~\ref{bem:errorterm}~.
  Further, it follows from (\ref{step4:c}) and (\ref{step4:g}) that
  $\eta(j,n) \leq \#([-(n-j),-1] \cap \Omega_\infty) \leq
  \frac{n-j}{10}$ (see (\ref{cond:hfunctions2})). Finally
  $n=l^-_{p(m)} \geq 4p(m)$ by (\ref{eq:lpmestimates}), such that
  $\alpha^{-\viertel n} \leq \epsilon$.\label{foot:step4f}}
\begin{equation} \label{step4:h}
       |\xi_{m+t}(\beta,l^-_p) - \xi_0(\beta^*,l^-_{q'})| \ \leq \
        \epsilon \cdot (6+K\cdot S_\infty(\alpha^{\viertel})) \ .
\end{equation}
As $d(\omega_{m+t},0) \geq \frac{3}{4} \cdot \frac{\epsilon}{L_2}$ (see
(\ref{eq:domegat})), it follows from (\ref{cond:Fcontraction}) and
(\ref{cond:sharppeak}) that
\begin{eqnarray}
      \lefteqn{\xi_{m+t+1}(\beta,l^-_p) \ \geq } \nonumber \\ & \geq &
      \xi_1(\beta^*,l^-_{q'}) +
      \frac{3\epsilon}{4}  - \epsilon \cdot \frac{6 +
      K \cdot S_\infty(\alpha^{\viertel})}{2\walpha} 
       \ \stackrel{(\ref{cond:alpha3})}{\geq} \
      \xi_1(\beta^*,l^-_{q'}) + \frac{\epsilon}{2} \ . \label{step4:crucial}
\end{eqnarray}
Now first assume $p(m) \geq 2$, such that $\epsilon \leq \alphtel$.
(The case $p(m)=1$ has to be treated separately, see below.) Then we can
apply Lemma \ref{lem:orbitsthrowout}(b) to compare the orbits
\begin{equation}
  \label{step4:i}
x^1_1 \ld x^1_n \ := \ \xi_1(\beta^*,l^-_{q'}) \ld
\xi_{l^+_{p(m)}}(\beta^*,l^-_{q'})
\end{equation} 
and
\begin{equation}
  \label{step4:j}
x^2_1 \ld x^2_n \ := \ \xi_{m+t+1}(\beta,l^-_p)
\ld \xi_{m+t+l^+_{p(m)}}(\beta,l^-_p)
\end{equation}
to conclude that%
\footnote{We have $q=p(m)-1$ and $\err(\ldots) \leq K\cdot\epsilon$ as
  before (see Footnote~\ref{foot:step4f}). (\ref{step4:d}) yields that
  $x^1_{n+1} \in \Balphcl$ (note that $l^+_q+1 \in R_m$ by
  (\ref{eq:lpminrn})). Further, we have
\[
\tau(n) \ \leq \ \#\left([1,l^+_q] \setminus R_m\right) \ \leq \ 
\max\left\{0,\frac{2q-3}{4}\right\}
\]
and $\tau(j) \leq \#([1,j] \setminus R_m) \leq \frac{j}{8}$ by
Lemma~\ref{lem:regsets}(a). $\tau(n)-\tau(j) \leq \frac{n-j}{6}$
follows again from (\ref{eq:etajnestimate}), and finally $n =
l^+_{p(m)} \geq 5(p(m)-1)$ .\label{foot:step4g}}
\begin{equation} \label{step4:k}
    \xi_{m+t+l^+_q+1}(\beta,l^-_p) \ \geq \ \frac{2}{\alpha} \ . 
\end{equation}
As $J_r = J(m)$ is a maximal interval in $\Gamma_{l^+_{p_{r-1}}}$ we 
have $\lambda^+(m)+1 \in R_{l^+_{p_{r-1}}}$.  Therefore 
$\lambda^+(m)+1+t \in R_{l^+_p+1}$ follows from the recursive
structure of this set. Consequently, we can choose $\tilde{k} =
m+t+l^+_{p(m)}+1$.

\ \\ Finally, suppose $p(m)=1$. In this case we still have
$\xi_{m+t+1}(\beta,l^-_p) \geq \xi_1(\beta^*,l^-_{q'}) +
\frac{\epsilon}{2}$ by (\ref{step4:crucial}). There are two
possibilities: Either $\xi_{m+t+1}(\beta,l^-_p) \geq
\frac{2}{\alpha}$. As $m+1 \in R_{l_{p_{r-1}}}$ (see
Remark~\ref{bem:regularsets}(c)), we have $m+t+1 \in R_{l^+_p+1}$ due
to the recursive structure of this set. Thus, we can choose $\tilde{k}
= m+t+1$. On the other hand, if $\xi_{m+t+1}(\beta,l^-_p) \in
B_{\frac{2}{\alpha}}(0)$ then we can apply (\ref{cond:Fcontraction})
again and obtain
\[
\xi_{m+t+2}(\beta,l^-_p) \ \geq \
\xi_2(\beta^*,l^-_{q'}) + 2\walpha\cdot\epsilon - K \epsilon \
\stackrel{(\ref{cond:alpha1})}{\geq} \ 
\xi_2(\beta^*,l^-_{q'}) + \walpha\cdot\epsilon \geq \frac{1}{\alpha}
\]
by (\ref{cond:alphagamma0}), as $\xi_2(\beta^*,l^-_q) \in \Balphcl$ by
(\ref{step4:d}) and $\epsilon \geq \alphtel$.  Thus, we can choose
$\tilde{k} = m+t+2$ in this case. Note that $m+t+2$ is contained in $
R_{l^+_p+1}$ for the same reasons as $m+t+1$.

\qed

\solidqed

\ \\ Now we can show that $\xi_{l^+_p+1}(\beta,l^-_p) \in \Balphcl$
implies $\beta \in \left[1+\walphtel,1+\frac{3}{\walpha}\right]$ and
thus complete the proof of Part I of the induction statement: Suppose
$\xi_{l^+_p+1}(\beta,l^-_p) \in \Balphcl$. By Steps~2 and 3 we know
that $\xi_0(\beta,l^-_p) \geq \gamma$. This implies that
\[
    \xi_1(\beta,l^-_p) \ \in \
    \left[1+\frac{3}{2\walpha}-\beta,1+\frac{3}{\walpha}-\alphtel-\beta\right]
\]
(see assumptions (\ref{cond:Ffixedpoints}) and
(\ref{cond:Fcontraction})). As Step 4 implies that $\xi_1(\beta,l^-_p)
\in \Balphcl$, and $\frac{1}{2\walpha} \geq \alphtel$ by
(\ref{cond:alpha1}), this gives the required estimate.

\ \\


\ \\
\underline{\textbf{Step 5:}} \ \ \ \textit{Part II of the induction
  statement implies Part III}

\ \\Actually, the situation is a little bit
more complicated than the headline above may suggest. In fact,
the both remaining parts of the induction statement have to be proved
simultaneously by induction on $n$. However, in each step of the
induction Part II will imply Part III.

In order to make this more precise, assume that Part II with $q=p$ holds for
all $n\leq N$, with $N \in [l^+_p+1,\nu(p+1)]$. What we will now show is that
in this case Part III(a) holds as well whenever $n_1,n_2 \leq N$, and
similarly Part III(b) holds whenever $n_2 \leq N$.

\ \\
Suppose that $N \in [l^+_p+1,\nu(p+1)]$ and Part II with $q=p$ holds
for all $n \leq N$. Further, let $n_2 \leq N$ and $n_1 \in R_{n_2}$.
Then the Part II of the induction statement applied to $q=p$ and
$n=n_2$ yields that $\xi_{n_1}(\beta^+_{q,n_2},l^-_p) \in \Balphcl$,
and for $n=n_1$ we obtain that $\xi_{n_1}(\beta,l^-_p)\in\Balphcl$
implies $\xi_j(\beta,l^-_p) \in \Balphcl \ \forall j \in R_{n_1}$.
Thus all the assumption of Lemma~\ref{lem:parameterestimates} (with
$n=n_1$, $\beta_1=\beta^+_{p,n_1}$ and $\beta_2=\beta^+_{p,n_2}$) are
satisfied, such that
\[
   |\beta^+_{p,n_1} - \beta^+_{p,n_2}| \ \leq \ 2\alpha^{-\frac{n_1}{4}}
\] 
as required.

\ \\
For Part III(b) let $q_1 < p$, $n_1 \in [l^+_{q_1}+1,\nu(q_1+1)]$ and 
$n_2 \in [l^+_p+1,N]$. First suppose $q_1 < p-1$. Then Part III(b) of 
the induction statement (with $q=p-1$ and $n_2 = \nu(p)$) yields
\[ 
    |\beta^+_{q_1,n_1} - \beta^+_{p-1,\nu(p)}| \ \leq \ 3 \cdot
    \sum_{i=q_1+1}^{p-1}  \alpha^{-i} \ . 
\]
Further, Part II of the induction statement (with $q=p$ and $n=n_2$) 
yields that $\xi_{l^+_p+1}(\beta^+_{p,n_2},l^-_p) \in \Balphcl$ (note 
that $l^+_p+1 \in R_{n_2}$ by (\ref{eq:lpminrn})), and consequently
\[
   | \beta^+_{p-1,\nu(p)} - \beta^+_{p,n_2}| \leq \  \alpha^{-p}
\]
by Step III. Altogether, we obtain
\[
   |\beta^+_{q_1,n_1} - \beta^+_{p,n_2}| \ \leq \ |\beta^+_{q_1,n_1} -
    \beta^+_{p-1,\nu(p)}| + | \beta^+_{p-1,\nu(p)} - \beta^+_{p,n_2}|
    \ \leq \ 3\cdot \sum_{i=q_1+1}^p \alpha^{-i} \ .
\]
On the other hand, if $q_1=p-1$ then Part III(a) (with $q=q_1=p-1$ and
$n_2=\nu(p)$) in 
combination with $n_1 \geq l^+_{p-1} \geq 4p$ (see (\ref{eq:lpmestimates}))
yields 
\[
|\beta^+_{q_1,n_1} - \beta^+_{p-1,\nu(p)}| \ \leq \ 
2\alpha^{-frac{n_1}{4}} \ \leq \ 2\alpha^{-p} \ ,
\]
such that 
\[
|\beta^+_{q_1,n_1} - \beta^+_{p,n_2}| \ \leq \ |\beta^+_{q_1,n_1} -
    \beta^+_{p-1,\nu(p)}| + | \beta^+_{p-1,\nu(p)} - \beta^+_{p,n_2}|
    \ \leq \ 
    3\alpha^{-p} 
\]
as required. Finally, note that
\[
         3\cdot \sum_{i=q_1+1}^p \alpha^{-i} \ \leq \
         \frac{3S_\infty(\alpha)}{\alpha} \cdot \alpha^{-q_1} \ \leq \
         \alpha^{-q_1} 
\]
by (\ref{cond:alpha4}).

\solidqed

\ \\
Now we can already use the parameter estimates up to $N$ (in the way
mentioned above) during the induction step $N \ra N+1$ in the proof of
Part II.


\ \\
\underline{\textbf{Step 6:}} \ \ \ \textit{Proof of Part II for $q=p$.} 

\ \\ In order to prove Part II of the induction statement for $q=p$,
we will proceed by induction on $n$. Steps 2--4 show that the
statement holds for $n=l^+_p+1$. Suppose now that it holds for all
$n\leq N$, where $N\in[l^+_p+1,\nu(p+1)-1]$. We have to show that it
then holds for $N+1$ as well. In order to do so, we distinguish three
different cases: First, if $N+1$ is not admissible there is nothing to
prove. Secondly, if both $N$ and $N+1$ are admissible then necessarily
$p(N)=0$, otherwise $N+1$ would be contained in $J(N)$. Thus
$d(\omega_N,0) \geq \frac{4\gamma}{L_2}$, and in addition Part II of
the induction statement with $q=p$ and $n=N$ implies that
$\beta^\pm_{p,N} \in \left[1+\walphtel,1+\frac{3}{\walpha}\right]$.
Therefore Lemma~\ref{lem:basicestimate} yields that
\[
    \xi_{N+1}(\beta^+_{p,N},l^-_p) \ > \ \alphtel
\]
and
\[
    \xi_{N+1}(\beta^-_{p,N},l^-_p) \ < \ -\alphtel \ .
\]
Consequently $\xi_{N+1}(\beta,l^-_p) \in \Balphcl$ implies that
$\xi_N(\beta,l^-_p) \in \Balphcl$, and everything else follows from Part II of
the induction statement for $n=N$.

Thus, it remains to treat the case where $N+1$ is admissible but $N \notin
A_N$. By (\ref{eq:aninequality}) this also means that $N \notin A_{N+1}$.
Consequently there exists an interval $J \in \J_{N+1}$ which contains $N$,
such that $J = [t,N]$ where $t:=\lambda^-(m_J)$. Note that $t-1,t,m_J \in
A_{m_J}$ by Lemma \ref{lem:endpoints}(a). In particular $m_J$ and $t-1$ are
admissible. First of all, we will prove the following claim.
\begin{claim}
  \label{claim:step6a}
$\xi_{N+1}(\beta,l^-_p) \in \Balphcl$ implies $\xi_{t-1}(\beta,l^-_p) \in
\Balphcl$.
\end{claim}
\proof\ 
It suffices to show that 
\begin{equation}
  \label{step6:a}
\xi_{N+1}(\beta^+_{p,t-1},l^-_p) \ > \ \alphtel 
\end{equation}
(see (\ref{eq:betapluspndef})--(\ref{eq:betaminpn})). Let $m:=m_J$, $\beta^+:=
\beta^+_{p,t-1}$ and $\beta^*:=\beta^+_{p,m}$. Using Part II of the induction
statement (with $q=p$ and $n=m$) we obtain
\begin{equation}
  \label{step6:b}
    \xi_j(\beta^*,l^-_p) \ \in \Balphcl \ \ \ \forall j \in R_m 
\end{equation}
and $\beta^* \in \left[1+\walphtel,1+\frac{3}{\walpha}\right]$. Further, the same
statement with $n=t-1$ implies
\begin{equation}
  \label{step6:c}
    \xi_j(\beta^+,l^-_p) \ \geq \ \gamma \ \ \ \forall j \in [-l^-_p,0]
    \setminus \Omega_\infty 
\end{equation}
and
\begin{equation}
  \label{step6:d}
    \xi_j(\beta^+,l^-_p) \ \in \Balphcl \ \ \ \forall j \in R_{t-1} \ . 
\end{equation}
Finally Part III(a) of
the induction statement (with $q=p$, $n_1=t-1$ and $n_2=m$) yields that
\begin{equation}
  \label{step6:e}
  |\beta^+ - \beta^*| \ \leq \ \alpha^{-\frac{(t-1)}{4}} \ \leq \
   \alpha^{-\frac{l^+_p}{4}} \ \stackrel{(\ref{eq:lpmestimates})}{\leq} \
   \alpha^{-(p+1)} \ .
\end{equation}
Note that $t-1$ is contained in $\Omega_\infty(m)$ by (\ref{eq:jminominfty}),
and as $l^+_p+1 \notin \Omega_\infty$ this interval must be to the right of
$l^+_p+1$. Therefore $t-1 > l^+_p+1$.  Now all the assumptions for the
application of Lemma~\ref{lem:approaching} are satisfied%
\footnote{With $\beta^*$ and $m$ as above, $q= p(m) \ (\leq p)$,
  $l=l^*=l^-_p$, $k=0$ and $\beta=\beta^+$. Note that $p(t-1) = 0$ as
  $t-1$ is admissible, and $\xi_{t-1}(\beta^+,l^-_p) = \alphtel$ by
  definition of $\beta^+ = \beta^+_{p,t-1}$. Therefore
  Lemma~\ref{lem:basicestimate} implies
  $\xi_{m-l^-_{p(m)}}(\beta^+,l^-_p) = \xi_{t}(\beta^+,l^-_p) \geq
  \gamma$. \label{foot:step6a} }
and we obtain
\begin{equation}
  \label{step6:f}
 \{ j \in [t,m] \mid \xi_j(\beta^+,l^-_p) < \gamma \} \ \ssq \ \Omtil_{p(m)-2}
 \ .
\end{equation}
Using Lemma~\ref{lem:omegatransition}(b) this further means that 
\begin{equation}
  \label{step6:g}
 \{ j \in [-l^-_{p(m)},0] \mid \xi_{j+m}(\beta^+,l^-_p) < \gamma \} \ \ssq \
 \Omega_\infty 
 \ .
\end{equation}
Now we compare the orbits 
\begin{equation}
  \label{step6:referenceone}
x_1^1 \ld x^1_n \ := \ \xi_{-l^-_{p(m)}}(\beta^+,l^-_p) \ld
\xi_{-1}(\beta^+,l^-_p)
\end{equation}
and 
\begin{equation}
  \label{step6:newone}
x^2_1 \ld x^2_n \ := \ \xi_t(\beta^+,l^-_p) \ld
\xi_{m-1}(\beta^+,l^-_p) \ ,
\end{equation}
using Lemma~\ref{lem:orbitscontraction} with $\epsilon := L_2 \cdot
d(\omega_m,0) \in [\alpha^{-p(m)},\alpha^{-(p(m)-1)})$, to conclude that%
\footnote{As $\beta_1=\beta_2=\beta^+$ we
  have $\err(\ldots) \leq K\epsilon$, see
  Remark~\ref{bem:errorterm}. (\ref{step6:c}) and (\ref{step6:g}) imply that 
 \[
     \eta(j,n) \ \leq \ \#([-(n-j),-1]\cap \Omega_\infty) \ \leq \
  \frac{n-j}{10} 
 \]
 by (\ref{cond:hfunctions2}). Finally $n=l^-_{p(m)} \geq 4(p(m)+1)$ by
  (\ref{eq:lpmestimates}), such that $\alpha^{-\frac{n}{4}} \leq
  \alpha^{-(p(m)+1)} \leq \epsilon$.  \label{foot:step6b} }
\[
    |\xi_m(\beta^+,l^-_p) - \xi_0(\beta^+,l^-_p)| \ \leq \ \epsilon
     \cdot (6+K\cdot S_\infty(\alpha^{\viertel})) \ .
\]
As (\ref{step6:c}) and (\ref{step6:g}) in particular imply that
$\xi_m(\beta^+,l^-_p),\xi_0(\beta^+,l^-_p) \geq \gamma$, it follows from
(\ref{cond:Fcontraction}) and (\ref{cond:sharppeak}) that
\begin{eqnarray} \nonumber
\lefteqn{\xi_{m + 1}(\beta^+,l^-_p) \ \geq} \\ & \geq & \xi_1(\beta^+,l^-_p)
+ \epsilon - \frac{\epsilon \cdot (6+K\cdot
S_\infty(\alpha^{\viertel}))}{2\walpha} \ \stackrel{(\ref{cond:alpha3})}{\geq}
\xi_1(\beta^+,l^-_p) + \frac{\epsilon}{2} \ . 
\end{eqnarray}
Now first assume $p(m) \geq 2$, such that $d(\omm,0) \leq
\frac{\alpha^{-1}}{L_2}$. Then we can apply Lemma~\ref{lem:orbitsthrowout}(b)%
\footnote{\label{foot:step6c}With $\epsilon = L_2\cdot d(\omm,0)$ as
  above, such that $q=p(m)-1$. $\err(\ldots) \leq K\cdot \epsilon$
  follows again from Remark~\ref{bem:errorterm}.
  $x^1_{n+1}\in\Balphcl$ follows from (\ref{step6:d}) as $l^+_{p(m)}+1
  \in R_{t-1}$ by (\ref{eq:lpminrn}). (\ref{step6:d}) also implies
  $\tau(n) \leq \frac{2p(m)-5}{4}$ and $\tau(j) \leq \frac{j}{8} \ 
  \forall j \in [1,n]$ by Lemma~\ref{lem:regsets}(a) and
  $\tau(n)-\tau(j) \leq \frac{n-j}{6}$ by (\ref{eq:etajnestimate}).
  Finally $n=l^+_{p(m)} \geq 5(p(m)-1)$ by (\ref{eq:lpmestimates}). }
to the sequences
\begin{equation}
  \label{step6:referencetwo}
x_1^1 \ld x^1_n \ := \ \xi_1(\beta^+,l^-_p) \ld \xi_{l^+_{p(m)}}(\beta^+,l^-_p)
\end{equation}
and 
\begin{equation}
  \label{step6:newtwo}
x^2_1 \ld x^2_n \ := \ \xi_{m+1}(\beta^+,l^-_p) \ld \xi_{N}(\beta^+,l^-_p) \ ,
\end{equation}
which yields that $\xi_{N+1}(\beta^+,l^-_p) = x^2_{n+1}\geq \frac{2}{\alpha}$ as required
for (\ref{step6:a}). 

\ \\
It remains to address the case $p(m) = 1$. Note that in this case $p(j) = 0 \ \forall j
\in [m+1,N]$ (see Lemma~\ref{lem:endpoints}(a)). There are two
possibilities: Either $\xi_{m+1}(\beta^+,l^-_p) \geq \alphtel$, in which case
$\xi_{N+1}(\beta^+,l^-_p) \geq \gamma > \alphtel$ follows from the repeated
application of Lemma~\ref{lem:basicestimate}~. Otherwise,
$\xi_{m+1}(\beta^+,l^-_p) \in \Balphcl$. As $1,2\in R_{t-1}$ (see
(\ref{eq:lpminrn})) it follows from (\ref{step6:d}) that
$\xi_1(\beta^+,l^-_p), \xi_2(\beta^+,l^-_p) \in \Balphcl$ as well. Therefore
(\ref{cond:Fexpansion}) implies that 
\begin{eqnarray} \nonumber
  \xi_{m+2}(\beta^+,l^-_p) & \geq & \xi_2(\beta^+,l^-_p) +
  2\walpha\cdot \epsilon - K \cdot \epsilon \\
  &\stackrel{(\ref{cond:alpha1})}{\geq}& \xi_2(\beta^+,l^-_p) + \walpha\cdot
  \epsilon \ \geq \ \frac{2}{\alpha}  \label{step6:h}
\end{eqnarray}
as $\epsilon \geq \alphtel$ in this case.  Again, we obtain
$\xi_{N+1}(\beta^+,l^-_p) \geq \gamma > \alphtel$ by repeated
application of Lemma~\ref{lem:basicestimate}~.

\qed

\ \\
Now suppose $\xi_{N+1}(\beta,l^-_p) \in \Balphcl$. Then by
Claim~\ref{claim:step6a} there holds $\xi_{t-1}(\beta,l^-_p) \in
\Balphcl$. As we can already apply Part II of the induction statement with
$q=p$ and $n=t-1$, this further implies (\ref{eq:indstatementI}),
$\beta \in \left[1+\walphtel,1+\frac{3}{\walpha}\right]$  and 
\[
  \xi_j(\beta,l^-_p) \ \in \ \Balphcl \ \ \ \ \  \forall j \in R_{t-1} \ . 
\]
Note that $R_{N+1} \cap [1,t-1] = R_{t-1}$ (see
Remark~\ref{bem:regularsets}(b)), such that $R_{N+1} = R_{t-1} \cup
R(J) \cup\{N+1\}$. Therefore, in order to complete this step and
thereby the proof of Induction scheme~\ref{thm:indscheme}, it only
remains to show that
\begin{claim}
 \label{claim:step6b}
$\xi_{N+1}(\beta,l^-_p) \in \Balphcl$ implies $\xi_j(\beta,l^-_p) \in \Balphcl
  \ \forall j \in R(J)$.
\end{claim}
\proof\  The proof of this statement is very similar to the proof of Step 4,
and likewise we will use two further claims, namely Claim~\ref{claim:step6c}
and Claim~\ref{claim:step6d} below, which are the analogues of
Claim~\ref{claim:step4a} and Claim~\ref{claim:step4b}~. Suppose
$\xi_j(\beta,l^-_p) > \alphtel$ for some $j \in R(J)$. We have to distinguish
two cases (note that $R(J) \cup \Gamma^+(J) = J^+$): Either $j+1 \in R(J)$. As
$d(\omj,0) \geq \frac{3\gamma}{L_2}$ by Lemma~\ref{lem:regularsetsbasic}(b),
Lemma~\ref{lem:basicestimate} implies $\xi_{j+1}(\beta,l^-_p) \geq \gamma \geq
\frac{2}{\alpha}$. Therefore we can apply Claim~\ref{claim:step6c} with
$k=j+1$. On the other hand, if $j+1\notin R(J)$, then Claim \ref{claim:step6d}
(with $k=j$) yields the existence of a suitable $\tilde{k}$ and we can again
apply Claim~\ref{claim:step6c}, this time with $k=\tilde{k}$. In both cases we
obtain that $\xi_j(\beta,l^-_p) \geq \alphtel$ implies $\xi_{N+1}(\beta,l^-_p) >
\alphtel$. As we are in the case of one-sided forcing, the fact that
$\xi_j(\beta,l^-_p) \leq -\alphtel$ implies $\xi_{N+1}(\beta,l^-_p) < -\alphtel$
is obvious. This proves the claim.

\qed

\begin{claim}
 \label{claim:step6c}
Suppose $\xi_k(\beta,l^-_p) \geq \frac{2}{\alpha}$ for some $k \in
R(J)$. Then $\xi_{N+1}(\beta,l^-_p) \ > \ \alphtel$.
\end{claim} 
\proof\  First of all, if $p_J=1$ then $p(j)=0 \ \forall j \in J^+$ by
Lemma~\ref{lem:endpoints}(a), and the claim follows from the repeated
application of Lemma~\ref{lem:basicestimate}~. Thus we can assume $p_J \geq
2$. Claim~\ref{claim:step6a} together with Part II of the induction statement
with $q=p$ and $n=t-1$ imply that
\begin{equation}
  \label{step6:i}
  \xi_j(\beta,l^-_p) \ \in \ \Balphcl \ \ \ \ \ \forall j \in R_{t-1} \supseteq
  R_{l^+_{p_J}+1} \supseteq R_{l^+_{p_J}}
\end{equation}
(see (\ref{eq:lplusstabilized}) for the inclusions). Consequently, we can
apply Lemma~\ref{lem:orbitsthrowout}(a)%
\footnote{We choose $\epsilon = L_2\cdot d(\omega_{m_{p_J}},0) \in
  [\alpha^{-p_J},\alpha^{-(p_J-1)}]$, such that $q=p_J-1$ and
  $\err(\ldots) \leq K \cdot \epsilon$. Note that $k\in
  R(J)$ implies $k-m_J \in R_{l^+_{p_J}}$ by (\ref{eq:recursivedef1}), and
  further $l^+_{p_J}+1 \in R_{l^+_{p_J}+1}$ by (\ref{eq:lpminrn}). Therefore
  $x^1_1,x^1_{n+1} \in \Balphcl$ by (\ref{step6:i}). Finally $\tau(n) \leq
  \min\{0,\frac{2p(m)-5}{4}\}$ by Lemma~\ref{lem:regsets}(a) and
  $\tau(n)-\tau(j) \leq \frac{n-j}{6}$ by (\ref{eq:etajnestimate}).\label{foot:step6d}}
to the sequences
\begin{equation}
  \label{step6:referencethree}
x^1_1 \ld x^1_n \ := \ \xi_{k-m_J}(\beta,l^-_p) \ld \xi_{l^+_{p_J}}(\beta,l^-_p)
\end{equation}
and  
\begin{equation}
  \label{step4:newthree}
  x^2_1 \ld x^2_n \ := \ \xi_k(\beta,l^-_p) \ld \xi_N(\beta,l^-_p) \ .
\end{equation}
to obtain that $\xi_{N+1}(\beta,l^-_p) = x^2_{n+1} \geq
\frac{2}{\alpha}$.

\qed
 
\begin{claim}
 \label{claim:step6d}
Suppose $k\in R(J),\ k+1\in \Gamma^+(J)$ and $\xi_k(\beta,l^-_p)
  \geq \alphtel$. Then there exists some $\tilde{k} \in R(J)$ with
  $\xi_{\tilde{k}}(\beta,l^-_p) \geq \frac{2}{\alpha}$.
\end{claim}
\proof\  Let $J_1:=J$, $m_1:=m_J$ and $p_1:=p_J$. As in the proof of
Claim~\ref{claim:step4b} we can find sequences $p_1 > \ldots > p_r
\geq 0$ and $m_1 > \ldots > m_r \in [1,l^+_{p_{r-1}}]$ with $p_i =
p(m_i) \leq p_{i-1}-3$, such that $k-\sum_{i=1}^{r-1} m_i + 1 =
\lambda^-(m_r)$ for some $r \in \N$. Let $m := m_r$ and $t :=
\sum_{i=1}^{r-1} m_i$. (The only difference to Claim~\ref{claim:step4b}
is that $r=1$ is not possible.) Likewise, we have
\begin{equation} \label{step6:j} 
    d(\omega_t,0) \ \leq \  \viertel \cdot
    \frac{\alpha^{-(p(m)+1)}}{L_2} \ \leq \ \frac{d(\omm,0)}{4} \ .
\end{equation}
Again, we choose some $q' \geq p(m)$ such that $l^+_{q'}+1 \leq m \leq
\nu(q'+1)$. As $m \leq l^+_{p_{r-1}} \leq l^+_p < \nu(p-2)$ (see
(\ref{eq:lpmestimates})) we can assume that $q' \leq p-2$. 

We now want to apply Lemma~\ref{lem:approaching} with $\beta^* :=
\beta^+_{q',m}$ $q=p(m)$, $l^* = l^-_{q'}$ $l=l^-_p$ and $k=t$. In
order to check the assumptions, note that we can apply Part II of the 
induction statement (with $q=q'$ and $n=m$) to %
$\beta^* := \beta_{q',m}$ and obtain that %
$\beta^* \in \left[1+\walphtel,1+\frac{3}{\walpha}\right]$,
\begin{equation}
  \label{step6:k} \xi_j(\beta^*,l^-_{q'}) \ \geq \ \gamma \ \ \ \ \ \forall j
  \in [l^-_{q'},0] \setminus \Omega_\infty 
\end{equation}
and
\begin{equation}
  \label{step6:l} \xi_j(\beta^*,l^-_{q'}) \ \in \Balphcl \ \ \ \ \ \forall j
  \in R_m \setminus \Omega_\infty \ .
\end{equation}
In addition, Step 3%
\footnote{Note that $\xi_{N+1}(\beta,l^-_p) \in \Balphcl$ implies
  $\xi_{t-1}(\beta,l^-_p) \in \Balphcl$ by Claim~\ref{claim:step6a},
  which in turn implies $\xi_{l^+_p+1}(\beta,l^-_p) \in \Balphcl$ as
  $l^+_p+1 \in R_{t-1}$, see (\ref{eq:lpminrn}).}
together with Part III of the induction statement%
\footnote{With $q_1 = q'$, $q=p-1$, $n_1=m$ and $n_2=\nu(p-1)$. }
imply that
\begin{equation} \label{step6:m}
|\beta-\beta^*| \ \leq \ |\beta - \beta_{p-1,\nu(p)}| +
  |\beta_{p-1,\nu(p)} - \beta^*| \ \leq \ \alpha^{-p} +  
\alpha^{-q'} \ \leq \ 2\alpha^{-p(m)} \  
\end{equation}
Finally,  $\xi_{k+1}(\beta,l^-_p) \geq
\gamma$ if $\xi_k(\beta,l^-_p) \geq \alphtel$ by Lemma
\ref{lem:basicestimate}%
\footnote{Note that $k+1=\lambda^-(m)$ in the claim above
  corresponds to $m+k-l^-_{p(m)}$ in Lemma~\ref{lem:approaching}~.}
. Thus Lemma \ref{lem:approaching} yields 
\begin{equation} \label{step6:n}
    \{ j \in [\lambda^-(m),m] \mid \xi_{j+t}(\beta,l^-_p) < \gamma \} \ \ssq \
    \Omtil_{p(m)-2} \ .
\end{equation}
Consequently (Lemma~\ref{lem:omegatransition}(b))
\begin{equation} \label{step6:o}
    \{ j \in [-l^-_{p(m)},0] \mid \xi_{j+m+t}(\beta,l^-_p) < \gamma \}
    \ \ssq \ 
    \Omega_\infty \ .
\end{equation}
This means that we can compare the two sequences 
\begin{equation}
  \label{step6:referencefour}
x^1_1 \ld x^1_n \ := \ \xi_{-l^-_p(m)}(\beta^*,l^-_{q'}) \ld 
\xi_{-1}(\beta^*,l^-_{q'}) 
\end{equation}
and 
\begin{equation}
  \label{step4:newfour}
x^2_1 \ld x^2_n \ := \ \xi_{m+t-l^-_{p(m)}}(\beta,l^-_p) \ld 
\xi_{m+t-1}(\beta,l^-_p)
\end{equation}
 via Lemma~\ref{lem:orbitscontraction} with %
$\epsilon \ := \ L_2 \cdot d(\omm,0) \ \in $   
$(\alpha^{-p(m)},\alpha^{-(p(m)-1)}]$
to obtain that%
\footnote{ Note that $d(\omega_{m+t},0) \leq \frac{2\epsilon}{L_2}$ (see   
  (\ref{step6:j})) and %
  $|\beta-\beta^*| \leq 2\alpha^{-p(m)} \leq 2\epsilon$ by  
  (\ref{step6:m}), such that $\err(\ldots) \leq K \cdot \epsilon$ by
  Remark~\ref{bem:errorterm}~. Further, it follows from
  (\ref{step6:k}) and (\ref{step6:o}) that $\eta(j,n) = \#([-(n-j),-1]
  \cap \Omega_\infty) \leq \frac{n-j}{10}$ (see
  (\ref{cond:hfunctions2})). Finally $n=l^-_{p(m)} \geq 4p(m)$ by
  (\ref{eq:lpmestimates}), such that $\alpha^{-\frac{n}{4}} \leq
  \epsilon$.\label{foot:step6f}}
\begin{equation} \label{step6:p}
       |\xi_{m+t}(\beta,l^-_p) - \xi_0(\beta^*,l^-_{q'})| \ \leq \
        \epsilon \cdot (6+K\cdot S_\infty(\alpha^{\viertel})) \ .
\end{equation}
Note that (\ref{step6:k}) and (\ref{step6:o}) in particular imply that
$\xi_0(\beta^*,l^-_{q'})\geq \gamma$ and $\xi_{m+t}(\beta,l^-_p) \geq
\gamma$. As $d(\omega_{m+t},0) \geq \frac{3}{4} \cdot
\frac{\epsilon}{L_2}$ (see (\ref{step6:j})), (\ref{cond:sharppeak}) in
combination with (\ref{cond:Fcontraction}) therefore implies
\begin{eqnarray}
      \lefteqn{\xi_{m+t+1}(\beta,l^-_p) \ \geq } \nonumber \\ & \geq &
      \xi_1(\beta^*,l^-_{q'}) +
      \frac{3\epsilon}{4}  - \epsilon \cdot \frac{6 +
      K \cdot S_\infty(\alpha^{\viertel})}{2\walpha} 
       \ \stackrel{(\ref{cond:alpha3})}{\geq} \
      \xi_1(\beta^*,l^-_{q'}) + \frac{\epsilon}{2} \ . \label{step6:crucial}
\end{eqnarray}
Now first assume $p(m) \geq 2$, such that $\epsilon \leq
\frac{\alpha^{-1}}{L_2}$. (The case $p(m) =1$ has to be treated separately,
see below.) Then we can apply Lemma \ref{lem:orbitsthrowout}(b), with
$\epsilon$ as above, to compare the
orbits
\begin{equation}
  \label{step6:q}
x^1_1 \ld x^1_n \ := \ \xi_1(\beta^*,l^-_{q'}) \ld
\xi_{l^+_{p(m)}}(\beta^*,l^-_{q'})
\end{equation} 
and
\begin{equation}
  \label{step6:r}
x^2_1 \ld x^2_n \ := \ \xi_{m+t+1}(\beta,l^-_p)
\ld \xi_{m+t+l^+_{p(m)}}(\beta,l^-_p)
\end{equation}
to conclude that%
\footnote{We have $q=p(m)-1$ and $\err(\ldots) \leq K\epsilon$ before (see
  Footnote~\ref{foot:step6f}). (\ref{step6:l}) yields that $x^1_{n+1} \in \Balphcl$
  (note that $l^+_q+1 \in R_m$ by (\ref{eq:lpminrn})). Further we have
\[
\tau(n) \ \leq \ \#([1,l^+_{p(m)}] \setminus R_m \ \leq \ 
\max\left\{0,\frac{2p(m)-5}{4}\right\}
\]
as well as $\tau(j) \leq \#([1,j] \setminus R_m) \leq \frac{j}{8}$ by
Lemma~\ref{lem:regsets}(a). $\tau(n)-\tau(j) \leq \frac{n-j}{6}$
follows again from (\ref{eq:etajnestimate}), and finally $n =
l^+_{p(m)} \geq 5(p(m)-1)$ by (\ref{eq:lpmestimates}).\label{foot:step6g}}
\begin{equation} \label{step6:s}
    \xi_{m+t+l^+_{p(m)}+1}(\beta,l^-_p) \ \geq \ \frac{2}{\alpha} \ .
\end{equation}
As $J_r = J(m)$ is a maximal interval in $\Gamma_{l^+_{p_{r-1}}}$ we
have $\lambda^+(m)+1 \in R_{l^+_{p_{r-1}}}$.  Therefore
$\lambda^+(m)+1+t \in R(J)$ follows from the recursive structure of
the regular sets.  Consequently, we can choose $\tilde{k} =
\lambda^+(m)+1+t = m+l^+_{p(m)}+t+1$.

\ \\ Finally, suppose $p(m)=1$. In this case we still have
$\xi_{m+t+1}(\beta,l^-_p) \geq \xi_1(\beta^*,l^-_{q'}) +
\frac{\epsilon}{2}$ by (\ref{step6:crucial}). There are two
possibilities: Either $\xi_{m+t+1}(\beta,l^-_p) \geq
\frac{2}{\alpha}$. As $m+1 \in R_{l_{p_{r-1}}}$ (see
Remark~\ref{bem:regularsets}(c)), we have $m+t+1 \in R(J)$ due to the
recursive structure of this set. Thus, we can choose $\tilde{k} =
m+t+1$. On the other hand, if $\xi_{m+t+1}(\beta,l^-_p) \in
B_{\frac{2}{\alpha}}(0)$ then we can apply (\ref{cond:Fcontraction})
again and obtain
\[
\xi_{m+t+2}(\beta,l^-_p) \ \geq \
\xi_2(\beta^*,l^-_{q'}) + 2\walpha\cdot\epsilon - K \epsilon \
\stackrel{(\ref{cond:alpha1})}{\geq} \ 
\xi_2(\beta^*,l^-_{q'}) + \walpha\cdot\epsilon \geq \frac{1}{\alpha}
\]
by (\ref{cond:alphagamma0}), as
$\xi_1(\beta^*,l^-_{q'}),\xi_2(\beta^*,l^-_{q'}) \in \Balphcl$ by
(\ref{step6:l}) and $\epsilon \geq \alphtel$.  Thus, we can choose
$\tilde{k} = m+t+2$ in this case. Note that $m+t+2$ is contained in $
R(J)$ for the same reasons as $m+t+1$.

\qed

\solidqed


\section{Construction of the sink-source-orbits: Symmetric forcing}
\label{Symmetricsetting}
For the symmetric setting, we will use two additional assumptions on the
parameters, namely 
\begin{eqnarray} \label{cond:alphagammasymmetric}
 \frac{4\gamma}{L_2} + \frac{S_\infty(\alpha)}{\alpha \cdot L_2} & < & \halb \
 ; \\ g_{| B_{\frac{4\gamma}{L_2}}(0)} \ \geq 0 & \textrm{ and } & g_{|
 B_{\frac{4\gamma}{L_2}}(\halb)} \ \leq 0 \ . \label{cond:gvalues}
\end{eqnarray}
Due to the Lipschitz-continuity of $g$ by assumption (\ref{eq:g-lipschitz}), condition
 (\ref{cond:gvalues}) can of course be ensured by choosing $\gamma \leq
 L_2/(4L_1)$. \medskip

Further, we remark that the symmetry condition (\ref{eq:systemsymmetry})
reduces the possible alternatives in Theorem~\ref{thm:schwarzian} and leads to
the following corollary:
\begin{cor}[Corollary 4.3 in \cite{jaeger:2003}]
  \label{cor:symmetricschwarzian}
Suppose $T$ satisfies all assertions of
Theorem~\ref{thm:schwarzian} and has the symmetry given by
(\ref{eq:systemsymmetry}). Then one of the following holds:
\romanlist
\item There exists one invariant graph $\varphi$ with
  $\lambda(\varphi) \leq 0$. If $\varphi$ has a negative Lyapunov
  exponent, it is always continuous. Otherwise the equivalence class
  contains at least an upper and a lower semi-continuous
  representative.
\item There exist three invariant graphs $\varphi^- \leq \psi \leq
  \varphi^+$ with $\lambda(\varphi^-) = \lambda(\varphi^+) < 0$ and
  $\lambda(\psi) > 0$. $\varphi^-$ is always lower semi-continuous and
  $\varphi^+$ is always upper semi-continuous. Further, if one
  of the three graphs is continuous then so are the other two, if none
  of them is continuous there holds
  \[\esscl{\Phi^-} =
  \esscl{\Psi} = \esscl{\Phi^+} \ .
  \] 
  In addition, there holds 
  \[ \textstyle
  \varphi^-(\theta) = -\varphi^+(\theta+\halb)
  \]
  and 
  \[ \textstyle
  \psi(\theta) \ = \  -\psi(\theta+\halb) \ .
  \]
\listend
\end{cor}
Consequently, if we can show that there exists an SNA in a system of this
kind, then we are automatically in situation (ii). Thus there will be two
symmetric strange non-chaotic attractors which embrace a self-symmetric
strange non-chaotic repellor, as claimed in Theorem~\ref{thm:symmetricsna}~.
\medskip

In order to repeat the construction from Section~\ref{Construction} for the
case of symmetric forcing, we have to define admissible times and the
sets $R_N$ again. However, this time there are two critical intervals instead
of one, namely $B_{\frac{4\gamma}{L_2}}(0)$ and
$B_{\frac{4\gamma}{L_2}}(\halb)$, corresponding to the maximum and minimum of
the forcing function $g$. Therefore, we have to modify
Definition~\ref{def:omegasets} in the following way:
\begin{definition}
   \label{def:modifiedomegasets}
For $p \in \N_0 \cup \{\infty\}$ let $Q_p : \Z \ra \N_0$ be defined by
\[
    Q_p(j) := \left\{ \begin{array}{cll} q  & \mathrm{if} &
    d(\omj,\{0,\halb\}) \in \left[ S_{p-q+1}(\alpha) \cdot
    \frac{\alpha^{-q}}{L_2},S_{p-q+2}(\alpha)\cdot \frac{\alpha^{-(q-1)}}{L_2}
    \right) \ \textrm{ for } q \geq 2
    \\ \\
     1 & \mathrm{if} & d(\omj,\{0,\halb\}) \in \left[ S_p(\alpha) \cdot 
    \frac{\alpha^{-1}}{L_2},  \frac{4\gamma}{L_2} +
    S_{p}(\alpha)\cdot \frac{\alpha^{-1}}{L_2} \cdot (1-\ind_{\{0\}}(p))  \right)
    \ \ \ . 
     \\ \\
     0 & \mathrm{if} & d(\omj,\{0,\halb\}) \geq \frac{4\gamma}{L_2} +
    S_{p}(\alpha)\cdot \frac{\alpha^{-1}}{L_2} \cdot (1-\ind_{\{0\}}(p)) 
     \end{array} \right. 
\]
if $j \in \Z \setminus \{ 0 \}$ and $Q_p(0) := 0$. Again,
let $p(j) := Q_0(j)$. Further let
\[
    \begin{array}{llll} \textstyle
        \nutil(q) & := & 
    \min\left\{ j \in \N \mid d(\omj,\{0,\halb\}) < 
    3S_\infty(\alpha)\cdot \frac{\alpha^{-(q-1)}}{L_2} \right\} & \textrm{if } q \geq 2
    \ \textrm{ and } \\ \\
    \nutil(1) & := & \min\left\{ j
    \in \N \mid d(\omj,\{0,\halb\}) < 3\left(\frac{4\gamma}{L_2} +
    S_\infty(\alpha)\cdot\frac{\alpha^{-1}}{L_2} \right) \right\} \ .
\end{array}
\]
Apart from this, we define all the quantities
$\Omega^{(\pm)}_p(j),\Omega^{(\pm)},\Omtil^{(\pm)}$ and $\nu$
exactly in the same way as in Definition~\ref{def:omegasets}, only
using the altered definitions of the functions $Q_p$. Finally, we let
\[
s(j) \ := \ \left\{ \begin{array}{rll} 1 & \textit{if } \ \ d(\omj,0)
    \ \leq \ \frac{4\gamma}{L_2} +
    \frac{S_\infty(\alpha)}{\alpha \cdot L_2} \\ \\
    -1 & \textit{if } \ \ d(\omj,\halb) \ \leq \ \frac{4\gamma}{L_2} +
    \frac{S_\infty(\alpha)}{\alpha \cdot L_2} \\ \\
    0 & \textit{otherwise} \end{array} \right. \ .
\]

\end{definition}
In other words, we have just replaced $d(\omj,0)$ by
$d(\omj,\{0,\halb\})$ and introduced the function $s$ in order to tell
whether $\omj$ is close to $0$ or to $\halb$. However, if we let
$\tilde{\omega} := 2 \omega \bmod 1$ there holds
\[ \textstyle
       d(\omj,\{0,\halb\}) \ = \ \halb d(\tilde{\omega}_j,0) \ .
\]
This means that Definition~\ref{def:omegasets} with $\tilde{\omega}$
and $\tilde{L}_2 := \halb L_2$ yields exactly the same objects as
Definition~\ref{def:modifiedomegasets} with $\omega$ and $L_2$.
Therefore, if we define all the quantities $l^\pm_q,\ J(m),$ $A_N,\ 
\Lambda_N,$ $R_N,\ \Gamma_N$, ect.\ exactly in the same way as in
Section~\ref{Tools}, only with respect to
Definition~\ref{def:modifiedomegasets} instead of
Definition~\ref{def:omegasets}, then all the results from
Sections~\ref{Approximatingsets}--\ref{Regulartimes} will literally
stay true. The only exception is Lemma~\ref{lem:regularsetsbasic}(b),
where we can even replace $d(\omj,0)$ by $d(\omj,\{0,\halb\})$.
Further, in Section~\ref{Comparingorbits} we did not use any specific
assumption on $g$ apart from the Lipschitz-continuity. Thus, we have
all the tools from Section~\ref{Tools} available again.

Therefore, the only difference to the preceding section is the fact that the
mapping $\beta \mapsto \xi_n(\beta,l)$ is not necessarily monotone
anymore (where the $\xi_n(\beta,l)$ are defined exactly as before, see
Definition~\ref{def:parameterfamily}).  Hence, instead of
considering arbitrary $\beta$ as in Induction
statement~\ref{thm:indscheme} we  have to restrict to certain intervals %
$I^q_n = [\beta^+_{q,n},\beta^-_{q,n}]$ %
($q\in\N_0,\ n \in [l^+_q+1,\nu(q+1)]$ admissible) on which the
dependence of $\xi_n(\beta,l^-_q)$ on $\beta$ is monotone.  The parameters
$\beta^\pm_{q,n}$ will again satisfy
\begin{equation}
  \label{eq:betapluspndefS} 
  \xi_n(\beta_{q,n}^+,l^-_q) \ = \ \alphtel 
\end{equation}
and
\begin{equation} 
  \label{eq:betaminpndefS}
  \xi_n(\beta_{q,n}^-,l^-_q) \ = \ -\alphtel \ ,
\end{equation}
but they cannot be uniquely defined by these equations anymore. 

The fact which makes up for the lack of monotonicity, and for the
existence of the second critical region $B_{\frac{4\gamma}{L_2}}(\halb)$, is
that by deriving information about the orbits $\xi_n(\beta,l)$ we get
another set of reference orbits for free: It follows directly from
(\ref{eq:systemsymmetry}) that
\begin{equation} \label{eq:zetadef}
    \zeta_n(\beta,l) \ := \ T_{\beta,-\omega_l+\halb,n+l}(-3) \ = \
    -\xi_n(\beta,l) 
\end{equation}
(Similar as in Definition~\ref{def:parameterfamily}, the
$\zeta_n(\beta,l)$ correspond to the forward orbit of the points
$(\omega_{-l}+\halb,-3)$, where we suppress the $\theta$-coordinates
again). Consequently, we have
\begin{equation}
  \label{eq:xizetaequivalenceI}
  \xi_n(\beta,l) \in \Balphcl \ \equi \ \zeta_n(\beta,l) \in \Balphcl
\end{equation}
and 
\begin{equation}
  \label{eq:xizetaequivalenceII}
  \xi_n(\beta,l) \geq \gamma \ \equi \ \zeta_n(\beta,l) \leq - \gamma
  \ .
\end{equation}
In the case of symmetric forcing the induction statement reads as
follows:
\begin{indscheme}
  \label{thm:indschemeS}
Suppose the assumptions of Theorem~\ref{thm:symmetricsna} are satisfied and in
addition (\ref{cond:alpha1}), (\ref{cond:hfunctions1}),
(\ref{cond:hfunctions2}), (\ref{cond:u})--(\ref{cond:alpha3}),
(\ref{cond:alphagammasymmetric}) and (\ref{cond:gvalues}) hold.

Then for any $q \in \N_0$ and all admissible $n \in [l^+_q+1,\nu(q+1)]$ there
exists an interval $I^q_n = [\beta^+_{q,n},\beta^-_{q,n}]$, such that
$\beta^\pm_{q,n}$ satisfy (\ref{eq:betapluspndefS}) and
(\ref{eq:betaminpndefS}) and in addition
\begin{list}{\textbf{\Roman{enumi}.} \ \ }{\usecounter{enumi}}
\item $\beta \in I^q_{l^+_q+1}$ implies
  \begin{equation} \label{eq:indstatementIS}
           \xi_j(\beta,l^-_q) \ \geq  \ \gamma \ \ \ \forall j \in
           [-l^-_q,0] \setminus \Omega_\infty \ .
  \end{equation}
  Further $I^q_{l^+_q+1} \ssq \left[1+\walphtel,1+\frac{3}{\walpha}\right]$.
\item For each admissible $n \in [l^+_q+1,\nu(q+1)]$ the mapping $\beta
  \mapsto \xi_n(\beta,l^-_q)$ is strictly monotonically decreasing on
  $I^q_n$, (\ref{eq:indstatementIS}) holds for all $\beta\in I^q_n$ and  
  \begin{equation} \label{eq:indstatementIIaS} 
       I^q_n \ \ssq \ I^q_j \ \ssq \
       \left[1+\walphtel,1+\frac{3}{\walpha}\right]\ \ \ \ \forall j \in A_n \cap 
       [l^+_q+1,n] \ . 
  \end{equation}
  Further, for any $\beta \in I^q_n$ there holds
  \begin{equation} \label{eq:indstatementIIcS}
       \xi_j(\beta,l^-_q) \ \in \ \Balphcl \ \ \ \ \ \forall j \in R_n  \ .
  \end{equation}
  \item 
  \alphlist
  \item If $n_1 \in [l^+_q+1,\nu(q+1)]$ for some $q \geq 1$ there
    holds
    \begin{equation} \label{eq:indstatementIIIaS}
      |\beta^+_{q,n_1} - \beta^-_{q,n_1}| \ \leq \
      2\alpha^{-\frac{n_1}{4}} \ .
    \end{equation}
    In particular, in combination with (\ref{eq:indstatementIIaS})
    this implies that 
    \begin{equation} \label{eq:indstatementIIIbS}
      |\beta^\pm_{q,n_1} - \beta| \ \leq \
      2\alpha^{-\frac{n_1}{4}}  \ \ \ \ \ \forall \beta \in I^q_{n_2}
    \end{equation}
    whenever $n_2 \in [l^+_q+1,\nu(q+1)]$ and $n_1 \in A_{n_2}$ (as
    $I^q_{n_2} \ssq I^q_{n_1}$ in this case).
    \item Let $1 \leq q_1 < q$, $n_1 \in [l^+_{q_1}+1,\nu(q_1+1)]$  and $n_2
      \in [l^+_{q}+1,\nu(q)+1]$. Then 
      \begin{equation} \label{eq:indstatementIIIcS}
              |\beta^+_{q_1,n_1} - \beta^+_{q,n_2}| \ \leq \
              3 \cdot \hspace{-0.7eM}\sum_{i=q_1+1}^q \alpha^{-i}\ \leq \ \alpha^{-q_1} \ .
      \end{equation} 
  \listend
\end{list}
\end{indscheme}

Theorem \ref{thm:symmetricsna} now follows in exactly the same way as
Theorem~\ref{thm:snaexistence} from Induction scheme~\ref{thm:indscheme} (we
do not repeat the argument here). The additional statements about the symmetry
follow from Corollary~\ref{cor:symmetricschwarzian}~.

However, due to the lack of monotonicity we are not able to derive any further
information about the sink-source-orbit or the bifurcation scenario as in the
case of one-sided forcing. In particular, we have to leave open here whether
$\beta_0$ is the only parameter at which an SNA occurs, or if this does indeed
happen over a small parameter interval as the numerical observations suggest
(compare Section~\ref{Pitchfork}).
\medskip

%


\subsection{Proof of the induction scheme}

\textbf{Standing assumption:} In this whole subsection, we always assume that
the assumptions of Induction scheme~\ref{thm:indschemeS} are satisfied. 
\medskip

In order to start the induction we will need the following equivalent to
Lemma~\ref{lem:basicestimate}, which can be proved in exactly in the same way
(using that $d(\omj,0) \geq \frac{3\gamma}{L_2}$ implies $g(\omj) \leq
1-3\gamma$ by (\ref{cond:symmetricpeak}) and (\ref{cond:gvalues}), and
similarly $d(\omj,\halb) \geq \frac{3\gamma}{L_2}$ implies $g(\omj)
\geq - (1-3\gamma)$). 
\begin{lem}
  \label{lem:basicestimateS} Suppose that $\beta \leq 1 + \frac{4}{\walpha}$ and $j \geq
  -l$. If $d(\omj,0) \geq \frac{3\gamma}{L_2}$, then $\xi_j(\beta,l) \geq
  \alphtel$ implies $\xi_{j+1}(\beta,l) \geq \gamma$. Similarly, if
  $d(\omj,\halb) \geq \frac{3\gamma}{L_2}$ then $\xi_j(\beta,l) \leq -
  \alphtel$ implies $\xi_{j+1}(\beta,l) \leq -\gamma$.  Consequently,
  $\xi_{j+1}(\beta,l) \in \Balphcl$ implies $\xi_j(\beta,l) \in \Balphcl$
  whenever $d(\omj,\{0,\halb\}) \geq \frac{3\gamma}{L_2}$.
\end{lem}

Further, the following lemma replaces Lemma~\ref{lem:parameterestimates}~. It
will be needed to derive the required estimates on the parameters
$\beta^\pm_{q,n}$ as well as the monotonicity of $\beta \mapsto
\xi_n(\beta,l^-_q)$ on $I^q_n$.
\begin{lem}
  \label{lem:parameterestimatesS} Let $q \in \N$ and let $n
  \in [l^+_q+1,\nu(q+1)]$ be admissible. Further, assume
\begin{equation}
  \label{eq:pmesttimeabove}
  \xi_j(\beta,l^-_q) \ \geq \ \gamma \ \ \ \ \ \forall j \in
  [-l^-_q,0] \setminus \Omega_\infty
\end{equation}
and
\begin{equation}
  \label{eq:pmesttimebelow}
  \xi_j(\beta,l^-_q) \ \in \ \Balphcl \ \ \ \ \ \forall j \in R_n
  \setminus \{n\} \ .
\end{equation}
Then
\begin{equation}
  \label{eq:derivativeestimate}
  \frac{\partial}{\partial \beta} \xi_n(\beta,l^-_q) \ \leq \ -
  \alpha^{\frac{n-1}{4}} \ .
\end{equation}
\end{lem}
\proof\ 
We have
\begin{equation}
  \label{eq:derivativerule}
  \frac{\partial}{\partial\beta} \xi_{j+1}(\beta,l^-_q) \ = \
  F'(\xi_j(\beta,l^-_q)) \cdot \frac{\partial}{\partial\beta}
  \xi_j(\beta,l^-_q) - g(\omj)
\end{equation}
(compare (\ref{eq:xiderivatives})). In order to prove
(\ref{eq:derivativeestimate}) we first have to obtain a suitable
upper bound on
$|\frac{\partial}{\partial\beta}\xi_0(\beta,l^-_q)|$. Let
\[
    \eta(j) \ := \ \# \left( [-j,-1] \cap \Omega_\infty\right) \ .
\] 
We claim that under assumption (\ref{eq:pmesttimeabove}) and for any $l
\in [0,l^-_q]$ there holds
\begin{equation}
  \label{eq:destinduction}
  \left|\frac{\partial}{\partial\delta} \xi_0(\beta,l^-_q)\right| \ \leq \
  \left|\frac{\partial}{\partial\delta} \xi_{-l}(\beta,l^-_q)\right| \cdot
  \alpha^{-\halb(l-5\eta(l))} + \sum_{j=0}^{l-1}
  \alpha^{-\halb(j-5\eta(j))} \ .
\end{equation}
As $\eta(j) \leq \frac{j}{10}$ by
(\ref{cond:hfunctions2}) and %
$\frac{\partial}{\partial\delta} \xi_{-l^-_q}(\beta,l^-_q) = 0$ by
definition, this implies
\[
  \left|\frac{\partial}{\partial\delta} \xi_0(\beta,l^-_q)\right| \ \leq \
  S_\infty(\alpha^{\viertel}) \ .
\]
Using the fact that $\xi_0(\beta,l^-_q) \geq \gamma$ by assumption, 
this further yields
\begin{eqnarray}  \label{eq:xizeroderivative} \nonumber
     \lefteqn{\frac{\partial}{\partial\delta} \xi_1 (\beta,l^-_q) \ =} \\  & = &  
     F'(\xi_0(\beta,l^-_q)) \cdot \frac{\partial}{\partial\delta}
     \xi_0(\beta,l^-_q) - 1 \
     \stackrel{(\ref{cond:Fcontraction})}{\leq}  \ - 1 +
     \frac{S_\infty(\alpha^{\viertel})}{2\walpha} \
     \stackrel{(\ref{cond:alpha3})}{\leq} \ -\halb \ . 
\end{eqnarray}
We prove (\ref{eq:destinduction}) by induction on $l$. For $l = 0$ the
statement is obvious. In order to prove the induction step $l \ra
l+1$, first suppose that $-(l+1) \notin \Omega_\infty$, such that
$\eta(l+1) = \eta(l)$ and $\xi_{-(l+1)}(\beta,l^-_q) \geq \gamma$.
Then, using (\ref{eq:derivativerule}) we obtain
\begin{eqnarray*}
  \lefteqn{ \left|\frac{\partial}{\partial\delta}
      \xi_0(\beta,l^-_q)\right| \ \leq \ 
  \left|\frac{\partial}{\partial\delta} \xi_{-l}(\beta,l^-_q)\right|
  \cdot \alpha^{-\halb(l-5\eta(l))} + \sum_{j=0}^{l-1}
  \alpha^{-\halb(j-5\eta(j))} } \\
  & = & \left|F'(\xi_{-(l+1)}(\beta,l^-_q)) \cdot
  \frac{\partial}{\partial\delta} \xi_{-(l+1)}(\beta,l^-_q) -
  g(\omega_{-(l+1)})\right| \cdot
  \alpha^{-\halb(l-5\eta(l))} \\ && + \sum_{j=0}^{l-1}
  \alpha^{-\halb(j-5\eta(j))} \\
  & \stackrel{(\ref{cond:Fcontraction})}{\leq} & \left(\alpha^{-\halb} 
      \cdot \left|\frac{\partial}{\partial\delta}
    \xi_{-(l+1)}(\beta,l^-_q)\right| + 
  1\right) \cdot \alpha^{-\halb(l-5\eta(l))} + \sum_{j=0}^{l-1}
  \alpha^{-\halb(j-5\eta(j))}   \\
  & = & \left|\frac{\partial}{\partial\delta}
  \xi_{-(l+1)}(\beta,l^-_q)\right| \cdot
\alpha^{-\halb(l+1-5\eta(l+1))} + \sum_{j=0}^{l} 
  \alpha^{-\halb(j-5\eta(j))} \ .
\end{eqnarray*}
The case $\eta(l+1) = \eta(l)+1$ is treated similarly, using
(\ref{cond:Funiformbounds}) instead of (\ref{cond:Fcontraction})
(compare with the proof of Lemma~\ref{lem:orbitscontraction}). This
proves (\ref{eq:destinduction}), such that (\ref{eq:xizeroderivative}) 
holds.

\ \\
Now we can turn to prove (\ref{eq:derivativeestimate}). For any $k \in 
\N$ let
\[
  \tau(k) \ : = \ \#([1,k-1]\setminus R_n)  \ . 
\]
We will show the following statement by induction on $k$:
\begin{equation}
  \label{eq:destinductionII}
  \frac{\partial}{\partial\delta} \xi_k(\beta,l^-_q) \ \leq \ -\halb
  \cdot \left(\frac{3\walpha}{2} \right)^{k-1-5\tau(k)} \ \ \ \ \ \forall k 
  \in [1,n] \ .
\end{equation}
As $\tau(n) \leq \frac{n-1}{10}$ by Lemma~\ref{lem:regsets}(a), this
implies (\ref{eq:derivativeestimate}) whenever $n \geq l^+_q+1$. Note
that $l^+_q \geq 3$ by (\ref{eq:lpmestimates}) and $\tau(n) = 0$ for
all $n \leq 10$.

For $k = 1$ the statement is true by
(\ref{eq:xizeroderivative}). Suppose that (\ref{eq:destinductionII})
holds for some $k \geq 1$ and first assume that $\tau(k+1) = \tau(k)$.
Then
\begin{eqnarray*}
  \lefteqn{\frac{\partial}{\partial\delta} \xi_{k+1}(\beta,l^-_q) \ =
    \ F'(\xi_k(\beta,l^-_q)) \cdot
    \frac{\partial}{\partial\delta}\xi_k(\beta,l^-_q) - g(\omk)  } \\
  & \stackrel{(\ref{cond:Fexpansion})}{\leq} & - 2\walpha \cdot \halb
  \cdot 
  \left(\frac{3\walpha}{2} \right)^{k-1-5\tau(k)} + 1 \\
  & \stackrel{(*)}{\leq} & - 
  (2\walpha - 2) \cdot \halb \cdot \left(\frac{3\walpha}{2}
  \right)^{k-1-5\tau(k)} \\
  & \stackrel{(\ref{cond:alphagamma0})}{\leq} &  -\halb
  \cdot \left(\frac{3\walpha}{2} \right)^{k-5\tau(k+1)} \ ,
\end{eqnarray*}
($*$) where $\tau(k) \leq \frac{k-1}{10}$ by Lemma~\ref{lem:regsets}(a)
ensures that $\left(\frac{3\walpha}{2} \right)^{k-1-5\tau(k)}$ is
always larger than $1$.  The case $\tau(k+1) = \tau(k)+1$ is treated
similar again, using (\ref{cond:Funiformbounds}) instead of
(\ref{cond:Fexpansion}) (compare with the proof of
Lemma~\ref{lem:orbitsexpansion}). Thus we have proved
(\ref{eq:destinductionII}) and thereby the lemma.

\qed
\bigskip

As in Section~\ref{Construction}, in order to prove Induction
scheme~\ref{thm:indschemeS} we proceed in six steps. The overall strategy
needs some slight modifications in comparison to the case of one-sided
forcing, but in many cases the proofs of the required estimates stay literally
the same. In such situations we will not repeat all the details, but refer to
the corresponding passages of the previous section instead.

\ \\
\underline{\textbf{Step 1:}} \ \ \ \textit{Proof of the statement for
  $q=0$} \\ \ \\
Part I: Recall that $l^-_0 = l^+_0 = 0$ and note that $\xi_0(\beta,0)
= 3 \geq \gamma$ by definition, such that (\ref{eq:indstatementIS})
holds automatically. As $\frac{\partial}{\partial \beta}
\xi_1(\beta,0) = -1$, we can construct the interval $I^0_1$ by
uniquely defining $\beta^\pm_{0,1}$ via (\ref{eq:betapluspndefS}) and
(\ref{eq:betaminpndefS}). Further, we have
$\xi_1(\beta,0) = F(3) - \beta$. Using (\ref{cond:Ffixedpoints}) and
(\ref{cond:Fcontraction}), it is easy to check that 
\[ \textstyle
     F(3) \ \in \
     [x_\alpha,x_\alpha+\frac{2-\frac{2}{\walpha}}{2\walpha}] \ \ssq \ 
     [1+\frac{1}{\walpha}+\alphtel,1+\frac{3}{\walpha}-\alphtel] \ .
\]
Therefore %
$I^1_0 = [\beta^+_{0,1},\beta^-_{0,1}]$ must be contained in
$[1+\frac{1}{\walpha},1+\frac{3}{\walpha}]$.

\ \\
Parts II: We proceed by induction on $n$. For $n=1$ the statement
follows from the above. Suppose we have defined the intervals $I^0_n
\ssq [1+\frac{1}{\walpha},1+\frac{3}{\walpha}]$ with the stated
properties for all $n \leq N$, $N \in [1,\nu(1)-1]$. As $p(N) = 0$,
Lemma~\ref{lem:basicestimateS} yields that
\[
    \xi_{N+1}(\beta^+_{0,N},0)  >  \alphtel \ \ \ \  \textrm{ and }
    \ \ \ \ \xi_{N+1}(\beta^-_{0,N},0)  < -\alphtel \ .
\]
This means that we can find $\beta^\pm_{0,N+1}$ in $I^0_N$ which
satisfy (\ref{eq:betapluspndefS}) and
(\ref{eq:betaminpndefS}). Consequently %
$I^0_{N+1} := [\beta^+_{0,N+1},\beta^-_{0,N+1}] \ssq I^0_N$. It then
follows from Part II of the induction statement for $N$, that
$I^0_{N+1} \ssq I^0_j \ \forall j \in [1,N]=R_N$, in particular
$I^0_{N+1} \ssq I^0_1 \ssq [1+\walphtel, 1+\frac{3}{\walpha}]$ (note
that $A_N = [1,N]$ as $N \leq \nu(1)$). This proves
(\ref{eq:indstatementIS}) and (\ref{eq:indstatementIIaS}).
 
In order to see (\ref{eq:indstatementIIcS}) suppose that $\beta \in 
I^0_{N+1}$. Then $\xi_N(\beta,0) \in \Balphcl$ by the definition of 
$I^0_{N+1}\ssq I^0_N$ above, and therefore %
$\xi_j(\beta,0) \in \Balphcl \ \forall j \in [1,N] = R_N $ follows
from Part II of the induction statement for $N$. Finally, we can now
use Lemma~\ref{lem:parameterestimatesS} to see that
\begin{equation} \label{eq:qzeroderivative}
\frac{\partial}{\partial \beta} \xi_{N+1}(\beta,0) \ \leq \ -
\alpha^{\frac{N}{4}} \ .
\end{equation} 
This ensures the monotonicity of %
$\beta \mapsto \xi_{N+1}(\beta,0)$.

\ \\
As Part III of the induction statement is void for $q=0$, this
completes Step I.

\solidqed


\ \\
It remains to prove the induction step.  Assume that the statement of
Induction scheme \ref{thm:indschemeS} holds for all $q \leq p-1$. As in
Section~\ref{Onesidedproof}, the next two steps will prove Part I of
the induction statement for $p$. Further, we can again assume in Step
2 and 3 that 
\begin{equation}
  \label{eq:pgeqoneS} p \ \geq \ 2 \ .
\end{equation}
For the case $p=1$ note that the analogue of Lemma~\ref{lem:inductionstart}
holds again in the case of symmetric forcing, with $d(\omj,0)$ being replaced
by $d(\omj,\{0,\halb\})$, and this already shows Part I for $p=1$.


\ \\
\underline{\textbf{Step 2:}}  \ \ \ 
\textit{If
  $|\beta-\beta^+_{p-1,\nu(p)}| \leq \alpha^{-p}$, then
  $\xi_j(\beta,l^-_p) \geq \gamma \ \forall j \in [-l^-_p,0] \setminus
  \Omega_\infty$.}

\ \\
Actually, this follows in exactly the same way as Step 2 in
Section~\ref{Onesidedproof}~. The crucial observation is the fact that
Lemma~\ref{lem:approaching} literally stays true in the
situation of this section. As we will also need the statement for the reversed
inequalities in the later steps, we restate it here:
\begin{lem}
  \label{lem:approachingS}
  Let $q \geq 1$ $,l^*,l \geq 0$, $\beta^* \in
  \left[1+\walphtel,1+\frac{3}{\walpha}\right]$ and $|\beta - \beta^*| \leq 2\alpha^{-q}$.
  Suppose that $m$ is admissible, $p(m) \geq q$ and either $k=0$ or
  $p(k) \geq q$. Further, suppose
\[
      \xi_j(\beta^*,l^*) \ \in \ \Balphcl \ \ \ \forall j \in R_m \
\] 
and $\xi_{m+k-l^-_q}(\beta,l) \geq \gamma$. Then 
\[
        \{ j \in [m-l^-_q,m] \mid \xi_{j+k}(\beta,l) < \gamma \} \ \ssq \
        \Omtil_{q-2} \ .
\]
Similarly, if $\xi_{m+k-l^-_q}(\beta,l) \leq -\gamma$ then 
\[
        \{ j \in [m-l^-_q,m] \mid \xi_{j+k}(\beta,l) > -\gamma \} \ \ssq \
        \Omtil_{q-2} \ .
\]
\end{lem}
The application of this lemma in order to show the statement of Step 2 is 
exactly the same as in Section~\ref{Onesidedproof}~. The
proof of the lemma is the same as for Lemma~\ref{lem:approaching}, apart from two slight
modifications: First of all, Lemma~\ref{lem:basicestimateS} has to be
used instead of Lemma~\ref{lem:basicestimate}~. Secondly, in order to
show (\ref{eq:approaching}) two cases have to be distinguished. If
$s(k) = 1$ nothing changes at all. For the second case $s(k) =
-1$ it suffices just to replace the reference orbit
\[
   x^1_1 \ld x^1_n \ := \ \xi_{j^-_i-1}(\beta^*,l^*) \ld
   \xi_{j^+_i}(\beta^*,l^*) 
\]   
which is used for the application of Lemma~\ref{lem:orbitsthrowout}(a) 
by
\[
 x^1_1 \ld x^1_n \ := \ \zeta_{j^-_i-1}(\beta^*,l^*) \ld
   \zeta_{j^+_i}(\beta^*,l^*) \ .
\]
Then the reference orbit starts on the fibre $\omega_{j^-_i-1}+\halb$,
and is therefore $\frac{\alpha^{-(q-1)}}{L_2}$-close to the first fibre
$\omega_{j^-_i-1+k}$ of the second orbit
\[
   x^2_1 \ld x^2_n \ := \ \xi_{j^-_i-1+k}(\beta,l) \ld
\xi_{j^+_i+k}(\beta,l)\ ,
\]
such that the error term is sufficiently small again.  Due to
(\ref{eq:xizetaequivalenceI}) and (\ref{eq:xizetaequivalenceII}), all further
details then stay exactly the same as in the case $s(m)=1$. The reader should
be aware that even though the reference orbit changed, the set of times $R_m$
at which it stays in the expanding region is the same as before. This is all
which is needed in order to verify the assumptions of
Lemma~\ref{lem:orbitsthrowout}(a), which completes the proof of the
lemma. Finally, the additional statement for the reversed inequalities can be
shown similarly.

\solidqed


\ \\ \underline{\textbf{Step 3:}} \ \ \ \textit{Construction of
$I^p_{l^+_p+1} \ssq B_{\alpha^{-p}}(\beta^+_{p-1,\nu(p)})
$ \ .}

\ \\
Similar as in Step 3 of Section~\ref{Onesidedproof}, we define
$\beta^* := \beta^+_{p-1,\nu(p)}$, $\beta^+:=\beta^*-\alpha^{-p}$ and
$\beta^- := \beta^*+ \alpha^{-p}$. It then follows that 
\begin{equation}
  \label{step3S:a}
\xi_{l^+_p+1}(\beta^+,l^-_p) \ > \ \alphtel \ \ \ \ \textit{ and }
\ \ \ \ \ \xi_{l^+_p+1}(\beta^-,l^-_p) \ < \ - \alphtel \ .
\end{equation}
The proof is exactly the same as for Claim~\ref{claim:step3}, with reversed
inequalities for the case of $\beta^-$. This means that we can define the
parameters $\beta^\pm_{p,l^+_p+1}$ by
\begin{equation}
  \label{step3S:c}
  \beta^-_{p,l^+_p+1} \ := \ \min \left\{ \beta \in
  B_{\alpha^{-p}}(\beta^*) \mid \xi_{l^+_p+1}(\beta,l^-_p) =
  -\alphtel \right\} 
\end{equation}
and
\begin{equation}
  \label{step3S:b}
  \beta^+_{p,l^+_p+1} \ := \ \max \left\{ \beta \in
  B_{\alpha^{-p}}(\beta^*) \mid \beta < \beta^-_{p,l^+_p+1}, \
  \xi_{l^+_p+1}(\beta,l^-_p) =   \alphtel \right\} \ .
\end{equation}

Step 2 then implies
that (\ref{eq:indstatementIS}) is satisfied for 
\[
I^p_{l^+_p+1} \ := \
\left[\beta^+_{p,l^+_p+1},\beta^-_{p,l^+_p+1}\right] 
\]
and as $\beta^* \in \left[1+\walphtel,1+\frac{3}{\walpha}\right]$
there holds
\begin{equation} \label{step3S:d}  \textstyle
I^p_{l^+_p+1} \ \ssq \
\left[1+\walphtel-\alpha^{-p},1+\frac{3}{\walpha}+\alpha^{-p}\right] \ .
\end{equation}
$I^p_{l^+_p+1} \ \ssq \
\left[1+\walphtel,1+\frac{3}{\walpha}\right]$ will be shown after
Step 4. Apart from this the proof of Part I for
$q=p$ is complete.

\solidqed


\ \\
The next three steps will prove Part II and III of the induction
statement for $q=p$, proceeding by induction on %
$n \in [l^+_p+1,\nu(p)]$. Again we start the induction with
$n=l^+_p+1$.

\ \\ \underline{\textbf{Step 4:}} \ \ \
\textit{Proof of Part II for $q=p$ and $n=l^+_p+1$.}

\ \\
Let $\beta^* := \beta^+_{p-1,\nu(p)}$ again. We will prove the
following claim:
\begin{claim}
  \label{claim:step4Sa}
Suppose $\beta \in B_{\alpha^{-p}}(\beta^*)$ and
$\xi_{l^+_p+1}(\beta,l^-_p) \in \Balphcl$. Then 
\begin{equation}
  \label{step4S:a} \xi_j(\beta,l^-_p) \ \in \ \Balphcl \ \ \ \ \ \forall j \in
  R_{l^+_p+1} \ .
\end{equation}
\end{claim}
This follows more or less in the same way as Step 4 in 
Section~\ref{Onesidedproof}. Before we give the details, let us see
how this implies the statement of Part II for $q=p$ and $n=l^+_p+1$:

In Step 3 we have constructed $I^p_{l^+_p+1} \ssq
B_{\alpha^{-p}}(\beta^*)$.  Suppose $\beta \in
B_{\alpha^{-p}}(\beta^*)$ and $\xi_{l^+_p+1}(\beta,l^-_p) \in
\Balphcl$. Then Step 2 and the above claim ensure that the
assumptions (\ref{eq:pmesttimeabove}) and (\ref{eq:pmesttimebelow}) of
Lemma~\ref{lem:parameterestimatesS} are satisfied, and we obtain that
$\xi_{l^+_p+1}(\beta,l^-_p)$ is decreasing in $\beta$. In particular,
this applies to $\beta^+_{p,l^+_p+1}$. Consequently, if we increase
$\beta$ starting at $\beta^+_{p,l^+_p+1}$, then
$\xi_{l^+_p+1}(\beta,l^-_p)$ will decrease until it leaves the
interval $\Balphcl$. Due to the definition in (\ref{step3S:c}) this is
exactly the case when $\beta^-_{p,l^+_p+1}$ is reached. This yields
the required monotonicity on $I^p_{l^+_p+1}$, and
(\ref{eq:indstatementIIcS}) then follows from the claim. Note that
(\ref{eq:indstatementIS}) is already ensured by Step 2.

\ \\
The proof of Claim \ref{claim:step4Sa} is completely analogous to that 
of Step 4 in Section~\ref{Onesidedproof}: It will
follow in the same way from the the two claims below, which correspond 
to Claims~\ref{claim:step4a} and \ref{claim:step4b}~. 
\begin{claim}
   \label{claim:step4Sb}
Suppose $\beta \in B_{\alpha^{-p}}(\beta^*)$ and $\xi_{l^+_p+1}(\beta,l^-_p)
\in \Balphcl$. If $\xi_k(\beta,l^-_p) \geq \frac{2}{\alpha}$ for some $k \in
R_{l^+_p+1}$ then $\xi_{l^+_p+1}(\beta,l^-_p) \ > \ \alphtel$. Similarly, if
$\xi_k(\beta,l^-_p) \leq -\frac{2}{\alpha}$ then $\xi_{l^+_p+1}(\beta,l^-_p) \
< \ -\alphtel$.
\end{claim} 
For $\xi_k(\beta,l^-_p) \geq \frac{2}{\alpha}$ this can be shown exactly as
Claim~\ref{claim:step4a}. In the case $\xi_k(\beta,l^-_p) \leq
-\frac{2}{\alpha}$ it suffices just to reverse all inequalities. The analogue
to Claim~\ref{claim:step4b} holds as well:
\begin{claim}
  \label{claim:step4Sc}
  Suppose $\beta \in B_{\alpha^{-p}}(\beta^*)$ and
  $\xi_{l^+_p+1}(\beta,l^-_p) \in \Balphcl$. If $k\in R_{l^+_p+1},\ 
  k+1\in \Gamma_{l^+_p+1}$ and $\xi_k(\beta,l^-_p) \geq \alphtel$,
  then there exists some $\tilde{k} \in R_{l^+_p+1}$ with
  $\xi_{\tilde{k}}(\beta,l^-_p) \geq \frac{2}{\alpha}$. Similarly, if
  $\xi_k(\beta,l^-_p) \leq -\alphtel$ then there exists some
  $\tilde{k} \in R_{l^+_p+1}$ with $\xi_{\tilde{k}}(\beta,l^-_p) \leq
  -\frac{2}{\alpha}$.
\end{claim}
\proof\  In order to prove this, we can proceed as in the proof of
Claim \ref{claim:step4b}: Suppose first that $\xi_k(\beta,l^-_p) \geq
\alphtel$ and define $m,\ t$ and $q'$ in exactly the same way. As
these definitions only depend on the set $R_{l^+_p+1}$, which is the
same as before, there is no difference so far. Only instead of
(\ref{eq:domegat}) we obtain
\begin{equation}
  \label{step4S:domegat}
 d(\omega_t,\{0,{\textstyle \halb}\}) \ \leq \ \viertel \cdot
 \frac{\alpha^{-(p(m)+1)}}{L_2} 
\end{equation}
Now we can apply Lemma~\ref{lem:approachingS}, in the same way as
Lemma~\ref{lem:approaching} was applied in order to obtain
(\ref{step4:g}), to conclude that
\begin{equation} \label{step4S:b}
    \{ j \in [-l^-_{p(m)},0] \mid \xi_{j+m+t}(\beta,l^-_{p}) < \gamma \} \ \ssq \
    \Omega_\infty \ .
\end{equation}
For the further argument we have to distinguish two cases. If
$s(m+t)=1$, then we can use exactly the same comparison arguments as
in Section~\ref{Onesidedproof} to show that
$\xi_{m+t+l^+_{p(m)}+1}(\beta,l^-_p) \geq \frac{2}{\alpha}$ if
$p(m)\geq 2$. The details all remain exactly the same. Thus, we can
choose $\tilde{k} = m+t+l^+_{p(m)}+1$ if $p(m) \geq 2$ and again $\tilde{k}
= m+t+1$ or $m+t+2$ if $p(m)=1$.

On the other hand, suppose $s(m+t) = -1$. Then
$d(\omega_{m+t},0)\geq\frac{3\gamma}{L_2}$, and in addition (\ref{step4S:b})
implies that $\xi_{m+t}(\beta,l^-_p) \ \geq \
\gamma$. Lemma~\ref{lem:basicestimateS} therefore yields that
$\xi_{m+t+1}(\beta,l^-_p) \geq \gamma \geq \frac{2}{\alpha}$, such that we can
choose $\tilde{k} = m+t+1$.

\ \\ The case $\xi_k(\beta,l^-_p) \leq -\alphtel$ is then treated analogously:
First of all, application of Lemma~\ref{lem:approachingS} yields
\begin{equation} \label{step4S:c}
    \{ j \in [-l^-_{p(m)},0] \mid \xi_{j+m+t}(\beta,l^-_p) > -\gamma \} \ \ssq \
    \Omega_\infty \ ,
\end{equation}
in particular $\xi_{m+t}(\beta,l^-_p) \ \leq \ -\gamma$. If $s(m+t) =
1$, such that $d(\omega_{m+t},\halb)\geq\frac{3\gamma}{L_2}$, then
Lemma~\ref{lem:basicestimateS} yields that $\xi_{m+t+1}(\beta,l^-_p)
\leq -\gamma \leq -\frac{2}{\alpha}$ and we can choose $\tilde{k} =
m+t+1$.

On the other hand, if $s(m+t) = -1$, then we can again apply similar
comparison arguments as in the proof of Claim~\ref{claim:step4b} to
conclude that $\xi_{m+t+l^+_{p(m)}+1}(\beta,l^-_p) \leq
-\frac{2}{\alpha}$ if $p(m)\geq 2$ (and $\xi_{m+t+1}(\beta,l^-_p) \leq
-\frac{2}{\alpha}$ or $\xi_{m+t+2}(\beta,l^-_p) \leq
-\frac{2}{\alpha}$ if $p(m)=1$). Apart from the reversed inequalities,
the only difference now is that the reference orbits
$\xi_{-l^-_{p(m)}}(\beta^*,l^-_{q'}),$ $\ldots,$
$\xi_{-1}(\beta^*,l^-_{q'})$ and $\xi_1(\beta^*,l^-_{q'}) \ld
\xi_{l^+_{p(m)}}(\beta^*,l^-_{q'})$ in (\ref{step4:referencetwo}) and
(\ref{step4:i}) have to be replaced by
$\zeta_{-l^-_{p(m)}}(\beta^*,l^-_{q'}),$ $\ldots$,
$\zeta_{-1}(\beta^*,l^-_{q'})$ and $\zeta_1(\beta^*,l^-_{q'}),$
$\ldots,$ $\zeta_{l^+_{p(m)}}(\beta^*,l^-_{q'})$, respectively. Due to
(\ref{eq:xizetaequivalenceI}) and (\ref{eq:xizetaequivalenceII}), all
other details remain exactly the same as before, with
(\ref{cond:sharppeak}) being replaced by (\ref{cond:symmetricpeak}).

\qed

\solidqed

\ \\
Now we can also show that $I^p_{l^+_p+1} \ssq
\left[1+\walphtel,1+\frac{3}{\walpha}\right]$, which completes the
proof of Part I of the induction statement for $p$. Suppose that
$\beta \in I^p_{l^+_p+1}$. Then, due to Step~2 and the construction of
$I^p_{l^+_p+1} \ssq B_{\alpha^{-p}}(\beta^+_{p-1,\nu(p)})$ in Step 3,
(\ref{eq:indstatementIS}) holds, such that in particular
$\xi_0(\beta,l^-_p) \geq \gamma$. Thus, it follows from
(\ref{cond:Ffixedpoints}) and (\ref{cond:Fcontraction}) that
\[
    \xi_1(\beta,l^-_p) \ \in \
    \left[1+\frac{3}{2\walpha}-\beta,1+\frac{3}{\walpha}-\alphtel-\beta\right] .
\]
As Step 4 yields that $\xi_1(\beta,l^-_p) \in \Balphcl$ and
$\frac{1}{2\walpha} \geq \alphtel$, this implies $\beta \in
\left[1+\walphtel,1+\frac{3}{\walpha}\right]$ as required.

\ \\


\ \\
\underline{\textbf{Step 5:}} \ \ \ \textit{Part II of the induction
  statement implies Part III}

\ \\ As in Section~\ref{Onesidedproof}, we suppose that Part II with
$q=p$ holds for all $n\leq N$, with $N \in [l^+_p+1,\nu(p+1)]$, and
show that in this case Part III(a) holds as well whenever $n_1,n_2
\leq N$ and similarly Part III(b) holds whenever $n_2 \leq N$.

\ \\
Let $n_1 \leq N$ be admissible. As we assume that Part II of the
induction statement with $q=p$ holds for $n=n_1$, we can use
Lemma~\ref{lem:parameterestimatesS} to see that
$\frac{\partial}{\partial \beta} \xi_{n_1}(\beta,l^-_p) \ \leq \ 
-\alpha^{\frac{n_1}{4}}$ for all $\beta \in I^p_{n_1}$, which implies
(\ref{eq:indstatementIIIaS}). Then (\ref{eq:indstatementIIIbS}) is a
direct consequence of (\ref{eq:indstatementIIaS}). This proves Part
III(a).  Part III(b) follows in the same way as in Step 3 of
Section~\ref{Onesidedproof}~.

\solidqed


\ \\
\underline{\textbf{Step 6:}} \ \ \ \textit{Proof of Part II for $q=p$.} 

\ \\
In order to prove Part II of the induction statement for $q=p$, we proceed by
induction on $n$. In Step 4 we already constructed $I^p_{l^+_p+1}$ with the
required properties. Now suppose that $I^p_n$ has been constructed for
all admissible $n \in [l^+_p+1,N]$, where $N \in [l^+_p+1,\nu(p+1)-1]$. We now
have to construct $I^p_{N+1}$ with the required properties, provided $N+1$ is
admissible. Again, the case where $N$ is admissible as well is rather easy: In
this case $p(N)=0$, otherwise $N+1$ would be contained in $J(N)$. Therefore
Lemma~\ref{lem:basicestimateS} yields that 
\begin{equation}
  \label{step6S:a} \xi_{N+1}(\beta^+_{p,N},l^-_p) \ > \ \alphtel
\end{equation}
and
\begin{equation}
  \label{step6S:b} \xi_{N+1}(\beta^-_{p,N},l^-_p) \ < \ -\alphtel \ .
\end{equation}
Consequently, we can find $\beta^\pm_{p,N+1} \in I^p_N$ which satisfy 
(\ref{eq:betapluspndefS}) and (\ref{eq:betaminpndefS}), such that %
$I^p_{N+1} = [\beta^+_{p,N+1},\beta^-_{p,N+1}] \ssq I^p_N$. Note that
$R_N = R_{N+1}\setminus\{N+1\}$ by (\ref{eq:regularstabilized}).
Therefore Part II of the induction statement for $n=N$ implies that we
can apply Lemma~\ref{lem:parameterestimatesS} to any $\beta \in
I^p_{N+1}$, and this yields the monotonicity of
$\xi_{N+1}(\beta,l^-_p)$ on $I^p_{N+1}$. All other required statements 
for $n=N+1$ then follow directly from Part II of the induction 
statement for $n=N$.

It remains to treat the case where $N+1$ is admissible but $N$ is not
admissible. As in Step 6 of Section~\ref{Onesidedproof} we have to
consider the interval $J \in \J_{N+1}$ which contains $N$, i.e.\ %
$J = [t,N]$ with $t:=\lambda^-(m_J)$. In order to construct
$I^p_{N+1}$ inside of $I^p_{t-1}$ we prove the following claim
(compare Claim~\ref{claim:step6a}):
\begin{claim}
  \label{claim:step6Sa}
$\xi_{N+1}(\beta^+_{p,t-1},l^-_p) > \alphtel$ and
  $\xi_{N+1}(\beta^-_{p,t-1},l^-_p) < -\alphtel$. 
\end{claim}
\proof\ 
We only give an outline here, the details can be checked exactly as in
the proof of Claim~\ref{claim:step6a}~. Note that it sufficed there to 
show~(\ref{step6:a}), such that the problem is analogous. 

Let $\beta^+:=\beta^+_{p,t-1}$ and $m:=m_J$. First, we can apply
Lemma~\ref{lem:approachingS} with $q=p(m)$, $l=l^*=l^-_p$,
$\beta^*=\beta^+_{p,m}$, $m$ as above, $k=0$ and $\beta=\beta^+$ to
obtain that
\begin{equation}
  \label{step6S:c}
 \{ j \in [-l^-_{p(m)},0] \mid \xi_{j+m}(\beta^+,l^-_p) < \gamma \} \ \ssq \
 \Omega_\infty 
 \end{equation}
 (compare (\ref{step6:b})--(\ref{step6:g})).  Then we have to
 distinguish two cases. If $s(m)=1$, we can proceed as in the proof
 of \ref{claim:step6a} to show that $\xi_{N+1}(\beta^+,l^-_p) \geq
 \frac{2}{\alpha}$. On the other hand suppose $s(m)=-1$, such that
 $d(\omega_{m_J},0) \geq \frac{4\gamma}{L_2}$. In this case
 (\ref{step6S:c}) implies in particular that $\xi_{m}(\beta^+,l^-_p)
 \geq \gamma$, and Lemma~\ref{lem:basicestimateS} therefore yields
 $\xi_{m+1}(\beta^+,l^-_p) \geq \gamma \geq \frac{2}{\alpha}$.
 Similar to the case $s(m)=1$ we can now compare the orbits
\begin{equation}
  \label{step6S:referenceone}
x_1^1 \ld x^1_n \ := \ \zeta_{1}(\beta^+,l^-_p) \ld
\zeta_{l^+_{p(m)}}(\beta^+,l^-_p)
\end{equation}
and 
\begin{equation}
  \label{step6S:newone}
x^2_1 \ld x^2_n \ := \ \xi_{m+1}(\beta^+,l^-_p) \ld
\xi_N(\beta^+,l^-_p) \ ,
\end{equation}
(see (\ref{step6:referencetwo}) and (\ref{step6:newtwo})), with the
difference that it suffices to use Lemma~\ref{lem:orbitsthrowout}(a)
instead of (b). Note that the information we have about the orbit
(\ref{step6S:referenceone}) is exactly the same as for the orbit
(\ref{step6:referencetwo}) (see (\ref{eq:xizetaequivalenceI})). Thus,
we also obtain $\xi_{N+1}(\beta^+,l^-_p) > \alphtel$ in this case.

The proof for $\xi_{N+1}(\beta^-,l^-_p) < -\alphtel$ is then
analogous. This time, it suffices to use
Lemma~\ref{lem:orbitsthrowout}(a) for the case $s(m)=1$, whereas
Lemma~\ref{lem:orbitsthrowout}(b) has to be invoked in order to
compare the orbits $x_1^1 \ld x^1_n$ $:=$ $\zeta_{1}(\beta^+,l^-_p),$
$\ldots,$ $\zeta_{l^+_{p(m)}}(\beta^+,l^-_p)$ and $x^2_1 \ld x^2_n$
$:=$ $\xi_{m+1}(\beta^+,l^-_p),$ $\ldots,$ $\xi_N(\beta^+,l^-_p)$ in
case $s(m)=-1$, but the details for the application are again the same 
as before.

\qed

\ \\ Using the above claim, we see that 
\begin{equation}
  \label{step6S:betamindef}
 \beta^-_{p,N+1} \ := \ \min\left\{\beta\in I^p_{t-1} \mid
   \xi_{N+1}(\beta,l^-_p) = -\alphtel \right\}
\end{equation}
and
\begin{equation}
  \label{step6S:betaplusdef}
 \beta^+_{p,N+1} \ := \ \max\left\{\beta\in I^p_{t-1} \mid
   \beta < \beta^-_{p,N+1},\ \xi_{N+1}(\beta,l^-_p) = \alphtel \right\}
\end{equation}
are well defined, such that $I^p_{N+1} := [\beta^+_{p,N+1},\beta^-_{p,N+1}]
\ssq I^p_{t-1}$. Then, due to Part II of the induction statement for $n=t-1$,
(\ref{eq:indstatementIS}) holds for all $\beta \in I^p_{N+1}$ and similarly
\begin{equation}
  \label{step6S:d}
\xi_j(\beta,l^-_p) \ \in \ \Balphcl \ \ \ \ \ \forall j \in R_{t-1}
\end{equation}
whenever $\beta \in I^p_{N+1}$. As $R_{N+1} = R_{t-1} \cup R(J) \cup
\{N+1\}$, it remains to obtain information about $R(J)$. Thus, in
order to complete this step we need the following claim, which is the
analog of Claim~\ref{claim:step6b}:
\begin{claim}
  \label{claim:step6Sb}
Suppose $\beta\in I^p_{N+1}$ and $\xi_{N+1}(\beta,l^-_p) \in
\Balphcl$. Then $\xi_j(\beta,l^-_p) \in \Balphcl$ $\forall j \in
R(J)$. 
\end{claim}
Similar to Claim~\ref{claim:step6b}, this follows from two further
claims, which are the analogues of Claims~\ref{claim:step6c} and
\ref{claim:step6d}~. Before we state them, let us see how we can use
Claim~\ref{claim:step6Sb} in order to complete the induction step $N
\ra N+1$ and thereby the proof of Step 6:

Suppose that $\beta \in I^p_{N+1}$ and $\xi_{N+1}(\beta,l^-_p) \in
\Balphcl$. Then (\ref{step6S:d}) together with the claim imply that 
\begin{equation}
  \label{step6S:e}
  \xi_j(\beta,l^-_p)  \ \in \ \Balphcl \ \ \ \ \ \forall j \in R_{N+1} 
  \ .
\end{equation}
In addition (\ref{eq:indstatementIS}) holds, as mentioned before
(\ref{step6S:d}). Consequently, Lemma~\ref{lem:parameterestimatesS}
(with $q=p$ and $n=N+1$) implies that
\[
    \frac{\partial}{\partial \beta} \xi_{N+1}(\beta,l^-_p) \ \leq \ -
    \alpha^{\frac{N}{4}} \ .
\]
In particular, this is true for $\beta=\beta^+_{p,N+1}$, and when
$\beta$ is increased it will remain true until
$\xi_{N+1}(\beta,l^-_p)$ leaves $\Balphcl$, i.e.\ all up to
$\beta^-_{p,N+1}$. This proves the required monotonicity of $\beta
\mapsto \xi_{N+1}(\beta,l^-_p)$ on $I^p_{N+1}$, and thus Part II of the 
induction statement holds for $n=N+1$. 
\begin{claim}
  \label{claim:step6Sc}
  Suppose $\xi_k(\beta,l^-_p) \geq \frac{2}{\alpha}$ for some $k \in
  R(J)$. Then $\xi_{N+1}(\beta,l^-_p) > \alphtel$. Similarly, if
  $\xi_k(\beta,l^-_p) \leq -\frac{2}{\alpha}$ then
  $\xi_{N+1}(\beta,l^-_p) < -\alphtel$.
\end{claim}
This is proved exactly as Claim~\ref{claim:step6c}, with all
inequalities reversed for the case $\xi_k(\beta,l^-_p) \leq
-\frac{2}{\alpha}$. 
\begin{claim}
  \label{claim:step6Sd} Suppose $k\in R(J)$, $k+1 \in \Gamma^+(J)$ and 
  $\xi_k(\beta,l^-_p) \geq \alphtel$. Then there exists some
  $\tilde{k} \in R(J)$ with $\xi_{\tilde{k}}(\beta,l^-_p) \geq
  \frac{2}{\alpha}$. Similarly, if $\xi_k(\beta,l^-_p) \leq -\alphtel$
  there exists some $\tilde{k} \in R(J)$ with
  $\xi_{\tilde{k}}(\beta,l^-_p) \leq -\frac{2}{\alpha}$.
\end{claim}
\proof\  This can be shown in the same way as
Claim~\ref{claim:step6d}: Suppose first that $\xi_k(\beta,l^-_p) \geq
\alphtel$ and define $m,\ t$ and $q'$ as in the proof of
Claim~\ref{claim:step6d}. As these definitions only depend on the sets
of regular points, which are the same as before, there is no
difference so far. Only instead of (\ref{step6:j}) we obtain
\begin{equation}
  \label{step6S:f}
 d(\omega_t,\{0,{\textstyle\halb}\}) \ \leq \ \viertel \cdot
 \frac{\alpha^{-(p(m)+1)}}{L_2} 
\end{equation}
Nevertheless, we can apply Lemma~\ref{lem:approachingS}, in the same
way as Lemma~\ref{lem:approaching} was applied in order to obtain
(\ref{step6:o}), to conclude that
\begin{equation} \label{step6S:g}
    \{ j \in [-l^-_{p(m)},0] \mid \xi_{j+m+t}(\beta,l^-_p) < \gamma \} \ \ssq \
    \Omega_\infty 
\end{equation}
(compare (\ref{step6:k})--(\ref{step6:o})).  For the further argument
we have to distinguish two cases. If $s(m+t)=1$ and $p(m) \geq 2$,
then we can use exactly the same comparison arguments as for
Claim~\ref{claim:step6d} (compare
(\ref{step6:referencefour})--(\ref{step6:s})) to show that
$\xi_{m+t+l^+_{p(m)}+1}(\beta,l^-_p) \geq \frac{2}{\alpha}$.  The details
all remain exactly the same.  Thus, we can choose $\tilde{k} =
m+t+l^+_{p(m)}+1$ if $p(m)\geq 2$, and similarly $\tilde{k} = m+t+1$
or $m+t+2$ if $p(m)=1$.

On the other hand, suppose $s(m+t) = -1$. Then
$d(\omega_{m+t},0)\geq\frac{3\gamma}{L_2}$, and in addition (\ref{step6S:g})
implies that $\xi_{m+t}(\beta,l^-_p) \ \geq \
\gamma$. Lemma~\ref{lem:basicestimateS} therefore yields that
$\xi_{m+t+1}(\beta,l^-_p) \geq \gamma \geq \frac{2}{\alpha}$, such that we can
choose $\tilde{k} = m+t+1$.

\ \\ The case $\xi_k(\beta,l^-_p) \leq -\alphtel$ is then treated analogously:
First of all, application of Lemma~\ref{lem:approachingS} yields
\begin{equation} \label{step6S:h}
    \{ j \in [-l^-_{p(m)},0] \mid \xi_{j+m+t}(\beta,l^-_p) > -\gamma \} \ \ssq \
    \Omega_\infty \ ,
\end{equation}
in particular $\xi_{m+t}(\beta,l^-_p) \ \leq \ -\gamma$. If $s(m+t) =
1$, such that $d(\omega_{m+t},\halb)\geq\frac{3\gamma}{L_2}$, then
Lemma~\ref{lem:basicestimateS} yields that $\xi_{m+t+1}(\beta,l^-_p)
\leq -\gamma \leq -\frac{2}{\alpha}$ and we can choose $\tilde{k} =
m+t+1$.

On the other hand, if $s(m+t) = -1$, then we can again apply similar
comparison arguments as in the proof of Claim~\ref{claim:step6d}
(compare (\ref{step6:referencefour})--(\ref{step6:s}))to conclude that
$\xi_{m+t+l^+_{p(m)}+1}(\beta,l^-_p) \leq -\frac{2}{\alpha}$ if $p(m)\geq 2$
(and again $\xi_{m+t+1}(\beta,l^-_p) \leq -\frac{2}{\alpha}$ or
$\xi_{m+t+2}(\beta,l^-_p) \leq -\frac{2}{\alpha}$ if $p(m)=1$). Apart
from the reversed inequalities the only difference now is that the
reference orbits $\xi_{-l^-_{p(m)}}(\beta^*,l^-_{q'}),$ $\ldots,$
$\xi_{-1}(\beta^*,l^-_{q'})$ and $\xi_1(\beta^*,l^-_{q'}) \ld
\xi_{l^+_{p(m)}}(\beta^*,l^-_{q'})$ in (\ref{step6:referencefour}) and
(\ref{step6:q}) have to be replaced by
$\zeta_{-l^-_{p(m)}}(\beta^*,l^-_{q'}),$ $\ldots$,
$\zeta_{-1}(\beta^*,l^-_{q'})$ and $\zeta_1(\beta^*,l^-_{q'}),$
$\ldots,$ $\zeta_{l^+_{p(m)}}(\beta^*,l^-_{q'})$, respectively. Due to
(\ref{eq:xizetaequivalenceI}) and (\ref{eq:xizetaequivalenceII}), all
other details remain exactly the same as before.

\qed

\solidqed

%
%

\end{document}